\documentclass[12pt]{article}
 \usepackage[russian]{babel}
\input{xy}
\xyoption{all}
\usepackage{amsfonts,amsmath,amssymb,amsthm,graphicx,afterpage}

\usepackage{myrusadd,verstka}
\usepackage
{obstruct}
\usepackage{frac}
\usepackage{epsfig}

\def\Q{{\mathbb Q}}
\def\SO{{\mathop SO}}
\def\pro{\mbox{pr}}
\def\cyl{\mbox{cyl}}
\def\con{\mbox{con}}
\def\sgn{\mbox{sgn}}
\long\def\comment#1\endcomment{}

\newtheoremstyle{mydefinition}
  {3pt}
  {3pt}
  {\normalfont}
  {\parindent}
  {\bfseries}
  {.}
  { }
  {}

\theoremstyle{mydefinition}
\newtheorem{pr}{}[section]

\begin{document}

\comment

1. Introduction

Chapter 1. 2-dimensional manifolds.
2. Disks with strips.
3. 2-dimensional manifolds.
4. Thickenings of graphs.
5. Homology of 2-manifolds.
6. Involutions.

Chapter 2. Vector fields and their collections.
7. Vector fields in the plane.
8. Vector fields on 2-manifolds.
9. 3-manifolds.
10. Vector fields on higher-dimensional manifolds.
11. Classification of sections and integrable systems.
12. Collections of vector fields.
13. Non-immersibility and non-embeddability.
14. Fiber bundles and their applications.

Chapter 3. Homology of manifolds.
15. The simplest methods of homology computation.
16. Poincare duality and Alexander-Pontryagin duality.
17. Obstructions for cobordisms.

Chapter 4. Homotopy classification of maps.
18. Brouwer-Hopf theorem
19. Generalizations of the Hopf theorem.
20. Immersions, embeddings, knottedness, and homotopy spheres.

\endcomment

\bigskip
\newpage
\centerline{\bf\uppercase{Алгебраическая топология}}
\bigskip
\centerline{\bf\uppercase{с геометрической точки зрения}
\footnote{Изд-во МЦНМО, Москва, в печати.
Выложено на http://arxiv.org/ с разрешения издательства.
В этой версии не внесена правка в некоторые рисунки, которую в будущем внесет изд-во МЦНМО. 
Некоторые pcx рисунки (например, рис. 5 и 6) не компилируются архивным компилятором.
Их можно посмотреть в более часто обновляемой версии
http://www.mccme.ru/circles/oim/obstruct.pdf .}
}
\bigskip
 \centerline{А.\?Скопенков
\footnote{
Московский Физико-Технический Институт,
Независимый Московский Университет,
Инфо: www.mccme.ru/\~{ }skopenko.
Частично поддержан Российским Фондом Фундаментальных Исследований, Гранты номер 07-01-00648a,
06-01-72551-NCNILa и 12-01-00748-a, Грантом Президента РФ МД-4729.2007.1, стипендией П. Делиня,
основанной на его Премии Бальзана 2004 года, грантами фонда Саймонса 2011-2013 годов и стипендией фонда Династия 2014 года.}
}

\bigskip
{\bf Abstract.}
This book is purely expository and is in Russian.
It is shown how in the course of solution of interesting geometric problems close to applications
naturally appear main notions of algebraic topology
(such as homology groups, obstructions and invariants, characteristic classes).
Thus main ideas of algebraic topology are presented with minimum technicalities.
Familiarity of a reader with basic notions of topology (such as 2-dimensional manifolds and vector fields)
is desirable, although definitions are given at the beginning.
The book is accessible to undergraduates and could also be an interesting easy reading for professional mathematicians.

\bigskip
\hfill{\it Посвящается памяти Юрия Петровича Соловьева}

\tableofcontents

\newpage
\section{Введение}\label{0intr}

\hfill{\it Nothing was changed, but now it made sense.}

\hfill{\it U.\?K.\?Le Guin, The Beginning Place.
\footnote{\it Ничего не изменилось, но теперь все было понятно.
У. К. Ле Гуин, Изначальное место
(пер. автора).}
}


\subsection{Зачем эта книга}

Алгебраическая топология является фундаментальной частью математики и имеет применения за ее пределами.
Как и для любой фундаментальной части науки, ее основные мотивировки и идеи можно доступно изложить человеку,
не имеющему глубоких специальных познаний.
Такому изложению посвящена настоящая книга
(вместе с~\cite{ST04}, \cite{BE82}, \cite{Pr95}, \cite{An03} и, возможно, другими книгами).
Ее особенность --- возможности познакомиться с этими мотивировками и идеями при сведении к необходимому минимуму {\it языка}.
Демонстрируя основные идеи на простейших частных случаях, свободных от технических деталей,
я надеюсь сделать алгебраическую топологию более доступной неспециалистам ---
в первую очередь студентам и математикам, работающим в других областях.

Процесс появления полезных алгебраических понятий продемонстрирован на примере наиболее наглядных топологических задач:
о вложениях графов, о построении и классификации векторных полей, погружений и вложений, о гомотопической классификации непрерывных отображений.
\footnote{Замечу, что важнейшие геометрические проблемы, ради которых была создана алгебраическая топология, в свою очередь
были мотивированы предыдущим развитием математики (причем не только геометрии, но и анализа и алгебры).
Мотивировать эти геометрические задачи не входит в цели настоящей книги.
Я либо привожу ссылки, либо апеллирую к непосредственной геометрической любознательности читателя.}
Определение групп гомологий естественно появится при решении указанных проблем и потому его не обязательно знать заранее.
В то же время для тех, кто уже изучал теорию гомологий, ее применение к конкретным задачам обычно оказывается нетривиальным и интересным.

Такое изложение является не нововведением, а, напротив, возвращением к подходу первооткрывателей
(которое, впрочем, автору обычно приходилось сначала переоткрывать и лишь потом убеждаться, что классики рассуждали так же, ср. \cite{Hi95}).
Такое изложение было обычным в `добурбакистский' период~\cite{ST04}, \cite{CL95}.
Изложение `от простого к сложному' и в форме, близкой к форме рождения
материала, продолжает устную традицию, восходящую к Лао Цзы и Платону,
а в современном преподавании математики представленную, например, книгами Пойа и журналом `Квант'.

Вслед за классиками, я ориентируюсь на объекты, которые основательнее всего укореняются в памяти (или подсознании).
Это --- отнюдь не системы аксиом и не формально-логические схемы доказательств,
а естественные построения для решения интересных проблем или изящные
доказательства красивых теорем, формулировки которых ясны и доступны.
Более того, именно по таким построениям и доказательствам, при наличии
некоторой математической культуры, читатель сможет восстановить более
абстрактный теоретический материал.
Обратное же, как показывает опыт, практически невозможно.
(Изучение `от недостаточно мотивированного общего к частному' часто приводит к абсурдному эффекту: изучившие курс воспроизводят
громоздкое определение, но не могут по этому определению привести ни одного содержательного примера определяемого объекта.)

Алгебраическая топология основана на следующей простой идее, часто встречающейся при решении школьных (в частности олимпиадных) задач.
{\it Невозможность} проделать некоторую конструкцию можно доказывать путем построения алгебраического {\it препятствия}
(или {\it инварианта}).
Например, из соображений четности.
Точно так же {\it неэквивалентность} конструкций часто доказывается путем построения алгебраического {\it инварианта}, их различающего
(этот инвариант является {\it препятствием} к эквивалентности).
Многие непохожие друг на друга задачи топологии аналогичным образом естественно приводят к похожим друг на друга {\it препятствиям}.
Именно ветви алгебраической топологии, непосредственно связанной с {\it теорией гомологий} и {\it теорией препятствий}, посвящена настоящая книга.


Изложение алгебраической топологии, начинающееся с длительного освоения немотивированных абстрактных понятий и теорий делает малодоступными ее
замечательные применения.
Приведу лишь один пример из многих.
Еще в 19 веке был придуман очень простой, наглядный и полезный инвариант многообразий --- форма пересечений, т.е.
умножение в гомологиях многообразий \cite{Hi95}.
Замечательным открытием Колмогорова и Александера 1930-х годов явилось обобщение этого инварианта на произвольные полиэдры (умножение в когомологиях).
Умножение Колмогорова-Александера менее наглядно и определяется более громоздко, чем форма пересечений (но зато имеет более продвинутые применения).
Определение формы пересечений через умножение Колмогорова-Александера делает малодоступными ее замечательные применения (для получения которых не требуется умножения Колмогорова-Александера).
Поэтому форму пересечений иногда просто переоткрывают~\cite{Mo89}.

Надеюсь, принятый стиль изложения не только сделает материал более доступным, но позволит сильным студентам (для которых доступно даже
абстрактное изложение) приобрести математический вкус с тем, чтобы разумно выбирать проблемы для исследования, а также ясно излагать
собственные открытия, не скрывая ошибок (или известности полученного результата) за чрезмерным формализмом.
К сожалению, такое (бессознательное) сокрытие ошибок часто происходит с молодыми математиками, воспитанными на чрезмерно формальных курсах (происходило и с автором этих строк; к счастью, почти все мои ошибки исправлялись {\it перед} публикациями).

Чтение этой книги и решение задач потребуют от читателя усилий, однако эти усилия будут сполна оправданы тем, что вслед за великими математиками
20-го века в процессе изучения геометрических проблем читатель откроет некоторые основные понятия и теоремы алгебраической топологии.
Надеюсь, это поможет читателю совершить собственные настолько же полезные открытия (не обязательно в математике)!

 \subsection{Содержание и используемый материал}

В этой книге изучаются одни из важнейших
объектов математики: многообразия и векторные поля на них.
Для многообразий методы алгебраической топологии наиболее наглядны.
(Например, второй класс Штифеля-Уитни замкнутого трехмерного многообразия есть
$\Z_2$-гомологический класс объединения тех окружностей, на которых линейно
зависимы некоторые два касательных векторных поля общего положения.)
Это позволяет быстро добраться до по-настоящему интересных и сложных результатов.
В этой книге собраны некоторые результаты и методы, касающиеся именно многообразий, а не полиэдров.
Впрочем, для более глубокого изучения многообразий более полиэдры все-таки понадобятся.

Книга предназначена в первую очередь для читателей, не владеющих алгебраической топологией
(хотя возможно, что часть этого текста будет интересна и специалистам).
Все необходимые алгебраические объекты (со страшными названиями группы гомологий, характеристические классы и т. д.)
естественно возникают и строго {\it определяются} в процессе исследования геометрических проблем.
В начале книги знакомство с алгебраическим понятием группы не обязательно, и это слово можно воспринимать  как синоним слова `множество'
(кроме тех мест, где явно оговорено противное).

Для удобства читателя приведены определения и используемые свойства графов (\S\ref{0intr}), двумерных и трехмерных многообразий (\S\ref{02man},\S\ref{03man}).
Эти сведения более полно и подробно изложены, например, в \cite{Pr04},~\cite{BE82}, \cite{MF90}, \cite{Ma03}.

В книге сначала показаны те идеи, которые видны на двумерных многообразиях (\S\ref{0dile}-\S\ref{0ve2ma}).
Затем --- те идеи, которые видны на трехмерных многообразиях (\S\ref{03man}-\S\ref{03hom} и \S\ref{0sec}; \S\ref{0vehima}  интересен даже для частного случая 3-многообразий).
Только потом рассматриваются многомерные многообразия.
При этом мое отношение к двумерным трехмерным многообразиям весьма потребительское --- они все-таки интересны мне не сами по себе, а как простые объекты для демонстрации идей, приносящих наиболее значительные плоды для многомерного случая.

Наиболее блистательные применения некоторой теории --- теоремы, в формулировках которых {\it нет} понятий из этой теории,
но при доказательствах которых без данной теории не обойтись.
Такие применения есть и у алгебраической топологии.
Для удобства читателя такие теоремы выделены жирным шрифтом в тексте и собраны в начале параграфов.
Доказательства этих теорем разбиваются на два шага.
Первый и обычно более простой шаг --- получение необходимого условия на языке теории препятствий --- приводится в этой книге.
Второй и более сложный шаг --- вычисление появляющихся препятствий --- приводится лишь в виде наброска, цикла задач или просто ссылки
(поскольку, по мнению автора, второй шаг лучше описан в литературе, чем первый).
Замечу, что для большинства простейших геометрических проблем очевидно, что полученное необходимое алгебраическое условие является достаточным.
А вот для более сложных геометрических проблем, которые здесь не приводятся (например, о классификации многообразий или вложений),
наиболее трудно именно доказать {\it достаточность} полученного алгебраического условия.

Больш\'ая часть материала сформулирована в виде задач, обозначаемых жирными числами.
Следует подчеркнуть, что многие задачи не используются в остальном тексте.
Красивые наглядные задачи, для решения которых не нужно никаких знаний, приведены уже в самом начале.
Для понимания и решения большинства задач достаточно знакомства с настоящим
текстом (точнее, его части, предшествующей формулировке задачи).
Если используемые в условии термины не определены в этом тексте и вам незнакомы, то соответствующую задачу следует просто игнорировать.
Отметим, что для решения задач достаточно понимания их формулировок и {\it не
требуется} никаких дополнительных понятий и теорий (кроме, быть может, задач со звездочками).

Приводимые задачи являются примерами интересных и полезных фактов, и читателю
будет полезно ознакомиться с самими фактами, даже если он не сможет их самостоятельно доказать.
Например, в некоторых задачах изложен план доказательства теорем, который
полезно понимать, даже если детали этого плана останутся недоступными.
Поэтому приводимые формулировки задач могут быть путеводителем по другим учебникам по алгебраической топологии, позволяя намечать интересные
конечные цели и отбрасывать материал, не являющийся для этих целей необходимым.
Впрочем, полезнее всего было бы обсуждать со специалистом как ваши решения задач, так и возникающие при решении трудности.

Если условие задачи является формулировкой утверждения, то подразумевается, что это утверждение требуется доказать.

Через $|AB|$ или $|A,B|$ обозначается расстояние между точками $A$ и $B$ (из контекста ясно, в каком пространстве).


 \subsection{Для знатоков}

В \S\ref{0vesy3} приводится набросок простого доказательства теорем Штифеля о
параллелизуемости ориентируемых трехмерных многообразий, теоремы об алгебрах с делением
и о невложимости проективных пространств.
Изложение несколько упрощено по сравнению с~\cite{Ki89}.
По возможности приводятся ссылки на книги и обзоры, а не на оригинальные статьи.

Пункты \ref{0cobeul}, \ref{0cobsig} и \ref{0hosp} содержат наборы красивых важных задач по основами теории гомологий и поэтому могут быть использованы на семинарских занятиях по этой теме.

Большинство материала излагается нестандартным образом, более удобным (по мнению автора) для начинающего.
Стандартная терминология теории препятствий не используется там, где (по мнению автора) она неудобна для начинающего.
Приведем здесь сравнение обычной терминологии и принятой в книге.
{\it Расстановки элементов группы $G$ на $i$-симплексах триангуляции $T$} ---
то же, что {\it $i$-мерные цепи на $T$ с коэффициентами в $G$}.
Группа таких расстановок обычно обозначается $C_i(T;G)$.
Множество $\partial C_{i+1}(T;G)$ всех границ образует подгруппу группы
$C_i(K;G)$, обозначаемую $B_i(T;G)$.
Когда $G=\Z_2$, мы пропускаем коэффициенты в
обозначениях цепей, циклов, границ и гомологий.

Кроме того, в этой книге препятствия лежат в группах {\it гомологий}, а не в группах {\it когомологий}
(изоморфных группам гомологий для многообразий).
Именно эта точка зрения (двойственная принятой в учебниках, но обычная для первооткрывателей) позволяет {\it наглядно изображать} препятствия.

В некоторых параграфах мы работаем с многообразиями в кусочно-линейной категории, а в некоторых --- в гладкой.
Поскольку любое гладкое многообразие имеет триангуляцию, согласованную с гладкой структурой, часто мы эту триангуляцию и используем.

 \subsection{Благодарности}

Благодарю А. Н. Дранишникова, Д. Б. Фукса, А. Т. Фоменко и Е. В. Щепина:
я учился алгебраической топологии по книге~\cite{FF89} и на семинаре
Дранишникова-Щепина в Математическом институте Российской академии наук.

Настоящая книга основана на лекциях, прочитанных автором на мехмате
Московского государственного университета им. М. В. Ломоносова, в Независимом
московском университете, на ФИВТ Московского Физико-технического Института,
в Летней школе `Современная математика', а также в
Кировской и Петербургской летних математических школах в 1994-2013 гг.
(Материал статей~\cite{RS00},~\cite{RS02} содержится в настоящей книге в существенно доработанном и расширенном виде.)
Я признателен М.\?Н.\?Вялому, А.\?А.\?Заславскому, С.\?К.\?Ландо, С.В. Матвееву, В.\?В.\?Прасолову, А.Б. Сосинскому, В.\?В.\?Успенскому, В.\?В.\?Шувалову,
а также всем слушателям этих лекций за черновые записки лекций и многочисленные обсуждения, способствовавшие улучшению изложения.
Я признателен В.\?В.\?Прасолову и В.\?В.\?Шувалову за возможность использовать подготовленные ими компьютерные версии рисунков.

Эта книга посвящена памяти Юрия Петровича Соловьева --- замечательного математика, считавшего важным изложение математики на
конкретном, доступном (и в то же время строгом) языке, в отличие от `птичьего' языка излишней абстракции.

\comment

Главы 2, 3, 6 и 7 практически независимы друг от друга (впрочем, в главе 3
рассматриваются примеры, аналогичные главе 2, но в среднем более сложные).
Начинать изучение книги можно с любой из них (при необходимости предварительно
ознакомившись с необходимыми понятиями по главе 1).
При этом сложность материала (и количество используемых понятий) внутри каждой
главы растет.
Поэтому вполне разумно переходить к новой главе, отложив на потом изучение
конца старой.
Более того, разобрав по главам 2, 3, 6 или 7 несколько примеров, мотивирующих
понятие групп гомологий, можно ознакомиться с абстрактным изложением этого
понятия в главе 4.
После этого разумно возвращаться к отложенному материалу.
Формально глава 4 не зависит от предыдущих глав.
Но в ней нет ответа на вопрос 'зачем', важного для начала изучения любой
теории.



\endcomment

 \subsection{Основные понятия теории графов}

Вероятно, вводимые здесь понятия знакомы читателю, но мы приводим четкие
определения, чтобы фиксировать терминологию (которая бывает другой в других книгах).
В б\'ольшей части книги используется только начало этого пункта
поэтому остаток
можно пропустить и возвращаться к нему по мере необходимости.


{\it Графом} (без петель и кратных ребер) называется конечное множество, некоторые
двухэлементные подмножества (т.е. неупорядоченные пары) которого выделены.
Синоним: одномерный симплициальный комплекс (компактный).
Элементы данного множества называются {\it вершинами}.
Выделенные пары вершин называются {\it ребрами}.
Каждое ребро соединяет различные вершины (нет петель), и
любые две вершины соединены не более чем одним ребром (нет кратных ребер).

\begin{figure}[h]\centering
\begin{tabular}{@{}ccc@{}}
\includegraphics{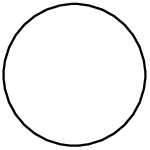}&
\includegraphics{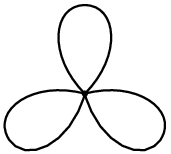}&
\includegraphics{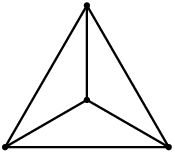}\\
$S^1$ & букет окружностей & $K_4$
\end{tabular}
\caption{Примеры графов (не все вершины отмечены!)}\label{RSo00-1.1}
\end{figure}

При работе с графами удобно пользоваться их изображениями.
Вершины изображаются точками (например, на плоскости или в пространстве).
Каждое ребро, соответствующее двухэлементному выделенному подмножеству,
изображается ломаной (или кривой), соединяющей соответствующие точки.
(Для простоты мы в основном будем рассматривать только такие рисунки, на
которых ребра изображаются ломаными, а не произвольными кривыми.)
На изображении ломаные могут пересекаться, но точки пересечения (кроме двух
концов ребра) не являются вершинами.
Важно, что граф и его изображение --- не одно и то же.


{\it Степенью} вершины графа называется число выходящих из нее ребер.

{\it Путем} в графе называется конечная последовательность вершин, в
которой любые две соседние вершины соединены ребром.
{\it Циклом} называется путь, в котором первая и последняя вершины совпадают.
Граф называется {\it связным}, если любые две его вершины можно соединить
путем.

Граф называется {\it деревом}, если он связен и не содержит
несамопересекающихся циклов.
Ясно, что в любом связном графе существует дерево, содержащее все его вершины.
Такое дерево называется {\it максимальным деревом}.

Вершина графа называется {\it висячей}, если из нее выходит ровно
одно ребро (т. е. если она имеет степень один).  Ребро графа называется
{\it висячим}, если одна из его концевых вершин висячая.

Грубо говоря, {\it подграф} данного графа --- это его часть.
Формально, граф $G$ называется {\it подграфом} графа $H$, если каждая
вершина графа $G$ является вершиной графа $H$, и каждое ребро графа $G$
является ребром графа $H$.
При этом две вершины графа $G$, соединенные ребром в графе $H$, не обязательно
соединены ребром в графе~$G$.

Количества вершин, ребер
и компонент
связности графа обозначаются $V$, $E$ и $C$, соответственно.


Графы $G_1$ и $G_2$ называются {\it изоморфными}, если существует
взаимно однозначное отображение $f$ множества $V_1$ вершин графа $G_1$
на множество $V_2$ вершин графа $G_2$, удовлетворяющее условию:
{\it вершины $A,B\in V_1$ соединены ребром в том и только в том случае, если
вершины $f(A),f(B)\in V_2$ соединены ребром.}

\begin{figure}
[h]\centering
\includegraphics{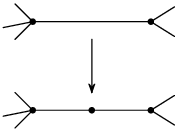}
\caption{Подразделение ребра}\label{podr}
\end{figure}

Операция {\it подразделения ребра} графа показана на рис.~\ref{podr}.
Два графа называются {\it гомеоморфными}, если от одного можно перейти
к другому при помощи операций подразделения ребра и обратных к ним.
Это эквивалентно существованию графа, который можно получить
из каждого из данных графов операциями подразделения ребра.

{\it Одномерным полиэдром (компактным)} (или {\it графом в топологическом
смысле}) называется класс эквивалентности некоторого графа по отношению гомеоморфности.
(Заметим, что одномерные полиэдры часто называют графами, опуская слова 'в топологическом смысле').
Будем сокращенно нызывать одномерный полиэдр {\it 1-полиэдром}.
Представляющий граф называется {\it триангуляцией} соответствующего 1-полиэдра.

Например, через $S^1$ обозначается граф, гомеоморфный циклу (или соответствующий 1-полиэдр).

При работе с 1-полиэдрами мы часто будем пользоваться представляющими их графами,
опуская проверку независимости от выбора представителя в классе эквивалентности.

В этом абзаце (не используемом в дальнейшем) понятия склейки, фигуры (или, научно, топологического пространства)
и то\-по\-ло\-ги\-чес\-кой эквивалентности (или, научно, гомеоморфности) не определяется, а считается
наглядно очевидными или известными читателю.
Каждому графу $G$ соответствует {\it тело графа} $|G|$ --- фигура, получающаяся из конечного числа отрезков
отождествлением некоторых их концов в соответствии с графом $G$.
1-полиэдры находятся во взаимно-однозначном соответствии с фигурами,
топологически эквивалентными телу некоторого графа.
Можно доказать следующий критерий: {\it фигуры $|G_1|$ и $|G_2|$ топологически эквивалентны
тогда и только тогда, когда графы $G_1$ и $G_2$ гомеоморфны}.

Для тополога 1-полиэдры интересны именно как фигуры, а не как классы эквивалентности графов;
графы же средства изучения и хранения в компьютере 1-полиэдров.
А вот дискретному геометру и комбинаторщику интересны в первую очередь графы, но оказываются полезными и 1-полиэдры.



\newpage
\section{Наглядные задачи о поверхностях}\label{0dile}



\hfill{\it Wissen war ein bisschen Schaum, der \"uber eine Woge tanzt.}

\hfill{\it Jeder Wind konnte ihn wegblasen, aber die Woge blieb.}

\hfill{\it E.\?M.\?Remarque, Die Nacht von Lissabon.}
\footnote{Знание здесь --- только пена, пляшущая на волне.
Одно дуновение ветра --- и пены нет. А волна есть и будет всегда.
Э. М. Ремарк, Ночь в Лиссабоне (Пер. Ю. Плашевского).}

\subsection{Графы на поверхностях}

{\it Сферой $S^2$} (стандартной) называется множество точек $(x,y,z)\in\R^3$, для которых
$x^2+y^2+z^2=1$ (или, что то же самое, множество всех точек $(x,y,z)$ вида
$(\cos\varphi\cos\psi,\sin\varphi\cos\psi,\sin\psi)$).
Если Вы не знакомы с декартовыми координатами в пространстве, то здесь и далее координатное определение можно опустить,
а работать с наглядным описанием и изображением на рисунке.

\begin{figure}[h]\centering
\includegraphics{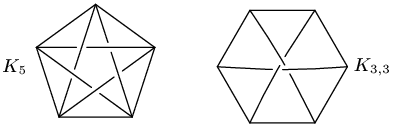}
\caption{Непланарные графы }\label{k5}
\end{figure}

\begin{figure}[h]\centering
\includegraphics{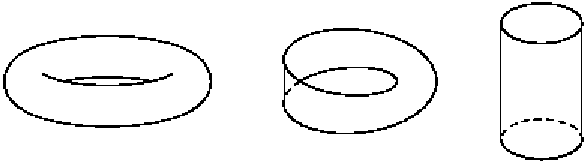}
\caption{Тор, лента Мебиуса и цилиндр}\label{tor}
\end{figure}

\begin{pr}\label{1-sph}
Можно ли нарисовать без самопересечений графы $K_5$ и $K_{3,3}$ (рис.~\ref{k5}) \quad

(a) на сфере? \quad

(b) на боковой поверхности $\{(x,y,z)\in\R^3\ |\ x^2+y^2=1, \ 0\le z\le1\}$ цилиндра (рис. \ref{tor})?
\end{pr}

{\it Тор} --- это поверхность бублика, рис.  \ref{tor}  слева.
Формально, {\it тором} $T^2$ (стандартным) называется фигура, образованная вращением окружности $(x-2)^2+y^2=1$ вокруг оси $Oy$.

\begin{figure}[h]
\hskip2cm\special{em:graph 9-8} \hskip5cm\special{em:graph 5-2}
\vskip2.8cm
\caption{Склейки прямоугольной полоски, дающие тор и ленту Мебиуса}
\label{glu}
\end{figure}

Эта фигура получена из (двумерного) квадрата склейкой  пар его противоположных сторон
`с одинаковыми направлениями', т.е. без поворота, рис.  \ref{glu} слева.

\begin{figure}[h]
\hskip4cm
\special{em:graph 5-1}
\vskip2.8cm
\caption{Кольцо; гомеоморфные фигуры (убрать правую часть)}\label{5-1}
\end{figure}

{\it Кольцом} называется любая фигура, полученная из прямоугольной полоски склейкой
ее двух противоположных сторон `с одинаковым направлением', рис. \ref{5-1}.
Например, боковая поверхность цилиндра является кольцом.

{\it Лентой Мебиуса} называется любая фигура, полученная из длинной прямоугольной полоски склейкой двух ее противоположных сторон
`с противоположным направлением', т.е. с поворотом на $180^\circ$, рис. \ref{glu} справа.
{\it Стандартной лентой Мебиуса} называется поверхность в $\R^3$, заметаемая стержнем длины 1,
равномерно вращающимся относительно своего центра, при равномерном движении этого центра
по окружности радиуса 9, при котором стержень делает пол-оборота, рис. \ref{tor}.

\begin{pr}\label{eul-neul}
(a) Нарисуйте на торе две замкнутые несамопересекающиеся кривые
так, чтобы при разрезаниях по первой и по второй тор распадался бы на разное количество кусков.

(a') То же для ленты Мебиуса.
\end{pr}

\begin{pr}\label{eul-bet}
(a) Нарисуйте две замкнутые кривые на торе, при разрезании по объединению которых тор не распадается на куски.


(b) Существуют ли три замкнутые кривые на торе, при разрезании по объединению которых тор не распадается на куски?

(b') Существуют ли две замкнутые кривые на ленте Мебиуса, при разрезании по объединению которых лента Мебиуса не распадается на куски?
\end{pr}

Тор, лента Мебиуса (и другие фигуры) предполагаются {\it прозрачными}, т.е. точка (или подмножество),
`лежащая на одной стороне поверхности', `лежит и на другой стороне'.
Это аналогично тому, что при изучении геометрии мы говорим, например, о треугольнике на плоскости,
а не о треугольнике на верхней (или нижней) стороне плоскости.

\begin{pr}\label{1-draw}
Можно ли нарисовать на торе без самопересечений граф

(a) $K_5$?  \quad (a') $K_{3,3}$? \quad
(b) $K_8$? \quad (b') $K_{5,4}$? \quad
(c) $K_6$? \quad (c') $K_{3,4}$? \quad
(d) $K_7$? \quad (d') $K_{4,4}$?
\end{pr}

\begin{pr}\label{1-drawmob}
Можно ли нарисовать на ленте Мебиуса без самопересечений граф

(a) $K_{3,3}$? \quad (a') $K_6$? \quad (b) $K_7$? \quad (b') $K_{4,4}$? \quad (b'') $K_{3,5}$? \quad  (c) $K_{3,4}$?
\end{pr}

\begin{pr}\label{map}
{\it Картой на торе} называется разбиение тора на (криволинейные и изогнутые) многоугольники.
Раскраска карты на торе называется {\it правильной}, если разные многоугольники, имеющие общую граничную кривую, имеют разные цвета.
Любую ли карту на торе можно правильно раскрасить в
\quad
(a) 5 цветов? \quad (a') 6 цветов? \quad (b) 7 цветов?
\end{pr}

\begin{pr}\label{intplan}
По планете Тополога, имеющей форму тора, текут реки Меридиан и Параллель.
Маленький принц и Тополог прошли по планете и вернулись в их исходные точки (которые различны).
Маленький принц переходил Меридиан 9 раз и Параллель 6 раз, а Тополог --- 8 раз и 7 раз, соответственно.
Докажите, что их пути
пересекались.
\footnote{
При переходе реки персонаж оказывается на {\it другой} ее стороне.
Более формально, пересечение реки и пути персонажа {\it трансверсально}, см. строгое определение в \S\ref{02hom}.
\newline
При решении пунктов (b), (b') и (b'')
задач \ref{eul-bet}-\ref{map}, а также задачи \ref{intplan}, можно использовать неравенства Эйлера (задача \ref{euler}).
}
\end{pr}

\begin{figure}[h]
\centering{\special{em:graph 7-8} \includegraphics{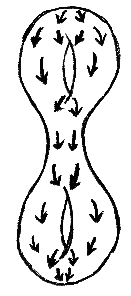}}
\caption{Сфера с тремя ручками и поле скоростей воды, стекающей по сфере с двумя ручками
}
\label{water}
\end{figure}


{\it Сферой с $g$ ручками} (стандартной) $S_g$ при $g\ge1$ называется поверхность, заданная в $\R^3$ уравнением $x^2+\Pi_{k=1}^g((z-4k)^2+y^2-4)=1$.
{\it Сферой с нулем ручек} (стандартной) $S_0$ называется сфера $S^2$.
Сферой с $g$ ручками изображена на рис. \ref{tor} справа и \ref{water} (для двух и трех ручек; для одной ручки это тор).

\begin{figure}[h]
\hskip5cm
\special{em:graph 9-12}
\vskip1.4cm
\caption{`Цепочка окружностей' на плоскости}
\label{9-12}
\end{figure}

Уравнение $\Pi_{k=1}^g((z-4k)^2+y^2-4)=0$ задает `цепочку окружностей' на плоскости $Oyz$ (рис. \ref{9-12}).

\begin{figure}[h]
\hskip5cm\special{em:graph 9-7}
\vskip2.1cm
\caption{Приклеивание ручки}\label{9-7}
\end{figure}

Поэтому сфера с $g$ ручками получена из сферы вырезанием $2g$ дисков и последующей заклейкой $g$ пар краевых окружностей этих дисков криволинейными боковыми поверхностями цилиндров, рис. \ref{9-7}.


\begin{figure}[h]\centering
\caption{Диск с $n$ лентами Мебиуса}\label{f-nmob}
\end{figure}

{\it Диском с $n$ лентами Мебиуса} называется объединение диска и $n$ `отделенных' ленточек,
при котором каждая ленточка приклеивается двумя отрезками к краевой окружности $S$ диска и направления на этих отрезках, задаваемые произвольным направлением на $S$, `сонаправлены вдоль ленточки', рис. \ref{f-nmob}.

\begin{figure}[h]\centering
\includegraphics{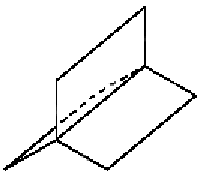}
\caption{Книжка с тремя листами}\label{TI}
 \end{figure}

Возьмем в трехмерном пространстве $n$ прямоугольников $XYB_kA_k$, $k=1,2,\dots,n$, любые два из которых пересекаются только по отрезку $XY$.
{\it Книжкой с $n$ листами} называется объединение этих прямоугольников, см. рис. \ref{TI} слева для $n=3$.

\begin{pr}\label{1-drany} Любой граф можно нарисовать без самопересечений
\quad

(a) в  пространстве.
\quad

(b) на сфере с некоторым количеством ручек, зависящим от графа.

(c) на диске с некоторым количеством лент Мебиуса, зависящим от графа.

(d) в книжке с некоторым количеством листов, зависящим от графа.
\quad

(e) в книжке с тремя листами (рис. \ref{TI}).
\end{pr}

\subsection{Применения неравенства Эйлера}

В этом параграфе слово `поверхность' понимается не как математический термин,
а как собирательное название сфер с ручками и дисков с лентами Мебиуса.

Пусть на поверхности нарисован без самопересечений граф.
Назовем {\it гранью} каждый из связных кусков, на которые распадается поверхность при разрезании по всем ребрам графа.
Задача \ref{eul-neul} показывает, что количество граней зависит от способа изображения графа.

\begin{pr}\label{euler}
{\bf Неравенства Эйлера.}
{\it Пусть на сфере с $g$ ручками или на диске с $m$ лентами Мебиуса нарисован без самопересечений  связный граф с $V$ вершинами и $E$ ребрами.
Для диска с лентами Мебиуса предположим, что граф не пересекает краевой окружности.
Обозначим через $F$ число граней.
Тогда соответственно}
$$V-E+F\ge 2-2g\quad\text{и}\quad V-E+F\ge2-m.$$
\end{pr}

Вряд ли у Вас получится доказать эти неравенства без знания дальнейшего материала.
Сначала порешайте следующие задачи на применение неравенств.

\begin{pr}\label{eul-betti} {\bf Теорема Бетти-Римана.}

(a) {\it Объединение любых $2g+1$ окружностей (возможно, пересекающихся) на сфере с $g$ ручками разбивает ее.}

(b) {\it Объединение любых $m+1$ окружностей (возможно, пересекающихся) на диске с $m$ лентами Мебиуса разбивает его.}
\end{pr}

\begin{pr}\label{1-nonrea}
(a) $K_8$ не вложим в тор.

(b) Для любой сферы с ручками найдется граф, не вложимый в нее.

(a') $K_7$ не вложим в ленту Мебиуса.

(b') Для любого диска с лентами Мебиуса найдется граф, не вложимый в него.
\end{pr}

{\it Ориентируемым родом} $g(G)$ графа $G$ называется наименьшее число $g$, для которого $G$ вложим в сферу с $g$ ручками.
Например, ориентируемый род графов $K_3$ и $K_4$ равен 0, графов $K_5,K_6$ и $K_7$ равен 1, а графа $K_8$ равен 2.

{\it Неориентируемым родом} $m(G)$ графа $G$ называется наименьшее число $m$, для которого $G$ вложим в диск
с  $m$ лентами Мебиуса.
Например, неориентируемый род графа $K_6$ равен 1, а графа $K_7$ равен 2.

\begin{pr}\label{eul-orge} (a) $g(K_n)\ge (n-3)(n-4)/12$; \qquad (b) $g(K_{m,n})\ge (m-2)(n-2)/4$.

На самом деле, в этих формулах неравенство можно заменить на равенство
\cite{Pr04}.

(a',b') Найдите аналогичные оценки на $m(K_n)$ и $m(K_{l,n})$.
\end{pr}

\begin{pr}\label{eul-heaw}
(a) Если $g>0$ и $m$ --- наибольшее целое число, для которого $g\ge(m-2)(m-3)/12$,
то любой граф на сфере с $g$ ручками имеет вершину степени не более $m$.

(b) {\bf Теорема Хивуда.} {\it Если $g>0$ и $n$ --- наибольшее целое число, для которого $g\ge(n-3)(n-4)/12$, то любой граф на сфере с $g$ ручками можно правильно раскрасить в $n$ цветов.}

(Ввиду результатов Рингеля о вложениях графа $K_n$ \cite{Pr04}, меньшего числа цветов не хватит для $K_n$.)

(с) Любой граф на ленте Мебиуса можно правильно раскрасить в 10 цветов.

(d) Сформулируйте и докажите аналог теоремы Хивуда для диска с $m$ лентами Мебиуса.
\end{pr}


\subsection{Наглядные задачи о разрезаниях}


\begin{pr}\label{1-mob}
(a) Можно ли ленту Мебиуса разрезать так, чтобы получилось кольцо?

(b) Можно ли из ленты Мебиуса вырезать непересекающиеся кольцо и ленту Мебиуса?

(c) Можно ли ленту Мебиуса разрезать на непересекающиеся кольцо и ленту Мебиуса?

(d)* Можно ли из ленты Мебиуса вырезать две непересекающиеся ленты Мебиуса?
(Вряд ли у Вас получится решить эту задачу без изучения дальнейшего материала.)
\end{pr}

\begin{figure}[h]
\hskip5cm\special{em:graph 9-4}
\vskip2.2cm
\caption{Бутылка Клейна: склейка квадрата (переместить с рис. \ref{9-5}) и изображение в $\R^3$}\label{kleinb}
\end{figure}

Рассмотрим в $\R^4$ окружность $x^2+y^2=1$ и семейство ее нормальных трехмерных плоскостей.
Рассмотрим в $\R^4$ окружность $x^2+y^2=1$ и семейство ее нормальных трехмерных плоскостей.
{\it Бутылкой Клейна} (стандартной) $K$ называется поверхность в $\R^4$, заметаемая окружностью $\omega$, центр
которой равномерно описывает окружность $x^2+y^2=1$, а окружность $\omega$ в то же время равномерно поворачивается
на угол $\pi$ (поворачивается в движущейся нормальной трехмерной плоскости относительно своего диаметра, движущегося
вместе с нормальной трехмерной плоскостью).
См. рис. \ref{kleinb}.
{\it Бутылкой Клейна} (нестандартной) называется любая фигура, полученная из прямоугольной полоски такой склейкой ее пар
противоположных сторон, при которой одна пара склеивается `с одинаковым направлением', а другая `с противоположным направлением', рис. \ref{kleinb}.

\begin{pr}\label{1-klein}
(a) Разрежьте бутылку Клейна на две ленты Мебиуса.

(b) Разрежьте бутылку Клейна так, чтобы получился (одна) лента Мебиуса.
\end{pr}

{\it Сферой с дыркой} (нестандартной) называется сфера, из которой удалена внутренность двумерного диска.
Аналогично определяются (нестандартные) тор с дыркой, сфера с ручками и дыркой, бутылка Клейна с дыркой и т.д.
Формально, (стандартной) {\it сферой с $g$ ручками и дыркой}
называется часть сферы с $g$ ручками,
лежащая не выше той плоскости, которая чуть ниже касательной плоскости в верхней точке (т.е. плоскости $z\le 4g+1$).

\begin{pr}\label{1-cut} Вырежьте из книжки с тремя листами (рис. \ref{TI})
\quad

(a) ленту Мебиуса. \quad
(b) тор с дыркой. \quad

(c) сферу с двумя ручками и одной дыркой. \quad
(d) бутылку Клейна с дыркой. \quad
\end{pr}

\subsection{Топологическая эквивалентность (гомеоморфность)}

\begin{figure}[h]
\hskip4cm
\special{em:graph 1-8}
\vskip1.6cm
\caption{Тор с дыркой, naoborot}\label{tordyr}
\end{figure}

В этом параграфе понятие {\it гомеоморфности} (то\-по\-ло\-ги\-чес\-кой эквивалентности) не определяется строго,
см. строгое определение в \S\ref{02man}.
Для `доказательства' гомеоморфности в этом параграфе нужно нарисовать цепочку картинок, аналогичную рис.~\ref{tordyr}.
При этом разрешается временно {\it разрезать} фигуру, а потом {\it склеить} берега разреза.
Например,

$\bullet$ сфера с одной ручкой (рис. \ref{9-7}, \ref{tor}) гомеоморфна тору (рис. \ref{tor}).

$\bullet$ 
диск с двумя ленточками на рис.~\ref{tordyr} справа гомеоморфен тору с дыркой (рис.\ref{tordyr} слева).

$\bullet$ три ленточки на рис. \ref{glu} гомеоморфны (здесь уже не обойтись без разрезания).

$\bullet$ две ленточки на рис. \ref{5-1} справа гомеоморфны (и здесь не обойтись без разрезания).

Ленточки на рис. \ref{glu} и на рис. \ref{5-1} справа не гомеоморфны.
Мы займемся {\it негомеоморфностью} в следующем параграфе, когда появится строгое определение.

\begin{figure}[h]\centering
\includegraphics{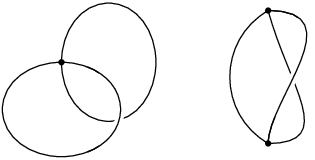}
\caption{Утолщения слов  $(abab)$, $(abacbc)$ и $(abcabc)$ (изменить)}\label{dilen}
\end{figure}

\begin{figure}[h]
 \special{em:graph 5-5} \hskip8cm \includegraphics{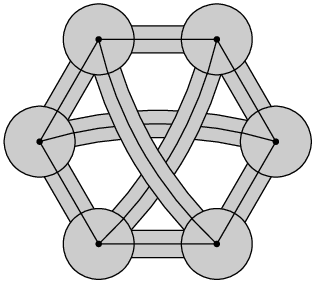}
\caption{Чему гомеоморфны эти фигуры?}\label{td}
\end{figure}

\begin{pr}\label{1-thi}
(a,b)
Фигуры на рис. рис.~\ref{dilen} в центре и справа гомеоморфны тору с двумя дырками.

(c) Фигура на рис.~\ref{td} слева гомеоморфна тору с дыркой.

(d) Гомеоморфна ли фигура на рис.~\ref{td} справа сфере с ручками и дырками? Если да, то чему равно их число?
\end{pr}

\begin{figure}[h]\centering
\begin{tabular}{@{}cccc@{}}
\includegraphics{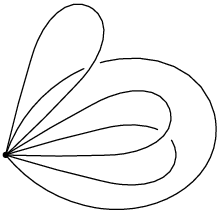}&
\includegraphics{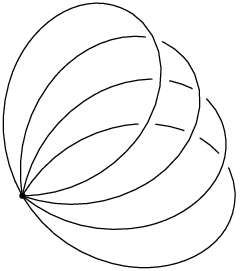}&
\includegraphics{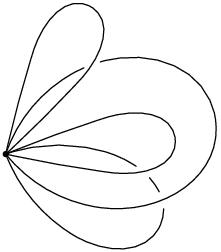}&
\includegraphics{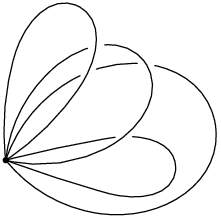}\\
\end{tabular}
\caption{Утолщения слов из четырех букв
(переделать!)}
\label{netor}
\end{figure}

\begin{pr}\label{1-thi4}
(a,b,c,d)
Фигуры на рис.~\ref{netor} гомеоморфны сфере с двумя ручками и одной дыркой.
\end{pr}

\begin{figure}[h]\centering
\caption{(a) Вывернутая ручка.
\newline
(b) Диск с двумя `перекрученными' `отделенными' ленточками.
\newline
(c) Утолщение слова $(aabcbc)$ с соответствием $w(a)=1$ и $w(b)=w(c)=0$}
\label{pereruchka}
\end{figure}

\begin{figure}[h]\centering
\caption{(a) Краевые окружности ленты Мебиуса с дыркой равноправны?
\newline
(b) Гомеоморфны ли кольца с двумя лентами Мебиуса?}
\label{f-brain}
\end{figure}

\begin{pr}\label{1-nohom}
(a) Лента Мебиуса  с ручкой гомеоморфна ленте Мебиуса с вывернутой ручкой, рис.~\ref{9-7}, \ref{pereruchka}.a.

(b) Диск с двумя `перекрученными' `отделенными' ленточками, рис. \ref{pereruchka}.b,
гомеоморфен бутылке Клейна (рис.~\ref{kleinb}) с дыркой.

(c) Фигура на рис. \ref{pereruchka}.c гомеоморфна диску с тремя лентами Мебиуса

(d) Фигуры на рис. \ref{f-brain}.a гомеоморфны.

(e) Кольцо с двумя `перекрученными' `отделенными' ленточками, приклеенными к одной краевой окружности кольца,
гомеоморфно кольцу с двумя `перекрученными' ленточками, приклеенными к разным краевым окружностям кольца, рис. \ref{f-brain}.b.
\end{pr}

\begin{figure}[h]\centering
\includegraphics{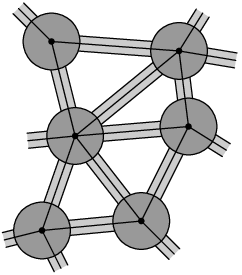}
\caption{Шапочки, ленточки и заплатки}\label{cap}
\end{figure}

{\it Регулярной окрестностью} графа в поверхности называется объединение шапочек и ленточек,
соответствующих вершинам и ребрам графа, рис. \ref{cap} (см. формальное определение в конце \S\ref{02man}).


\begin{figure}[h]\centering
\includegraphics{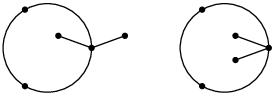}
\caption{Различные вложения графа в плоскость}\label{razvlo}
\end{figure}


\begin{pr}\label{intreg}
(a) Регулярные окрестности образов разных вложений графа в плоскость (рис.~\ref{razvlo}) гомеоморфны.

(b) Регулярные окрестности образов разных вложений некоторого графа в поверхность не обязательно гомеоморфны.

(с) Регулярные окрестности образов изотопных вложений данного графа в поверхность гомеоморфны.

Два вложения $f_0,f_1:G\to N$ называются {\it изотопными}, если существует семейство $U_t:N\to N$
гомеоморфизмов, непрерывно зависящее от параметра $t\in[0,1]$, для которого $U_1\circ f_0=f_1$.

(d) Регулярные окрестности образов гомотопных вложений данного графа в поверхность гомеоморфны.
Два непрерывных отображения $f_0,f_1:G\to N$ называются {\it гомотопными},
если существует семейство $f_t:G\to N$ непрерывных отображений, непрерывно зависящее от параметра $t\in[0,1]$.
Вряд ли у Вас получится доказать этот факт без знания дальнейшего материала.
\end{pr}

\subsection{Топологическая эквивалентность дисков с ленточками}

В этом пункте,
не используя даже понятия графа, мы покажем одну из основных идей доказательства теоремы классификации 2-многообразий (\S\ref{02man}).

\begin{figure}[h]\centering
\includegraphics{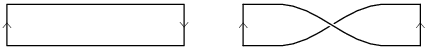}
\caption{Стрелки, противонаправленные вдоль ленточки}\label{lent}
\end{figure}

Пусть имеется слово длины $2n$ из $n$ букв, в котором каждая буква встречается дважды.
Возьмем двумерный диск (например, выпуклый многоугольник на плоскости).
Ориентируем его краевую (=граничную) окружность.
Отметим на ней непересекающиеся отрезки, отвечающие буквам данного слова, в том порядке, в котором буквы идут в слове.
Для каждой буквы соединим (не обязательно в плоскости) соответствующие ей два отрезка ленточкой-прямоугольником
(так, чтобы разные ленточки не пересекались).
 При этом стрелки на окружности должны быть {\it противонаправлены} вдоль ленточки, рис. \ref{lent}.
{\it Утолщением} данного слова, или {\it диском с $n$ неперекрученными ленточками} называется
объединение построенных диска и ленточек.

Примеры дисков с неперекрученными ленточками
приведены на рис.~\ref{dilen} и \ref{netor}.

В этом параграфе понятие {\it краевой окружности}  не определяется строго.
Неформально, {\it краевая окружность} диска с ленточками --- связный кусок множества его точек, к которым он подходит `с одной стороны'.
Краевые окружности дисков с ленточками на рис.~\ref{dilen} выделены жирным [сделать это].


\begin{pr}\label{1-boucir}
(a) Сколько краевых окружностей может быть у диска с двумя неперекрученными ленточками?

(b) По слову длины $2n$ из $n$ букв, в котором каждая буква встречается дважды, постройте граф, число компонент
связности которого равно числу краевых окружностей утолщения данного слова.
(Значит, это число можно находить на компьютере, не рисуя рисунка.)
\end{pr}

\begin{pr}\label{1-eul}
(a)  {\bf Формула Эйлера.} Диск с $n$ неперекрученными ленточками,
имеющий $F$ краевых окружностей, гомеоморфен сфере с $(n+1-F)/2$ ручками и $F$ дырками.

(b)* {\bf Формула Мохара.}
Пусть имеется диск с $n$ ленточками.
Построим матрицу $n\times n$ следующим образом.
Если $a\ne b$ и буквы $a$ и $b$ в соответствующем слове чередуются --- идут в порядке $abab$, а не $aabb$ (т.е. ленточки $a,b$ 'перекрещиваются'), то в клетке $a\times b$ поставим единицу.
В остальных клетках поставим  нули.
Обозначим через $r$ ранг над $\Z_2$ полученной матрицы.
Тогда $r$ четно и диск с ленточками гомеоморфен сфере с $r/2$ ручками и некоторым количеством дырок.
\end{pr}


Пусть имеется слово длины $2n$ из букв $1,2,\dots,n$, в котором каждая буква встречается дважды, и отображение $w:\{1,2,\dots,n\}\to\{0,1\}$.
Возьмем двумерный диск.
Ориентируем его краевую окружность.
Отметим на ней непересекающиеся отрезки, отвечающие буквам данного слова, в том порядке, в котором буквы идут в слове.
Для каждой буквы соединим (не обязательно в плоскости) соответствующие ей два отрезка ленточкой-прямоугольником (так, чтобы разные ленточки не пересекались).
При этом стрелки на окружности должны быть {\it противонаправлены} вдоль ленточки $k$, если $w(k)=0$, и
{\it сонаправлены} вдоль ленточки, $w(k)=1$.
{\it Утолщением} данного слова с данным отображением $w$, или {\it диском с $n$ ленточками} называется
 объединение построенных диска и ленточек.

На рис. \ref{pereruchka}.bc и \ref{f-nmob} изображены соответственно

$\bullet$ утолщение слова $(aabb)$ с соответствием $w(a)=w(b)=1$,

$\bullet$ утолщение слова $(aabcbc)$ с соответствием $w(a)=1$ и $w(b)=w(c)=0$,

$\bullet$ диск с $n$ лентами Мебиуса, т.е. утолщение слова $(1122\dots nn)$ с соответствием $w(1)=w(2)=\dots=w(n)=1$.

\begin{pr}\label{1-thino}
(a) Сколько краевых окружностей может быть у диска с двумя ленточками?

(b) Чему может быть гомеоморфен диск с двумя ленточками?
\end{pr}

\begin{pr}\label{1-eulno}
(a) К одной из краевых окружностей диска с $n$ лентами Мебиуса  и $k>0$ дырками приклеим перекрученную
(относительно этой краевой окружности) ленточку.
Полученная фигура гомеоморфна  диску с $n+1$ лентой Мебиуса и $k$ дырками.

(b) {\bf Формула Эйлера.}
Диск с $n$ ленточками, среди которых есть перекрученная, имеющий $F$ краевых окружностей,
гомеоморфен диску с  $n+1-F$ лентами Мебиуса и $F-1$ дыркой.

(c)* {\bf Формула Мохара.}
Пусть имеется слово длины $2n$ из букв $1,2,\dots,n$, в котором каждая буква встречается дважды, и ненулевое  отображение $w:\{1,2,\dots,n\}\to\{0,1\}$.
Построим матрицу $n\times n$ следующим образом.
Если $a\ne b$ и буквы $a$ и $b$ в соответствующем слове чередуются --- идут в порядке $abab$, а не $aabb$,
то в клетке $a\times b$ поставим единицу.
В диагональной клетке $a\times b$ поставим число $w(a)$.
В остальных клетках поставим  нули.
Обозначим через $r$ ранг над $\Z_2$ полученной матрицы.
Тогда соответствующий диск с ленточками гомеоморфен диску с $r$ лентами Мебиуса и некоторым количеством дырок.
\end{pr}

\subsection{Ответы и указания к некоторым задачам}

\smallskip
{\bf \ref{1-sph}.} (a), (b). Нельзя.

(a) Если бы  граф $K_5$ был нарисован на сфере без самопересечений, то выкинув из сферы точку,
не лежащую на $K_5$, мы получили бы плоскость, содержащую $K_5$.
Противоречие.

(b) Граф $K_5$ не планарен, а цилиндр можно спроектировать без самопересечений на плоскость.

\begin{figure}[h]\centering
\includegraphics{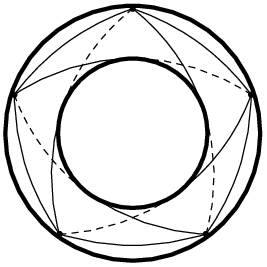}\qquad \includegraphics{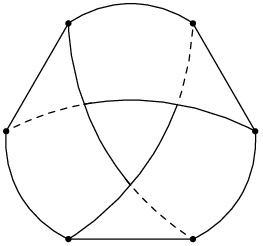}
\caption{Реализация непланарных графов}\label{real}
\end{figure}

\smallskip
{\bf \ref{1-draw}.} Ответы: (b), (b') нельзя, остальные пункты --- можно.

Красивая реализация графа $K_5$ на торе изображена на рис.~\ref{real}.

Имеются также другие решения.
Например, можно нарисовать граф $K_5$ на плоскости с {\it одним} самопересечением и...

\smallskip
{\bf \ref{1-drawmob}.} Ответы: (b), (b') нельзя, остальные пункты --- можно.

Красивая реализация графа $K_{3,3}$ на ленте Мебиуса изображена на рис.~\ref{real}.

\smallskip
{\bf \ref{1-drany}.}
(a) Нарисуем данный граф (возможно, с самопересечениями) на плоскости так, чтобы ребра не самопересекались.
Если образовались точки пересечения (кроме вершин) более чем двух ребер, то
подвинем некоторые ребра так, чтобы остались только двукратные точки пересечения.
Поднимем одно из каждых двух пересекающихся ребер в пространство так, чтобы каждое пересечение пропало.


(b,c) Используйте идею решения пункта (a).

(e) Можно считать, что точки самопересечения `хорошие' и лежат на одной прямой.
 Прилепим третий лист по этой прямой.
Теперь в малой окрестности каждой из точек пересечения ребер поднимем одно из ребер
`мостиком' над другим ребром на третий лист.
Так все точки пересечения будут ликвидированы.

\bigskip
{\bf \ref{eul-betti}.} \cite[\S11.4]{Pr04}

\smallskip
{\bf \ref{1-nonrea}.}
(a) Из неравенств Эйлера и $2E\ge 3F$ получается противоречие.

(b) Аналогично (a) граф $K_{g+15}$ не вложим в сферу с $g$ ручками.

(a',b') Аналогично (a).

\smallskip
{\bf \ref{eul-orge}.} Аналогично \ref{1-nonrea}.a.

\smallskip
{\bf \ref{eul-heaw}.} (a,b) \cite[\S13.2]{Pr04}

\begin{figure}[h]
\hskip4cm
\special{em:graph 9-5}
\vskip2.7cm
\caption{Разрезы бутылки Клейна (добавить разрез справа!) }\label{9-5}
\end{figure}

\smallskip
{\bf \ref{1-klein}.} (a) Разрежьте рис. \ref{kleinb} плоскостью симметрии. Или см. рис. \ref{9-5}.b (пока нет).

(b) См. рис. \ref{9-5}.a.

\begin{figure}[h]\centering
\caption{Утолщения слов $(abcabc)$ и $(abacbc)$ на торе}\label{abc}
\end{figure}

\smallskip
{\bf \ref{1-thi}.}
(a), (b) См. рис. \ref{abc}.

\begin{figure}[h]\centering
\includegraphics{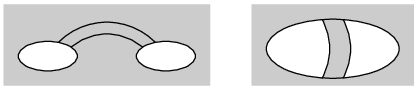}
\caption{Добавление неперекрученной ленточки (поправить!)}\label{strip}
\end{figure}

(a,b) Другое решение. Утолщение слова $(abab)$ (рис.~\ref{dilen} слева) гомеоморфно тору с дыркой.
При добавлении ленточки получается тор с двумя дырками, см. рис.~\ref{strip} справа.
Докажите, что добавляемая ленточка не перекручена относительно той краевой окружности, к которой она крепится.

(d) Выделите максимальное дерево и докажите, что эта фигура гомеоморфна диску с 4 ленточками.


\smallskip
{\bf \ref{1-thi4}.}
(b) Утолщение слова $(abab)$ (рис.~\ref{dilen} слева) гомеоморфно тору с дыркой.
Он лежит в утолщении слова $(abcdabcd)$.
При замене этого тора с дыркой на диск получится утолщение слова из 4 букв с одной краевой окружностью, т.е. тор с дыркой.
(Можно и непосредственно проверить, что получится утолщение слова $(cdcd)$.)
Значит, утолщение слова $(abcdabcd)$ получено из двух торов с дырками склейкой по отрезкам краевых окружностей.
Поэтому данное утолщение гомеоморфно сфере с 2 ручками с дыркой.

(b) Другое решение.  Выберем диск с тремя ленточками, рассмотренный в задаче \ref{1-thi}.b, рис. \ref{dilen}.
Он гомеоморфен тору с двумя дырками.
При добавлении четвертой ленточки добавляется ручка и удаляется дырка, см. рис.~\ref{strip} слева.
Докажите, что добавляемая ленточка не перекручена относительно `согласованных' ориентаций тех двух краевых окружностей, к которым она крепится.

\begin{figure}[h]\centering
\caption{Внутренняя и внешняя краевые окружности ленты Мебиуса с дыркой меняются местами}
\label{f-nonor}
\end{figure}

\smallskip
{\bf \ref{1-nohom}.} (b) Используйте \ref{1-klein}.a.

(с) Фигура на рис. \ref{pereruchka}.c гомеоморфна ленте Мебиуса с ручкой.
Диск с тремя лентами Мебиуса гомеоморфен ленте Мебиуса с вывернутой ручкой ввиду (b).
Используйте (a).

(d)  См. рис. \ref{f-nonor}.

(e) Используйте (d).

\bigskip
{\bf \ref{1-boucir}.} (b) 1 или 3.

\smallskip
Ленточки в диске с ленточками называются {\it перекрещивающимися}, если отрезки, по которым они приклеиваются к диску,  чередуются на его краевой окружности (т.е. идут в порядке  $(abab)$, а не $(aabb)$).

\smallskip
{\bf \ref{1-eul}.}
(a) Сначала докажите для случая, когда нет перекрещивающихся ленточек.
Если же есть две перекрещивающиеся ленточки, то остальные ленточки приклеиваются к краевой окружности получившегося тора с дыркой.
Они не перекручены относительно нее (докажите!).
Заменим тор с дыркой на диск.
При этом количество ленточек уменьшится на 2, а количество краевых окружностей не изменится.
(Расположение ленточек не обязательно получается из исходного вычеркиванием двух букв.)
Применим предположение индукции.

{\it Указание к другому решению.}
Индукция по $n$ (аналогично задачам \ref{1-thi} и \ref{1-thi4}).
Диск с нулем ленточек имеет одну краевую окружность и гомеоморфен диску.
При приклеивании новой ленточки возможны два случая на рис.~\ref{strip}.

\smallskip
{\bf \ref{1-thino}.} (a) 1, 2 или 3.

(b) Диску с 3 дырками, тору с дыркой, ленте Мебиуса с дыркой или бутылке Клейна с дыркой.

\smallskip
{\bf \ref{1-eulno}.}
(a) Аналогично задаче \ref{1-nohom}.e.

(b) Аналогично формуле Эйлера для неперекрученных ленточек.
Диск с одной перекрученной ленточкой гомеоморфен диску $D$ с лентой Мебиуса.
Пусть теперь есть другие ленточки.
Можно считать, что они крепятся к краевой окружности диска $D$ по отрезкам, не пересекающим того отрезка,
по которому к $D$ прикреплена лента Мебиуса.

Если среди других ленточек есть перекрученная (относительно краевой окружности диска $D$, а не исходного диска!), то
заменим диск $D$ с лентой Мебиуса на диск.
Применим предположение индукции.
Используя (a), получим формулу Эйлера.

Если среди других ленточек нет перекрученной и нет перекрещивающихся, формула Эйлера очевидна.

Если среди других ленточек нет перекрученной и есть перекрещивающиеся, то применим задачу \ref{1-nohom}.с.
Заменим полученный диск с тремя лентами Мебиуса на диск с одной лентой Мебиуса (= с одной перекрученной ленточкой).
Применим предположение индукции.
Получим формулу Эйлера.

\begin{figure}[h]\centering
\includegraphics{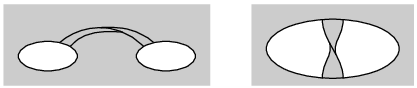}
\caption{Добавление перекрученной ленточки (изменить!)}
\label{strip1}
\end{figure}

{\it Указание к другому решению.}
Индукция по $n$ с использованием рис.~\ref{strip}, \ref{strip1} и задачи \ref{1-nohom}.

\newpage
\section{Двумерные многообразия}\label{02man}



\subsection{Определение двумерного симплициального комплекса}


Мы определим двумерные многообразия и полиэдры комбинаторно.
Это удобно как для теории, так и для хранения в памяти компьютера.

{\it Двумерным симплициальным комплексом} называется семейство двухэлементных и трехэлементных подмножеств конечного множества,
которое вместе с каждым трехэлементным множеством содержит все три его двухэлементные подмножества.
(Похожие объекты в комбинаторике называются гиперграфами.)
Будем сокращенно называть двумерный симплициальный комплекс просто {\it 2-комплексом}.

Элементы данного конечного множества называются {\it вершинами} 2-комплекса,
выделенные двухэлементные подмножества --- {\it ребрами} 2-комплекса, а
выделенные трехэлементные подмножества --- {\it гранями} 2-комплекса.

\begin{figure}[h]\centering
\includegraphics{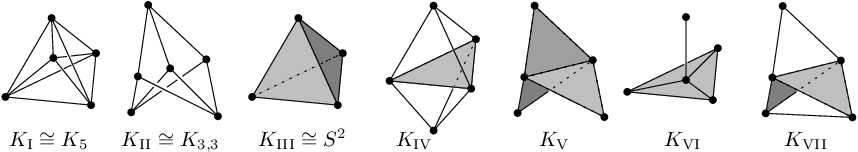}
\caption{Двумерные комплексы, не вложимые в плоскость}\label{hj}
\end{figure}


Примеры 2-комплексов приведены на рис. \ref{hj}.

{\it Кнопкой} называется 2-комплекс с вершинами $c,0,1,2,3$, в котором выделены трехэлементные подмножества $\{0,1,2\},\{0,1,3\}$ и $\{0,2,3\}$; выделены все их двухэлементные подмножества и $\{c,0\}$.
См. рис. \ref{hj}, $K_{VI}$.

{\it Книжкой с $n$ листами} называется 2-комплекс с вершинами $a,b,1,2,\dots,n$, в котором выделены трехэлементные подмножества $\{a,b,1\},\{a,b,2\},\dots,\{a,b,n\}$ и все их двухэлементные подмножества.
См. рис. \ref{hj}, $K_V$, для $n=3$, ср. рис. \ref{TI} слева.

\begin{figure}[h]\centering
\includegraphics{real-4.mps}
\caption{Полный 2-комплекс с 5 вершинами}\label{delta5}
\end{figure}

{\it Полным 2-комплексом с $n$ вершинами} (или двумерным остовом $n$-мерного симплекса) называется 2-комплекс с вершинами $0,1,2,\dots,n$, в котором все двухэлементные и трехэлементные подмножества выделены.
См. рис. \ref{hj}, $K_{III}$ для $n=4$ (этот 2-комплекс называется также сферой $S^2$) и рис. \ref{delta5} для $n=5$.


 \subsection{Важные неформальные замечания о 2-комплексах}

\begin{figure}[h]\centering
\includegraphics[width=\textwidth]{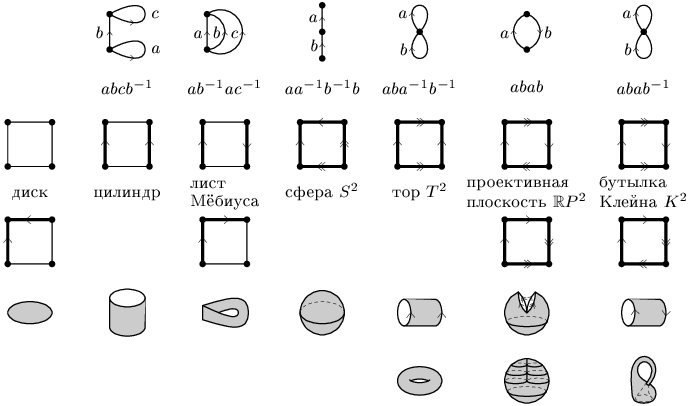}
\caption{Простейшие 2-комплексы (и соответствующие фигуры-тела); убрать 1-ю строку!}\label{exam}
\end{figure}

2-комплекс можно строить при помощи {\it `склейки'} сторон квадрата или даже многоугольника.
См. вторую и третью строки на рис. \ref{exam}.
{\it Шутовской колпак Зимана} получается склейкой сторон $\vec{AB}=\vec{AC}=\vec{BC}$ треугольника $ABC$.

Эта конструкция формализуется понятием {\it клеточного разбиения}, см.
далее;
здесь нам не понадобится эта формализация, достаточно интуитивного представления.

\begin{figure}[h]\centering
\includegraphics{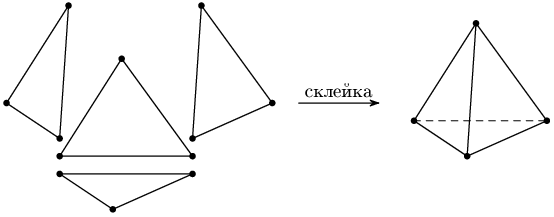}
\caption{Построение тела 2-комплекса}\label{razb}
\end{figure}

Как и по графу, по 2-комплексу естественно строится геометрическая фигура, называемая его {\it телом}.
Неформально, эта фигура получена склейкой треугольников и отрезков, соответствующих ребрам и граням 2-комплекса.
Склейка осуществляется не обязательно в трехмерном пространстве: либо в многомерном пространстве, либо даже абстрактно, независимо от объемлющего пространства.
Формально тело определяется в
\cite{Sk}, параграф `реализуемость двумерных комплексов'.

Например, на рис.~\ref{razb} изображено построение тела полного 2-комплекса с 4 вершинами.

Обратно, как и графы, 2-комплексы можно задавать {\it фигурами}, в т.ч. `гладкими' и самопересекающимися, т.е. их {\it телами}.
См. четвертую и пятую строки на рис. \ref{exam}.

2-комплекс --- комбинаторный объект.
Невозможно, например, взять точку на его грани.
Однако `взятие точки на грани тела 2-комплекса' формализуется `взятием новой вершины нового 2-комплекса, образовавшихся при подразделении этой грани', рис.~\ref{podra1} справа.
Как правило, мы не будем доводить наглядные рассуждения до такого формализма.

\subsection{Гомеоморфность двумерных комплексов}

\begin{figure}[h]\centering
\includegraphics{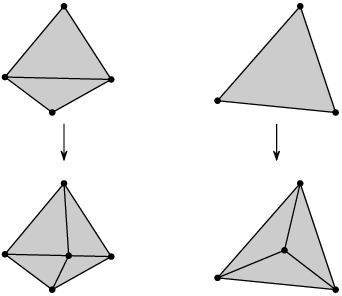}
\caption{Подразделения ребра и грани двумерного комплекса}\label{podra1}
\end{figure}

Операция {\it подразделения ребра} изображена на рис.~\ref{podra1} слева.

\begin{pr}\label{def-sdiv}
Операция {\it подразделения грани} на рис.~\ref{podra1} справа выражается через операцию подразделения ребра.
\end{pr}

Два 2-комплекса {\it гомеоморфны}, если от одного можно перейти к другому при помощи операций подразделения ребра и обратных к ним.

\begin{pr}\label{def-exam} 2-комплексы в каждой одной колонке на рис.~\ref{exam} гомеоморфны между собой, а из разных колонок --- нет.
(Указание: негомеоморфность можно доказывать по мере чтения следующих пунктов).
\end{pr}

{\it Двумерным полиэдром} называется класс гомеоморфности 2-комплекса.
Представители этих классов эквивалентности называются {\it триангуляциями} соответствующего двумерного полиэдра.
Именно полиэдры и многообразия
интересны топологу.
А комплексы и даже тела --- средства их изучения и хранения в компьютере.
Впрочем, комплексы и тела интересны дискретному геометру и комбинаторщику.

Два подмножества евклидова пространства называются {\it гомеоморфными}, если существуют взаимно-обратные непрерывные отображения между ними.

\smallskip
{\bf Теорема.}
{\it 2-комплексы гомеоморфны тогда и только тогда, когда их тела гомеоморфны.}

\smallskip
Заметим, что аналог этого результата для 5-комплексов неверен.

\subsection{Локально евклидовы двумерные комплексы}

2-комплекс называется {\it локально евклидовым}, если для любой его вершины $v$ все
грани, ее содержащие, образуют цепочку
$$\{v,a_1,a_2\},\ \{v,a_2,a_3\}\dots\{v,a_{n-1},a_n\}\quad\text{или}\quad
\{v,a_1,a_2\},\ \{v,a_2,a_3\}\dots\{v,a_{n-1},a_n\},\ \{v,a_n,a_1\}.$$
Если для всех $v$ имеет место второй случай, то локально евклидов 2-комплекс называется {\it замкнутым}.

Например, 2-комплексы на рис. \ref{exam} локально евклидовы.
Из них замкнуты только последние четыре.
2-комплекс, определенный диском с ленточками  (\S\ref{0dile}),
локально евклидов.

\begin{pr}\label{def-trian}
(a) Существует не локально евклидов 2-комплекс, к каждому ребру которого примыкает 2 грани.

(b) 2-комплекс, гомеоморфный локально евклидовому, сам локально евклидов.
\end{pr}

{\it Краем} (или границей) $\partial N$ локально евклидова 2-комплекса $N$ называется
объединение всех таких его ребер, которые содержатся только в одной грани.

\begin{pr}\label{4-bou}
(a) Край является несвязным объединением циклов.

(b) Количество краевых окружностей одинаково для гомеоморфных локально евклидовых 2-комплексов.

(c) 2-комплексы, `представляющие' цилиндр и ленту Мебиуса, не гомеоморфны.

(d) Диски с разными количествами дырок не гомеоморфны.
\end{pr}

`Заклеив' каждую краевую окружность диска с ленточками диском, получим замкнутый локально евклидов  2-комплекс.

{\it Двумерным многообразием} называется класс гомеоморфности локально евклидова 2-комплекса.
Этот 2-комплекс называется {\it триангуляцией} соответствующего 2-многообразия.
В математике вместо термина `локально евклидов 2-комплекс' используется термин `триангуляция 2-многообразия'.
Это неудобно для начинающего, поскольку при изучении 2-многообразий с кусочно-линейной точки зрения изначальным объектом являются 2-комплексы, и через них определяются 2-многообразия.
В этом параграфе мы используем термин `локально евклидов 2-комплекс', а в дальнейшем --- `триангуляция 2-многообразия' или даже `2-многообразие', если речь идет о свойстве 2-комплексов, инвариантном относительно гомеоморфности.

\subsection{Ориентируемость локально-евклидовых 2-комплексов}

{\it Ориентация} двумерного плоского треугольника --- упорядочение его вершин с точностью до четной перестановки.
Ориентацию удобно задавать окружностью со стрелкой, лежащей в треугольнике (или упорядоченной парой неколлинеарных векторов).


\begin{figure}[h]\centering
\includegraphics{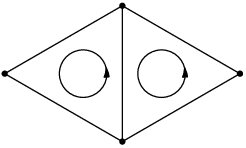}
\caption{Согласованные ориентации}\label{prot}
\end{figure}

{\it Ориентацией} локально евклидова 2-комплекса называется набор ориентаций на всех его гранях, которые {\it согласованы} вдоль каждого ребра, лежащего в двух гранях т.е. которые задают с двух сторон этого ребра {\it противоположные} направления (рис.~\ref{prot}).

Локально евклидов 2-комплекс называется {\it ориентируемым}, если у него существует ориентация.
2-Многообразие называется {\it ориентируемым}, если оно имеет ориентируемого представителя.

Понятие ориентируемости невозможно ввести для произвольных 2-комплексов (подумайте, почему).
Но можно ввести для 2-комплексов, каждое ребро которых содержится не более, чем в двух гранях.

\begin{pr}\label{4-oriex}
 (а) Гомеоморфные локально евклидовы 2-комплексы ориентируемы или нет одновременно.

(b) Сфера, тор и сфера с ручками ориентируемы, а лента Мебиуса, бутылка Кляйна и проективная плоскость (рис.~\ref{exam}) не ориентируемы.

(c) Тор не гомеоморфен бутылке Клейна.
\end{pr}

\begin{pr}\label{4-ori}
(a) Ориентируемость сохраняется при вырезании дырки.

(b) Диск с ленточками (\S\ref{0dile})
ориентируем тогда и только тогда, когда нет перекрученных ленточек.

(c)* Локально евклидов 2-комплекс ориентируем тогда и только тогда, когда он не содержит подкомплекса,
гомеоморфного ленте Мебиуса.
\end{pr}



\subsection{Эйлерова характеристика}\label{02maneul}

{\it Эйлеровой характеристикой} 2-комплекса $K$ с $V$ вершинами, $E$ ребрами и $F$ гранями называется число
$$\chi(K):=V-E+F.$$

\begin{pr}\label{4-eul}
(a) Найдите эйлеровы характеристики 2-комплексов с рис. \ref{exam}.

(b) Сформулируйте и докажите формулу включений-исключений для эйлеровой характеристики.

(c) Эйлеровы характеристики гомеоморфных 2-комплексов одинаковы.

(d) $\chi(P\times Q)=\chi(P)\times\chi(Q)$ для графов $P$ и $Q$.
\end{pr}

\begin{pr}\label{4-eulexa}
(a) Эйлерова характеристика сферы с $g$ ручками равна $2-2g$.

(b) Сферы с разными количествами ручек не гомеоморфны.

(c) Эйлерова характеристика сферы с $g$ ручками, $m$ пленками Мебиуса и $h$ дырками равна $2-2g-m-h$.

(d) Эйлерова характеристика равна $2-2g$, $1-2g$ и $-2g$ для сферы, проективной
плоскости или бутылки Кляйна с $g$ ручками, соответственно.

(e) Какие 2-комплексы из (d) гомеоморфны, а какие --- нет?

(f) {\bf Неравенство Эйлера.} {\it Пусть $G\subset K$ подкомплекс в 2-комплексе и к любому ребру в $K-G$ примыкает не более двух граней из $K-G$.
Обозначим через $F(K-G)$ количество связных кусков дополнения $F-G$.
Тогда $\chi(G)+F(K-G)\le\chi(K)$. }
\end{pr}

 Вложение подкомплекса (например, графа) в 2-комплекс называется {\it клеточным}, если
каждый связный кусок дополнения гомеоморфен открытому диску.
Приведем формальное комбинаторное определение.

Пусть в 2-комплексе задан подкомплекс (например, граф).
Рассмотрим {\it дополнение до  открытой регулярной окрестности подкомплекса}, т.е. объединение симплексов второго барицентрического подразбиения 2-комплекса, не пересекающих подкомплекс.
(Регулярная окрестность графа --- объединение шапочек и ленточек с рис. \ref{cap}.)
Вложение подкомплекса в 2-комплекс называется {\it клеточным}, если каждая связная компонента этого дополнения
гомеоморфна диску (или, эквивалентно, разбивается любой ломаной с концами на границе этой компоненты).
Например, вложение точки в сферу клеточно, а в тор --- нет.

\begin{pr}\label{4-euler}
 {\bf Формула Эйлера.}
{\it Пусть $G\subset K$ клеточно вложенный подкомплекс в 2-комплексе.
Обозначим через $F(K-G)$ количество связных кусков дополнения $F-G$.
Тогда $\chi(G)+F(K-G)=\chi(K)$. }
\end{pr}


\subsection{Классификация двумерных многообразий}

\begin{figure}[h]\centering
\includegraphics{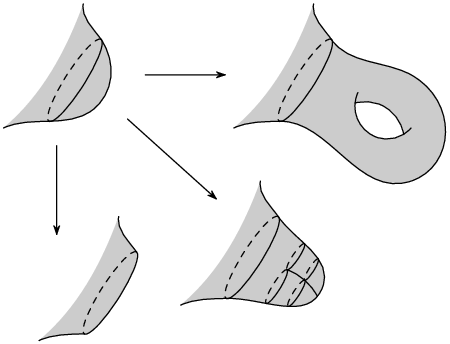}
\caption{Приклеивание ручки и пленки Мебиуса; вырезание дырки}
\label{hanhol}
\end{figure}

Новые локально евклидовы 2-комплексы получаются операциями {\it приклеивания ручки}
$N\to N\# T^2$, {\it приклеивания пленки Мебиуса} $N\to N\#\R P^2$ и {\it вырезания дырки}, рис.~\ref{hanhol}.

\begin{pr}\label{4-homnor}
(a) $\R P^2$ с дыркой гомеоморфно ленте Мебиуса.

(b) $\R P^2\#\R P^2\cong K^2$.
\qquad

(c) $\R P^2\# K^2\cong\R P^2\#T^2$.
\end{pr}





 {\bf Теорема классификации.}
{\it (1) Любой связный ориентируемый локально евклидов 2-комплекс гомеоморфен сфере с ручками и дырками.
Число $g$ ручек (оно называется родом) и число $h$ дырок однозначно определяются по 2-комплексу ($2-2g=\chi(N)$),
т.е. сферы с $g$ ручками и $h$ дырками не гомеоморфны для различных пар $(g,h)$.

(2) Любой связный неориентируемый локально евклидов 2-комплекс гомеоморфен сфере с пленками Мебиуса и дырками.
Число $m$ пленок Мебиуса (оно называется неориентируемым родом) и число $h$ дырок однозначно определяются по 2-комплексу ($2-m=\chi(N)$), т.е. сферы с $m$ пленками Мебиуса и $h$ дырками не гомеоморфны для различных пар $(m,h)$.}

\smallskip
{\it Набросок доказательства.}
Негомеоморфность доказывается путем подсчета количества связных компонент края и эйлеровой
характеристики, см. задачи \ref{4-bou}.b, \ref{4-eul}.c и \ref{4-eulexa}.c.

Докажем гомеоморфность.
Построим {\it шапочки, ленточки и заплатки}, как на рис.~\ref{cap}.
Приведем более аккуратное построение (которое читатель легко сможет совсем формализовать в соответствии
с определением 2-комплекса).

\begin{figure}[h]\centering
\includegraphics{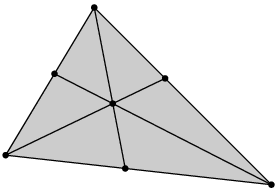}
\caption{Барицентрическое подразделение}\label{bar}
\end{figure}

Назовем {\it барицентрическим подразделением грани 2-комплекса} замену ее на шесть новых граней, полученных
проведением `медиан' в треугольнике, представляющем эту грань (рис.~\ref{bar}).
{\it Барицентрическим подразделением 2-комплекса} назовем результат барицентрического подразбиения всех его граней.
Так как барицентрическое подразделение можно осуществить с помощью конечного числа подразделений ребер, то этот результат гомеоморфен исходному 2-комплексу,

Пусть $T''$ --- 2-комплекс, полученный из данного 2-комплекса $T$ двукратным применением операции барицентрического подразделения. Назовем

$\bullet$ {\it шапочкой} объединение всех граней 2-комплекса $T''$, содержащих некоторую вершину 2-комплекса $T$.

$\bullet$ {\it ленточкой} объединение всех граней 2-комплекса $T''$, пересекающих некоторое ребро 2-комплекса $T$,
но не содержащих никакой вершины 2-комплекса $T$.

$\bullet$ {\it заплаткой} компоненту связности объединения оставшихся граней 2-комплекса $T''$.

Итак, мы `разбили' $T$ на шапочки, ленточки и заплатки.

Рассмотрим граф $G$, являющийся объединением всех ребер 2-комплекса $T$.
Пусть $G_0$ --- максимальное дерево графа $G$.
Через $T_0$ обозначим объединение шапочек и ленточек, пересекающихся с $G_0$.
Легко видеть, что $T_0$ гомеоморфно `двумерному диску'.
Если мы будем последовательно добавлять ленточки и заплатки к $T_0$,
то получим последовательность 2-комплексов
$T_0\subset T_1\subset T_2\subset\dots\subset T_p=T$.
Таким образом, $T$ получается из несвязного объединения дисков (возможно, из
одного диска) приклеиванием ленточек и заклеиванием дырок.
Поэтому гомеоморфность вытекает из задач \ref{1-eul}.a, \ref{1-eulno}.b и \ref{4-ori}.ab.
QED

\begin{pr}\label{4-cell}
(a) Если связный граф вложим в сферу с $g$ ручками, то он клеточно вложим в сферу с не более, чем $g$ ручками.

(b) Если связный граф вложим в диск с $m$ лентами Мебиуса, то он клеточно вложим в диск с не более, чем $m$ лентами Мебиуса.
\end{pr}

Было бы интересно решить эти задачи без использования теоремы классификации 2-многообразий.

\subsection{Ответы и указания к некоторым задачам}

 {\bf \ref{4-eulexa}.} (f)
Индукция по количеству ребер в $K-G$.
Если таких ребер нет, то $\chi(G)+F(K-G)=\chi(K)$.
Если есть, то обозначим через $G_+\subset K$ подкомплекс, полученный из $G$ добавлением одного из таких ребер.
Так как к любому ребру в $K-G$ примыкает не более двух граней комплекса $K-G$, то $F(K-G_+)\le F(K-G)+1$.
Поэтому и по предположению индукции $\chi(G)+F(K-G)\le \chi(G_+)+F(K-G_+)\le\chi(K)$.

{\it Другой способ (для локально-евклидовых 2-комплексов).}
Если, например, $K$ --- сфера с ручками, то аналогично доказательству теоремы классификации $K$ гомеоморфен сфере с $(2-V+E-F)/2$ ручками.

\comment

\smallskip
{\bf \ref{4-euler}.}
Докажите, что $V-E+F=V'-E'+F'$ для количеств $V',E',F'$ вершин, ребер и граней 2-комплекса $K$.

{\bf Дополнение: другое определение двумерных многообразий.}

{\it Ориентированным циклом} называется циклическая последовательность
ориентированных ребер, в которой начало каждого следующего ребра совпадает с
концом предыдущего
(ориентированные циклы, отвечающие тору и проективной плоскости на
рис.~\ref{exam}, различны).
{\it 2-Схемой} называется граф вместе с некоторым количеством выделенных в
нем ориентированных циклов (возможно, самопересекающимися и пересекающиеся друг с другом).
Выделенные циклы называются {\it гранями} 2-схемы.

Примеры схем, у которых ровно один выделенный цикл, имеющий длину 4, изображены
на рис.~\ref{exam} в верхней строке.

Из данной 2-схемы можно получать новые 2-схемы склеиванием ребер,
которое обозначается на рисунках проставлением на склеиваемых
ребрах одинаковых стрелочек (рис.~\ref{TI}).
Заметим, что эти стрелочки не имеют ничего общего с направлением ребер в
соответствующем направленном цикле (последнее обычно ясно из контекста и
поэтому не указывается нигде, кроме верхней строки на рис.~\ref{exam}).
Так из схемы 'квадрат' (первая схема во втором ряду на рис.~\ref{exam})
получаются схемы из верхней строки на рис.~\ref{exam}: под каждой схемой
первой строки нарисовано соответствующее отождествление ребер схемы 'квадрат'.
Напомним, что в верхней строке стрелки на ребрах обозначают их ориентацию в
выделенном цикле, а во второй строке совсем другое --- способ склеивания ребер.
В каждом рисунке второй строки ориентация в выделенном цикле (длины 4) задается
его обходом против часовой стрелки (при данном изображении на плоскости нашего
графа, являющегося циклом длины 4).

\begin{figure}[h]\centering
\includegraphics{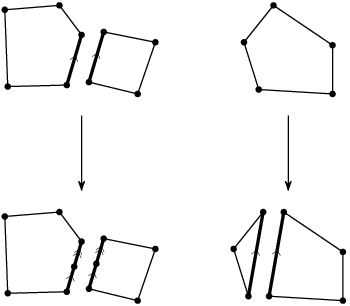}
\caption{Подразделение ребра и подразделение грани}\label{podrasc}
\end{figure}

Операции {\it подразделения ребра} и {\it подразделения грани
(проведения диагонали)} изображены на рис.~\ref{podrasc} (читатель
легко восстановит формальное комбинаторное определение этих
операций по рисункам). Две 2-схемы называются {\it гомеоморфными
(топологически эквивалентными)}, если от одной можно перейти к
другой при помощи этих операций и обратных к ним. Обозначение:
$T_1\cong T_2$.
{\it Компактным двумерным полиэдром} называется класс эквивалентности 2-схем
с точностью до гомеоморфности.
Представители этого класса эквивалентности называются {\it 2-схемами
соответствующего 2-полиэдра}.
Мы будем сокращенно называть компактные двумерные полиэдры 2-полиэдрами.

2-полиэдры и 2-схемы соотносятся так же, как 1-полиэдры (т.е. графы с
точностью до гомеоморфизма, или  графы в топологическом смысле) и графы
(с точностью до изоморфизма, или графы в комбинаторном смыслах).
Однако 2-полиэдры и 2-схемы уже важно различать.

По каждой 2-схеме естественно строится {\it фигура} (тело 2-схемы),
полученная из графа 2-схемы {\it склейкой} двумерных многоугольников,
соответствующих выделенным циклам 2-схемы (склейка осуществляется не в
трехмерном пространстве, а абстрактно).
Например, на рис.~\ref{razb} изображено построение тела схемы, состоящей
из графа $K_4$, в котором выделены все циклы длины 3.
Мы не даем строгого определения фигуры и склейки (т.е. топологического
пространства и операции склейки топологических пространств), поскольку
формально их не используем (все рассуждения можно перевести на язык классов
эквивалентности схем).
Однако важно, что на самом деле мы изучаем именно {\it фигуры} посредством
представления их 2-схемами.

Склеивание ребер 2-схемы порождает {\it склеивание ребер} соответствующего
2-полиэдра (рис.~\ref{TI},~\ref{exam}).
Примеры 2-полиэдров, полученных склейкой сторон квадрата, изображены в
третьем ряду на рис.~\ref{exam}.

\begin{figure}[h]\centering
\includegraphics{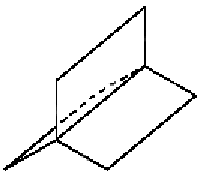}
\caption{Склеивание ребер}\label{TI}
\end{figure}

Вершина 2-схемы  называется {\it плоской}, если выходящие из этой вершины
ребра можно так занумеровать $e_1,e_2,\dots,e_n,e_{n+1}=e_1$, что выделенные
циклы проходят через эту вершину только вдоль ребер $e_i,e_{i+1}$ и только по
одному разу для каждого $i=1,2,\dots,n-1$ или $i=1,2,\dots,n$.
(Это равносильно тому, что линк вершины гомеоморфен окружности или отрезку.)
Компактный двумерный полиэдр называется {\it компактным двумерным
многообразием} (сокращенно: 2-многообразием) или поверхностью, если для
некоторой (или, эквивалентно, для любой
\footnote{Доказательство этой эквивалентности является (несложным) упражнением}
) его схемы каждая вершина этой схемы плоская.
Это свойство неформально означает, что тело схемы локально устроено как
плоскость или полуплоскость.
В этом пункте, а также в тех, где участвуют только 2-многообразия (т.е. не
участвуют $n$-многообразия с $n\ge3$) мы будем называть их просто
многообразиями.

{\it Краем} (или границей) схемы многообразия называется объединение таких
ребер схемы, которые содержатся только в одной грани и причем однократно.
Многообразие называется {\it замкнутым}, если край некоторой
(или, эквивалентно, любой) его схемы пуст.

\endcomment

\newpage
\section{Гомологии двумерных многообразий}\label{02hom}

\hfill{\it I should say it meant something simple}

\hfill{\it and obvious, but then I am no philosopher!}

\hfill{\it I.\?Murdoch, The Sea, the Sea.
\footnote{\it Должен бы сказать, что это означает нечто простое и очевидное,
однако я не философ! А. Мэрдок, Море, море, пер автора.}
}

 \subsection{Критерий ориентируемости и определение цикла}\label{02homcy}

Определение 2-многообразия и его ориентируемости приведено в \S\ref{02man}.
Существует красивый и простой критерий ориентируемости 2-многообразия: `не содержит ленты Мебиуса'.
Существует простой алгоритм распознавания ориентируемости связного 2-многообразия: сначала ориентируем
произвольно одну грань, затем на каждом шаге будем ориентировать грань,
соседнюю с некоторой уже ориентированной, пока не ориентируем все грани или не
получим несогласованности ориентаций вдоль некоторого ребра.

В этом параграфе мы приведем алгебраический критерий ориентируемости,
который по сути  является лишь переформулировкой определения ориентируемости на алгебраический язык.
Но он важен не сам по себе, а как иллюстрация теории препятствий.
(Например,
{\it классификации утолщений}, см. \cite{Sk}, конец параграфа `утолщения графов'.)

\smallskip
{\bf Теорема ориентируемости.}
{\it Замкнутое 2-многообразие $N$ ориентируемо тогда и только тогда, когда его}
первый класс Штифеля-Уитни $w_1(N)\in H_1(N)$ {\it нулевой.}

\smallskip
Группа $H_1(N)$ и класс $w_1(N)$ определены позже.
Они естественно возникают и строго определяются в процессе {\it придумывания} теоремы ориентируемости,
к которому мы сейчас перейдем.

В этом пункте слово `группа' можно рассматривать как синоним слова `множество' (кроме задач \ref{origraph} и \ref{orisurf}).
И теорема ориентируемости, и приводимые построения останутся интересными.

{\it Клеточным разбиением} 2-комплекса называется клеточное вложение (см. конец \S\ref{02man}) некоторого графа этот 2-комплекс.
Возможно, начинающему будет более понятен термин {\it разбиение на многоугольники}.
Например, 2-комплекс является клеточным разбиением себя.
См. примеры `склейки из многоугольников' в начале \S\ref{02man}.
Многие построения (в том числе приводимые здесь) удобно проделывать не для 2-комплексов а для клеточных разбиений.
Ибо у `интересных' 2-комплексов `много' граней, но можно найти их `экономные' клеточные разбиения.

\smallskip
{\it Определение препятствующего цикла.}
Возьмем произвольное клеточное разбиение $T$ данного 2-многообразия $N$.
Возьмем набор $o$ ориентаций на гранях разбиения $T$.
Покрасим ребро разбиения $T$ в красный цвет, если ориентации граней, примыкающих к ребру с двух сторон, {\it рассогласованы вдоль этого ребра},
т.е. задают на этом ребре {\it одинаковые} направления.

\begin{figure}[h]\centering
\includegraphics{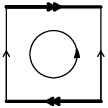}
\caption{Набор $o$ ориентаций и препятствующий цикл $\omega(o)$}\label{cycle}
\end{figure}

Например, на рис.~\ref{cycle} бутылка Кляйна представлена в виде склейки сторон квадрата, т.е. разбита на один многоугольник.
Грани, примыкающие к горизонтальному ребру  с двух сторон, совпадают.
Но их (точнее, ее) ориентации рассогласованы вдоль этого ребра.
Кроме того, ориентация единственной грани согласована с собой вдоль вертикального ребра.

Итак, если разбиение не является 2-комплексом (`триангуляцией'), то даже если две грани, примыкающие к ребру, совпадают, их ориентации могут быть рассогласованы вдоль этого ребра.
Кроме того, для одной пары граней (совпадающих или нет) их ориентации могут быть согласованы вдоль одного ребра и расогласованы вдоль другого.

\begin{pr}\label{oriobcy}
Для каждого ребра клеточного разбиения  проективной плоскости на один многоугольник (т.е. ее представления в виде склейки сторон квадрата, рис.~\ref{exam}) выясните, согласована ли вдоль этого ребра с собой ориентация единственной грани.
\end{pr}

Объединение красных ребер называется {\it препятствующим циклом} $\omega(o)$.

На рис.~\ref{cycle}  препятствующий цикл состоит из горизонтального ребра (жирные линии).

\begin{pr}\label{oricycle}
(a) Изобразите препятствующий цикл для клеточного разбиения проективной плоскости на один многоугольник, рис.~\ref{exam}.

(b) Набор $o$ ориентаций граней определяет ориентацию клеточного разбиения тогда и только тогда, когда $\omega(o)=\emptyset$.

(c) Из каждой вершины выходит четное число красных ребер.

(d) 2-Многообразие с краем, получающееся вырезанием `малой окрестности' препятствующего цикла $\omega(o)$ из $N$, ориентируемо.
\end{pr}

{\it Гомологическим циклом} в графе называется набор его ребер, для которого
из каждой вершины выходит четное число ребер набора.
 Слово `гомологический' будем опускать, а циклы в смысле теории графов будем называть `замнутыми кривыми'.
Например, объединение ребер клеточного разбиения бутылки Кляйна на один многоугольник (рис.~\ref{cycle}) является `восьмеркой', поэтому в этом графе четыре цикла.

\begin{pr}\label{orinumber} (a) Сумма по модулю 2 циклов есть цикл.

(b) Найдите число циклов в связном графе с $V$ вершинами и $E$ ребрами.
\end{pr}

\begin{pr}\label{origraph}
(Эта задача не используется в дальнейшем.)
Группа $H_1(G)$ всех циклов в графе $G$ с операцией суммы по модулю 2
называется {\it  группой гомологий графа $G$} (одномерной с коэффициентами  $\Z_2$).

(a) Одномерные группы гомологий гомеоморфных графов изоморфны.

(b) $H_1(G)\cong\Z_2^{E-V+1}$ для связного графа $G$ с $V$ вершинами и $E$ ребрами.

(c) Несамопересекающиеся циклы порождают $H_1(G)$.
\end{pr}

 \subsection{Ориентируемость: гомологичность циклов}\label{02homhocy}

{\it Изменение препятствующего цикла и определение границы грани.}
Если $\omega(o)\ne\emptyset$, то $o$ не определяет ориентации клеточного разбиения $T$.
Но еще не все потеряно: можно попытаться изменить $o$ так, чтобы препятствующий цикл стал пустым.
Для этого выясним, как $\omega(o)$ зависит от $o$.
Ответ удобно сформулировать с использованием суммы по модулю 2 (т.е. симметрической разности) наборов ребер в произвольном графе.

Если $T$ --- 2-комплекс, то назовем {\it гомологической границей $\partial a$ грани $a$} набор ребер геометрической границы этой грани.

\begin{figure}[h]\centering
\includegraphics{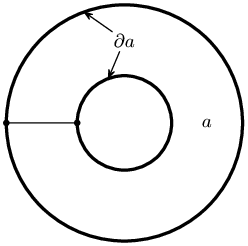}
\caption{Гомологическая (алгебраическая) граница сложной грани}\label{part}
\end{figure}

Для произвольного клеточного разбиения $T$ определение границы более сложное.
Назовем {\it гомологической границей $\partial a$ грани $a$} набор всех тех ребер геометрической границы этой грани,
к которым грань примыкает с нечетного числа сторон (рис.~\ref{part}).

Как и для циклов, слово `гомологическая' будем опускать.
Для клеточного разбиения бутылки Кляйна на один многоугольник (рис.~\ref{cycle}) граница единственной грани пуста.

\begin{pr}\label{oribou}
(a) Что будет границей единственной грани клеточного разбиения проективной плоскости на один многоугольник, рис.~\ref{exam}?

(b) Граница грани является циклом.

(c) При изменении ориентации одной грани $a$ цикл $\omega(o)$ изменяется на сумму по модулю 2
с границей этой грани: $\omega(o')-\omega(o)=\partial a$.

(d) При изменении ориентации {\it нескольких граней} $a_1,\dots,a_k$ цикл
$\omega(o)$ изменяется на сумму по модулю 2 границ этих граней:
$$
\omega(o')-\omega(o)=\partial a_1+\dots+\partial a_k.
$$
\end{pr}

Назовем {\it границей} сумму границ нескольких граней.
Назовем циклы {\it гомологичными} (или сравнимыми по модулю границ), если их разность есть граница.

\begin{pr}\label{orichange}
(a) При изменении набора $o$ ориентаций препятствующий цикл $\omega(o)$ заменяется на гомологичный цикл.

(b) Если $\omega(o)$ является границей, то можно изменить $o$ на $o'$ так, чтобы получилось $\omega(o')=0$.
\end{pr}

Определение цикла и гомологичности осмыслено для клеточного разбиения произвольного 2-комплекса
(не обязательно замкнутого локально евклидова).

\begin{pr}\label{orihom}
(a) Любые два цикла для клеточного разбиения сферы на один многоугольник, рис.~\ref{exam}, гомологичны.


(b) Краевые окружности на торе с двумя дырками гомологичны (для любого разбиения).

(с) Краевая окружность ленты Мебиуса гомологична пустому циклу (для любого клеточного разбиения).

(d) Любые два цикла для клеточного разбиения шутовского колпака Зимана на один многоугольник, рис.~\ref{exam}, гомологичны.
\end{pr}

\begin{pr}\label{oritor} Для  клеточного разбиения тора на один многоугольник, рис.~\ref{exam}

(a) цикл `меридиан' не гомологичен пустому.

(b) только одинаковые циклы гомологичны.

(c) имеется ровно 4 попарно не гомологичных цикла.
\end{pr}

\begin{pr}\label{orirp2} Для клеточное разбиения проективной плоскости на один многоугольник, рис.~\ref{exam}

(a) цикл `$\R P^1$' не гомологичен пустому.

(b) имеется ровно 2 попарно не гомологичных цикла.
\end{pr}

\begin{pr}\label{orikl}
(a) Для клеточного разбиения бутылки Кляйна на один многоугольник (рис.~\ref{cycle}) имеется ровно 4 попарно не гомологичных цикла.

(b) Гомологичность является отношением эквивалентности на множестве циклов.

(c) Любой цикл в 2-комплексе гомологичен некоторому циклу в произвольном клеточном разбиении этого 2-комплекса
 (т.е. в произвольном графе, клеточно вложенном в 2-комплекс).

(d) Если два цикла в клеточном разбиении 2-комплекса гомологичны в 2-комплексе, то они  гомологичны и в этом
клеточном разбиении.

(e) Любой цикл становится гомологичным `несамопересекающемуся циклу' после некоторого измельчения клеточного разбиения.
\end{pr}

В рассуждениях  классами гомологичности циклов удобно сначала работать с представляющими их циклами, а потом
доказывать независимость от выбора представляющих циклов.

 \subsection{Ориентируемость: гомологии и первый класс Штифеля-Уитни}\label{02homhom}

 {\it Определение группы $H_1(T)$ клеточного разбиения $T$ 2-комплекса.}
Набор ребер клеточного разбиения $T$ называется {\it циклом по модулю 2}, если из каждой вершины выходит четное число ребер набора.
{\it Границей $\partial a$ грани $a$ по модулю 2\ } называется набор всех тех ребер геометрической границы этой грани, к которым грань примыкает с нечетного числа сторон (рис.~\ref{part}).
Назовем {\it границей по модулю 2\ } сумму границ по модулю 2 нескольких граней.
Два цикла по модулю 2 называются {\it гомологичными по модулю 2}, если их разность есть граница по модулю 2.
{\it Одномерной группой гомологий $H_1(T)$ с коэффициентами $\Z_2$ клеточного разбиения $T$}\ называется
группа циклов по модулю 2 с точностью до гомологичности по модулю 2.
Слова `по модулю 2' далее опускаются.

\smallskip
Например, $|H_1(T)|=0,0,4,4,2$ для клеточных разбиений на один многоугольник сферы, шутовского колпака Зимана, тора, бутылки Клейна и проективной плоскости.

Важно, что эта группа, возникающая при решении конкретной задачи об ориентируемости
(и других задач!) коротко определяется независимо от этой задачи.
Для определения не нужны ни замкнутость, ни локальная евклидовость.


\begin{pr}\label{orih1}
(a) Если $T$ --- клеточное разбиение 2-комплекса $T'$, то $|H_1(T)|=|H_1(T')|$.

(b) Одномерные группы гомологий гомеоморфных 2-комплексов находятся во взаимно-однозначном соответствии друг с другом.
\end{pr}

{\it Одномерной группой гомологий $H_1(N)$ 2-многообразия (или 2-комплекса) $N$ с коэффициентами $\Z_2$} называется группа $H_1(T)$ для любого клеточного разбиения любого 2-комплекса, представляющего $N$.
Это определение корректно ввиду задачи \ref{orih1}.ab (и \ref{orisurf}.cd).

{\it Первым классом Штифеля-Уитни} клеточного разбиения $T$ называется класс гомологичности препятствующего цикла:
$$
w_1(T):=[\omega(o)]\in H_1(T).
$$
Это определение корректно ввиду утверждения \ref{orihom}.a.

{\it Первым классом Штифеля-Уитни} 2-многообразия $N$ называется первый класс Штифеля-Уитни любого клеточного разбиения $T$ любого 2-комплекса, представляющего $N$: $w_1(N):=w_1(T)$.
Это определение корректно в следующем смысле.

\begin{pr}\label{orimap}
При `естественных' биекциях, построенных Вами при решении задач \ref{orih1}.ab групп гомологий  клеточных разбиений гомеоморфных 2-комплексов первый класс Штифеля-Уитни одного клеточного разбиения переходит в первый класс Штифеля-Уитни другого.
\end{pr}

{\it Набросок доказательства теоремы ориентируемости.}
Ясно, что условие $w_1(N)=0$ необходимо для ориентируемости.
Обратно, пусть $w_1(N)=0$.
Возьмем 2-комплекс $T$, представляющий $N$.
Тогда $w_1(T)=0$.
Значит, существуют такой набор $o$ ориентаций граней 2-комплекса $T$, что $\omega(o)$ есть граница.
Поэтому можно изменить $о$ на $o'$ так, чтобы $\omega(o')=0$.
Тогда $T$ ориентируем.
Следовательно, $N$ ориентируемо.
QED

 \subsection{Вычисления гомологий и обобщение}\label{02homcal}

В следующих задачах `найти группу гомологий' означает `найти количество элементов в ней'.
Если Вы знакомы с понятием группы, то --- `найти известную группу, изоморфную данной'.
Группа гомологий конечно порожденная абелева, поэтому данную задачу легко формализовать.

\begin{pr}\label{orihocl} Найдите группу гомологий

(a)  сферы с $g$ ручками. \quad
(b) сферы с $g$ ручками и $h$ дырками.

(с) сферы с $m$ пленками Мебиуса. \quad
(d) сферы с $m$ пленками Мебиуса и $h$ дырками.
\end{pr}

\begin{pr}\label{orihosum}
(a) $H_1(M\#N)\cong H_1(M)\oplus H_1(N)$ для замкнутых 2-многообразий $M$ и $N$.

(b) Верна ли эта формула для незамкнутых 2-многообразий $M$ и $N$?
\end{pr}

\begin{pr}\label{orisurf}
Пусть $T$ --- клеточное разбиение замкнутого 2-многообразия
(например, связный локально евклидов замкнутый 2-комплекс).

(a) Сумма по модулю 2 границ всех граней разбиения $T$ пуста.

(b) Сумма по модулю 2 границ любого собственного подмножества множества всех граней разбиения $T$ непуста.

(c) $|H_1(T)|=2^{2-\chi(T)}(=2^{2-V+E-F})$.
\end{pr}

\begin{pr}\label{orisurft}
(Эта задача не используется в данном параграфе.)
Пусть $T$ --- клеточное разбиение замкнутого 2-многообразия
(например, связный локально евклидов замкнутый 2-комплекс).

(a) В множестве $H_1(T)$ корректно определена операция суммы формулой $[\alpha]+[\beta]=[\alpha+\beta]$.

(b) Множество $H_1(T)$ с этой операцией является группой.

(c) Если $T$ --- клеточное разбиение 2-комплекса $T'$, то $H_1(T)\cong H_1(T')$.

(d) Одномерные группы гомологий гомеоморфных 2-комплексов изоморфны.

(e) $H_1(T)\cong\Z_2^{2-\chi(T)}$.
\end{pr}

\begin{pr}\label{oribd} Придумайте определение препятствующего цикла для клеточного разбиения локально евклидова 2-комплекса, не являющегося замкнутым.
Нарисуйте препятствующий цикл для разбиения ленты Мебиуса на один многоугольник.
\end{pr}

Пусть $T$ --- клеточное разбиение локально евклидова 2-комплекса (не обязательно замкнутого).
{\it Циклом по модулю края} в называется набор ребер в $T$, из каждой внутренней вершины которого выходит четное число ребер.
Границы определяются дословно так же, как и для выше.
Два цикла по модулю края называются {\it гомологичными по модулю края}, если их разность есть сумма границы и некоторого набора ребер края.
Группа $H_1(T,\partial)$ гомологий по модулю края и класс $w_1(T)\in H_1(T,\partial)$ определяются аналогично предыдущему.

\begin{pr}\label{orihobd} Найдите группу гомологий по модулю края

(a) сферы с $g$ ручками и $h$ дырками. \quad

(b) сферы с $m$ пленками Мебиуса и $h$ дырками.

(c) {\bf Теорема ориентируемости.}
{\it Локально евклидов 2-комплекс $T$ ориентируем тогда и только тогда, когда его
первый класс Штифеля-Уитни $w_1(T)\in H_1(T,\partial)$ нулевой.}
\end{pr}

 \subsection{Форма пересечений}\label{02homint}

Форма пересечений --- один из важнейших инструментов и объектов исследования для топологии и ее приложений.
См.~\cite{DZ93}.
Она естественно появляется, например, при решении задач \ref{intplan}, \ref{intreg}.d и \ref{intrans}.
(Если Вам непонятны используемые в них слова или не получается решить, то пропустите их и читайте дальше.)
Задачи \ref{intreg}.abc можно решить и без формы пересечений.
См. также формулу Мохара --- задачи \ref{1-eul}.b и \ref{1-eulno}.c.

\begin{figure}[h]\centering
\includegraphics{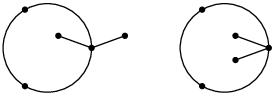}
\caption{Трансверсальное пересечение и нетрансверсальное пересечение (изменить рисунок!)}
\label{trans}
\end{figure}

 Точка $x$ пересечения двух несамопересекающихся ломаных
на 2-многообразии называется {\it трансверсальной}, если любая достаточно малая окружность $S_x$ с центром в $x$ пересекает ломаные по парам точек, {\it чередующимся} вдоль окружности (т.е. если обозначить через $A_1,B_1$ точки пересечения первой ломаной с $S_x$ и через $A_2,B_2$ точки пересечения второй ломаной с $S_x$, то эти точки пересечения расположены на окружности в порядке $A_1A_2B_1B_2$).
Иными словами, если два звена одной ломаной, выходящие из точки пересечения, находятся по разные стороны от другой ломаной в малой окрестности точки пересечения, рис. \ref{trans}.

\begin{pr}\label{intrans}
Пусть $N$ --- 2-многообразие и $a,b$ --- замкнутые кривые на нем.
Будем считать, что $a$ и $b$

$\bullet$ являются подграфами некоторого 2-комплекса, представляющего $N$;

$\bullet$ {\it общего положения}, т.е. трансверсально (рис. \ref{trans}) пересекаются в конечном числе точек,
никакая из которых не является ни точкой самопересечения $a$, ни точкой самопересечения $b$.

(a) $|a\cap b|\mod2$ не меняется при замене  $a$ и $b$ на гомологичные им кривые, удовлетворяющие тем же условиям
(подграфы, соответствующие кривым, являются гомологическими циклами; это определяет `гомологичность' кривых).

(b) Если $b$ `близко' к $a$, то $|a\cap b|$ нечетно, если при обходе вдоль $a$ `меняется ориентация' 2-многообразия, и четно, если не меняется.
\end{pr}

Построение из задачи \ref{intrans} легко доработать до определения формы пересечений, использующего понятие трансверсальности.
Но мы приведем другое определение.
В нем вместо понятия трансверсальности используется более простое понятие двойственного клеточного разбиения, см.
ниже.

Возьмем некоторое клеточное разбиение $U$ 2-многообразия $N$ (точнее, некоторого представляющего его 2-комплекса).

\begin{figure}[h]\centering
\includegraphics{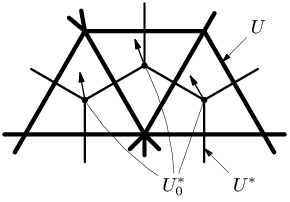}
\caption{Двойственное клеточное разбиение. Убрать стрелочки, сделать дв-е пунктирным}
\label{dual}
\end{figure}

\smallskip
{\it Определение двойственного клеточного разбиения (рис.~\ref{dual}).}
Выберем в каждой грани клеточного разбиения $U$ точку и обозначим полученное множество точек через $U_0^*$.
Для каждого ребра $a$ разбиения $U$ соединим ребром $a^*$ точки множества $U_0^*$, соответствующие соседним вдоль ребра $a$ граням.
(Могут появиться петли.)
Ребро $a^*$ называется {\it двойственным} к $a$.
Полученный граф обозначим через $U_1^*$.
Этот граф клеточно вложен в 2-многообразие $N$ (точнее, в представляющий его 2-комплекс).
Полученное клеточное разбиение $U^*$ называется {\it двойственным} к $U$.

\smallskip
{\it Определение формы пересечений.}
Возьмем двойственное клеточное разбиение $U^*$.
Для наборов ребер $X$ разбиения $U$ и $Y$ разбиения $U^*$ положим $X\cap Y$ равным количеству их точек пересечения по модулю 2.
Для цикла $X$ разбиения $U$ и цикла $Y$ разбиения $U^*$ определим $[X]\cap[Y]:=X\cap Y$.
Это определение корректно задает {\it форму пересечений} $\cap:H_1(U)\times H_1(U^*)\to\Z_2$ (см. задачу \ref{intwell}.b ниже).
Так как имеются канонические изоморфизмы $H_1(U)\cong H_1(U^*)\cong H_1(N)$, ввиду аналога задач \ref{orisurft}.cd получается {\it форма пересечений} $\cap:H_1(N)\times H_1(N)\to\Z_2$.


\begin{pr}\label{intwell}
(a) Пересечение цикла и границы равно нулю.

(b) Определение формы пересечений корректно.

(с) Форма пересечений симметрична и билинейна: $\alpha\cap\beta=\beta\cap\alpha$ и
$(\alpha+\beta)\cap\gamma=\alpha\cap\gamma+\beta\cap\gamma$.

(d) Найдите форму пересечений сферы с $g$ ручками.
(Т.е. найдите матрицу этой формы в некотором базисе группы гомологий.)

(e) Найдите форму пересечений сферы с $m$ пленками Мебиуса.

(f) Если $f:H_1(T)\to H_1(T_1)$  и $f^*:H_1(T^*)\to H_1(T_1^*)$ --- изоморфизмы из задачи \ref{orih1}, то
$\alpha\cap\beta=f(\alpha)\cap f^*(\beta)$.
\end{pr}

\begin{pr}\label{intstwh} Пусть $N$ --- замкнутое 2-многообразие.
Определение первого класса Штифеля-Уитни $w_1(N)\in H_1(N)$ приведено в предыдущем пункте.

(a) $w_1(N)\cap a=a\cap a$ для любого $a\in H_1(N)$.
\qquad
(b) $w_1(N)\cap w_1(N)\equiv\chi(N)\mod2$.
\end{pr}

\begin{pr}\label{intpoin}
 {\it Двойственность Пуанкаре.}
(a) Форма пересечений любого замкнутого 2-многообразия $N$ невырождена,
т.е. для любого $\alpha\in H_1(N)-\{0\}$ существует такое $\beta\in H_1(N)$, что $\alpha\cap\beta=1$.

(b) Для любого 2-многообразия $N$ умножение $H_1(N,\partial)\times H_1(N)\to\Z_2$ из \ref{intbd}.f невырождено.
\end{pr}

\begin{pr}\label{intbd}
  Пусть $N$ --- 2-многообразие с непустым краем.

(a) Определите аналогично умножение $H_1(N)\times H_1(N)\to\Z_2$.

(b) Сформулируйте и докажите аналоги задач \ref{intwell} для умножения из (a).

(c) Умножение из (a) может быть вырождено.

(d) Найдите умножение из (a) и его ранг для   сферы с $g$ ручками и $h$ дырками.

(e) Найдите умножение из (a) и его ранг для сферы с $m$ пленками Мебиуса и $h$ дырками.

(f) Определите аналогично умножение $H_1(N,\partial)\times H_1(N)\to\Z_2$.
(Не забудьте доказать корректность определения!)

(g) Сформулируйте и докажите аналоги задач \ref{intwell}, \ref{intstwh}.a и для умножения из (f).

(h) (загадка)
Почему не получается аналогично определить умножение $H_1(N,\partial)\times H_1(N,\partial)\to\Z_2$?
\end{pr}



\smallskip
\subsection{Ответы, указания и решения к некоторым задачам}

{\bf \ref{orinumber}.} (b) Выделим в графе максимальное дерево, тогда вне этого дерева $E-V+1$ ребер.

Ответ: $2^{E-V+1}$.

\smallskip
{\bf \ref{origraph}.} (a) Это достаточно доказать для случая, когда один граф получается из другого подразделением ребра.

(b) Осторожно, это не следует из $|H_1(G)|=2^{E-V+1}$.
Но следует из этого и того, что сумма любого элемента с собой равна нулю.

\smallskip
{\bf \ref{orichange}.} (b) Пусть $\omega(o)\hm=\partial a_1+\dots+\partial a_s$.
Возьмем набор $o'$ ориентаций граней, отличающийся от $o$ в точности на гранях $a_1,\dots,a_s$.
Тогда по задаче \ref{oribou}.b получим $\omega(o')=0$.

\smallskip
{\bf \ref{orihom}.} (a) Имеется только пустой цикл.

(b,c) Сумма границ всех граней равна разности циклов.

(d) Граница единственной грани равна единственному циклу.

\smallskip
{\bf \ref{oritor}, \ref{orirp2}, \ref{orikl}.a.} Граница единственной грани пуста.

\smallskip
{\bf \ref{orikl}.}
(с) Следует из того, что для любого цикла в 2-комплексе грани 2-комплекса можно так раскрасить в 2 цвета, чтобы

$\bullet$ при переходе через ребро цикла, не лежащее на ребре разбиения, цвет менялся, и

$\bullet$ при переходе через ребро, не лежащее ни на цикле, ни на ребре разбиения, цвет сохранялся.

(d) Граница в 2-комплексе между циклами разбиения есть сумма границ некоторых граней разбиения.

\smallskip
{\bf \ref{orih1}.} (a) Следует из \ref{orikl}.cd

(b) Это достаточно доказать для случая, когда один 2-комплекс получается из другого подразделением ребра.

\smallskip
{\bf \ref{orihocl}.} Ответы: \quad (a) $\Z_2^{2g}$; \quad (b) $\Z_2^{2g+h-1}$ при $h>0$; \quad (c) $\Z_2^m$; \quad (d) $\Z_2^{m+h-1}$ при $h>0$.

(a) Постройте клеточное разбиение на один $4g$-угольник.

(с) Постройте клеточное разбиение на один $2m$-угольник.

(b) {\it Первый способ.} Пусть $h>0$ и $N_h$ --- сфера с $g$ ручками и $h$ дырками.
Определите гомоморфизм $H_1(N_h)\to H_1(N_{h-1})$ `индуцированный' включением $N_h\to N_{h-1}$.
Докажите, что он сюръективен.
Докажите, что его ядро нулевое при $h=1$ и изоморфно $\Z_2$ при $h>1$.

(b) {\it Второй способ.}
Если тело 2-комплекса $X$ есть {\it деформационный ретракт} тела 2-комплекса $Y$ (ср. с задачей \ref{vfann})
или $X$ {\it сдавливается} на $Y$, то $H_1(X)\cong H_1(Y)$.

(d) Аналогично (b).

\smallskip
{\bf \ref{orihobd}.} Ответы: \quad (a) $\Z_2^{2g}$; \quad (b) $\Z_2^m$.

Докажите, что группа гомологий 2-многообразия по модулю края изоморфна группе гомологий замкнутого 2-многообразия, полученного из данного 2-многообразия заклейкой каждой краевой окружности диском.

\smallskip
{\bf \ref{intreg}.} (a) Используйте формулу Эйлера.

(d) Обозначим через $R_1$ и $R_2$ регулярные окрестности образов графа $G$ при двух гомотопных вложениях.
Далее $k=1,2$.
Определите $H_1(R_k)$ и постройте изоморфизм $i_k:H_1(G)\to H_1(R_k)$.
Определите форму пересечений $\cap_k:H_1(R_k)\times H_1(R_k)\to\Z_2$ (см. задачу \ref{intbd}.a).
Докажите, что изоморфизм $i_1i_2^{-1}$ сохраняет форму пересечений.
Используйте, что форма пересечений однозначно определяет 2-многообразие с непустым краем (см. задачу \ref{intbd}.d).

\smallskip
{\bf \ref{intwell}.}
(d) Постройте клеточное разбиение на один $4g$-угольник.
В соответствующем базисе матрица стандартная симплектическая, т.е. составлена из $g$ диагональных блоков
$\left(\begin{matrix} 0&1 \\ 1&0 \end{matrix}\right)$.

(e) Постройте клеточное разбиение на один $2m$-угольник.
В соответствующем базисе матрица единичная.

\smallskip
{\bf \ref{intstwh}.} (b) Следует из теоремы классификации 2-многообразий.
Алгебраическое доказательство, не использующее классификации 2-многообразий, можно получить из рассмотрения
формы пересечений и (a).
Или из рассуждения в начале пункта `наборы нормальных полей' и погружаемости любого 2-многообразия в $\R^3$.

\smallskip
{\bf \ref{intpoin}.} (a) Ввиду задачи \ref{orikl}.d можно считать, что $\alpha$ представлен замкнутой несамопересекающейся ломаной на $N$.
Рассмотрите два случая: ломаная разбивает многообразие или нет.
Другой способ приведен в \S\ref{0dual}.

\smallskip
{\bf \ref{intbd}.} (d,e) Используйте базис, фактически построенный при вычилении группы гомологий.

(d) Ответ: в некотором базисе матрица имеет $g$ диагональных блоков
$\left(\begin{matrix} 0&1 \\ 1&0 \end{matrix}\right)$, остальные элементы нулевые.
Ранг равен $2g$.

(e) Ответ: в некотором базисе матрица имеет $m$ единиц на диагонали, остальные элементы нулевые.
Ранг равен $m$.


\newpage
\section{Инволюции}\label{0inv}

\subsection{Примеры инволюций}

Читатель наверняка знает, насколько полезны центральные, осевые и зеркальные симметрии в элементарной геометрии.
Настолько же полезны их обобщения --- инволюции --- в `высшей' математике.
В свою очередь, понятие инволюции является частным случаем более общего понятия действия группы на множестве.
Мы не затрагиваем здесь это понятие, поскольку идеи и методы, который мы хотим проиллюстрировать, достаточно
хорошо видны на примере инволюций.

{\it Инволюцией} на подмножестве $N\subset\R^m$ называется непрерывное
отображение $t:N\hm\to N$, для которого $t(t(x))=x$ для любого $x\in N$.

Инволюция $t$ назывется {\it свободной}, если она не имеет неподвижных точек, т.е. если $t(x)\ne x$ при любом $x$.
(Заметим, что центральные, осевые и зеркальные симметрии пространства не удовлетворяют этому условию.)

\begin{figure}[h]\centering
\includegraphics{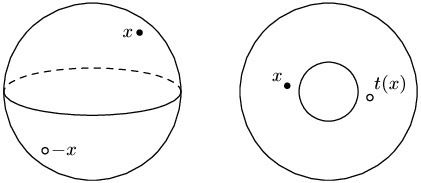}
\caption{Инволюции (1) и (2)}\label{dk1}
\end{figure}

{\it Примеры свободных инволюций.}

(1)
Отображение $\lessmskips{2mu}S^2\to S^2$ сферы $\lessmskips{2mu}S^2=\{(x,y,z)\in\R^3\ |\ x^2+y^2+z^2=1\}$,
заданное формулой $(x,y,z)\mapsto(-x,-y,-z)$. Это --- {\it антиподальная} инволюция (рис.~\ref{dk1}).

(2) Отображение $A^2\to A^2$ кольца $A^2=\{(r,\varphi)\in\R^2\ |\ 1\le r\le2\}$, заданное формулой
$(r,\varphi)\mapsto(r,\pi+\varphi)$ (рис.~\ref{dk1}).

\begin{figure}[h]\centering
\includegraphics{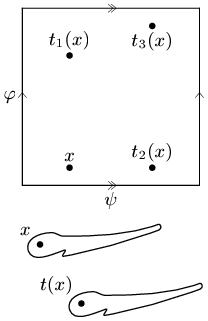}
\caption{Инволюции (3)--(6)}\label{noname}
\end{figure}

(3)--(5) Отображения $t_i:T^2\to T^2$, заданные формулами (рис.~\ref{noname})
$$t_1(\varphi,\psi)=(\pi+\varphi,\psi),
\quad t_2(\varphi,\psi)=(\varphi,\pi+\psi)\quad\text{и}
\quad t_3(\varphi,\psi)=(-\varphi,\pi+\psi).$$
(6) Возьмем несвязное объединение $X\sqcup X=X\times\{+1,-1\}$ двух копий
пространства $X$.
Определим инволюцию $t_X$ на $X\times\{+1,-1\}$ формулой
$(x,s)\mapsto(x,-s)$.
Это --- {\it тривиальная} инволюция (рис.~\ref{noname}).

(7) {\it Фазовым пространством} движения без столкновений двух различных частиц на множестве $N$
называется множество $\t N:=\{(x,y)\ |\ x,y\in N,\ x\ne y\}$.
Определим на этом множестве инволюцию формулой $(x,y)\mapsto(y,x)$.

\begin{pr}
(a) На сфере с любым числом ручек существует свободная инволюция.

(b)* Существует ли свободная инволюция на $\R P^2$?
\end{pr}

\begin{pr}
(a) Приведите пример свободной инволюции на окружности.

(b) На любом ли графе существует свободная инволюция?
 \end{pr}

Пункты (b) этих задач решите хотя бы  для {\it симплициальных} инволюций, определение которых см. ниже для $\R P^2$ или можно дать аналогично нижеприведенному для графа.

 Пусть $T$ --- триангуляция 2-многообразия $N$.
{\it Симплициальная относительно $T$ свободная инволюция} --- это разбиение вершин триангуляции $T$ на пары, удовлетворяющее следующим четырем условиям:

--- каждое ребро соединяет вершины из разных пар
(это условие принято, поскольку любое непрерывное отображение отрезка в себя имеет неподвижную точку);

--- если $\{a,a'\}$ --- пара, $\{b,b'\}$ --- пара и $ab$ --- ребро, то $a'b'$ ---ребро.

--- никакие две вершины в границе грани не лежат в одной паре
(это условие принято, поскольку любое непрерывное отображение диска в себя имеет неподвижную точку);

--- если $\{a_i,a'_i\}$ --- пары и $a_1\dots a_k$ --- грань, то $a_1'\dots a_k'$ --- грань.

Для 2-многообразия $N$, реализованного в качестве подмножества в $\R^3$,
инволюция $t$ на $N$ называется {\it симплициальной относительно разбиения
$T$}, если $t$ переводит вершины в вершины, каждое ребро $ab$ в ребро
$t(a)t(b)$ и каждую грань $a_1\dots a_k$ в грань $t(a_1)\dots t(a_k)$.

{\it Мы будем рассматривать только инволюции, симплициальные относительно некоторого разбиения.}

\begin{pr}
 (a) {\bf Теорема о четности.}
{\it Эйлерова характеристика 2-многообразия с симплициальной свободной инволюцией четна.}

(b) На проективной плоскости с ручками не существует свободной инволюции.

(c) Существует ли свободная инволюция на бутылке Клейна?

(d) На каких 2-многообразиях существует свободная инволюции?
 \end{pr}

 В 2-многообразии $N$ со свободной инволюцией $t$ отождествим точки $x$ и $t(x)$ для любого $x$.
Полученное в результате отождествления пространство будет 2-многообразием
(этим фактом можно пользоваться в дальнейшем без доказательства).
Оно называется {\it факторпространством} пространства $N$ по инволюции $t$ и обозначается $N/t$.
Наша инволюция называется инволюцией {\it на} пространстве $N$ {\it над} пространством $N/t$.

\begin{pr}
  В вышеприведенных примерах (1)--(6)

(a) инволюции симплициальны относительно некоторой триангуляции;

(b) факторпространства гомеоморфны пространствам
$\R P^2$, $A^2$, $T^2$, $T^2$, $K^2$, $X$.
 \end{pr}

\begin{pr}
 (a) Свободные инволюции над данным графом $N$ можно строить следующим образом.
Раскрасим вершины данного графа $N$ в $V$ разных цветов.
Для каждой вершины графа $N$ возьмем две вершины того же цвета.
Это будут вершины графа $N'$.
Для каждого ребра $ij$ графа $N$ возьмем два ребра, каждое из которых
соединяет разноцветные вершины цветов $i$ и $j$ графа $N'$.
Это будут ребра графа $N'$.
На полученном раскрашенном графе $N'$ определяется инволюция, меняющая
местами одноцветные вершины.

(b) Пусть $G$ --- граф с вершинами $a$ и $b$, двукратным ребром
$ab$ и двумя петлями $aa$ и $bb$.
Приведите пример свободной инволюции графа $G$ над графом 'восьмерка'.
 \end{pr}

Инволюция $t:N\to N$ называется {\it тривиальной}, если существует такое 2-многообразие $X$ и гомеоморфизм
$N\to X\sqcup X$, при которой $t$ переходит в тривиальную инволюцию $t_X$ на $X\sqcup X$ из примера (6).
Например, антиподальная инволюция на сфере нетривиальна.

\begin{pr}
Любая свободная инволюция одного из следующих типов является тривиальной:

(a) на графе над деревом.

(b) на 2-многообразии над диском.
\qquad

(c) на 2-многообразии над сферой.
\end{pr}

Эти и следующие задачи разумно решить сначала для {\it конкретных} выбранных Вами простых разбиений рассматриваемых 2-многообразий,
потом для {\it произвольных} разбиений.

\begin{pr}
 Для каких пар $(M,N)$ 2-многообразий существует свободная инволюция 2-многообразия $M$ над 2-многообразием $N$?
 \end{pr}

\subsection{Классификация инволюций}

Инволюции $t_1,t_2:N\to N$ называются {\it эквивалентными}, если существует триангуляция $T$ 2-многообразия $N$ и изоморфизм $f:T\to T$, для которого $t_2\circ f=f\circ t_1$
(т.е. изоморфизм $f$ переводит каждую пару $(x,t_1(x))$ в некоторую пару $(y,t_2(y))$.

\begin{pr}\label{inv-mn}
(a) Инволюции из примеров (3) и (4) эквивалентны.


(b) Приведите пример двух неэквивалентных свободных инволюций на одном и том же 2-многообразии.
{\it Указание.}
Используйте факторпространства.
Ответ: инволюции из примеров (4) и (5) не эквивалентны.

(c) Приведите пример двух не эквивалентных свободных инволюций на сфере с двумя
ручками, факторпространство по каждой из которых гомеоморфно тору.
 \end{pr}

С инволюциями можно работать на эквивалентном языке двулистных накрытий.
Неформально говоря, двулистное накрытие --- отображение, при котором каждая точка имеет два прообраза.
Формально, {\it (неразветвленным) двулистным накрытием} называется отображение $p:N\to N/t$ для некоторой
(свободной) инволюции $t$.
См. эквивалентное определение в \S\ref{03man}.

Инволюция называется {\it тривиальной}, если соответствующее двулистное накрытие тривиально.

Двулистные накрытия $p_1:X_1\to X$ и $p_2:X_2\to X$ называются {\it эквивалентными}, если существует изоморфизм
$f:X_1\to X_2$ комплексов, для которого $p_1=p_2\circ\varphi$.

Обозначим через $S(M)$ множество двулистных накрытий над $M$ ('инволюций над
$M$') с точностью до эквивалентности над $M$.

\smallskip
{\bf Теорема классификации двулистных накрытий над графом.} {\it Множество
$S(G)$ симплициальных двулистных накрытий над графом $G$ с точностью до
эквивалентности находится во взаимно однозначном соответствии с
$H_1(G)\cong\Z_2^{E-V+C}$.}

\smallskip
{\it Набросок доказательства.}
Рассмотрим симплициальное двулистное накрытие $p:\t G\to G$ над графом $G$.
Для каждой вершины $v$ графа $G$ одну из точек в $p^{-1}$ назовем верхней, а вторую нижней.
Этот выбор обозначим через $o$.
Поставим на ребре $e$ графа $G$ единицу, если соответствующие два ребра
$p^{-1}e$ в $\t G$ соединяют разноименные вершины, и ноль, если одноименные.
Назовем полученную расстановку {\it препятствующей} и обозначим $\omega(o)$.
(Пояснение к названию, не используемое в дальнейшем: эта расстановка препятствует тривиальности двулистного накрытия.)

{\it Одномерная группа гомологий} $H_1(G)$ графа $G$ определяется аналогично \S\ref{02hom} как группа циклов в графе $G$.
Для любого цикла $g$ сумма $\omega(o)\cdot g$ значений $\omega(o)$ по всем ребрам цикла $g$ не зависит от набора
$o$ стрелочек.
Поэтому формула $w_1(p)[g]=\omega(o)\cdot g$ корректно задает линейную функцию $w_1(p):H_1(G)\to\Z_2$.
Она называется {\it первым классом Штифеля-Уитни} двулистного накрытия.
Доказательство того, что соответствие $p\mapsto w_1(p)$ является биекцией, оставляем в качестве упражнения.
 QED

\begin{pr} (a) При этой биекции композиция инволюций переходит в сумму циклов.

(b) Сколько классов эквивалентности симплициальных двулистных накрытий над тором?

(c) Сколько классов эквивалентности симплициальных свободных инволюций на торе?
 \end{pr}

 Перейдем к классификации симплициальных двулистных накрытий над данным 2-многообразием $M$.
Пусть задана некоторая инволюция $t:N\to N$ над $M$.
Мы построим препятствие к тривиальности инволюции $t$.
С помощью этого препятствия нетрудно классифицировать симплициальные двулистные накрытия над $M$.
Рекомендуем читателю разобрать следующие рассуждения сначала для инволюции $t:T^2\to T^2$ из примера (2).

\begin{figure}[h]\centering
\includegraphics[scale=.24]{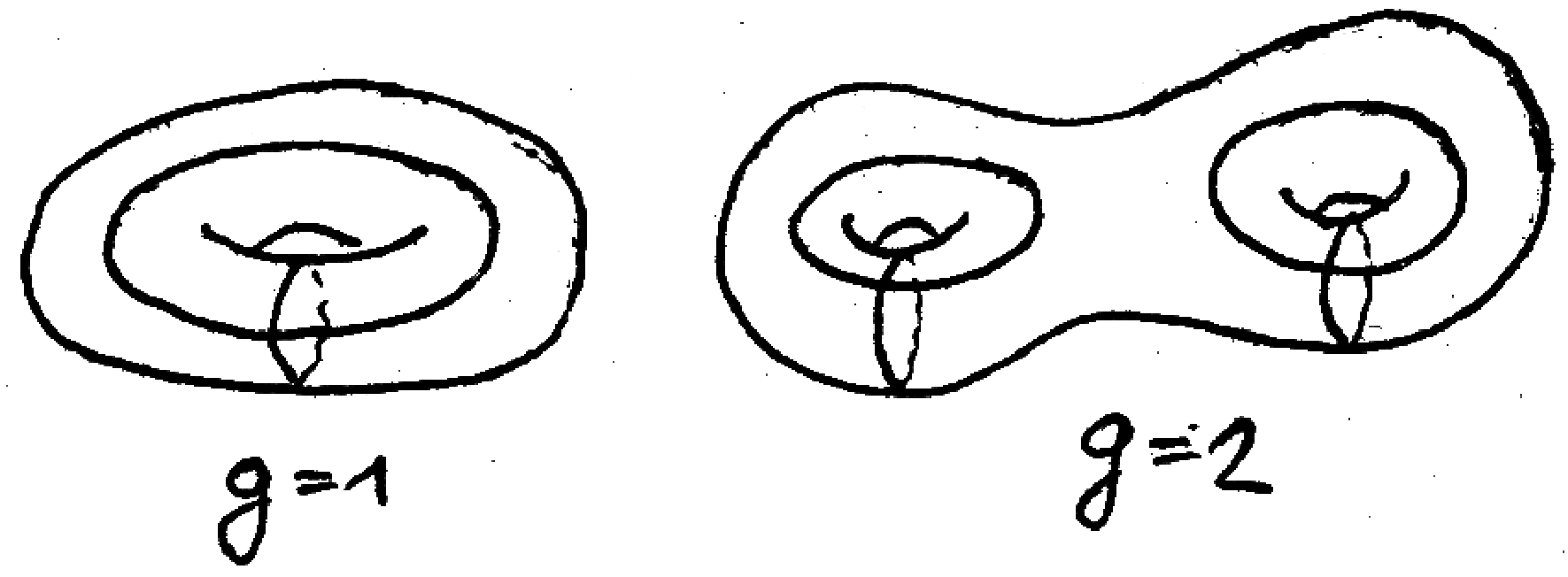}
\caption{Набор кривых для классификации инволюций}\label{2g}
\end{figure}

Выберем на сфере $M$ с $g$ ручками $2g$ кривых, как на рис.~\ref{2g}.
Для двулистного накрытия $p:\t M\to M$ и одной из этих кривых поставим на
кривой 0, если соответствующая инволюция над ней тривиальна (т.е. является
перестановкой соответствующих точек из двух копий окружности), и 1, если нетривиальна.
Получим набор $w_1(p)$ из $2g$ нулей и единиц, называемый {\it первым классом
Штифеля-Уитни} двулистного накрытия $p$ относительно данного выбора кривых~$L$.

Аналогичные наборы из $2g$ нулей и единиц можно получить,
рассматривая некоторые другие наборы замкнутых кривых на $M$.

\begin{pr}
Постройте аналогичный набор из $2-\chi(M)$ кривых на {\it неориентируемом}
многообразии $M$ (так, чтобы была выполнена следующая теорема).
\end{pr}

{\bf Теорема классификации  двулистных накрытий.}
{\it Для замкнутого 2-многообразия $M$ отображение
$w_1:S(M)\to\Z_2^{2-\chi(M)}$ (из множества $S(M)$ симплициальных двулистных
накрытий над графом $G$ с точностью до эквивалентности в группу
$H_1(M)\cong\Z_2^{2-\chi(M)}$ одномерных гомологий 2-многообразия
$M$ с коэффициентами в $\Z_2$) является взаимно-однозначным соответствием. }

\smallskip
Заметим, что группа гомологий $H_1(M)$ естественно появилась при классификации
двулистных накрытий независимо от ее появления при изучении ориентируемости.

\begin{pr}
 (a) При отображении $w_1$ тривиальная инволюция переходит в нулевой элемент.

(b) При отображении $w_1$ композиция инволюций переходит в
сумму циклов.

(c) Докажите теорему классификации  двулистных накрытий.

(d) Сформулируйте и докажите аналог теоремы классификации двулистных накрытий для многообразий с краем.

(e) То же для утолщений графов.

(f) Сформулируйте и докажите аналог теоремы классификации двулистных накрытий
для двулистных накрытий многообразий, сохраняющих ориентацию.
 \end{pr}

 \subsection{Другое доказательство теоремы классификации инволюций}

Это доказательство сложнее предыдущего.
Зато, в отличие от предыдущего, построения этого доказательства применимы для
классификации двулистных накрытий на 2-многообразиях, заданных триангуляцией (например, в компьютере).
Кроме того, это доказательство может быть обобщено для трехмерных (и многомерных) многообразий.

\smallskip
{\it Построение препятствующего цикла.}
Рассмотрим свободную симплициальную инволюцию $t:N\to N$, для которой $M=N/t$.
Возьмем некоторое разбиение $T$ многообразия $N$.
Построим двойственное разбиение $U$ (см. определение в \S\ref{02hom}).
Они порождают двойственные разбиения $T/t$ и $U/t$ многообразия $M=N/t$.

\begin{figure}[h]\centering
\includegraphics{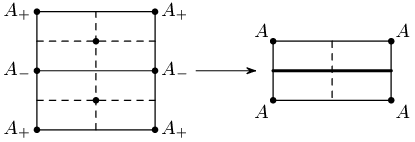}
\caption{Препятствующий цикл для инволюции на торе (над тором)}\label{dk4}
\end{figure}

Для каждой пары инволютивных вершин разбиения $T$ одну из них назовем верхней, а вторую нижней.
Этот выбор обозначим через $o$.
Покрасим в красный цвет те ребра разбиения $U/t$, для которых соответствующие два
ребра разбиения $T$ соединяют разноименные вершины.
Объединение красных ребер называется {\it препятствующим циклом} $\omega(o)$.
На рис.~\ref{dk4} жирной линией изображен препятствующий цикл для инволюции на
торе (т.е. для соответствующего двулистного накрытия над тором) из примера (3).

Ясно, что

(*) Данное распределение $o$ верхних и нижних вершин определяет
тривиализацию инволюции тогда и только тогда, когда $\omega(o)=\emptyset$.

(**) В графе $\omega(o)$ {\it нет изолированных вершин и из каждой вершины
выходит четное число красных ребер}, т.е. этот граф является (гомологическим) циклом.

\begin{pr}
 (a) Докажите эти утверждения.

(b) Вне графа $\omega(o)$ двулистное накрытие тривиально.
 \end{pr}

 {\it Изменение препятствующего цикла.}
Если $\omega(o)\ne\emptyset$, то $o$ не определяет тривиализации инволюции,
но еще не все потеряно: можно попытаться изменить $o$, так чтобы
препятствующий цикл стал пустым.
Для этого выясним, как $\omega(o)$ зависит от $o$.

Изменим названия двух инволютивных вершин разбиения $T$, соответствующих
вершине разбиения $T/t$, лежащей в грани $a$.
Тогда цикл $\omega(o)$ изменится на сумму по модулю 2 (т.е. симметрическую
разность) с границей этой грани: $\omega(o')=\omega(o)+\partial a$.
Поэтому при изменении названий нескольких инволютивных вершин разбиения $T$
цикл $\omega(o)$ изменяется на сумму по модулю 2 границ соответствующих граней разбиения $U/t$:
$$\omega(o')=\omega(o)+\partial a_1+\dots+\partial a_k.$$
Ясно, что\nopagebreak

(i) При изменении набора $o$ препятствующий цикл $\omega(o)$ заменяется на гомологичный цикл.

(ii) Если $\omega(o)$ является границей, то можно изменить $o$ на $o'$, так чтобы получилось $\omega(o')=0$.

\smallskip
{\it Завершение доказательства.}
Напомним, что {\it одномерной группой гомологий} разбиения $T$ (или многообразия $M$) называется группа $H_1(M)$
циклов с точностью до гомологичности.
Напомним также, что $H_1(M)\cong\Z_2^{2-\chi(M)}$.

{\it Первым классом Штифеля-Уитни} инволюции $t$ (или соответствующего
двулистного накрытия) называется класс гомологичности препятствующего цикла:
$$w_1(t):=[\omega(o)]\in H_1(M).$$
Ясно, что он является {\it полным} препятствием к
тривиальности инволюции $t$:

{\it инволюция $t$ тривиальна тогда и только тогда, когда $w_1(t)\ne0$.}

Теперь теорема классификации двулистных накрытий вытекает из (ii).
 QED

\newpage
\section{Векторные поля на плоскости}\label{0vepl}

 \hfill{\it 'You mean...' he would say, and then he would rephrase what I had said}

\hfill{\it in some completely simple and concrete way, which sometimes illuminated}

\hfill{\it  it enormously, and sometimes made nonsense of it completely.}

\hfill{\it I. Murdoch, Under the Net.
\footnote{\it 'Вы хотите сказать...' --- начинал он и пересказывал мои слова
конкретно и просто, после чего моя мысль либо оказывалась много понятнее и
глубже, либо оборачивалась полнейшей чепухой. А. Мердок, Под сетью, пер. М.
Лорие.}
}

\subsection{Интересные примеры и теоремы}

 Исследование векторных полей было начато Анри Пуанкаре в качественной
теории дифференциальных уравнений, возникшей при изучении небесной механики.
С тех пор векторные поля являются одним из важнейших объектов топологии и ее приложений.
Подробные мотивировки см. в `похвальном слове векторным полям'~\cite{An03}.
Проблема существования векторных полей и их наборов
определяет лицо теории препятствий.
Ее методы можно применять ко многим другим задачам.
(О вложениях  многообразий и графов см., например, \S\ref{0emb} и [S].)


\begin{pr}\label{vfcon}
(a) На плоскости дано $2n$ красных и $2n$ синих точек общего положения.
Докажите, что существует прямая,  по каждую сторону от которой находится по
$n$ красных и синих точек.

(b) Пусть дан (выпуклый) многоугольник на плоскости и точка $z$ вне него.
Тогда существует прямая, содержащая $z$ и делящая многоугольник на две части равной площади.



(c) Для любого (выпуклого) многоугольника на плоскости существует прямая,
делящая его на две части с равными площадями и периметрами.



(d) Любой ли бутерброд с маслом можно разрезать прямолинейным разрезом на две равноценные части?

(e) То же для бутерброда с маслом и сыром.

(f) На неровном полу табуретку всегда можно поставить на все четыре ножки.

(g) Вокруг любого выпуклого многоугольника на плоскости можно описать квадрат.

(h)* Любую плоскую выпуклую фигуру диаметра 1 можно поместить в правильном
шес\-ти\-у\-голь\-ни\-ке, расстояние между противоположными сторонами которого равно 1.
\end{pr}

Приведенные задачи решаются при помощи {\it соображений непрерывности,} т.е. теоремы о промежуточном значении непрерывной функции.
Эта теорема является `нульмерной версией' основной теоремы топологии (см. далее).

\begin{figure}[h]
\centering\hskip-6cm\special{em:graph 6-1}
\vskip1.7cm
\caption{Векторное поле на подмножестве плоскости}
\label{6-1}
\end{figure}

{\it Векторным полем} на подмножестве плоскости называется семейство векторов $v(x)$ на плоскости
в его точках $x$, непрерывно зависящих от точки $x$.

\begin{figure}[h]
\centering\hskip-6cm\special{em:graph 6-2}
\vskip2.7cm
\caption{Радиальное, центральное,  и седловое векторные поля на плоскости}
\label{6-2}
\end{figure}

{\it Примеры} векторных полей на плоскости:
$a(x,y)=(x,y)$ ({\it радиальное}), $b(x,y)=(y,-x)$ ({\it центральное}), $c(x,y)=(y,x)$ ({\it седловое}),
$u(x,y)=(1,0)$ ({\it постоянное}).

Пусть $N\subset\R^m$, $Y\subset\R^n$, причем $N$ замкнуто и ограничено.
Отображение $f:N\to Y$ называется {\it непрерывным}, если для любого $\varepsilon>0$ существует такое $\delta>0$,
что при любых $x,y\in N$ с условием $|x-y|<\delta$ выполнено $|f(x)-f(y)|<\varepsilon$.
{\it Далее все отображения считаются непрерывными и прилагательное `непрерывное' опускается.}

Векторное поле на подмножестве плоскости --- то же, что непрерывное отображение из этого подмножества в плоскость.

Обозначим окружность и диск через
$$S^1=\{(x,y)\in\R^2\ |\ x^2+y^2=1\}=\{z\in\C\ :\ |z|=1\}\quad\text{и}\quad
D^2=\{(x,y)\in\R^2\ |\ x^2+y^2\le1\}.$$
Векторное поле называется {\it ненулевым}, если все его векторы ненулевые.


\begin{pr}\label{vf-brower}
Следующие утверждения эквивалентны.

(1) {\bf Теорема непродолжаемости.}
{\it Радиальное векторное поле $a(x,y)=(x,y)$ на граничной окружности $S^1$ диска $D^2$ не продолжается до ненулевого
векторного поля на диске.}

(2) {\bf Несминаемость диска на окружность.}
{\it Не существует отображения диска в его граничную окружность, тождественного на этой окружности,
т.е. отображения $f:D^2\to S^1$, для которого $f(x)=x$ при $x\in S^1$.}
 (Или, говоря неформально, барабан нельзя смять на его обод).

(3) {\bf Теорема Брауэра о неподвижной точке.}
{\it Любое отображение $f:D^2\to D^2$ диска в себя имеет неподвижную точку, т.е. такую точку $x\in D^2$, что $f(x)=x$.}
\end{pr}

В этой и следующей задачах требуется именно доказать эквивалентность, доказывать сами утверждения не требуется.
Утверждение (1) задачи \ref{vf-brower} вытекает из задачи \ref{vfalg}.a ниже.
Утверждение (3) задачи \ref{vf-bu} вытекает из задач \ref{vfalg}.ab ниже.

\begin{pr}\label{vf-bu}
Обозначим сферу через $S^2=\{(x,y,z)\in\R^3\ |\ x^2+y^2+z^2=1\}.$
Отображение $f:S^2\to S^1$ называется {\it эквивариантным} (или {\it нечетным}),
если $f(-x)=-f(x)$ для любого $x\in S^2$.
Следующие утверждения эквивалентны.
(Каждое из них называется {\bf теоремой Борсука-Улама}.)

{\it (1) Для любого отображения $f:S^2\to\R^2$ существует $x\in S^2$ такое, что $f(x)=f(-x)$.}

(Т.е. в любой момент времени на Земле найдутся диаметрально
противоположные точки, в которых температура и давление совпадают.)

{\it (2) Не существует эквивариантного отображения $S^2\to S^1$.

(3) Никакое эквивариантное отображение $S^1\to S^1$ не продолжается на диск $D^2$.

(4) Если сфера $S^2$ является объединением трех замкнутых множеств, то одно из них содержит диаметрально
противоположные точки.}
 \end{pr}


\subsection{Гомотопность векторных полей и непрерывных отображений}

\begin{pr}\label{vfplex}
Для любых двух векторных полей существует непрерывная деформация одного в другое, т.е. семейство $v_t$ векторных
полей, непрерывно зависящее от параметра $t\in[0,1]$, для которого $v_0$ есть первое поле и $v_1$ есть второе поле.
(В отличие от дальнейшего, векторное поле $v_t$ не предполагается ненулевым.)
\end{pr}

Отсюда  видно, что при непрерывной деформации векторного поля его {\it траектории} меняются разрывно.
(Для ненулевых векторных полей см. задачу \ref{vfpl-}.)

Два ненулевых векторных поля называются {\it гомотопными}, если одно можно получить из другого непрерывной
деформацией, в процессе которой векторное поле остается ненулевым.
Такой деформацией называется семейство $v_t$ ненулевых векторных полей, непрерывно зависящее от параметра $t\in[0,1]$,
для которого $v_0$ есть первое поле и $v_1$ есть второе поле.

\begin{pr}\label{vfpl-}
(a) Любое ненулевое векторное поле $v$ на плоскости гомотопно полю $-v$.

(b) Радиальное векторное поле на окружности гомотопно центральному.
\end{pr}

\begin{pr}\label{vfany}
Любые два ненулевых векторных поля на $N$ гомотопны, если
\qquad

(a) $N=0\times[0,1]$; \qquad
(b) $N=D^2$; \qquad

(c) $N$ --- произвольное дерево на плоскости.
\qquad
(d) $N$ --- вся плоскость.
\end{pr}


\begin{pr}\label{vfd2}
(a) Каждое из утверждений задачи \ref{vf-brower} равносильно тому, что

\qquad (4) радиальное векторное поле на окружности $S^1$ не гомотопно постоянному.

(b) Ненулевое векторное поле на окружности $S^1$ гомотопно постоянному тогда и только тогда,
когда оно продолжается на диск $D^2$.

(с) Векторные поля $v(z):=9z^3-2z^2+z-1$ и $w(z):=z^3$ на окружности $S^1$ гомотопны.
\end{pr}


\begin{pr}\label{vfequi} Для любых ненулевых векторных полей $u,v,w$

(a) $v$ гомотопно $v$ ({\it рефлексивность}).

(b) если $u$ гомотопно $v$, то $v$ гомотопно $u$ {\it cимметричность}).

(c) если $u$ гомотопно $v$ и $v$ гомотопно $w$, то $u$ гомотопно $w$ ({\it транзитивность}).
\end{pr}

Векторное поле называется {\it единичным}, если все его векторы единичные.

\begin{pr}\label{vfunit}
(a) Любое ненулевое векторное поле гомотопно единичному.

(b) Два единичных векторных поля гомотопны тогда и только тогда, когда существует семейство $v_t$ {\it единичных}
векторных полей, непрерывно зависящее от параметра $t\in[0,1]$,
для которого $v_0$ есть первое поле и $v_1$ есть второе поле.
\end{pr}

Единичное векторное поле на подмножестве плоскости --- то же, что непрерывное отображение из этого
подмножества в окружность.

Два отображения $f_0,f_1:N\to S^1$ из  подмножества $N\subset\R^2$ называются {\it гомотопными}, если
существует семейство $h_t:N\to S^1$ отображений, непрерывно зависящее от параметра $t\in[0,1]$, для которого $h_0=f_0$ и $h_1=f_1$.
Это эквивалентно существованию такого непрерывного отображения $H:N\times[0,1]\to S^1$, что
$H(x,0)=f_0(x)$ и $H(x,1)=f_1(x)$.
Обозначение: $f_0\simeq f_1$.


Для гомотопности отображений справедливы свойства, аналогичные гомотопности векторных полей.
Например, гомотопность отображений --- отношение эквивалентности;
любые два отображения дерева в окружность гомотопны.

\begin{pr}\label{vfann} Обозначим кольцо $A:=\{(x,y)\in\R^2\ |\ 1\le x^2+y^2\le2\}$.

(a) Существует отображение $r:A\to S^1$, тождественное на $S^1$.

(b) Любое отображение $S^1\to S^1$ продолжается до отображения $A\to S^1$.

(c) Для любых двух отображений $A\to S^1$ любая гомотопия между их сужениями на $S^1$ продолжается до
гомотопии между исходными отображениями.
\end{pr}

\begin{pr}\label{vfbor} (a) Для любого отображения $f:[0,1]\to S^1$ любая гомотопия отображения $f|_{\{0,1\}}$ продолжается до некоторой гомотопии отображения $f$.

(b) Для любого отображения $f:D^2\to S^1$ любая гомотопия отображения $f|_{S^1}$ продолжается до некоторой гомотопии отображения $f$.

(c) {\it Теорема Борсука о продолжении гомотопии.}
Для любых подграфа $A$ плоского графа $N$ и отображения $f:N\to S^1$ любая гомотопия отображения $f|_A$
продолжается до некоторой гомотопии отображения $f$.

(d) Приведите пример подмножеств $A\subset N$ плоскости, отображения $f:N\to S^1$ и гомотопии отображения $f|_A$,
не  продолжаемой до гомотопии отображения $f$.
\end{pr}


\subsection{Число оборотов вектора и его применения}

Существует ли подмножество плоскости и негомотопные единичные векторные поля на нем?
Да. Например, седловое векторное поле на окружности не гомотопно радиальному (а, значит, и центральному).
Доказательство негомотопности непросто.
Для него, и для доказательства вышеприведенных интересных теорем, необходимо следующее понятие.

Неформально, {\it степенью} $\deg v$ ненулевого векторного поля $v$ на окружности $S^1$
называется число оборотов вектора $v(x)$ (вокруг своего неподвижного начала)
при однократном обходе точкой $x$ окружности $S^1$ против часовой стрелки.
Приведем строгое определение этого понятия.

Виду равномерной непрерывности отображения $S^1\to S^1$, соответствующего единичному векторному полю $v$,
существует $M>0$, для которого $|v(x)-v(y)|<1$ при $x-y<1/M$.
Для $k=1,2,\dots,M$ обозначим через $2\pi\Delta_k\in(-\pi,\pi)$ угол, на который нужно повернуть против часовой стрелки вектор $v(e^{2\pi i(k-1)/M})$, чтобы получить вектор $v(e^{2\pi ik/M})$.
Более точно, $\Delta_k\in(-1/2,1/2)$ однозначно определяется из условия
$v(e^{2\pi ik/M})=v(e^{2\pi i(k-1)/M})e^{2\pi i\Delta_k}$.
Определим
\footnote{Т.е. $\deg v:=\int_0^{2\pi} ds(u):=\lim\limits_{\max(u_{k+1}-u_k)\to0}
\{\sum\limits_{k=1}^M \Delta_k\ |\ 0=u_0<u_1<\dots<u_{M-1}<u_M=2\pi\},$ где
$\Delta_k\in(-1/2,1/2)$ при $u_k-u_{k-1}<1/M$ однозначно определяется из условия
$v(e^{iu_k})=v(e^{iu_{k-1}})e^{2\pi i\Delta_k}$.
Заметим, что $\Delta_k$ зависят также от $v$ и $M$; мы не указываем это в обозначениях.}
$$\deg v:=\sum\limits_{k=1}^M \Delta_k\in\Z.$$
Аналогично определяется степень {\it ненулевого} векторного поля на окружности или отображения окружности в себя.

Нужно доказать корректность определения степени.
Т.е. то, что $\deg v$ действительно целое число, не зависящее от выбора числа $M$.
Это можно доказать аналогично доказательству независимости определенного интеграла от выбора последовательности разбиений.
Более удобный способ --- задача \ref{vfcov}.cd.


\begin{pr}\label{vfsur} (a) Найдите, хотя бы для одного $M$, степень {\it стандартной $n$-намотки} $w_n$, т.е.,
отображения  $S^1\to S^1$ (или, что то же самое, векторного поля на $S^1$), заданного формулой
$$w_n(z)=z^n.$$
\quad
(b) Если $u$ и $v$ --- ненулевые векторные поля на окружности $S^1$ и для любой точки $z\in S^1$ векторы $u(z)$ и $v(z)$ симметричны относительно касательной к окружности в точке $z$.
Тогда $\deg u+\deg v=2$.
\end{pr}

Нам понадобится также непрерывная зависимость степени от векторного поля (точнее, от его гомотопии).
Это удобно доказать при помощи следующего обобщения определения степени.

{\it Поднятием} (или угловой функцией) отображения $s:N\to S^1$ называется отображение $\t s:N\to\R$ такое,
что $e^{i\t s}=s$.


\begin{pr}\label{vfcov}
(a) Найдите все поднятия пути (т.е. отображения отрезка) $s:[0,2\pi]\to S^1$, $s(t)=e^{it}$.

(b) Найдите все поднятия композиции пути $s$ из (a) и стандартной $n$-намотки $w_n$.

(с) По единичному векторному полю $v$ на окружности определим путь $s_v:[0,1]\to S^1$ формулой
$s_v(t)=v(e^{2\pi it})$.
Если $\t s_v$ --- поднятие пути $s_v$, то  $\deg v=\dfrac{\t s_v(1)-\t s_v(0)}{2\pi}.$

(d) Поднятие единственно, т.е. если $\t s,\t s':[0,1]\to\R$ --- два поднятия одного и того же
пути $[0,1]\to S^1$, причем $\t s(0)=\t s'(0)$, то $\t s=\t s'$.
\end{pr}

\begin{pr}\label{vfext}
(a) {\bf Лемма о поднятии пути.}
Любой путь $s:[0,1]\to S^1$ имеет поднятие $\t s:[0,1]\to\R$.

(a') Для любых пути $s:[0,1]\to S^1$ и точки $\t x\in\R$ с условием $e^{i\t x}=s(0)$
существует поднятие $\t s:[0,1]\to\R$ пути $s$, для которого $\t s(0)=\t x$.

(b) {\bf Лемма о поднятии гомотопии.}
Любая гомотопия $H:[0,1]^2\to S^1$ имеет поднятие.


(c) Степени гомотопных единичных векторных полей равны.

(d) $w_n$ не гомотопно $w_m$ при $m\ne n$.

(e) Существует отображение $S^1\to S^1$, не имеющее поднятия.
\end{pr}

\begin{pr}\label{vfalg}
(a) Ненулевое векторное поле на окружности не продолжается на диск, если
число оборотов вектора при обходе этой окружности (точнее, степень) не равно нулю.

(b) Степень любого нечетного отображеения $S^1\to S^1$ нечетна.

(c) Для любых трех выпуклых многогранников в пространстве существует плоскость, делящая каждый из них на две части равных объемов.

(d) Если $u$ и $v$ ненулевые векторные поля на окружности $S^1$, причем $|u(x)|>|v(x)|$ для любой точки $x\in S^1$,
то $\deg u=\deg(u+ v)$.

(e) Найдите степень векторного поля $v(z):=9z^3-2z^2+z-1$.

(f) {\bf Основная теорема алгебры.} {\it Любой непостоянный многочлен с
комплексными коэффициентами имеет комплексный корень.}
\end{pr}


\begin{figure}[h]
\centering\hskip-6cm\special{em:graph 7-9}
\vskip2.6cm
\caption{Диск с дырками}
\label{7-9}
\end{figure}

\begin{pr}\label{vfsum} (a) Для любого ненулевого векторного поля на кольце степени
его сужений на две граничные окружности кольца равны.

(b) Для любого ненулевого векторного поля на диске с дырками (рис. \ref{7-9}) степень
его сужения на внешнюю граничную окружность равна сумме степеней его сужений внутренние граничные окружности.
(Каждая окружность ориентирована против часовой стрелки.)
\end{pr}

\subsection{Классификация векторных полей на подмножествах плоскости}

\begin{pr}\label{vfbas}
(a) Любое единичное векторное поле степени $n$ на окружности  гомотопно $w_n$.

(b) Любые два ненулевые векторные поля на окружности равных степеней гомотопны.

(c) {\bf Теорема продолжаемости.}
{\it Ненулевое векторное поле продолжается с граничной окружности диска на диск
тогда и только тогда, когда число оборотов вектора при обходе этой окружности (точнее, степень) равно нулю.}
\end{pr}

Обозначим через $V(N)$ множество единичных векторных полей на подмножестве $N$ плоскости, с точностью до гомотопности
в классе единичных векторных полей (т.е., с точностью до непрерывной деформации, в процессе которой векторное поле остается единичным).
`Нормировка' определяет взаимно-однозначное соответствие между $V(N)$ и множеством {\it ненулевых} векторных полей на $N$ с точностью до гомотопности (задача \ref{vfunit}).
Множество $V(N)$ находится также во взаимно-однозначном соответствии с множеством отображений $N\to S^1$
с точностью до гомотопности.

Из задач \ref{vfsur}.a, \ref{vfcov}.cd, \ref{vfext}.a'd и \ref{vfbas}.a вытекает следующий результат.

\smallskip
{\bf Основная теорема топологии.}
{\it Любое единичное векторное поле $v$ на окружности гомотопно стандартной $\deg v$-намотке, и единичные
векторные поля разных степеней  на окружности не гомотопны.

Любое отображение $f:S^1\to S^1$ гомотопно стандартной $\deg f$-намотке, и отображения $S^1\to S^1$
разных степеней не гомотопны.

Иными словами, степень $\deg:V(S^1)\to\Z$ является взаимно-однозначным соответствием, переводящим класс
$d$-кратной намотки в число $d$. }

\smallskip
Этот результат назван основной теоремой топологии по аналогии с `Основной теоремой алгебры'
(из него Основная теорема алгебры и вытекает).
Суть дела, конечно, лучше бы отражали названия `основная теорема алгебры полиномов' и
`основная теорема одномерной топологии', но эти названия в литературе не используются.
Обобщения см. в \S\ref{0vehima}.

\begin{pr}\label{hopf-gra}
(a) Пусть $N$ --- несвязное объединение или букет $k$ окружностей (рис.~\ref{RSo00-1.1}).
Фиксируем произвольно направление на каждой из этих окружностей.
Для векторного поля $v$ на $N$ поставим на каждой из этих окружностей степень сужения поля $v$ на эту окружность.
Полученную расстановку $k$ целых чисел обозначим $\deg v$.
Тогда $\deg$ определяет биекцию $V(N)\to \Z^k$.

(b) Для плоского графа $N$ множество $V(N)$ `не меняется' при стягивании ребра (осуществляемом в плоскости).

(c) Для плоского графа $N$ постройте биекцию $V(N)\to\Z^{E-V+C}=\Z^{F-1}$.

(d) Для графа $N$ (не обязательно плоского) постройте биекцию между множеством отображений $N\to S^1$ с точностью до гомотопности и $\Z^{E-V+C}$.
\end{pr}


\begin{pr}\label{vfgr}
Опишите $V(N)$ для дискa с $n$ дырками $N$ (начните с $n=0,1$).
\end{pr}

Здесь `описать' означает построить `естественное' взаимно-однозначное
соответствие между $V(N)$ и некоторым `известным' множеством.
`Известность' множества означает как минимум описание количества его элементов, а как максимум ---
наличие `естественных' операций на множестве и на $V(N)$, сохраняемых соответствием.

\begin{pr}\label{vfhom}
{\bf Теорема гомотопности.}
{\it Любые два единичных векторных поля на квадрате $[0,1]^2$, совпадающие на его границе, гомотопны неподвижно на границе (т.е. гомотопны так, что $v_t(x)=v_0(x)$ для любых $x\in\partial[0,1]^2$ и $t\in[0,1]$).}
\end{pr}

Более подробное изложение приводится в~\cite[\S6]{Pr95},~\cite[\S\S4-6]{An03},~\cite[\S\S20,21]{BE82}.

\subsection{Указания к некоторым задачам}

\smallskip
{\bf \ref{vfcon}.} (b) Зафиксируем любую направленную прямую, проходящую через $z$.
Для каждого $\alpha\in[0,2\pi]$ обозначим через
$S_+(\alpha)$ (соответственно $S_-(\alpha)$) площадь пересечения многоугольника с правой (соответственно с левой)  полуплоскостью относительно направленной прямой, образующей угол  $\alpha$ с зафиксированной направленной прямой.
Обозначим $f(\alpha):=S_+(\alpha)-S_-(\alpha)$.
Тогда $f(0)=-f(\pi)$.

Обозначим через $R$ радиус круга с центром в $z$, содержащего многоугольник.
 Тогда $|S_+(\alpha)-S_+(\alpha')|<|\alpha-\alpha'|R^2/2$.
Поэтому функция $S_+$ непрерывна.
Аналогично функция $S_-$ непрерывна.
Значит, функция $f$ непрерывна.
Из этого и $f(0)=-f(\pi)$ следует, что существует $\alpha_0$, для которого $f(\alpha_0)=0$.

\smallskip
{\bf \ref{vf-brower}.} $(1)\Rightarrow(2)$, поскольку отображение $D^2\to S^1$ --- то же, что векторное поле.

Докажем, что $(2)\Rightarrow(3)$.
Пусть, напротив, отображение $f:D^2\to D^2$ не имеет неподвижных точек.
Для $z\in D^2$ обозначим через $\varphi(z)\in S^1$ точку пересечения с $S^1$ луча с началом в $f(z)$ и проходящего через $z$. Иными словами, определим отображение $\varphi:D^2\to S^1$ формулой $\varphi(z):=f(z)+(z-f(z))\dfrac{-f(z)\cdot\delta+\sqrt{(f(z)\cdot\delta)^2+|\delta|^2(1-|f(z)|^2)}}{|\delta|^2}$, где $\delta:=z-f(z)$.
Противоречие c (2).

Докажем, что $(3)\Rightarrow(2)$. Пусть, напротив, существует продолжение $f:D^2\to S^1$ тождественного отображения $S^1\to S^1$. Тогда $-f:D^2\to S^1\subset D^2$ --- отображение, не имеющее неподвижных точек.
Противоречие c (3).

Докажем, что $(1)\Rightarrow(3)$ (хоть это уже и не нужно).
Пусть, напротив, отображение $f:D^2\to D^2$ не имеет неподвижных точек.
Определим векторное поле $v$ на $2D^2$ формулой
$v(z)=\begin{cases}z-f(z) & |z|\le1 \\
z-(2-|z|)f(z/|z|) & 1\le|z|\le2\end{cases}$.
Противоречие c (1).

Докажем, что $(3)\Rightarrow(1)$ (хоть это уже и не нужно).
Пусть, напротив, существует ненулевое векторное поле $v$ на диске, сужение которого на краевую окружность диска является радиальным полем.
Определим отображение $f:2D^2\to 2D^2$ формулой
$f(z)=\begin{cases}z+v(z) & |z|\le1 \\
2ze^{i(|z|-1)}/|z| & 1\le|z|\le2\end{cases}$.
Противоречие c (3).

\smallskip
{\bf \ref{vf-bu}.} $(1)\Rightarrow(2)$ очевидно.

Докажем, что $(2)\Rightarrow(1)$.
Пусть, напротив, для отображения $f:S^2\to\R^2$ и любой $x\in S^2$ выполнено $f(x)\ne f(-x)$.
Определим отображение $\t f:S^2\to S^1$ формулой $\t f(x):=\dfrac{f(x)-f(-x)}{|f(x)-f(-x)|}$.
Это отображение нечетно.
Противоречие c (1).

Для доказательства $(2)\Rightarrow(3)$ используйте ортогональную проекцию верхней полусферы на диск.

Для доказательства $(4)\Rightarrow(2)$ предположите, что $f:S^2\to S^1$ --- эквивариантное отображение и рассмотрите
замкнутые множества $f^{-1}(e^{[0,2\pi/3]i})$, $f^{-1}(e^{[2\pi/3,4\pi/3]i})$, $f^{-1}(e^{[4\pi/3,2\pi]i})$.

Докажем, что $(1)\Rightarrow(4)$. Пусть, напротив, $S^2=A\cup B\cup C$ --- покрытие сферы замкнутыми множествами, ни одно из которых не содержит диаметрально противоположных точек.
Определим функцию $\varphi_A:S^2\to[-1,1]$ формулой $\varphi_A:=\dfrac{|x,A|-|x,-A|}{|x,A|+|x,-A|}$...

Тогда $\varphi_A^{-1}(\pm1)=\pm A$.
Аналогично определим функцию $\varphi_B:S^2\to[-1,1]$.
Применим (1) к отображению $\varphi_A\times\varphi_B:S^2\to[-1,1]^2$.
Получим $x\in S^2$, для которого $\varphi_A(x)=\varphi_A(-x)$ и $\varphi_B(x)=\varphi_B(-x)$.
Если $x\in A$, $-x\in-A$. Поэтому  $\varphi_A(x)-\varphi_A(-x)=2\ne0$ --- противоречие.
Если $-x\in A$, то $x\in -A$. Поэтому $\varphi_A(x)-\varphi_A(-x)=-2\ne0$ --- противоречие.
Поэтому $x\not\in\pm A$.
Аналогично $x\not\in\pm B$.
Значит, $x\in C$ и $-x\in C$.

Идея другого доказательства. Ясно, что $A\cap B\cap C=\emptyset$.
Построим отображение сферы в {\it нерв} покрытия $S^2=\overline A\cup\overline  B\cup\overline  C$, гомеоморфный $S^1$.

\bigskip
{\bf \ref{vfplex}.} Положим $v_t(z):=tv_1(z)+(1-t)v_0(z)$.

\smallskip
{\bf \ref{vfpl-}.}
(a) $v_t(z)$ получается из $v(z)$ поворотом на $\pi t$ в положительном направлении. Т.е. $v_t(z)=v(z)e^{\pi it}$.

(b) $v_t(z)$ получается из $v(z)$ поворотом на $\pi t/2$ в положительном направлении. Т.е. $v_t(z)=v(z)e^{\pi it/2}$.

\smallskip
{\bf \ref{vfany}.} (a,b) Формула  $v_t(x)=v(tx)$ определяет гомотопию между произвольным векторным полем $v$ и постоянным.

(c) Используйте индукцию. Шаг --- отбрасывание висячего ребра.

(d) Формула  $v_t(x)=v(tx)$ определяет гомотопию между произвольным векторным полем $v$ и постоянным.

\smallskip
{\bf \ref{vfd2}.}
(a) Докажем, что $(4)\Rightarrow(1)$.
Пусть, напротив, $v$ --- продолжение радиального поля до ненулевого поля на $D^2$.
Определим гомотопию $H_t:S^1\to S^1$ формулой $H_t(z):=v(tz)$, где $z\in S^1$ и $t\in[0,1]$.
Это гомотопия между радиальным $H_1$ и постоянным $H_0$ векторными полями на окружности.
 Противоречие c (4).

Докажем, что $(1)\Rightarrow(4)$.
Пусть, напротив, $H_t$ --- гомотопия между радиальным $H_1$ и постоянным $H_0$ векторными полями на окружности.
Определим векторное поле на $D^2$ формулой $v(tz):=H_t(z)$, где $z\in S^1$ и $t\in[0,1]$.
Это продолжение радиального поля.
Противоречие c (1).

(b) Аналогично (a).

\smallskip
{\bf \ref{vfequi}.} (c)
Положим
$H(x,t)=\begin{cases} H_1(x,2t)    &0\leqslant t\leqslant 1/2 \\
          H_2(x,2t-1)  &1/2\leqslant t\leqslant 1\end{cases},$
где $H_1$ и $H_2$ --- гомотопии, соединяющие $u$ и $v$, $v$ и $w$, соответственно.

\smallskip
{\bf \ref{vfunit}.} (a) $v_t(z):=\dfrac{v(z)}{t|v(z)|+1-t}$.

\smallskip
{\bf \ref{vfann}.} (a) $r(z):=z/|z|$.

(b) Следует из (a).

\smallskip
{\bf \ref{vfbor}.} (a,b) Нарисуйте стакан.

(c) Индукция по размеру подграфа.

(d) Можно взять $N=[0,1]\times0$, $A=(0,1]\times0$, $f(x):=0$ и $f_t(x):=t\sin(1/x)$.

\begin{figure}[h]
\centering\hskip-6cm\special{em:graph 7-12}
\vskip3.7cm
\caption{Векторные поля, симметричные относительно касательной}
\label{7-12}
\end{figure}

\bigskip
{\bf \ref{vfsur}.} Ответ: (a) $n$.

(b) Из рисунка \ref{7-12} получаем $\varphi_2+\varphi_1=2\alpha+\pi$.

\smallskip
{\bf \ref{vfcov}.} Ответы: (a) $\t s_k(t):=t+2\pi k$, $k\in\Z$. \quad (b) $\t s_k(t):=nt+2\pi k$, $k\in\Z$.

(с) $\Delta_k=\t s_v(k/m)-\t s_v((k-1)/m)$.

(d) Для каждого $x\in [0,1]$ имеем $\t s(x)-\t{s'}(x)=2\pi n(x)$ для некоторого целого $n(x)$.
Функция $n:[0,1]\to\R$ непрерывна как разность непрерывных функций.
Это и $n(0)=0$ влечет $n(x)=0$ для любого $x$.
Т.е. $s=\t{s'}$.

Другое решение: рассмотрите $\sup\{x\in[0,1]\ :\ \t s=\t s'\text{ на }[0,x]\}$.

\smallskip
{\bf \ref{vfext}.} (a) Аналогично определению степени.
Ввиду равномерной непрерывности отображения $s$ существуют $M>0$, для которого $|s(x)-s(y)|<1$ при $x-y<1/M$.
Для $t\in[0,1]$ и $k=1,2,\dots,n$ обозначим через $\Delta_k(t)\in(-1,1)$ угол, на который нужно повернуть против часовой стрелки точку $s(t(k-1)/n)$ окружности, чтобы получить точку $s(tk/n)$.
Иными словами, $\Delta_k(t)\in(-1,1)$ однозначно определяется из условия $s(tk/M)=s(t(k-1)/M)e^{i\Delta_k(t)}$.
Определим $\t s(t):=\t x+\sum\limits_{k=1}^M \Delta_k(t).$
Проверьте, что $\t s$ --- искомое поднятие.

Идея более короткого, но менее конструктивного решения: рассмотрите
$\sup\{x\in[0,1]\ :\text{ существует поднятие на }[0,x]\}$.
См. детали в \cite[стр. 60-61]{An03}.

Идея еще более короткого, но еще менее конструктивного решения: для каждого $z\in S^1$ постройте
{\it гомеоморфизм} $h:p^{-1}(S^1-z)\to (S^1-z)\times(0,1)$, для которого $p|_{p^{-1}(S^1-z)}=\pro_1\circ h$.
Этот подход особенно полезен в (b).

(b) Рассмотрите $\sup\{t\in[0,1]\ :\text{ существует поднятие на }[0,t]\}$.
Ср. \cite[стр. 61-63]{An03}.

(c) Следует из (b).

(d) Следует из (c).

\smallskip
{\bf \ref{vfalg}.}
(a) Следует из задачи \ref{vfext}.(c) аналогично \ref{vfd2}.a.

(c) Примените теорему Борсука-Улама.

(k) Для каждого $x\in S^2$ существует единственная плоскость $\alpha_1(x)$,
перпендикулярная вектору $x$ и делящая первый многогранник на две части равного объема.
Обозначим $\varphi_1(x):=x\cdot y$, где $y\in\alpha(x)$.
Получим непрерывное отображение $\varphi_1:S^2\to\R$.
Ясно, что $\varphi_1(-x)=-\varphi_1(x)$.
Аналогично, используя второй и третий многогранники, определим
плоскости $\alpha_2(x),\alpha_3(x)$ и
непрерывные отображения $\varphi_2,\varphi_3:S^2\to\R$.
Применим теорему Борсука-Улама \ref{vf-bu}.1 к отображению
$\varphi:=(\varphi_3-\varphi_1)\times(\varphi_2-\varphi_1):S^2\to\R^2$.
Получим $x_0\in S^2$, для которого $\varphi(x_0)=\varphi(-x_0)$.
Так как $\varphi(-x)=-\varphi(x)$, то $\varphi(x_0)=0$.
Т.е. $\varphi_1(x_0)=\varphi_2(x_0)=\varphi_3(x_0)$.
Значит, $\alpha_1(x_0)=\alpha_2(x_0)=\alpha_3(x_0)$.

(e) 3.

(f) Подробнее см. \cite[\S6]{Pr95}, \cite[\S6]{An03}, \cite[\S24]{BE82}.

\begin{figure}[h]
\centering\hskip-6cm\special{em:graph 6-18}
\vskip3.9cm
\caption{Гомотопия поднятий, неподвижная на концах}
\label{6-18}
\end{figure}

\smallskip
{\bf \ref{vfbas}.} (a,b) Докажите, что поднятия гомотопны неподвижно на концах.

(с) Выведите из (b) аналогично \ref{vfd2}.b.

\smallskip
{\bf \ref{vfsum}.} (a) \cite[Теорема 6.5]{Pr95}. \qquad
(b) \cite[Теорема 7.3]{Pr95}.

\smallskip
{\bf \ref{vfgr}.}
Ответ:
$\Z^{n-1}$.
Отображение сужения $V(A)\to V(S^1)$ является взаимно-однозначным соответствием по задаче \ref{vfann}.
Используйте двумерный аналог задачи \ref{vfbor}.

\smallskip
{\bf \ref{vfhom}.} По лемме о поднятии гомотопии поднятия существуют.
Докажите, что они гомотопны неподвижно на граничной окружности.

\smallskip
Если приводимые здесь указания к доказательствам недостаточны, обращайтесь к
\cite[\S\S2,6.1,6.3,8.8,8.9]{Pr04}.

\newpage
\section{Векторные поля на двумерных многообразиях}\label{0ve2ma}

\subsection{Векторные поля на сфере и торе}

{\it Касательным векторным полем} на подмножестве сферы $S^2$ называется семейство касательных
к нему векторов $v(x)$ в его точках $x$, непрерывно зависящих от точки $x$.

(Расстоянием между касательными векторами $AA_1$ и $BB_1$ к сфере $S^2$ в точках $A,B\in S^2$
называется $|AB|+|A_1B_1|$.
Теперь непрерывность отображения из $S^2$ в множество всех касательных векторов к $S^2$ определяется
аналогично началу \S\ref{0vepl}.)

\begin{figure}[h]
\centering{\hskip0cm\special{em:graph 7-1}\hskip4cm \special{em:graph 7-2}}
\vskip2.3cm
\caption{Касательные векторные поля на сфере}
\label{7-1i7-2}
\end{figure}

\begin{pr}\label{vfs2}
(a) Постройте касательное векторное поле на сфере $S^2$, у которого вектор нулевой только в одной точке.

(b) {\bf Теорема о еже.} {\it На стандарной сфере не существует ненулевого касательного векторного поля.}
\end{pr}

Два ненулевых касательных векторных поля называются {\it гомотопными}, если одно можно получить из другого
непрерывной деформацией, в процессе которой векторное поле остается ненулевым и касательным.
Такой деформацией называется семейство $v_t$ ненулевых касательных векторных полей, непрерывно зависящее
от параметра $t$, для которого $v_0$ есть первое поле и $v_1$ есть второе поле.

\begin{pr}\label{vfhemi}
Любые два ненулевых касательных векторных поля на верхней полусфере
$D^2_+:=\{(x,y,z)\in\R^3\ |\ x^2+y^2+z^2=1, \ z\ge0\}$ гомотопны.
\end{pr}

Понятие {\it касательного векторного поля} на рассматриваемых ниже поверхностях вводится
дословно аналогично случаю сферы.
{\it Будем называть единичное касательное векторное поле просто полем}.

{\it Гомотопность} полей на поверхности определяется дословно так же, как и на сфере.

Для приложений важно уметь описывать множество $V(N)$ {\it единичных касательных векторных полей с точностью до
гомотопности} на поверхности $N$.
Встречаются и более сложно формулируемые проблемы, для решения которых полезно сначала научиться описывать
множество $V(N)$.

\begin{pr}\label{vfclex}
Опишите $V(N)$ для \qquad

(a) сферического слоя $\{(x,y,z)\in\R^3\ |\ x^2+y^2+z^2=1, \ -\frac12\le z\le \frac12\}$.

(b) боковой поверхности $\{(x,y,z)\in\R^3\ |\ x^2+y^2=1, \ 0\le z\le1\}$ цилиндра.
\end{pr}

\begin{figure}[h]\centering
\includegraphics{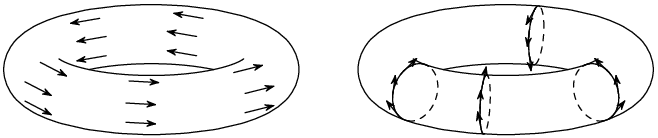}
\caption{`Параллельное' и `меридиональное' касательные векторные поля на торе}\label{vector}
\end{figure}

 Примеры полей на торе --- `параллельное' и `меридиональное' (рис.~\ref{vector}).

\begin{pr}\label{vftor}
(a) `Параллельное' поле на торе гомотопно `меридиональному'.

(b) Любое поле $v$ на на торе гомотопно полю $-v$.

(c,d) Сформулируйте и докажите аналог задачи \ref{vfunit} для тора.

(e) Приведите пример двух негомотопных полей на торе.
\end{pr}

Для доказательства негомотопности полей на торе необходим следующий инвариант.
Поставим в соответствие полю $v$ на торе
{\it число $r_p(v)$ оборотов вектора этого поля при обходе по параллели тора}.
Но количество оборотов определено только для обхода по замкнутой кривой на плоскости (а не в пространстве).
Поэтому для его определения нужно непрерывно отождествить касательные плоскости в разных точках тора
(подумайте, что это значит и как это сделать).
Вместо этого представим поле на торе в виде склейки единичного векторного поля на квадрате, лежащем в плоскости.
Точнее, воспользуемся без доказательства наличием взаимно-однозначного соответствия между $V(T^2)$
и множеством таких единичных векторных полей на квадрате $[0,1]^2\subset\R^2$, что $v(x,0)=v(x,1)$ и $v(0,x)=v(1,x)$
для любого $x$, с точностью до гомотопии в классе таких полей.
Для такого поля на квадрате определим $r_p(v):=\deg v|_{[0,1]\times0}$.
Ср. с задачами \ref{vfs2}.b, \ref{vfhemi} и \ref{vfclex}.
Наличием такого соответствия (и аналогичного соответствия для других поверхностей) можно и далее пользоваться без
доказательства.
(Заметим, что в \cite[\S7]{Pr95} оно используется даже без явной формулировки.)

Аналогично определяется {\it число $r_m(v)$ оборотов вектора поля при обходе по меридиану тора}.
Аналогичные инварианты можно получить, рассматривая другие замкнутые кривые на торе.


\begin{pr}\label{vfcltor}
(a) Обозначим через $T^2_0$ тор с дыркой, т.е. пересечение тора $T^2$ с полупространством $x\le 5/2$.
Отображение $r_p\times r_m:V(T^2_0)\to\Z\times\Z$ являются взаимно-однозначным соответствием.

(b) {\bf Теорема классификации векторных полей на торе.}
{\it Отображение $r_p\times r_m:V(T^2)\to\Z\times\Z$ является взаимно-однозначным соответствием.}

(c) Если на торе с дыркой заданы два поля, то их сужения на краевую окружность этой дырки гомотопны.
\end{pr}

Более подробное изложение приводится в~\cite[\S7]{Pr95}.

\subsection{Векторные поля на других поверхностях}

В каждой точке сферы с ручками (\S\ref{0dile}) имеется касательная плоскость.

Для решения следующей задачи полезно `изображение' поверхности $S_{g,0}$ в виде диска с неперекрученными ленточками
(\S\ref{0dile}; аналогично рис. \ref{dilen}, \ref{netor}).

\begin{pr}\label{vfhole}
(a) На $S_{g,0}$ существует поле.

(b) Опишите $V(S_{g,0})$.

(c) Если на $S_{g,0}$ заданы два поля, то их сужения на краевую окружность дырки гомотопны.
\end{pr}


\begin{pr}\label{vfepart}
{\bf Теорема Эйлера-Пуанкаре (частный случай).}
{\it Среди сфер с ручками только торе имеет ненулевое касательное векторное поле.}
\end{pr}


\begin{pr}\label{vfmob}
(a) Постройте единичное касательное векторное поле на листе Мебиуса $M$.

(b) Гомотопно ли построенное Вами поле полю $-v$?

(c) Опишите $V(M)$ для листа Мебиуса $M$. \qquad

(d) Верно ли, что сужения на граничную окружность листа Мебиуса $M$ любых двух единичных касательных векторных полей
на листе Мебиуса $M$ гомотопны?
\end{pr}

Для решения этой задачи представьте единичное касательное векторное поле на листе Мебиуса  в виде склейки
единичного векторного поля на квадрате, лежащем в плоскости.
Точнее, воспользуйтесь без доказательства наличием взаимно-однозначного соответствия между $V(M)$ и множеством
таких единичных векторных полей на квадрате $[0,1]^2\subset\R^2$, что $v(0,y)$ получено из
$v(1,1-y)$ симметрией относительно оси $Ox$ для любого $y$, с точностью до гомотопии в классе таких векторных полей.

\begin{pr}\label{vfklein}
(a) На бутылке Клейна $K$ существует единичное касательное векторное поле.

(b) Опишите $V(K)$ для бутылки Клейна $K$.
\end{pr}

Для решения этой задачи представьте единичное касательное векторное поле на бутылке Клейна  в виде склейки
единичного векторного поля на квадрате, лежащем в плоскости.
Точнее, воспользуйтесь без доказательства наличием взаимно-однозначного соответствия между $V(K)$ и множеством
таких единичных векторных полей на квадрате $[0,1]^2\subset\R^2$, что $v(x,0)=v(x,1)$ и $v(0,x)$ получено из
$v(1,1-x)$ симметрией относительно оси $Ox$ для любого $x$, с точностью до гомотопии в классе таких векторных полей.

\subsection{Поля направлений}

{\it Полем направлений} на поверхности называется семейство касательных к ней прямых $l(x)$ в точках
$x$, непрерывно зависящих от точки $x$.
Поля направлений связаны с {\it лоренцевыми метриками} на поверхностях, возникающими в физике~\cite{Ko01}.

\begin{pr}\label{vfdir}
(a) Если на поверхности существует ненулевое касательное векторное поле, то существует и поле направлений.

(b) Обратное тоже верно.
\end{pr}


\begin{pr}\label{vfcldir} Определите гомотопность полей направлений и классифицируйте поля направлений
с точностью до гомотопности на

(a) окружности в плоскости; \quad (b) торе; \quad (c) листе Мебиуса; \quad (d) бутылке Клейна.
\end{pr}

\subsection{Обобщение на двумерные многообразия}

Пусть $D$ --- декартово произведение интервалов (конечных или бесконечных) на прямой.
{\it Гладкой регулярной параметризованной двумерной поверхностью} называется бесконечно дифференцируемое отображение
$r:D\to\R^m$ (или, что то же самое, упорядоченный набор $m$ отображений $x_1,x_2\dots,x_m:D\to\R$), производная которого 
невырождена в любой точке (т.е. пара векторов $r_u(u,v),r_v(u,v)$ линейно независима при любых $(u,v)\in D$).
Далее слова {\it гладкая регулярная} опускаются.

{\it Двумерным гладким подмногообразием в $\R^m$} (или двумерной гладкой непараметризованной поверхностью) называется 
подмножество $N\subset\R^m$, для любой точки $x\in N$ которого существует такая ее замкнутая окрестность $Ox$ в $\R^m$, 
что $N\cap Ox$ является образом $r(D)$ инъективной параметризованной поверхности $r:D\to\R^m$.
Двумерные гладкие подмно\-го\-об\-ра\-зия мы будем коротко называть 2-многообразиями.
(Чтобы не было путаницы, то, что названо 2-многообразием в \S\ref{02man}, будем называть {\it кусочно-линейным
2-многообразием}.)

\begin{pr}\label{vfsubexa}
Нарисуйте следующие подмножества в $\R^3$.
Докажите, что они являются двумерными гладкими подмногообразиями в $\R^3$.

(a) Квадрат со стороной $a$ на плоскости.  

(b) Боковая поверхность прямого кругового цилиндра радиуса $R$ и высоты $h$ (рис. \ref{tor}).

(c) Боковая поверхность прямого кругового усеченного конуса с радиусами оснований $R,r$ и высотой $h$.

(d) Сфера радиуса $R$.

(e) Неформально, {\it тор} --- поверхность бублика, рис.  \ref{tor}  слева.
Формально, {\it стандартным тором с радиусами $R$ и $r$} называется фигура, образованная вращением окружности $(x-R)^2+y^2=r^2$ вокруг оси $Oy$.
Эта фигура получена из (двумерного) квадрата склейкой  пар его противоположных сторон
`с одинаковыми направлениями', т.е. без поворота.

(f) 
{\it Стандартной лентой Мебиуса} называется поверхность в $\R^3$, заметаемая стержнем длины 2,
равномерно вращающимся относительно своего центра, при равномерном движении этого центра
по окружности радиуса 9, при котором стержень делает пол-оборота, рис. \ref{tor}.
Эта фигура получена из длинной прямоугольной полоски склейкой двух ее противоположных сторон `с противоположным направлением', т.е. с поворотом на $180^\circ$.

(h) Поверхность вращения графика бесконечно дифференцируемой положительной функции $f:\R\to\R$.

(i) Седлообразная поверхность (гиперболический параболоид) $z=xy$.

(j) Двуполостный гиперболоид $z^2=x^2+y^2+1$.

(k) Однополостный гиперболоид $z^2=x^2+y^2-1$.

(l) График бесконечно дифференцируемой функции $f:D\to\R$.

(m) Прообраз нуля при бесконечно дифференцируемой функции  $f:\R^3\to\R$, производная (=градиент) которой ненулевая в каждой точке.
\end{pr}

\begin{pr}\label{vfsub34}
(a) Любое 2-многообразие в $\R^3$ становится 2-многообразием в $\R^4$, если рассмотреть $\R^3$ как подмножество в $\R^4$.

(b) Рассмотрим в $\R^4$ окружность $x^2+y^2=1$ и семейство ее нормальных трехмерных плоскостей.
{\it Бутылкой Клейна} (стандартной) $K$ называется поверхность в $\R^4$, заметаемая окружностью $\omega$, центр
которой равномерно описывает окружность $x^2+y^2=1$, а окружность $\omega$ в то же время равномерно поворачивается на 
угол $\pi$ (поворачивается в движущейся нормальной трехмерной плоскости относительно своего диаметра, движущегося 
вместе с нормальной трехмерной плоскостью).
См. рис. \ref{kleinb}.
Эта фигура получена из прямоугольной полоски такой склейкой ее пар противоположных сторон, при которой одна пара 
склеивается `с одинаковым направлением', а другая `с противоположным направлением', рис. \ref{kleinb}.

(c) Образ отображения $S^2\to\R^4$, заданного формулой $(x,y,z)\mapsto(x^2,xy,yz,zx)$, является 2-многообразием.
(Это `реализация' проективной плоскости $\R P^2$.)
\end{pr}

\begin{pr}\label{vfsubnot}
Не являются 2-многообразиями ни объединение двух (или трех) координатных плоскостей в $\R^3$,
ни конус $z^2=x^2+y^2$  в $\R^3$.
\end{pr}

{\it Клеточным разбиением 2-многообразия} называется такое разбиение 2-многообразия на криволинейные треугольники,
что любые два треугольника либо не пересекаются, либо пересекаются по вершине, либо пересекаются по ребру.

\smallskip
{\bf Теорема триангулируемости.}  \cite[Теорема 10.6 из дополнения]{MS74}
{\it (a) Для любых 2-многообразия и $\varepsilon>0$ существует
триангуляция 2-многообразия, диаметр каждого треугольника в которой меньше $\varepsilon$.

(b) Любые две триангуляции одного 2-многообразия кусочно-линейно гомеоморфны (в смысле \S\ref{02man}).}

\smallskip
Мы будем сокращенно называть клеточные разбиения 2-многообразия {\it разбиениями}.
{\it Эйлеровой характеристикой} $\chi(N)$ 2-многообразия $N$ называется эйлерова характеристика произвольного его
разбиения.
Это определение корректно ввиду теоремы триангулируемости (и аналога задачи \ref{4-eul}.b).

{\it Касательной плоскостью} $T_xN$ к 2-многообразию $N$ в точке $x$ называется образ плоскости $\R^2$ при 
производной в точке $r^{-1}(x)$ некоторого отображения $r:D\to\R^m$ из определения 2-многообразия.

\begin{pr}\label{vfweltan}
(a) Этот образ есть плоскость, проходящая через $P$ и содержащая векторы $r_u(r^{-1}(x))$ и $r_v(r^{-1}(x))$.

(b) Это определение корректно, т.е. не зависит от выбора отображения $r$.
\end{pr}

Теперь касательные векторные поля на 2-многообразиях определяются
аналогично случаю сферы.
Касательное векторное поле называется {\it ненулевым}, если все его векторы ненулевые.

Если для любой точки $x$ 2-многообразия $N$ существуют $Ox$ и $r:[0,1]^2\to\R^m$ (из определения 2-многообразия),
для которых $x\in r((0,1)^2)$, то $N$ называется {\it замкнутым}.
Множество всех тех точек $x$ 2-многообразия $N$, для которых таких $Ox$ и $r$ нет, называется {\it краем}
мно\-го\-об\-ра\-зия $N$.

\textit{Примеры} замкнутых 2-многообразий: вышеописанные стандартные сфера, тор и сфера с $g$ ручками в $\R^3$, а также
бутылка Клейна в $\R^4$.
 
\textit{Примеры} 2-многообразий с непустым краем: кольцо, цилиндр, лента Мебиуса, тор с дыркой.

\begin{pr}\label{vfexi}
На любом связном 2-многообразии
с непустым краем существует ненулевое касательное векторное поле.
\end{pr}

{\bf Теорема Эйлера-Пуанкаре.} {\it На замкнутом 2-многообразии имеется ненулевое касательное векторное поле
тогда и только тогда, когда эйлерова характеристика этого 2-многообразия нулевая.}

\smallskip
Этот результат доказан в следующих двух пунктах.

\smallskip
{\bf Замечание.}
Взаимно-однозначное соответствие (т.е. биекция) $f:N\to M$ между замкнутыми 2-многообразиями $N\subset\R^n$ и
$M\subset\R^m$ называется {\it диффеоморфизмом}, если для любой точки $x\in N$ найдутся окрестности $Ox\subset\R^n$
точки $x$ и $Of(x)\subset\R^m$ точки $f(x)$, а также отображения $r:[0,1]^2\to\R^n$ и $q:[0,1]^2\to\R^m$
из определения замкнутого 2-многообразия, для которых
$fr=q$.
Замкнутые 2-многообразия $N\subset\R^n$ и $M\subset\R^m$ называются {\it диффеоморфными}, если существует диффеоморфизм
$f:N\to M$.
Соответствие между классами диффеоморфности 2-многообразий и классами кусочно-линейной гомеоморфности кусочно-линейных
2-многообразий, определяемое триангулируемостью, корректно определено и является взаимно однозначным.
Поэтому эйлерова характеристика 2-многообразия нулевая тогда и только тогда, когда это  2-многообразие диффеоморфно
тору или бутылке  Клейна.
Впрочем, проверять эту диффеоморфность проще всего именно через эйлерову характеристику.


\comment

Подмножество $N\subset\R^m$  называется \emph{двумерным гладким подмно\-го\-об\-ра\-зием},
если для любой точки $x\in N$ найдутся ее окрестность $U$ и диффеоморфизм $U\to\R^2\hm\times \R^{m-2}$, переводящий
пересечение $U\cap N$ в $\R^2\times0$ или в $\R\times[0,+\infty)\times0$.
Если для всех $x$ имеет место первый случай, то $N$ называется {\it замкнутым}.
Множество всех тех точек $x$, для которых имеет место второй случай, называется {\it краем} подмно\-го\-об\-ра\-зия $N$.

{\it Касательной плоскостью} к 2-многообразию называется прообраз плоскости $\R^2\times0$ при производной
в точке $x$ диффеоморфизма из определения 2-многообразия.

Пусть $n<m$ и $N\subset\R^n$, $M\subset\R^m$ --- замкнутые 2-многообразия.
Обозначим через $i:\R^{n-2}=\R^{n-2}\times0\subset \R^{n-2}\times\R^{m-n}=\R^{m-2}$ включение.
Взаимно-однозначное соответствие (т.е. биекция) $f:N\to M$ называется {\it диффеоморфизмом}, если для любой точки
$x\in N$ найдутся окрестности $U\subset\R^n$ точки $x$ и $V\subset \R^m$ точки $f(x)$, а также диффеоморфизмы
$p:U\to\R^2\times \R^{n-2}$ и $q:V\to\R^2\times \R^{m-2}$ из определения замкнутого 2-многообразия, для которых
$f(U\cap N)\subset V$ и $qf=(\id \R^2\times i)p$ на $V$.

\endcomment

\subsection{Касательные векторные поля: общее положение}


\smallskip
{\it Набросок определения числа Эйлера $e(N)$ для замкнутого 2-многообразия $N\subset\R^m$.}
Возьмем на $N$ {\it гладкое} касательное векторное поле {\it общего положения}, у которого допускаются нули.
{\it Гладкость} означает существование производной поля в любой точке $x\in N$.
Производная в точке $x\in N$ является линейным оператором $T_xN\to T_xN$ из касательного пространства $T_xN$ в точке $x$ к $N$ в себя; производная  представляется матрицей $2\times2$.
{\it Общность положения} означает, что в каждой точке, где вектор нулевой, производная невырождена.
Существование такого поля фактически доказано в \cite[Ч. II, \S13]{DNF79}, \cite{Pr04}.
Из общности положения вытекает конечность нулей поля.
{\it Числом Эйлера} $e(N)$ называется сумма знаков определителей производной в нулях поля.
(Для определения знака нужна ориентация в $T_xN$, но при обращении ориентации получается такой же знак.)

\begin{figure}[h]\centering
\includegraphics{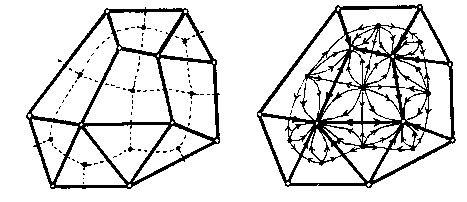}
\caption{Построение векторного поля общего положения по разбиению на многоугольники}
\label{eul}
\end{figure}

 \begin{figure}[h]\centering
\centering\hskip-6cm\special{em:graph 7-6}
\vskip2.7cm
\caption{Векторное поле в окрестности вершины, особых точек внутри грани  и на ребре}
\label{7-6}
\end{figure}

\smallskip
{\it Примеры полей общего положения:}

$\bullet$ ненулевое касательное векторное поле (если оно существует, то $e(N)=0$).

$\bullet$ векторное поле на рис.~\ref{water} (получаем $e(S_g)=2-2g$).

$\bullet$ векторное поле на рис.~\ref{eul},~\ref{7-6} (получаем $e(N)=E-V+F=\chi(N)$).

\smallskip
{\it Набросок доказательства корректности этого определения.}
Возьмем {\it гладкую} гомотопию {\it общего положения} между двумя касательными векторными полями $v$ и $v'$
общего положения на $N$.
Ее можно представлять себе как касательное векторное поле общего положения на $N\times I\subset\R^m\times I$, векторы которого параллельны гиперплоскости $\R^3\times0$.
{\it Гладкость} означает существование производной поля в любой точке $(x,t)\in N\times I$.
Производная в точке $(x,t)\in N\times I$ является линейным оператором $T_xN\times\R\to T_xN$; производная  представляется матрицей $3\times2$.
{\it Общность положения} означает, что в любой точке, в которой вектор нулевой, производная поля имеет ранг 2.
Существование такой гомотопии фактически доказано в \cite[Ч. II, \S13]{DNF79}, \cite{Pr04}.
Из общности положения вытекает, что множество нулей является несвязным объединением отрезков и окружностей.
На этих отрезках и окружностях можно ввести ориентацию (подумайте, как, и сравните с \cite[18.1]{Pr04},
\cite[Ч. II, \S13]{DNF79}).
Если один из этих отрезков соединяет {\it разные} основания $N\times0$ и $N\times1$, то в соответствующих нулях
полей $v$ и $v'$ определители производных имеют одинаковый знак, а если {\it одинаковые}, то разный.
Поэтому $e(v)=e(v')$.

\smallskip
Необходимость в теореме Эйлера-Пуанкаре и равенство $\chi(N)=e(N)$ доказаны выше.
Для доказательства достаточности нужно при $e(N)=0$ получить ненулевое касательное векторное поле из касательного
векторного поля общего положения `сокращением' точек разных знаков аналогично \cite[\S18.3]{Pr04}.

\subsection{Касательные векторные поля: триангуляции}

\smallskip
Приведем определение числа Эйлера через триангуляции и соответствующее доказательство теоремы Эйлера-Пуанкаре
для ориентируемых 2-многообразий.
Оно похоже на \cite[\S14]{BE82}, ср. \cite[\S18]{Pr04}.
Оно может быть обобщено на случай нормальных полей на 2-многообразиях и на случай многомерных многообразий.

Вообще, роль теории препятствий состоит в {\it сведении} топологических задач на
{\it произвольном} многообразии к похожим задачам для {\it простейших, модельных} многообразий.
Важно, что для применения теории препятствий можно воспользоваться {\it результатом}
решения этих простейших задач, не вникая в его доказательство.
Указанные простейшие задачи могут решаться, в частности, средствами теории препятствий.
Мы будем использовать доказанные ранее теоремы продолжаемости и гомотопности.

\begin{figure}[h]\centering
\includegraphics{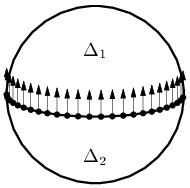}
\caption{Начало построения векторного поля на 2-многообразии}\label{vesph}
\end{figure}

Слово `касательное поле' будет означать `ненулевое касательное векторное поле'.!!!
Для сферы следующие рассуждения иллюстрируются рисунком~\ref{vesph}.


\smallskip
{\it Начало доказательства теоремы Эйлера-Пуанкаре.}
Возьмем некоторое разбиение $U$ 2-многообразия  $N$ на многоугольники.
Обозначим через $U^*$ двойственное разбиение (см. определение в \S\ref{02hom}).
Выберем $U$ настолько мелким, чтобы касательные плоскости в любых двух точках
любой грани разбиения $U^*$ не были бы ортогональны.

Очевидно, что можно построить поле на $U_0^*$.
Ясно, что это поле единственно (с точностью до непрерывной деформации в классе полей на $U_0^*$).
Поэтому существование поля на $N$ равносильно продолжаемости построенного поля с $U_0^*$ на $N$.

\smallskip
{\it Построение препятствующей расстановки.}
Возьмем произвольное ребро $e$ разбиения $U^*$.
Ввиду мелкости разбиения $U$ ортогональная проекция касательной плоскости в произвольной
точке этого ребра на касательную плоскость $\alpha_e$ в некоторой фиксированной точке этого ребра
переводит ненулевые векторы в ненулевые.
Значит, касательные плоскости в разных точках этого ребра можно отождествить с одной плоскостью $\alpha_e$.
Если на плоскости лежит отрезок, и в его концах заданы ненулевые вектора (лежащие в плоскости),
то это поле из двух векторов можно продолжить до поля ненулевых векторов на всем отрезке.
Поэтому построенное на $U_0^*$ поле можно продолжить на $U_1^*$.
Заметим, что такое продолжение неоднозначно.
Обозначим полученное на $U_1^*$ поле через $v$.
Пример такого поля для сферы приведен на рис.~\ref{vesph}.

Попробуем теперь продолжить поле $v$ с $U_1^*$ на все $N$.
Фиксируем ориентацию 2-многообразия $N$.
Возьмем произвольную грань $\Delta$ разбиения $U^*$.
Ввиду мелкости разбиения $U$ ортогональная проекция касательной плоскости в произвольной
точке этой грани на касательную плоскость $\alpha_\Delta$ в некоторой фиксированной точке этой грани
переводит ненулевые векторы в ненулевые.

Ориентация грани $\Delta$ порождает направление на ее граничной окружности.
При обходе этой окружности вдоль этого направления ортогональная проекция
вектора нашего поля на $\alpha_\Delta$ повернется на некоторое целое число оборотов.
Поставим полученное число оборотов в вершину исходного разбиения $U$, лежащую в диске $\Delta$.
Для случая сферы (рис.~\ref{vesph}) в обоих вершинах  будут стоять единицы.
Полученная расстановка целых чисел в вершинах называется {\it препятствующей} и
обозначается $\varepsilon(v)$.
По теореме о продолжаемости

{\it продолжение поля с $U_1^*$ на $N$ возможно тогда и только тогда, когда $\varepsilon(v)=0$.}

\smallskip
{\it Определение различающей расстановки.}
Если $\varepsilon(v)\ne0$, то поле $v$ не продолжается на $N$, но еще не все потеряно:
можно попытаться так изменить поле $v$ на $U_1^*$,
чтобы препятствующая расстановка стала равной нулю.
Для этого выясним, как $\varepsilon(v)$ зависит от $v$.

\begin{figure}[h]\centering
\includegraphics{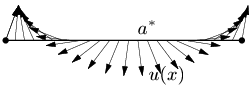}
\caption{Подкручивание векторного поля на один оборот}\label{vecdif}
\end{figure}

Различие между полями $v$ и $u$ на $U_1^*$, совпадающим на $U_0^*$, можно измерять (и задавать) так.
Фиксируем ориентацию каждого ребра разбиения $U$.
Ориентируем каждое ребро $a^*$ разбиения $U^*$ так, чтобы стрелки на $a$ и на $a^*$ (в этом порядке)
образовали `положительный базис на $N$'.
Пусть точка $x$ движется по ребру $a^*$ вдоль направления, а потом обратно.
При движении 'туда' будем рассматривать вектор $u(x)$, а при движении 'обратно' --- вектор $v(x)$.
Поставим на ориентированном ребре $a$ разбиения $U$ число оборотов
ортогональной проекции рассматриваемого вектора на $\alpha_a$ при этом движении.
(Например, для поля $u$ на рис.~\ref{vecdif} и поля $v$, направленного
вертикально вверх, на ребре стоит $-1$.)
Полученную расстановку целых чисел на ориентированных ребрах разбиения $U$ назовем {\it различающей} и обозначим $d(v,u)$.

\smallskip
{\it Изменение препятствующей расстановки и определение границы ребра.}
При изменении поля на ребре $a^*$ 'на $+1$ оборот' к $\varepsilon(v)$
прибавляется расстановка $+1$ в начале ребра $a$ и $-1$ в его конце
(и 0 на всех остальных вершинах).
Эта расстановка называется {\it границей ребра $a$} и обозначается $\partial a$.
Если
$$
d(v,v')=n_1a_1+\dots+n_sa_s,\quad\text{то}\quad
\varepsilon(v)-\varepsilon(v')=n_1\partial a_1+\dots+n_s\partial a_s.
$$

{\it Завершение доказательства.}
Из последней формулы вытекает, что {\it сумма} чисел препятствующей расстановки не зависит от $v$.
А из рис.~\ref{eul} можно понять, что эта сумма равна $\chi(N)$.
QED

\smallskip
Приведем другое завершение доказательства.
Оно длиннее предыдущего, но является иллюстрацией общего подхода теории препятствий.

\smallskip
{\it Определение группы $H_0(U;\Z)$, класса $e(U)$ и другое завершение доказательства.}
Назовем {\it границей} сумму $n_1\partial a_1+\dots+n_s\partial a_s$ границ
нескольких ребер с целыми коэффициентами.
Назовем расстановки $\varepsilon_1$ и $\varepsilon_2$ {\it гомологичными},
если $\varepsilon_1-\varepsilon_2$ есть граница.
Ясно, что

(i) При изменении поля $v$ на $U_1^*$ препятствующая расстановка $\varepsilon(v)$
заменяется на гомологичную расстановку.

(ii) Если $\varepsilon(v)$ является границей, то можно так изменить $v$ на $v'$
(на $U_1^*$, не меняя на $U_0^*$), чтобы получилось $\varepsilon(v')=0$.

(iii) Гомологичность является отношением эквивалентности на множестве расстановок.

(Для доказательства утверждения (ii) заметим, что если
$\varepsilon(v)\hm=n_1\partial a_1+\dots+n_s\partial a_s$, то
можно взять поле $v'$ на $U_1^*$, для которого
$d(v,v')=n_1a_1+\dots+n_sa_s$, тогда получим $\varepsilon(v')=0$.)

{\it Нульмерной группой гомологий $H_0(U;\Z)$ разбиения $U$ (с коэффициентами в $\Z$)}
называется группа расстановок с точностью до гомологичности.

{\it Классом Эйлера} разбиения $U$ называется класс гомологичности препятствующего цикла:
$$e(U):=[\varepsilon(v)]\in H_0(U;\Z).$$
Это определение корректно ввиду утверждения (i).

Теперь теорема Эйлера-Пуанкаре вытекает из утверждений (i) и (ii) вместе с нижеследующей задачей \ref{vfta}.cd.
 QED

\begin{pr}\label{vfta}
Возьмем произвольное разбиение $U$ замкнутого ориентируемого 2-многообразия $N$.

(a) Определите группу $H_0(U;\Z)$ независимо от рассуждений, в которых она появилась.

(b) Группа $H_0(U;\Z)$ и класс $e(U)$ не зависят от выбора $U$ для фиксированной $N$.

(c) $H_0(U;\Z)\cong\Z$.

(d) $\lessmskips{3mu}e(U)\hm=\chi(N)$.
\end{pr}

\begin{pr}\label{vfnonor} Докажите теорему Эйлера-Пуанкаре для неориентируемых 2-многообразий.
\end{pr}

\subsection{Нормальные векторные поля для двумерных многообразий}

\smallskip
{\it Единичным нормальным векторным полем} к окружности $S^1\subset\R^2\subset\R^3$ в трехмерном пространстве
называется семейство нормальных к ней единичных векторов $v(x)$ в точках $x$
окружности, непрерывно зависящих от точки $x$.
Будем называть единичное нормальное векторное поле просто {\it нормальным полем}.

Понятие {\it гомотопности} нормальных полей вводится дословно так же, как и для
ненулевых касательных векторных полей.

\begin{pr}\label{vfnos1}
(a) Постройте нормальное поле.

(b) Постройте нормальное поле, не гомотопное уже построенному.

(c) Опишите нормальные поля с точностью до гомотопности.

(d) Опишите нормальные поля с точностью до гомотопности для {\it заузленной} гладкой окружности в $\R^3$
(т.е. на связном замкнутом 1-подмногообразии в $\R^3$; его определение аналогично вышеприведенному
определению замкнутого 2-подмногообразия в $\R^3$).
\end{pr}

Для 2-многообразия $N\subset\R^m$ {\it единичным нормальным векторным полем} на $N$ называется семейство
единичных нормальных к $N$ (т.е. к касательной плоскости к $N$) векторов $v(x)$ в точках $x\in N$,
непрерывно зависящих от точки $x\in N$.
Будем называть единичное нормальное векторное поле просто {\it нормальным полем}.

\begin{pr}\label{vfnome}
Постройте нормальные поля для

(a) стандартной сферы в $\R^3$;

(b) стандартного листа Мебиуса, рассматриваемого как подмножество в $\R^4$.

(c) стандартной бутылки Кляйна в $\R^4$.
\end{pr}

2-многообразие $N\subset\R^m$ называется {\it ориентируемым}, если существует семейство ориентаций касательных  плоскостей к $N$ в точках $x\in N$, непрерывно зависящих от точки $x\in N$.


\begin{pr}\label{vfnoconcl}
(a) Любое ориентируемое 2-подмногообразие в $\R^3$ имеет нормальное поле.

(b) Никакое неориентируемое 2-подмногообразие в $\R^3$ не имеет нормального поля.

(c) Если ориентируемое 2-многообразие в $\R^4$ имеет нормальное поле,
то оно имеет пару линейно-независимых нормальных полей.

(d) Верно ли (c) без предположения ориентируемости?
\end{pr}


\begin{pr}\label{vfnobou}
Любое 2-многообразие с непустым краем в $\R^4$ имеет нормальное поле.
\end{pr}


{\bf Теорема о нормальных полях.}
{\it Любое ориентируемое 2-многообразие в $\R^m$ имеет нормальное поле.}

\begin{pr}\label{vfnoth}
(a) Докажите теорему о нормальных полях для $m=4$.

(b) Никакое 2-многообразие нечетной эйлеровой характеристики в $\R^4$  не имеет нормального поля.

(c) Существует 2-многообразие четной эйлеровой характеристики в $\R^4$, не имеющее нормального поля.
\end{pr}


Вряд ли у Вас получится решить эту задачу без решения следующих!

\begin{pr}\label{vfnoeu}
Пусть $N\subset\R^4$ --- замкнутое 2-многообразие.

(a) Определите {\it нормальное число Эйлера} $\overline e(N)\in\Z$ для 2-многообразия $N\subset\R^4$
как препятствие к существованию нормального поля.

Указание. Аналогично предыдущему (используя общее положение или триангуляцию), заменяя касательные поля на нормальные.
Корректность этого определения также доказывается аналогично случаю касательных векторных полей.

(b) Пусть $N'\subset\R^4$ --- 2-многообразие в {\it общем положении} с $N$, т.е. $N\cap N'$
есть конечное число точек пересечения, в каждой из которых касательные плоскости к $N$ и к $N'$ порождают все $\R^4$.
Фиксируем ориентацию на $N$ и в $\R^4$.
Для каждой из точек пересечения возьмем положительные базисы $e_1,e_2$ и $e_1',e_2'$
в касательных плоскостях к $N$ и к $N'$.
Если $e_1,e_2,e_1',e_2'$ --- положительный базис пространства $\R^4$, то назовем точку пересечения
положительной, иначе отрицательной.

Докажите, что сумма знаков точек пересечения равна $\overline e(N)$.

(c) Если $N$ ориентируемо, то $\overline e(N)=0$.

(d) $\overline e(N)$ четно.

(e)* $\overline e(N)\equiv2\chi(N)\mod4$.
\end{pr}

\begin{pr}\label{vfnorp2}
(a) Существует замкнутое 2-многообразие в $\R^4$ эйлеровой характеристики 1 (т.е., диффеоморфное $\R P^2$).

(b) Найдите его нормальное число Эйлера.

(с) Для любого $k$ с условиями $-2\mu\le k\le 2\mu$, $k\equiv2\mu\mod4$, существует такое
замкнутое неориентируемое 2-многообразие $N_{\mu,k}\subset\R^4$ эйлеровой характеристики $\mu$ (т.е., диффеоморфное сфере с $\mu$ пленками Мебиуса $N_\mu$), что $\overline e(N_{\mu,k})=k$.
(Заметим, что для $|k|>2\mu$ такого 2-многообразия не существует~\cite{Ma69}.
Удивительно, что доказательство этого факта сложно.)
\end{pr}


\begin{pr}\label{vfnocl}
Опишите нормальные поля в $\R^4$ с точностью гомотопности на стандартных

(a) сфере; \quad (b) торе; \quad (с) сфере с ручками; \quad (d) листе Мебиуса; \quad (e)* бутылке Кляйна.
\end{pr}


\begin{pr}\label{vfnohi}
 Докажите теорему о нормальных полях для $m\ge5$.
\end{pr}

\begin{pr}\label{vfnohic}
 (a)-(e) То же, что в задаче \ref{vfnocl}, с заменой $\R^4$ на $\R^5$.

(f) Любые два нормальных поля на 2-многообразии в $\R^6$ гомотопны.
\end{pr}

\subsection{Препятствие Уитни к вложимости}

Невложимость бутылки Кляйна в $\R^3$  можно доказать  с использованием многомерного аналога теоремы Жордана:
{\it если бы бутылка Кляйна была бы вложена в $\R^3$, то разбивала бы $\R^3$.}
Так доказывается невложимость неориентируемых $n$-многообразий в $\R^{n+1}$.
Но так не получается доказать невложимость $n$-многообразий в $\R^m$ для $m>n+1$.
Для последнего нужны препятствия к `старшей ориентируемости', т.е. препятствия Уитни.
Мы продемонстрируем идею построения препятствия Уитни на примере {\it другого} доказательства невложимости бутылки Кляйна в $\R^3$

{\it Кусочно-линейным вложением $f:K\to\R^m$} называется инъективное
отображение, линейное на каждом треугольнике некоторой триангуляции.

\smallskip
{\it Доказательство несуществования кусочно-линейного вложения бутылки Кляйна $K$ в $\R^3$.}
Возьмем {\it кусочно-линейное} отображение $f:K\to\R^3$ {\it общего положения}.
{\it Кусочно-линейность} означает линейность на каждом треугольнике некоторой триангуляции.
{\it Общность положения} означает, что образы никаких четырех вершин этой триангуляции не лежат в одной плоскости.
Определим {\it множество самопересечений}
$$\Sigma(f):=\Cl\{x\in K\ |\ \# f^{-1}fx >1\}.$$

\begin{pr}\label{vfwhs}
(a) Если $f$ --- вложение, то $\Sigma(f)=\emptyset$.

(b) Нарисуйте $\Sigma(f)$ для стандартного отображения $f:K\to\R^3$ (рис. \ref{kleinb}).

(c) $\Sigma(f)$ --- подграф данной триангуляции.

(d) $\Sigma(f)$ --- цикл, то есть степень любой его вершины четна.
 \end{pr}

Возьмем теперь окружность $S$ в $K$, которая при обычном представлении
бутылки Клейна в виде квадрата со склеенными сторонами представляется
каждой из тех параллельных сторон квадрата, которые склеиваются без перекрутки.
Эту окружность можно считать подграфом двойственного разбиения.
Поэтому $S$ (`трансверсально') пересекает $\Sigma (f)$ в конечном числе точек.
Определим {\it $S$-препятствие Уитни} как
$$w(f):=|S \cap \Sigma (f)|\mod 2.$$
Невложимость бутылки Клейна в $\R^3$ вытекает теперь из следующих задач \ref{vfwhw} и \ref{vfwhs}.d.
 QED

\begin{pr}\label{vfwhw}
(a) Если $f:K\to\R^3$ --- вложение, то $w(f)=0$.

(b) Для стандартного отображения $f:K\to\R^3$ (рис. \ref{kleinb}) выполнено $w(f)=1$.
 \end{pr}

\begin{pr}\label{vfwhspl}
Рассмотрим кусочно-линейную гомотопию $F:K\times I\to\R^3\times I$ общего положения между
кусочно-линейными отображениями $f_0,f_1:K\to\R^3$ общего положения.
{\it Кусочно-линейность} означает линейность на каждом тетраэдре некоторой триангуляции произведения $K\times I$.
{\it Общность положения} означает, что образы никаких пяти вершин этой триангуляции  не лежат в одной трехмерной гиперплоскости.
Определим $\Sigma (F)$ аналогично предыдущему:
$$\Sigma(F):=\Cl\{x\in K\times I\ |\ \# F^{-1}Fx >1\}.$$
\quad
(a) $\Sigma (f_k)\cap S=\Sigma (F)\cap(S\times k)$.

(b) Цилиндр $S\times I$ можно считать подграфом двойственного разбиения и $\Sigma(F)\cap(S\times I)$ --- граф.

(c) В этом графе вершины нечетной степени встречаются только на $S\times\{0,1\}$.

(d) $w(f_0)=w(f_1)$.
 \end{pr}

Заметим, что $\Sigma (F)$ --- 2-полиэдр и $\partial[\Sigma (F)]=[\Sigma(f_0)]-[\Sigma (f_1)]$.

\begin{pr}\label{vfwhrp}
Постройте препятствие Уитни к вложимости $\R P^2$ в $\R^3$ и докажите невложимость $\R P^2$ в $\R^3$.
\end{pr}

 \subsection{Гомотопическая классификация отображений 2-многообразий в окружность}


\begin{pr}\label{hopf-s2s1}
Любое отображение \quad (a) $S^2\to S^1$ \quad (b) $\R P^2\to S^1$ \quad (c) $S^n\to S^1$ при $n\ge 2$ \quad

гомотопно отображению в точку.
\end{pr}

Подмножество $A\subset X\subset\R^m$ называется {\it ретрактом} множества $X$, если существует отображение $X\to A$, тождественное на $A$.

\begin{pr}\label{hopf-ret}
(a) {\it Теорема Хопфа.}
Окружность $S\subset N$ в сфере с ручками $N$ является ее ретрактом тогда и только тогда, когда $N-S$ связно.

(d) При каком условии окружность в компактном 2-многообразии является его ретрактом?
\end{pr}

\begin{pr}\label{hopf-ts1}
Опишите множество отображений $N\to S^1$ с точностью до гомотопности для

(a) диска с $n$ ленточками $N$ (т.е. связного 2-многообразия с непустым краем);
\qquad (b) тора $N$; \qquad

(c) сферы с $g$ ручками $N$; \qquad
(d) сферы с $m$ пленками Мебиуса $N$.
\end{pr}

\subsection{Указания к некоторым задачам}

{\bf \ref{vfs2}.} (b) По векторному полю на верхней полусфере $D^2_+$ постройте векторное поле на диске $D^2$
при помощи центральной проекции $D^2_+\to D^2$ из точки $(-1,0,0)$.
Аналогично для нижней полусферы.
Посмотрите на рисунок \ref{7-3} и используйте задачу \ref{vfsur}.

\begin{figure}[h]
\centering\hskip-6cm\special{em:graph 7-3}
\vskip2.2cm
\caption{Векторное поле на экваторе с точек зрения Белого Медведя и Пингвина}
\label{7-3}
\end{figure}

Идея другого доказательства изображена на рисунке \ref{7-4}.

\begin{figure}[h]
\centering\hskip-6cm\special{em:graph 7-4}
\vskip2.7cm
\caption{Векторное поле на малой окружности, обходящей вокруг северного полюса}
\label{7-4}
\end{figure}

\smallskip
{\bf \ref{vfhemi}.} При помощи центральной проекции $D^2_+\to D^2$ из точки $(-1,0,0)$ постройте
взаимно-однозначное соответствие $V(D^2_+)\to V(D^2)$.

\smallskip
{\bf \ref{vfclex}.}  Ответы: $\Z$.

\smallskip
{\bf \ref{vftor}.} (a,b,c,d) Аналогично задачам \ref{vfpl-}.a,b и \ref{vfunit}.

(e) `Параллельное' поле и поле, соответствующее полю $v(x,y)=(\cos2\pi x,\sin2\pi x)$ на квадрате.

\smallskip
{\bf \ref{vfcltor}.} (a) Для доказательства сюръективности рассмотрите поле $v(x,y):=e^{2\pi i(mx+ny)}$.
Для доказательства инъективности используйте основную теорему топологии и лемму о продолжении гомотопии,
аналогичную задаче \ref{vfann}.c.

(b) К решению пункта (a) нужно добавить использование теорем о продолжении гомотопии и о гомотопности (задачи  \ref{vfbor}.b и \ref{vfhom}.a).

(c) Каждое из них гомотопно полю степени $a+b-a-b=0$.

\comment

{\bf \ref{vfepart}.}
Аналогично задаче \ref{vfs2}.
Обозначим через $C$ `цепочку окружностей' в плоскости $x=1$, заданную уравнением $\Pi_{k=1}^g((z-4k)^2+y^2-4)=0$    (рис. \ref{9-12}).
Определим отображение $f$ из `ближней части' ($x\ge0$) $S_{g,+}$ сферы с ручками в $C$.
...
Определим отображение $p$ из `ближней части' ($x\ge0$) $S_{g,+}$ сферы с ручками в диск с дырками
$\Delta$ в плоскости $x=0$, заданную уравнением $\Pi_{k=1}^g((z-4k)^2+y^2-4)\le1$ (рис. \ref{7-9}).
Для каждой точки $A\in S_{g,+}$ обозначим через $p(A)\in\Delta$ центральную
проекцию точки $A$ на плоскость $x=0$ из точки $f(A)$.
Центральная проекция на плоскость $x=0$ из точки $f(A)$ переводит касательную к $S_g$ плоскость в точке $A$
в плоскость $x=0$.
При этом ненулевые векторы переходят в ненулевые.
Получается векторное поле $v_+$ на $\Delta$.

Аналогично получается векторное поле $v_-$ на $\Delta$ из векторного поля на `дальней части' ($x\le0$)
сферы с ручками.
В каждой точке $B\in\partial\Delta$ векторы $v_+(B)$ и $v_-(B)$ не являются симметричными относительно прямой, перпендикулярной касательной к $\partial\Delta$ в точке $B$.
Поэтому $v_-|_{\partial\Delta}$ гомотопно полю, полученному из $v_+|_{\partial\Delta}$ симметрией относительно касательной.
Продолжим эту гомотопию до гомотопии всего поля $v_-$ (аналогично задаче \ref{vfbor}).
Теперь из задач \ref{vfsum} и \ref{vfsur} получается противоречие.

\endcomment

\smallskip
{\bf \ref{vfhole}.} (a) Отобразите (сомните) $S_{g,0}$ с самопересечениями в плоскость.

(b) Ответ: $\Z^{2g}$.

(c) Каждое из них гомотопно полю степени $a_1+b_1-a_1-b_1+\dots+a_g+b_g-a_g-b_g=0$.

\smallskip
{\bf \ref{vfmob}.} (с) Ответ:  $\Z_2$.

\smallskip
{\bf \ref{vfklein}.} (b) Ответ:  $\Z\oplus\Z_2$.

\bigskip
{\bf \ref{vfta}.} (d) Для поля $v$, изображенного на рис.~\ref{vesph}, $e(v)=1+1=2$;
разберитесь, почему не $1-1=0$; см. также рис.~\ref{water},~\ref{eul}.

\smallskip
{\bf \ref{vfnonor}.} Аналогично предыдущему --- определите класс Эйлера $e(N)$ как препятствие к построению
такого поля, а потом докажите его полноту и вычислите его.

 \bigskip
{\bf \ref{vfnos1}.} (c,d) Ответы: $\Z$.


\smallskip
{\bf \ref{vfnobou}, \ref{vfnohi}.} Сначала постройте нормальное поле на вершинах, затем продолжите его
на ребра и потом продолжите его на грани.

\smallskip
{\bf \ref{vfnoth}.}
(a) Следует из \ref{vfnoeu}.c и полноты препятствия, построенного в \ref{vfnoeu}.a.

(b) Следует из \ref{vfnoeu}.ac.

(c) Следует из \ref{vfnoeu}.a и \ref{vfnorp2}.a.

\smallskip
{\bf \ref{vfnoeu}.} (b) Постройте препятствие $e_1(N)$ к существованию замкнутого 2-многообразия $N'\subset\R^4$, близкого к $N$ и не пересекающего $N$.
Затем докажите, что $N\cap N'=e_1(N)=\overline e(N)$.

(с) Возьмем гомотопию $F:N'\times I\to\R^m$ общего положения между включением $N'\to\R^4$ и отображением, образ
которого не пересекает $N$.
Тогда $N\cap F(N'\times I)$ есть одномерное ориентированное подмногообразие с краем $N\cap N'$.

\smallskip
{\bf \ref{vfnorp2}.}
(a) См. \ref{vfsub34}.

(b) Ответ: 1.

(c) Возьмите связную сумму попарно непересекающихся 2-многообразий из (a) и полученных из него переносами.

\smallskip
{\bf \ref{vfnocl}.} Ответы: (a) $0$,  (b) $\Z\oplus\Z$.

\smallskip
{\bf \ref{vfnohic}.} Ответы: (a,b) $\Z$, (c) 0.

\smallskip
{\bf \ref{hopf-s2s1}.}
Аналогично теореме гомотопности (задача \ref{vfhom}).
Представьте отображение $S^2\to S^1$ в виде композиции $D^2\to S^2\to S^1$.
Выведите отсюда, что у любого отображения $f:S^2\to S^1$ существует поднятие, т.е. такое отображение
$\t f:S^2\to\R$, что $f(x)=\overline f(x)\mod2\pi$ для любого $x\in S^2$.

 \newpage
\section{Трехмерные многообразия}\label{03man}

\subsection{Определение трехмерного симплициального комплекса}

{\it Трехмерным симплициальным комплексом} называется семейство двухэлементных, трехэлементных и
четырехэлементных подмножеств конечного множества, которое вместе с каждым множеством содержит все его подмножества.
Мы будем сокращенно называть трехмерный симплициальный комплекс просто {\it 3-комплексом}.
Элементы данного конечного множества называются {\it вершинами} комплекса, выделенные

$\bullet$ двухэлементные подмножества --- {\it ребрами} комплекса,

$\bullet$ трехэлементные подмножества --- {\it гранями} комплекса,

$\bullet$ четырехэлементные подмножества --- {\it тетраэдрами (трехмерными гранями)} комплекса.

Они называются также {\it симплексами} размерностей 0, 1, 2 и 3, соответственно.

Определение {\it тела 3-комплекса} и {\it гомеоморфности тел} аналогично случаю 2-комплексов.

\smallskip
{\it Примеры 3-комплексов.}

(1) Полный 3-комплекс с 4 вершинами.
Его тело гомеоморфно шару
$D^3=\{(x,y,z)\in\R^3\ |
x^2+y^2+z^2\le1\}$.

(2) Полный 3-комплекс с 5 вершинами. Его тело гомеоморфно сфере $S^3=\{(x,y,z,t)\in\R^4\ |\ x^2+y^2+z^2+t^2=1\}$.

Для 2-комплекса $P$ можно определить 3-комплексы $P\times I$, $P\times S^1$ и $P\times G$--- {\it произведения}  2-комплекса $P$ на отрезок $I$, на окружность $S^1$ и на граф $G$ аналогично произведению графов.

\subsection{Гомеоморфность трехмерных комплексов}

\begin{figure}[h]\centering
\includegraphics{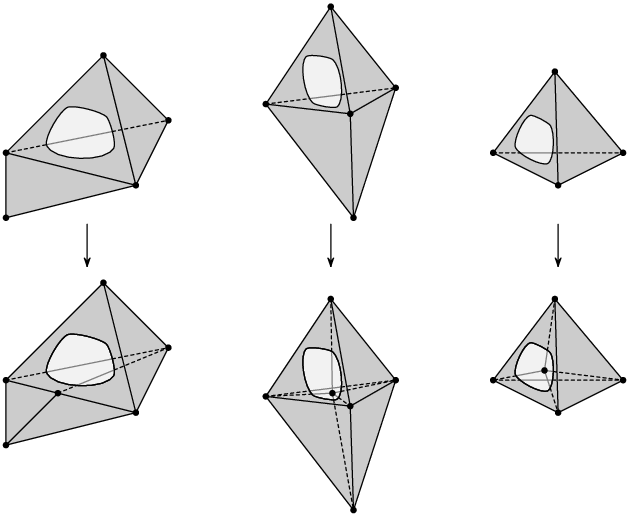}
\caption{Подразделение одномерной, двумерной и трехмерной граней}
\label{podra2}
\end{figure}

Операция {\it подразделения ребра} изображена на рис.~\ref{podra2} слева.
(Определение аналогично случаю 2-комплексов;
читатель легко восстановит формальное комбинаторное определение этой операции по рисунку.)

\begin{pr}\label{3-subd}
Операции {\it подразделения грани и тетраэдра} на рис.~\ref{podra2} в центре и справа выражаются через операцию подразделения ребра.
\end{pr}

Два 3-комплекса называются {\it гомеоморфными}, если от одного можно перейти к другому при помощи
операций  подразделения ребра и обратных к ним.
Обозначение: $P\cong Q$.

\begin{pr}\label{3-bouho}
(a) Трехмерные шар и сфера не гомеоморфны.

(b) Если произведение сферы с $g$ ручками на отрезок гомеоморфно произведению сферы с $h$
ручками на отрезок, то $g=h$.

(c) Придумайте две негомеоморфных триангуляции $P$ и $Q$ 2-многообразий, для которых $P\times I\cong Q\times I$.

(d) Для каких триангуляций $P$ и $Q$ 2-многообразий (возможно, неориентируемых и имеющих край)
$P\times I\cong Q\times I$?
\end{pr}

\subsection{Локально евклидовы трехмерные комплексы}

Неформально, трехмерным многообразием называется фигура (тело 3-комплекса), любая точка которой
имеет малую окрестность (звезду), топологически гомеоморфную трехмерному шару.

3-комплекс называется {\it триангуляцией 3-многообразия}, если для любой его вершины объединение симплексов, ее содержащих,
гомеоморфно полному 3-комплексу с 4 вершинами.
Пока `триангуляция 3-многообразия' --- единый термин, понятий 3-многообразия и его триангуляции мы еще не определили.
Такие 3-комплексы следовало бы называть `локально-евклидовыми', ср. \S\ref{02man}.

\begin{pr}\label{3-equi}
(a) Это равносильно тому, что линк
любой вершины 3-комплекса гомеоморфен некоторому 2-комплексу, представляющему сферу $S^2$ или диск $D^2$.

(b) Это равносильно тому, что линк любой вершины связен, является триангуляцией 2-многообразийя (т.е. локально-евклидов), ориентируем и имеет эйлерову характеристику 2 или 1.

(c) В триангуляции 3-многообразия к каждой грани примыкает 1 или 2 тетраэдра.

(d) Придумайте 3-комплекс, не являющийся триангуляцией 3-многообразия, линк любой вершины которого связен и для любого ребра $\{u,v\}$ которого все тетраэдры, содержащие это ребро, образуют одну из цепочек
$$\{u,v,a_1,a_2\},\ \{u,v,a_2,a_3\}\dots\{u,v,a_{n-1},a_n\},\ \{u,v,a_n,a_1\}
\quad
\text{ или }
$$
$$
\{u,v,a_1,a_2\},\ \{u,v,a_2,a_3\}\dots\{u,v,a_{n-1},a_n\}.$$
\end{pr}

(6) {\it Регулярная окрестность} 2-комплекса (в частности, графа) в некоторой триангуляции 3-многообразия
(в частности, в $S^3$).

{\it Регулярная окрестность} подкомплекса в комплексе --- симплициальная окрестность во втором барицентрическом подразбиении. Ср. \S\ref{02man}.
Таким окрестностям гомеоморфны, например, $D^3$ для точки и $F\times I$ для 2-комплекса $F$.


\begin{pr}\label{3-exa}
Следующие фигуры гомеоморфны (топологически, а не кусочно-линейно) телам некоторых триангуляций 3-многообразий.

(3) Проективное пространство $\R P^3$: оно получается из сферы $S^3$ отождествлением диаметрально противоположных точек (или, что то же самое, из диска $D^3$ отождествлением диаметрально противоположных точек на его граничной сфере).

(4)
Пусть $p,q$ --- взаимно простые целые положительные числа.
Определим {\it линзовое пространство}
$$L(p,q):=S^3/(z_1,z_2)\sim(z_1e^{2\pi i/p},z_2e^{2\pi iq/p})_{z_1,z_2\in\C,\ |z_1|^2+|z_2|^2=1}.$$
Оно получается склейкой граней объединения двух $p$-угольных пирамид по их общему основанию.
Каждая верхняя грань $A$ склеивается с нижней, полученной из $A$ композицией вращения на $2\pi q/p$ относительно прямой, соединяющей вершины пирамид, и симметрии относительно плоскости основания пирамид.

(5) `Трехмерная лента Мебиуса' $D^2\widetilde\times S^1$ получается из трехмерного цилиндра $D^2\times I$ склейкой   точек $(x,0)$ и $(\sigma(x),1)$ для любого $x\in D^2$.
Здесь $\sigma:D^2\to D^2$ --- осевая симметрия.
\end{pr}

\begin{pr}\label{3-cone}
(a) $L(1,1)=S^3$ и $L(2,1)\cong\R P^3$.

(b) $L(5,2)\cong L(5,3)$. \qquad

(c) $L(7,2)\cong L(7,3)\cong L(7,4)\cong L(7,5)$.

(d) Если $q_1\equiv \pm q_1^{\pm1}\mod p$ (знаки $\pm$ не обязательно согласованы),
то $L(p,q_1)\cong L(p,q_2)$.
\end{pr}

Примеры (3)--(5) можно представить при помощи {\it склейки граней многогранников}.
Так удобно задавать и многие другие 3-многообразия.
Приведем формализацию этой идеи.

Пусть в 3-комплексе задан 2-подкомплекс.
Рассмотрим {\it дополнение до  открытой регулярной окрестности 2-подкомплекс}, т.е. объединение симплексов второго барицентрического подразбиения 2-комплекса, не пересекающих 2-подкомплекс.
Вложение 2-подкомплекса в 3-комплекс называется {\it клеточным}, если каждая связная компонента этого дополнения
гомеоморфна 3-диску.
Например, вложение точки в 3-сферу клеточно.

{\it Клеточным разбиением} 3-комплекса называется клеточное вложение некоторого 2-комплекса в этот 3-комплекс, вместе с клеточным разбиением этого 2-комплекса.
Возможно, начинающему будет более понятен термин {\it разбиение на многогранники}.
Например, 3-комплекс является клеточным разбиением себя.

{\it Краем} (или границей) $\partial T$ триангуляции $T$ 3-многообразия называется
объединение тех (двумерных) граней, которые содержатся только в одном тетраэдре.

\begin{pr}\label{3-bou}
Чему гомеоморфен край

(a) Шара $D^3$?

(b) Произведения на окружность сферы с $g$ ручками и $h$ дырками?

(c) Произведения на отрезок сферы с $g$ ручками и $h$ дырками?

(d) Произведения на отрезок сферы с $m$ пленками Мебиуса и $h$ дырками?

(e) `Трехмерного листа Мебиуса' $D^2\widetilde\times S^1$?
\end{pr}

\begin{pr}\label{3-bougen}
(a) Край любой триангуляции 3-многообразия является триангуляцией замкнутого 2-многообразия (не обязательно связного).

(b) Края гомеоморфных триангуляций 3-многообразий гомеоморфны.
\end{pr}

{\bf Теорема.}
{\it Существует алгоритм распознавания гомеоморфности 3-комплексов сфере $S^3$ (т.е. полному 3-комплексу с 5 вершинами).}

\smallskip
{\bf Теорема.}
{\it Существует алгоритм распознавания гомеоморфности достаточно больших триангуляций 3-многообразий.}

\smallskip
Определение достаточно  больших (триангуляций) 3-многообразий
и доказательство этих теорем можно найти, например, в \cite{Ma03}.
Aлгоритмы достаточно сложны, они основаны на теории нормальных поверхностей Хакена, дополненой идеями Вальдхаузена,
Йохансена, Хемиона, Рубинштейна, Томпсон и Матвеева.

Существование алгоритма распознавания гомеоморфности произвольных (не обязательно достаточно больших) трехмерных многообразий, видимо,  следует из  справедливости геометризационной  гипотезы Терстона.
Автору неизвестно, появилось ли полное доказательство.

 {\it 3-полиэдром} называется класс эквивалентности 3-комплексов с точностью до гомеоморфизма.
Представители этого класса эквивалентности называются {\it триангуляциями} соответствующего 3-полиэдра.
{\it Трехмерным кусочно-линейным многообразием} называется класс гомеоморфности триангуляции 3-многообразия.
Мы будем сокращенно называть трехмерное кусочно-линейное многообразие просто {\it 3-многообразием}.

3-Многообразие называется {\it связным}, {\it ориентируемым} и т.д., если некоторый (или, эквивалентно, любой) представляющий его 3-комплекс связен, ориентируем (см. определение ниже) и т.д.

{\it Связная  сумма} трехмерных многообразий определяется аналогично двумерному случаю.

\subsection{Ориентируемость трехмерных многообразий}\label{03manor}

\begin{pr}\label{3-oriho}
Произведения тора и бутылки Клейна на окружность не гомеоморфны.
\end{pr}

{\it Ориентацией} грани комплекса называется упорядочение ее вершин с точностью до четной перестановки.
Ориентация $(1234)$ тетраэдра в комплексе {\it порождает} ориентации $(123),(234),(134),(124)$ (двумерных) граней.
Двумерная грань комплекса называется {\it внутренней}, если она лежит по крайней мере в двух тетраэдрах.
Триангуляция 3-многообразия называется {\it ориентируемой}, если можно так ввести ориентации на всех ее тетраэдрах,  чтобы с двух сторон каждой внутренней грани порождались бы {\it противоположные} ориентации (ср. рис.~\ref{prot}).
Указанный набор ориентаций на гранях называется {\it ориентацией} 3-комплекса.

\begin{pr}\label{3-ori}
(а) Гомеоморфные триангуляции 3-многообразий ориентируемы или нет одновременно.

(b) Произведение триангуляции  $P$ 2-многообразия на отрезок (или на окружность) ориентируемо тогда и только тогда, когда $P$ ориентируем.

(c) Триангуляция 3-многообразия  ориентируема тогда и только тогда, когда никакой гомеоморфный ей комплекс не
содержит подкомплекса, представляющего `трехмерную ленту Мебиуса' $D^2\t\times S^1$.
\end{pr}

\begin{pr}\label{3-oriemb}
(a) Бутылка Клейна и проективная плоскость;
\qquad

(b) Любое 2-многообразие

вложимы в некоторые ориентируемые триангуляции 3-многообразий (т.е. гомеоморфны некоторым их подкомплексам).
\end{pr}

\begin{pr}\label{ex3-ori}
Определите $H_2(N)$ и $w_1(N)\in H_2(N)$ так, чтобы выполнялся следующий результат.

 {\bf Теорема ориентируемости.}
{\it Замкнутое 3-многообразие $N$ ориентируемо тогда и только тогда, когда его}
первый класс Штифеля-Уитни $w_1(N)\in H_2(N)$ {\it нулевой.}
\end{pr}


Произведения 2-многообразий на отрезок обобщаются до {\it $I$-расслоений над 2-мно\-го\-об\-ра\-зи\-я\-ми}.
Пусть дано замкнутое 2-многообразие $F$.
Возьмем семейство $M$ непересекающихся окружностей на $F$.
Разрежем $F$ по этому семейству.
Получим 2-многообразие $F'$ с краем и инволюцией
$\sigma:\partial F'\to\partial F'$, не имеющей неподвижных точек.
Назовем {\it утолщением} 2-многообразия $F$ (или {\it пространством $I$-расслоения над $F$})
3-многообразие
$$F\widetilde\times_M I:=F'\times I/(x,t)\sim(\sigma(x),1-t)_{x\in\partial F', t\in I}.$$

\begin{pr}\label{3-iskew}
(a) Если $M=\emptyset$, то $F\widetilde\times_M I=F\times S^1$.

(b) Если на $F-M$ имеется ориентация, меняющаяся при пересечении каждой окружности из $M$
(т.е. если $M$ представляет $w_1(F)$, см.  \S\ref{02hom}), то $F\widetilde\times_M I$ ориентируемо.
\end{pr}

 \subsection{Эйлерова характеристика трехмерных комплексов}\label{03maneul}

{\it Эйлеровой характеристикой} 3-комплекса $K$ называется знакопеременная сумма количества
симплексов размерностей 0,1,2 и 3:
$$\chi(K):=V-E+F-P.$$

\begin{pr}\label{3-eul}
(a) Найдите эйлеровы характеристики 3-комплексов, соответствующих примерам из \ref{3-exa}.

(b) $\chi(P\times Q)=\chi(P)\times\chi(Q)$ для графа $P$ и 2-комплекса $Q$.

(c) Эйлеровы характеристики гомеоморфных 3-комплексов равны.

(d) Если $T'$ --- клеточное разбиение 3-комплекса $T$, то $\chi(T)=\chi(T')$, где
{\it эйлерова характеристика} $\chi(T')$ клеточного разбиения $T'$ есть знакопеременная сумма количества
граней разбиения размерностей 0,1,2 и 3: $\chi(T')=V'-E'+F'-P'$.

(e) Эйлерова характеристика триангуляции замкнутого 3-многообразия равна нулю.
\end{pr}

\begin{figure}[h]\centering
\includegraphics{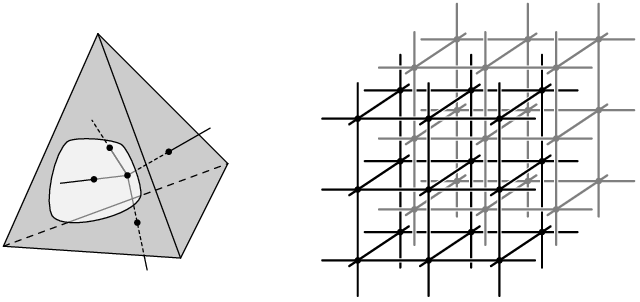}
\caption{Двойственные клеточные разбиения}\label{dual3}
\end{figure}

\smallskip
{\it Определение двойственного клеточного разбиения (рис.~\ref{dual3}).}
Возьмем некоторое клеточное разбиение $U$ триангуляции замкнутого 3-многообразия.
В каждом многограннике $x$ разбиения выберем по точке $x^*$.
Для каждой грани $f$ разбиения соединим ребром $f^*$ две построенные точки,
лежащие в соседних по этой грани многогранниках разбиения.
(Это ребро протыкает грань $f$ разбиения; рис.~\ref{dual3}.)
Для каждого ребра $a$ разбиения натянем двумерный многоугольник $a^*$ на
построенные ребра, соответствующие тем граням разбиения, которые содержат ребро $a$.
(Этот многоугольник протыкается ребром $a$; рис.~\ref{dual3}.)
Объединение полученных многоугольников является 2-комплексом, клеточно вложенным в $U$.
Полученное клеточное разбиение $U^*$ называется {\it двойственным} к $U$.

\subsection{Указания и решения к некоторым задачам}

\smallskip
{\bf \ref{3-equi}.} (c) $D^3\vee_{(1,0,0)} D^3$ или конус над тором.

\smallskip
{\bf \ref{3-bouho}.} Используйте край.

(c) Тор с дыркой и сфера с 3 дырками.

\smallskip
{\bf \ref{3-bou}.} (e) Бутылке Клейна.

\smallskip
{\bf \ref{3-eul}.} (c) Достаточно доказать, что эйлерова характеристика не меняется при подразделении ребра.

(d) Аналогично формуле Эйлера для графов на поверхностях (задача \ref{1-nonrea}.a).

(e) Используйте двойственное клеточное разбиение и (c,d).

\smallskip
{\bf \ref{3-oriho}.} Используйте ориентируемость.

\smallskip
{\bf \ref{3-cone}.}
(c,d) Это частные случаи пункта (e).

(e) Вытекает из $L(p,q)\cong L(p,p-q)$ и $L(p,q_1)\cong L(p,q_2)$ при $q_1q_2\equiv1\mod p$.
Для доказательства второго разрежьте бипирамиду $SNA_1A_2\dots A_p$, склейкой которой является линза,
на тетраэдры $SNA_kA_{k+1}$.
Склейте их в новую бипирамиду с `осью' $A_1A_2=A_2A_3=\dots$ и ребрами основания,
отвечающими размноженному отрезку $SN$.

\smallskip
{\bf \ref{3-oriemb}.}
(a)
Пусть
$S^1=\{z\in\Bbb C\ :
\ |z|=1\}$. Тогда 3-многообразие
$\dfrac{S^1\times [-1,1]\times [0,1]}{(z,t,0)\sim(\overline z,-t,1)}$ ориентируемо и содержит
бутылку Клейна $\dfrac{S^1\times0\times[0,1]}{(z,0,0)\sim(\overline z,0,1)}$.

Для проективной плоскости $\R P^2$, стандартно вложенной в $\R P^3$,
можно взять ее `регулярную' окрестность в $\R P^3$.

(b) Достаточно доказать это утверждение для замкнутых многообразий.
Для ориентируемых оно очевидно.
Для неориентируемых используйте классификацию и (a).

 \newpage
\section{Векторные поля для многомерных многообразий}\label{0vehima}

 \subsection{Векторные поля на подмножествах в $\R^m$}

\begin{pr}\label{hopf-bou}
Любое отображение в $S^2$ из сферы с ручками, пленками Мебиуса и хотя бы одной дыркой
гомотопно отображению в точку.
\end{pr}

Обозначим через $\pi^n(N)$ множество отображений $N\to S^n$ с точностью до гомотопности.

\begin{pr}\label{hopf-deg}
(a) Для отображения $f:S^2\to S^2$ определите {\it степень} $\deg f$ как сумму знаков прообразов значения общего положения ({\it регулярного} значения) при гладкой (или кусочно-линейной) аппроксимации.
(См. детали определения в \cite[\S18.1]{Pr04}.)

(b) Указанная сумма одинакова для двух гладко гомотопных гладких отображений и их общего регулярного значения.

(с) Степень корректно определена, т.е. не зависит от выбора регулярного значения и гладкой аппроксимации.

(d) Степень не меняется при гомотопии.

(e) Степень $\deg:\pi^2(S^2)\to\Z$ является взаимно-однозначным соответствием.

(f) Тождественное отображение сферы $S^2$ не гомотопно антиподальному.
\end{pr}

\begin{pr}\label{hopf-ts2}
(a-e) Сформулируйте и докажите аналоги предыдущей задачи для отображений $S^1\times S^1\to S^2$.
\end{pr}

\begin{pr}\label{hopf-n2s2}
(a) Для сферы с ручками (т.е. замкнутого ориентируемого 2-многообразия) $N$ существует биекция $\deg:\pi^2(N)\to\Z$.

(b) $|\pi^2(\R P^2)|=2$.
\qquad

(c) Опишите $\pi^2(N)$ для любого замкнутого 2-многообразия $N$.
\end{pr}

\begin{pr}\label{hopf-tos2}
(a) Любое отображение $S^k\to S^n$ гомотопно отображению в точку для $k<n$.

Иными словами, в $n$-мерном пространстве ненулевое векторное поле (векторов $n$-мерного пространства), заданное на границе $k$-мерного шара, продолжается до ненулевого векторного поля на шаре при $k<n$.

(b) 
В $n$-мерном пространстве ненулевое векторное поле (векторов $n$-мерного
пространства), заданное на границе $n$-мерного шара, продолжается до ненулевого
векторного поля на шаре тогда и только тогда, когда алгебраическое число векторов
на границе, смотрящих в данное направление общего положения, равно нулю.

(c) {\bf Многомерная основная теорема топологии.} {\it Существует биекция $\deg:\pi^n(S^n)\to\Z$.}
\end{pr}

Обозначим $S^n:=\{(x_1,\dots,x_{n+1})\in\R^{n+1}\ :\ x_1^2+\dots+x_{n+1}^2=1\}.$

\begin{pr}\label{vehexa}
Определите следующие стандартные подмножества в $\R^m$.
Касательные векторные поля на них определяются аналогично началу \S\ref{0ve2ma}.
Докажите, что на каждом из них существует ненулевое касательное векторное поле.

(a) $S^1\times S^1\times S^1$; \quad (b) $S^2\times S^1$; \quad (c) $S^3$;
\quad

(d) Декартово произведение сферы с ручками и $S^1$;
\quad
(e) $S^{2k-1}$; \quad
(f) $S^{2k-1}\times S^q$.
\end{pr}

Определения $n$-мерных {\it гладких} подмногообразий, касательных и нормальных векторных полей и полей направлений на них, а также гомотопности векторных полей, аналогичны случаю 2-многообразий (\S\ref{0ve2ma}).
(Те, кому трудно работать с многообразиями размерности больше 3, могут решать приводимые задачи для трехмерных многообразий --- даже этот случай интересен.)
Аналогично формулируются и доказываются теоремы о триангулируемости и о кусочно-линейной гомеоморфности любых двух  триангуляций одного многообразия.

Определение $n$-мерных {\it кусочно-линейных} многообразий для $n=3$ приведено в \S\ref{03man}
(для произвольного $n$ определение аналогично).

 \subsection{О существовании ненулевых касательных векторных полей}

\begin{pr}\label{vehgen}
(a) Ненулевое касательное векторное поле, заданное на вершинах триангуляции 3-многообразия, можно продолжить на его ребра.

(b) Ненулевое касательное векторное поле, заданное на объединении ребер триангуляции 3-многообразия, можно продолжить на объединение его двумерных граней.

(c) На любом многообразии с краем существует ненулевое касательное векторное поле.
\end{pr}

\begin{pr}\label{veheul}
 Эйлерова характеристика $\chi$ триангуляции многообразия определяется как знакопеременная сумма по $k$
количеств $k$-мерных граней.

(a) Эйлеровы характеристики кусочно-линейно гомеоморфных триангуляций многообразий равны.

(b) $\chi(P\times Q)=\chi(P)\times\chi(Q)$
\end{pr}

 Эйлерова характеристика $\chi(N)$ многообразия $N$ определяется как знакопеременная сумма по $k$
количеств $k$-мерных граней в клеточном разбиении этого многообразия.
Аналогично случаю 2-многообразий доказывается, что эйлерова характеристика не зависит от выбора такого разбиения данного многообразия.

Существуют простые способы вычислять эйлерову характеристику.

\begin{pr}\label{vehopf}
(a) {\bf Теорема Хопфа.} На замкнутом многообразии существует ненулевое
касательное векторное поле тогда и только тогда, когда эйлерова характеристика этого многообразия нулевая.

(b) {\bf Следствие.} На любом нечетномерном многообразии
существует ненулевое касательное векторное поле.

(c) {\bf Следствие.} На $n$-мерной сфере существует ненулевое касательное векторное поле тогда и только тогда, когда $n$ нечетно.
\end{pr}


 \subsection{О  существовании  ненулевых нормальных векторных  полей}


\begin{pr}\label{} На $n$-многообразии в $\R^{n+1}$ существует ненулевое нормальное поле тогда и только тогда, когда
подмногообразие ориентируемо.
 \footnote{Из этого вытекает, что не существует замкнутого неориентируемого $n$-многообразия в $\R^{n+1}$.
Действительно, пусть существует.
Аналогично теореме Жордана доказывается, что многообразие делит $\R^{n+1}$ на две части.
Тогда существует ненулевое нормальное поле, направленное во внешнюю часть.
Значит, многообразие ориентируемо.}
\end{pr}

 {\bf Теорема о нормальных полях.}
{\it Если $m\ge2n$ или $m\le n+2$, то для любого замкнутого ориентируемого $n$-подмногообразия в $\R^m$ существует ненулевое нормальное векторное поле.}

\begin{pr}\label{vehnor}
(a) Докажите теорему о нормальных полях для $m\ge 2n+1$.

(b)  На любом $n$-подмногообразии с непустым краем в $\R^{2n}$ имеется нормальное поле.

(с) Докажите теорему о нормальных полях для $m=2n$.

(d) То же для $m=n+2$.
\end{pr}

 При $n+2<m<2n$, или при $n+2=m$ для подмногообразий с краем,
ненулевого нормального векторного поля может не быть.

\begin{pr}\label{vehnexi}
(a) Существует ориентируемое 3-многообразие с краем в $\R^5$, не имеющее нормального поля.

(b) Существует 4-подмногообразие в $\R^7$, не имеющее нормального поля.
\end{pr}


\begin{pr}\label{veh35}
  Для любого ориентируемого 3-многообразия $N\subset\R^5$ с непустым краем

(a) постройте группу $H_1(N,\partial N;\Z)$ и препятствие $\overline e(f)\in H_1(N,\partial N;\Z)$
({\it нормальный класс Эйлера}) к существованию нормального поля.

(b)* если в $H_1(N,\partial N;\Z)$ нет элементов порядка 2, то $\overline e(f)$ четно, т.е.
$\overline e(f)=2x$ для некоторого $x\in H_1(N,\partial N;\Z)$.
\end{pr}

Доказывать утверждения \ref{vehnexi}.b и \ref{veh35}.b лучше при изучении дальнейшего материала.

Полный ответ на следующий вопрос (М. Хирша) неизвестны:
{\it при каких $(m,n)$ для любого $n$-подмногообразия в $\R^m$
существует ненулевое нормальное векторное поле?}


\subsection{О классификации ненулевых касательных полей}

\begin{pr}\label{vehinfty}
Постройте бесконечное количество попарно негомотопных ненулевых касательных векторных полей на

(a) $S^1\times S^1\times S^1$; \quad
(b) $S^1\times S^2$; \quad
(c) $S^3$.  \quad
\end{pr}


\begin{pr}\label{vehinvar}
(a) $V(S^3)\cong\Z$ (умножение на $V(S^3)$ задается умножением кватернионов).


(b)* Классифицируйте поля направлений на
3-многообразиях.
По поводу $n$-мерного случая и связи с лоренцевыми метриками см.~\cite{Ko01}.
\end{pr}

{\bf Теорема классификации векторных полей.}
{\it Если для замкнутого ориентируемого $n$-многообразия $N$ выполнено
$\chi(N)=0$, то существует сюръекция $V(N)\to H_1(N;\Z)$ на} группу одномерных гомологий многообразия $N$
(с целыми коэффициентами).

{\it Более того, существует такое отображение $D:V(N)\times
V(N)\hm\to H_1(N;\Z)$, что $D(\cdot,v)$ сюръективно и
$D(u,v)+D(v,w)=D(u,w)$ для любых $u,v,w\in V(N)$.

При $n=2$ отображение $D(\cdot,v)$ является взаимно однозначным соответствием.

При $n=3$ число прообразов класса $x\in H_1(N,\Z)$ при
отображении $D(\cdot,v)$ равно $\lessmskips{2mu}2x\cap H_2(N;\Z)$
(группа $H_2(N;\Z)$ и умножение $\lessmskips{2mu}\cap:H_1(N;\Z)\hm\times H_2(N;\Z)\to\Z$ определены в задачах в \S\S\ref{0vesy3},\ref{0homol}
аналогично \S\ref{02hom}).

При $n\ge4$ при отображении $D(\cdot,v)$ у каждого класса ровно два прообраза.}


\smallskip
Для $n=2$ эта теорема является фольклорным результатом начала 20-го века.
Для $n=3$ она является следствием теорем Штифеля
(\S\ref{0vesy3}) и Понтрягина (\S\ref{0hopf},~\cite{CRS}).
Для $n\ge4$ эта теорема, видимо, является фольклорным результатом середины 20-го века~\cite[Theorem 18.2]{Ko81}.
Этот результат получен с помощью теории, основы которой здесь излагаются.

По задаче \ref{vehinfty}.с указанное отображение $D(\cdot,v)$ не инъективно для $N=S^3$.

В следующем пункте мы определяем группу $H_1(N,\Z)$ и доказываем теорему классификации векторных полей для $n=2$.
Приводимый способ, более сложен, чем намеченный в задаче \ref{vfcltor}.b, зато
он позволяет определить инвариант $D$ для общего случая.

\begin{pr}\label{vehopfn}
 Опишите $\pi^n(N)$ для компактного $n$-многообразия $N$.
\end{pr}

 \subsection{Построение гомологического инварианта векторных полей}


\smallskip
{\it Набросок определения класса $D(u,v)$, использующего общее положение.}
Рассмотрим произвольную гомотопию $u_t$ общего положения между ненулевыми векторными полями
$u=u_0$ и $v=u_1$ на многообразии $N$ (векторы гомотопии не предполагаются ненулевыми).
Общность положения означает, что множество точек многообразия, в которых $u_t=0$ для некоторого $t$,
является несвязным объединением окружностей.
Используя нашу гомотопию, на этих окружностях можно ввести ориентацию.
Это объединение ориентированных окружностей с точностью до {\it гомологичности} и называется классом $D(u,v)$.
Два набора непересекающихся ориентированных окружностей называются
{\it гомологичными}, если существует ориентированное 2-многообразие в $N\times I$, ориентированный край
которого есть объединение первого набора в $N\times0$ и второго набора в $N\times 1$.
Доказательство корректности этого определения аналогично случаю $\chi(N)$.

\smallskip
Далее мы приводим аккуратное определение класса $D(u,v)$, использующее клеточное разбиение.

\smallskip
{\it Определение различающего цикла.}
Как и в доказательстве теоремы Эйлера-Пуанкаре через триангуляции, возьмем клеточное разбиение $U$ многообразия
$N$ и двойственное к нему клеточное разбиение $U^*$, а также фиксируем ориентации.
Рассмотрим два поля $u$ и $v$ на многообразии $N$.
Поле $u$ гомотопно такому, которое совпадает с $v$ на $U_0^*$.
Поэтому будем считать, что $u=v$ на $U_0^*$ (заметим, что наши построения
могут зависеть от выбора совмещения полей на $U_0^*$).
Тогда можно определить различающую расстановку $d(u,v)$, как и в доказательстве теоремы Эйлера-Пуанкаре.
Ясно, что $d(u,u)=0$ и что $d(u,v)+d(v,w)=d(u,w)$.
Поскольку поле $u$ продолжается на грани разбиения $U^*$, то выполнено следующее {\it правило Кирхгофа:}

{\it для любой вершины $A$ разбиения $U$ сумма чисел расстановки $d(u,v)$ на входящих в $A$ ребрах
равна сумме чисел расстановки $d(u,v)$ на выходящих из $A$ ребрах.}

Такие расстановки называются {\it циклами}.

\smallskip
{\it Изменение различающего цикла.}
Ясно, что если $d(u,v)=0$, то поля $u$ и $v$ гомотопны.
Обратное неверно, как показывает пример следующей гомотопии.
Для вершины $f^*$ изменим поле $u$ так, чтобы вектор в $f^*$ сделал один
оборот против часовой стрелки, вектора в маленькой окрестности вершины $f^*$
'потянулись' за вектором в $f^*$, а вне этой маленькой окрестности поле
осталось прежним.
Обозначим полученное поле через $u'$.
Понятно, что $d(u',u)$ есть расстановка $\pm1$ (в зависимости от ориентации)
на ребрах, ограничивающих грань $f$, и нулей на всех остальных ребрах.
Эта расстановка называется {\it границей грани $f$} и обозначается $\partial f$.

B окрестностях вершин $f^*_1,\dots,f^*_s$ сделаем описанную выше гомотопию
поля $u$, поворачивая векторы в этих гранях на $n_1,\dots,n_s$ оборотов, соответственно.
Обозначим полученное поле через $u(n_1,\dots,n_s)$ и назовем его {\it подкруткой} поля $u$.
Тогда
$d(u(n_1,\dots,n_s),u)=n_1\partial f_1+\dots+n_s\partial f_s.$

Теперь рассмотрим гомотопию $u_t$ между полями $u_0$ и $u_1$, совпадающими
на $U_0^*$ (правда, эта гомотопия может менять поле даже на $U_0^*$).
Для каждой грани $f_i$ обозначим через $n_i$ число оборотов при изменении $t$
от 0 до 1 вектора $u_t(f^*)$.
Легко проверить, что
$$
d(u_0,v)-d(u_1,v)=n_1\partial f_1+\dots+n_s\partial f_s.
$$

\smallskip
{\it Определение группы $H_1(N;\Z)$ и завершение доказательства.}
Назовем {\it границей} сумму $n_1\partial f_1+\dots+n_s\partial f_s$ границ
нескольких граней с целыми коэффициентами.
Назовем циклы $\gamma_1$ и $\gamma_2$ {\it гомологичными}, если
$\gamma_1-\gamma_2$ есть граница.

{\it Одномерной группой гомологий $H_1(U,\Z)$ разбиения $U$ (с коэффициентами в
$\Z$)} называется группа циклов с точностью до гомологичности.
Обозначим класс гомологичности различающего цикла через
$$D(u,v):=[d(u,v)]\in H_1(U,\Z).$$
Ясно, что отображение $D:V(N)\times V(N)\to H_1(U,\Z)$ уже не зависит от
выбора совмещений полей на $U_0^*$ и поэтому корректно определено.

Так как $d(u,v)+d(v,w)=d(u,w)$, то $D(u,v)+D(v,w)=D(u,w)$.

Чтобы доказать сюръективность отображения $D(\cdot,v)$, возьмем произвольный цикл $\gamma$.
Возьмем поле $u$ равным $v$ на $U_0^*$.
На каждом ребре $a^*$ определим $u$ как $\gamma(a)$-кратную подкрутку поля $v$.
Тогда $d(u,v)=\gamma$, поэтому $D(u,v)=[\gamma]$.
 QED

\smallskip
{\it Доказательство инъективности отображения $D(\cdot,v)$ для $n=2$.}
Из теоремы продолжаемости (см. пункт 'касательные векторные поля: триангуляции')
вытекает следующее утверждение (докажите!).

{\it Два ненулевых касательных к плоскости векторных поля, заданных на отрезке
(лежащем в плоскости) и совпадающие на его концах, гомотопны тогда и только
тогда, когда количество поворотов вектора первого поля при движении от
начала отрезка к его концу равно аналогичному количеству для второго поля.}

По этому утверждению и по теореме гомотопности если $u=v$ на $U_0^*$ и $d(u,v)=0$,
то поля $u$ и $v$ гомотопны неподвижно на $U_0^*$.

Если $D(u_0,v)=D(u_1,v)$ для некоторых полей $u_0$ и $u_1$, совпадающих на
$U_0^*$, то $d(u_0,v)-d(u_1,v)=n_1\partial f_1+\dots+n_s\partial f_s$ для
некоторых целых чисел $n_1,\dots,n_s$.
Тогда $u_0\simeq u_0(n_1,\dots,n_s)\simeq u_1$.
 QED

 \begin{pr}\label{vehomol}
Возьмем произвольное клеточное разбиение $U$ замкнутого связного 2-многообразия $N$.

(a) Дайте определение группы $H_1(U;\Z)$, независимое от рассуждений, в которых она появилась.

(b) Группа $H_1(U;\Z)$ и класс $D(u,v)$ `одинаковы' для клеточных разбиений гомеоморфных комплексов
и пары векторных полей $u,v$.
Более точно, для таких разбиений $U$ и $U'$ существует изоморфизм
$H_1(U;\Z)\to H_1(U';\Z)$, переводящий $D(U,u,v)$ в $D(U',u,v)$.

(Поэтому будем обозначать группу через $H_1(N,\Z)$.)

(c) $H_1(U;\Z)\cong\Z^{2g}$ для разбиения $U$ сферы с $g$ ручками.
\end{pr}

 \begin{pr}\label{}
 (a) Постройте отображение $D:V(N)\to H_1(N;\Z)$ для 3-многообразия $N$,
обобщая построение различающей расстановки.
(Определение группа  $H_1(N;\Z)$ естественно появится в процессе доказательства, поэтому
для решения данной задачи не обязательно знать, что это такое.)

(b) Докажите сюръективность отображения $D$.
\end{pr}

\subsection{Указания к некоторым задачам}

{\bf \ref{hopf-deg}, \ref{hopf-ts2}.} \cite[\S18.1,3]{Pr04}.

\smallskip
{\bf \ref{hopf-n2s2}.} Ответ для связного $N$: $\Z$, если $N$ ориентируемо и $\Z_2$ иначе.
Доказательство см. в \cite[\S18.3]{Pr04}.

\smallskip
{\bf \ref{hopf-tos2}.} Аналогично двумерному случаю.

\smallskip
{\bf \ref{veheul}.} (b) (d) Касательное векторное поле на $P\times Q$ общего положения, если оба его `проекции' на сомножители общего положения.
Пусть $e$ и $f$ --- касательные векторные поля общего положения на $P$ и на $Q$.
Тогда $e+f$ --- касательное векторное поле общего положения на $P\times Q$.

\smallskip
{\bf \ref{vehopf}.} (b) Препятствие $e(-v)$ к продолжению поля $-v$ противоположно по знаку препятствию $e(v)$.
 С другой стороны $e(-v)=e(v)$.

(c) Постройте векторное поле скоростей воды, текущей от северного полюса к южному.
Докажите, что $\chi(S^{2k})=2$.


\smallskip
{\bf \ref{vehnor}.}
(c,d) Пусть $f:N\to\R^m$ --- вложение замкнутого ориентируемого
$n$-многообразия $N$.
Постройте препятствие $\overline e(f)\in H_{2n-m}(N;\Z)$ ({\it нормальное
число Эйлера}) к существованию нормального поля.
Докажите, что $\overline e(f)=0$.

(d) Используйте, что любые два отображения  $S^k\to S^1$ гомотопны при $k>1$.

\smallskip
{\bf \ref{vehnexi}.} (a) Пример получается, если рассмотреть композицию $S^2\times D^1\to S^2\times D^3\to\R^5$ вложения, заданного формулой $(x,t)\to (x,tx)$ и стандартного вложения.
По нормальному полю к такому вложению легко построить касательное поле на $S^2$.

(b) Возьмем образ отображения $\C P^2\to\R^7$, определенного формулой
$(x:y:z)\mapsto(x\overline y, y\overline z, z\overline x,2|x|^2+|y|^2)$, где $|x|^2+|y|^2+|z|^2=1.$
В обозначениях \S\ref{0vesy} имеем $\overline w_2(\C P^2)=w_2(\C P^2)\ne0$.

\smallskip
{\bf \ref{veh35}.} (a) Аналогично предыдущему.

(b) Любой элемент группы $H_2(N;\Z)$ реализуется некоторым замкнутым ориентированным 2-подмногообразием $F\subset N$, см. \S\ref{0homol}.
Поэтому и ввиду двойственности Пуанкаре (\S\ref{0dual}) достаточно доказать, что $\overline e(f)\cap[F]\in \Z$ четно.
Это число есть препятствие к построению нормального к $f(N)$ поля на $F$.
Для его вычета по модулю 2 ввиду формулы Уитни-Ву (\S\ref{0vesy}) имеем $\rho_2(\overline e(f)\cap[F])+w_2(F)=0$, поскольку $5\varepsilon_F=\tau_F\oplus\nu_{F\subset N}\oplus\nu_{N\subset\R^5}|_F$.

\newpage
\section{Наборы векторных полей на 3-многообразиях} \label{0vesy3}

\subsection{О существовании наборов касательных полей}

Ученик Хопфа Эдуард Штифель рассмотрел проблему существования {\it пары, тройки, etc.}
{\it линейно независимых} касательных векторных полей на данном многообразии.

\begin{pr}\label{vesetri}
(a) Если $n$-многообразие ориентируемо, то любой набор из $n-1$ линейно
независимых касательных векторных полей на нем можно дополнить до набора из $n$ таких полей.

(b) Если на многообразии существует набор из $k$ линейно независимых касательных векторных полей, то
 существует набор из $k$ ортонормированных касательных векторных полей.
\end{pr}

 Развивая идеи Хопфа, Штифель около 1934 пришел к определению характеристических классов.
(Любопытно, что Штифель начал с частного случая ориентируемых 3-многообразий
и пытался построить пример такого многообразия, на котором не существует тройки
линейно независимых касательных векторных полей.)
Формализация была завершена Норманом Стинродом в 1940е гг.
При помощи построенной теории были доказаны следующий факт и результаты в начале \S\ref{0vesy}.

\smallskip
{\bf Теорема Штифеля.}
{\it На любом ориентируемом 3-многообразии существует тройка линейно
независимых касательных векторных полей.}

\smallskip
Для доказательства следующего элементарно формулируемого факта также нужны характеристические классы.
Далее в соответствующих местах написано, что можно вернуться к решению этой задачи.

\begin{pr}\label{vesepro}
Для связного 2-многообразия $F$ следующие условия равносильны:

$\bullet$ На $F\times S^1$ существует пара линейно независимых касательных векторных полей;

$\bullet$ На $F\times I$ существует пара линейно независимых касательных векторных полей;

$\bullet$ $F$ имеет непустой край или четную эйлерову характеристику.
\end{pr}

\subsection{Параллелизуемость некоторых трехмерных многообразий}


Далее в этом параграфе пара (тройка, четверка, ...) ортонормированных касательных векторных полей называется просто парой (тройкой, четверкой, ...) полей.

В этом пункте можно использовать без доказательства следующий факт
(впрочем, его тоже можете попытаться доказать).
Обозначим через $SO_3$
пространство положительных ортонормированных реперов в $\R^3$.
(Это то же, что пространство ортонормированных пар векторов в $\R^3$.)

\begin{pr}\label{veseso3}
Любое отображение граничной сферы трехмерного шара в $SO_3$ можно продолжить на весь шар.

Указание. Докажите, что $SO_3$ гомеоморфно

$\bullet$ пространству движений пространства $\R^3$, оставляющих начало координат неподвижным;

$\bullet$ пространству вращений пространства $\R^3$ относительно прямых, проходящих через начало координат;

$\bullet$ пространству прямых в $\R^4$, проходящих через фиксированную точку;

$\bullet$ замкнутому трехмерному шару, на границе которого склеены диаметрально противоположные точки.

Используйте последнюю из перечисленных 'моделей' пространства $SO_3$.
\end{pr}

Мы будем использовать без доказательства то, что на диффеоморфных многообразиях тройки
полей существуют или нет одновременно.

\begin{pr}\label{veses1s2}
Тройка полей существует на

(a) $S^3$; \quad (b) $\R P^3$; \quad (c) $D^1\times S^2$; \quad (d) $S^1\times S^2$;\quad

(e) произведении сферы с ручками и отрезка.

(f)* произведении сферы с ручками и окружности.
\end{pr}

\begin{pr}\label{vesehint}
(a) Если тройка полей существует на дополнении замкнутого 3-многообразия
до некоторого трехмерного шара, то такая тройка существует и на самом 3-многообразии.

(b) Если 3-многообразие с краем погружается в $\R^3$, то на нем есть тройка полей.

(c) Если на сфере с ручками и дыркой $F_0$, лежащей в произведении $F\times I$ сферы с ручками и отрезка,
есть тройка ортонормированных касательных (к $F\times I$) векторных полей, то эта тройка продолжается на $F$ (а значит, и на $F\times I$).
\end{pr}

При решении последнего пункта нужно либо построить клеточное разбиение сферы с ручками на многоугольники,
имеющие общую вершину, либо использовать то, что множество гомотопических классов отображений $S^1\to SO_3$
образует абелеву группу, либо то, что это множество состоит из двух элементов.

\begin{pr}\label{vesekle}
(a) Бутылка Клейна вложима в некоторое ориентируемое 3-многообразие с краем.

(b) На одном из таких 3-многообразий есть тройка полей.
\end{pr}

Мы докажем теорему Штифеля, сведя ее к следующему примеру, обобщающему предыдущую задачу.
Пусть дано замкнутое 2-многообразие $F$.
На $F$ существует семейство непересекающихся окружностей, вне которого имеется ориентация, меняющаяся при пересечении каждой окружности этого семейства.
(Это семейство представляет первый класс Штифеля-Уитни, см. \S\ref{02hom}.)
Разрежем $F$ по этому семейству.
Получим 2-многообразие $F'$ с краем и инволюцией
$\sigma:\partial F'\to\partial F'$, не имеющей неподвижных точек.
Назовем {\it ориентируемым утолщением} 2-многообразия $F$
3-многообразие $F'\times D^1/(x,t)\sim(\sigma(x),-t)_{x\in\partial F', t\in D^1}$.
(Это пространство того $I$-расслоения над $F$, которое является ориентируемым, см. определение в начале \S\ref{0vebun}.)

Если $F$ ориентируемо, то его ориентируемое утолщение гомеоморфно $F\times I$.

\begin{pr}*\label{vesethi}
(a) На любом ориентируемом утолщении любого замкнутого  2-многообразия
есть тройка полей.

(b) Ориентируемое утолщение определено корректно, т.е. не зависит от выбора семейства окружностей.

(c) Если замкнутое 2-многообразие лежит в ориентируемом 3-многообразии,
то некоторая окрестность 2-многообразия гомеоморфна ориентируемому утолщению.
\end{pr}

\begin{pr}\label{veseimm}
(a)* $\R P^2$ допускает погружение в $\R^3$.

(b) Любое 2-многообразие $F$ допускает погружение в $\R^3$.
\end{pr}

\subsection{Характеристические классы для 3-многообразий}\label{0vesy3ch}

 При построении характеристических классов возникают все новые и новые гомологические группы (в т.ч. в задачах).
Читателю полезно продумать эти примеры перед прочтением общего определения в начале \S\ref{0homol}
(куда можно заглянуть, чтобы проверить себя).

Следующий результат является важнейшим шагом в доказательстве теоремы Штифеля,
а для неориентируемых 3-многообразий интересен и сам по себе.
Определения группы $H_1(N)$ и класса $w_2(N)$ естественно возникают далее при попытке построить пару полей по остовам.

\smallskip
{\bf Лемма Штифеля.}
{\it На замкнутом 3-многообразии $N$ существует пара ортонормированных касательных векторных полей
тогда и только тогда, когда} второй класс Штифеля-Уитни $w_2(N)\in H_1(N)$ {\it нулевой.}

\smallskip
{\it Неформальное описание класса $w_2(N)$, использующее общее положение.}
(Это описание не используется в дальнейшем.)
Обозначим через $\Sigma$ подмножество пространства $(\R^3)^2$,
состоящее из пар линейно зависимых векторов.
Подмножество пространства $(\R^3)^2$, состоящее из пар векторов с
непропорциональными первыми координатами, шестимерно (имеет коразмерность 0).
Пересечение $\Sigma$ с этим подмножеством выделяется двумя независимыми уравнениями
(два определителя должны равняться нулю).
Поэтому оно четырехмерно (имеет коразмерность 2).
Аналогично рассматривая вторую и третью координату, получаем, что $\Sigma$
является объединением трех четырехмерных множеств и потому четырехмерно (имеет коразмерность 2).

Пара векторных полей на $\R^3$ --- то же, что отображение $\R^3\to(\R^3)^2$.
Множество точек пространства $\R^3$, в которых векторы этой
пары линейно зависимы, является прообразом подмножества $\Sigma$.
Поэтому для пары векторных полей {\it общего положения} оно является подмногообразием коразмерности 2, т.е. несвязным объединением окружностей (ориентация на них не важна).

Поэтому для пары касательных векторных полей общего положения на 3-многообразии
$N$ множество точек многообразия, в которых векторы этой
пары линейно зависимы, является несвязным объединением окружностей.
Оно называется {\it препятствующим циклом}.
Его класс {\it гомологичности} (определенной далее или в \S13) и называется классом $w_2(N)$.

Надеемся, что читателю интуитивно ясно понятие общего положения и без следующего определения.
Формально, пара касательных векторных полей {\it находится в общем положении}, если соответствующие этим полям сечения касательного расслоения находятся в общем положении друг с другом и с нулевым сечением. Тройка отображений многообразий в многообразие {\it находится в общем положении}, если все их попарные и тройные пересечения локально диффеоморфны соответствующим пересечениям линейных пространств общего положения. (Последнего `общего положения' мы не определяем, уж оно-то понятно читателю.)


\smallskip
{\it Начало рассуждения, приводящее к определению группы $H_1(N)$ и класса $w_2(N)$.}
Возьмем клеточное разбиение $U$ 3-многообразия $N$ на маленькие многогранники.
Обозначим двойственное к нему разбиение через $U^*$.
(Многогранники должна быть настолько мелкими, что угол между касательными пространствами в любых двух точках одного многогранника разбиения $U^*$ меньше $\pi/2$.)

Сначала мы построим пару полей в вершинах двойственного разбиения.
Затем будем пытаться продолжить эту пару на объединение ребер двойственного разбиения.
Потом --- на объединение его граней.
Далее --- на объединение его многогранников.

Попробуем продолжить пару полей на ребро двойственного разбиения.
Ввиду мелкости многогранников касательные пространства в точках ребра можно отождествить.
Поэтому пара полей на ребре --- то же, что отображение ребра в пространство пар полей в $\R^3$, т.е. в $SO_3$.
Очевидно, что это пространство связно.
Поэтому построенная в вершинах двойственного разбиения пара полей продолжается на объединение ребер двойственного  разбиения.

\begin{figure}[h]\centering
\includegraphics{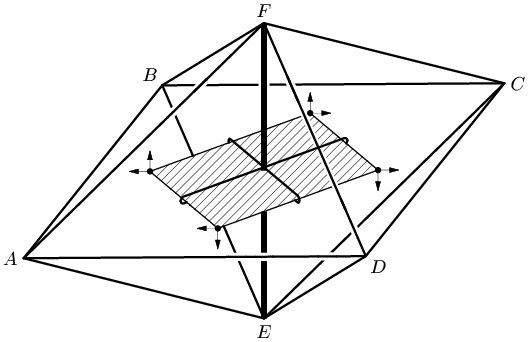}
\caption{Продолжение пары полей на грань двойственного рабиения}\label{nagran}
\end{figure}

Попробуем продолжить пару полей на грань $a^*$ двойственного разбиения (рис.~\ref{nagran}).
Ввиду мелкости многогранников касательные пространства в точках грани можно отождествить.
Поэтому пара полей на грани --- то же, что отображение грани в $SO_3$.
Если это отображение не продолжается на грань $a^*$, то покрасим ребро $a$
исходного разбиения (протыкающее грань $a^*$), в красный цвет.
Таким образом, паре $w$ ортонормированных векторных полей на объединении ребер двойственного разбиения
соответствует набор красных ребер $\varepsilon(w)$ исходного разбиения.
Этот набор ребер называется {\it препятствующим циклом}.

\begin{pr}\label{veseobcy}
(a) Вне препятствующего цикла пара полей существует.

(b) Число красных ребер, выходящих из каждой вершины, четно.

(c) При изменение пары полей на одном ребре двойственного разбиения к препятствующему циклу добавляется граница той двумерной грани исходного разбиения, которая соответствует выбранному ребру.
\end{pr}

\begin{pr}\label{veseobcyex}
(a) Найдите какой-нибудь препятствующий цикл для $F\times S^1$, где $F$ --- 2-многообразие.

(b)* Для любой триангуляции 3-многообразия объединение ребер ее барицентрического подразделения является препятствующим циклом для некоторой пары полей $w$.
\end{pr}

{\it Определение группы $H_1(N)$ и класса $w_2(N)$.}
Напомним, что наборы ребер с условием из задачи \ref{veseobcy}.b называются {\it (гомологическими) циклами}.
Точно так же, как в \S\ref{02hom}, определяются {\it границы}, {\it гомологичность} циклов и группа $H_1(N)$ классов гомологичности циклов (несмотря на то, что эти понятия возникли при решении {\it другой} задачи!).
{\it Вторым классом Штифеля-Уитни} 3-многообразия $N$ называется
$$
w_2(N):=[\varepsilon(w)]\in H_1(N).
$$
Аналогично определениям первого класса Штифеля-Уитни и числа Эйлера проверяется, что это определение корректно.

\begin{pr}\label{vesesti}
 (a) $w_2(N)=0$ тогда и только тогда, когда пару полей $w$ можно продолжить на объединение многоугольников двойственного разбиения.

Указание. Используйте без доказательства следующий факт (впрочем, его тоже
можете попытаться доказать):
существует ровно два гомотопических класса отображений окружности в $SO_3$.
Нетривиальный гомотопический класс представляется диаметром трехмерного
шара, из которого $SO_3$ получается склейкой.

(b) Докажите лемму Штифеля.

(c) Сформулируйте и докажите аналог леммы Штифеля для 3-многообразий с непустым краем.
\end{pr}

Теперь Вы уже можете решить задачу \ref{vesepro}.

\subsection{Простое доказательство теоремы Штифеля}

 {\it Набросок простого доказательства теоремы Штифеля. }
Можно считать, что данное ориентируемое 3-многообразие $N$ замкнуто.
Достаточно доказать, что

{\it Пересечение любого 2-многообразия $F$ в $N$, являющегося объединением некоторого количества граней некоторой триангуляции многообразия $N$,
и любого представителя класса $w_2(N)$, являющегося набором ребер двойственной триангуляции, состоит из четного числа точек.}

(Этого достаточно ввиду двойственности Пуанкаре и реализуемости элементов группы $H_2(N)$ вложениями, см. \S\S\ref{0homol},\ref{0dual};
группа $H_2(N)$ определена там.)

Четность числа точек такого пересечения является препятствием к {\it существованию
пары линейно независимых векторных полей на $F$, касательных к $N$} (а не к $F$!).

(Это препятствие $w_2(\tau_N|_F)\in\Z_2$ строится аналогично классу $w_2(N)$ с использованием общего положения или разбиения поверхности $F$ на многоугольники.
Свойство $w_2(\tau_N|_F)=w_2(N)\cap F$ проверяется очевидно.)

Это препятствие нулевое, поскольку такая пара полей существует.
Ее достаточно построить на трубчатой окрестности поверхности $F$ в $N$.
Так как $N$ ориентируемо, то эта окрестность гомеоморфна ориентируемому утолщению поверхности $F$   (задача \ref{vesethi}.c).
Значит, на этой окрестности
имеется пара полей (задача \ref{vesethi}.a).
QED

\smallskip
{\it Другое простое доказательство теоремы Штифеля.
\footnote{Это доказательство получено из обычно приводимого в книгах отбрасыванием ненужных обозначений.
Хотя оно немного сложнее вышеприведенного, оно интересно тем, что допускает обобщения.}
}
Аналогично предыдущему доказательству достаточно установить тривиальность препятствия $w_2(N)\cap F\in\Z_2$ к
{\it существованию пары линейно независимых векторных полей на $F$, касательных к $N$}.
Она следует из
$$w_2(N)\cap F=\rho_2\chi(F)+w_1(F)^2=0.$$
Здесь второе равенство доказано в задаче \ref{intstwh}.b.
Докажем первое равенство (аналогично формуле для класса Штифеля-Уитни произведения, задача \ref{vesnww}.a ниже).
Возьмем разбиение 2-многообразия $F$.
Рассмотрим такое касательное к $F$ поле $u$ на $F$, что $u\ne0$ вне некоторой точки $p\in F$.
Возьмем такую пару $v,v'$ касательных к $F$ полей на $F$, что

$\bullet$ $u\perp v$, $u\perp v'$,

$\bullet$ $v=0$ на объединении ребер некоторого представителя $\omega$ класса $w_1(F)$,
лежащего в исходном разбиении, причем при переходе через $\omega$ ориентация базиса $(u,v)$ меняется;

$\bullet$ $v'=0$ на объединении ребер некоторого представителя $\omega'$ класса $w_1(F)$,
лежащего в двойственном разбиении,
причем при переходе через $\omega'$ ориентация базиса $(u,v')$ меняется;

$\bullet$ $p\not\in\omega\cup\omega'$.

Так как $N$ ориентируемо, то определено векторное поле $u\times v'$, касательное к $N$ и нормальное к $F$.
Множество линейной зависимости пары полей $(u,v+u\times v')$ на $F$, касательных к $N$, есть объединение
множеств тех точек, для которых $u=0$, и тех точек, для которых обе пары $(u,v)$ и $(u,v')$ линейно зависимы.
Первое множество есть $p$, второе --- $\omega\cap\omega'$.
Точка $p$ входит в построенный представитель класса $w_2(N)\cap F$ с коэффициентом $\rho_2\chi(F)$,
а каждая точка из $\omega\cap\omega'$ --- с коэффициентом 1 (докажите!). Это доказывает первое равенство.
QED


\begin{pr}*\label{vesesimp} Используя задачу \ref{veseobcyex}.b, докажите, что
$w_2(N)=0$ для любого ориентируемого 3-многообразия $N$. (Ср. с задачей \ref{vesew2w1})
 \end{pr}

\subsection{Указания и решения к некоторым задачам}


\smallskip
{\bf \ref{vesetri}.}
(b) Ввиду каноничности процесса (Грама-Шмидта) ортогонализации.

\smallskip
{\bf \ref{vesepro}.}
Следует из леммы Штифеля и задачи \ref{veseobcyex}.a.
Класс $[*\times S^1]\in H_1(F\times S^1)$ ненулевой, поскольку его пересечение с классом
$[F\times*]\in H_2(F\times S^1)$ ненулевое.

(Для решения этой задачи группу  $H_1(F\times S^1)$ вычислять не нужно!)

\bigskip
{\bf \ref{veseso3}.}
Обозначим через $\pi_i(X)$ множество гомотопических классов непрерывных
отображений $S^i\to X$ (см. определение в начале \S\ref{0hopf}).
Тогда приведенный факт равносилен равенству $\pi_2(SO_3)=0$.
Указание для специалистов: $\pi_2(SO_3)=\pi_2(\R P^3)=\pi_2(S^3)=0$.

\smallskip
{\bf \ref{veses1s2}.} (a) {\it Одно решение} получается, если вспомнить, что $S^3$ есть группа
единичных кватернионов.

{\it Другое решение} получается, если сначала построить
тройку полей на нижней полусфере, а потом продолжить ее на верхнюю, используя задачу \ref{vesehint}.a.

(b) {\it Одно решение} получается, если вспомнить, что $\R P^3$
есть факторгруппа группы единичных кватернионов.

{\it Другое решение} получается, если сначала построить тройку полей
на дополнении до трехмерного шара, которое является ориентируемым
утолщением пространства $\R P^2$ (см. задачу \ref{vesethi}), а потом
продолжить ее на $\R P^3$, используя задачу \ref{vesehint}.a.

(с) Используйте вложение $D^1\times S^2\to\R^3$.

(d) {\it Одно решение.} Дополнение до трехмерного шара есть
$D^1_+\times S^2\cup_{\partial D^1_+\times D^2_+}D^1_\times D^2_+$.
 Постройте погружение этого дополнения в $\R^3$.
Используя его, можно построить тройку полей на этом дополнении (ср. с задачей \ref{vesehint}.b).
Потом можно продолжить ее на $\R P^3$, используя  задачу \ref{vesehint}.a.

{\it Другое решение} Аналогично (f).

(e) Аналогично (c).

(f) Обозначим данную сферу с ручками через $F$.
Так как на $F_0$ есть пара ортонормированных касательных (к $F_0$) векторных полей, то
на $F_0\times S^1$ есть тройка ортонормированных касательных векторных полей.
(Наличие такой тройки вытекает также из задачи \ref{vesehint}.b.)
По задаче \ref{vesehint}.c эта тройка продолжается на  $F_0\times S^1\cup F\times I$.
Значит, по задаче \ref{vesehint}.a эта тройка продолжается и на  $F\times S^1$.

\smallskip
{\bf \ref{vesehint}.} (a) Следует из задачи \ref{veseso3}.

(b) Нужная тройка полей приходит из тройки ортонормированных полей на $\R^3$.

(c) Препятствие в $\pi_1(SO_3)\cong\Z_2$ к существованию тройки ортонормированных векторных полей,
касательных к $F\times I$ (а не к $F$!) строится аналогично классу Эйлера.
По задаче \ref{veses1s2} такая тройка существует.
Значит, это препятствие нулевое.
Поэтому такая тройка, заданная на $F_0$, продолжается на $F$.

\smallskip
{\bf \ref{vesekle}.} (a) Пусть $S^1=\{z\in\Bbb C\ :\ |z|=1\}$.
Тогда 3-многообразие $\dfrac{S^1\times [-1,1]\times [0,1]}{(z,t,0)\sim(\bar z,-t,1)}$ ориентируемо и содержит
бутылку Клейна $S^1\times0\times[0,1]/(z,0,0)\sim(\bar z,0,1)$.

(b) Так как $\R^3$ ориентируемо, то погружение бутылки Клейна в $\R^3$ продолжается до погружения
(т.е. локального вложения) в $\R^3$ ориентируемого 3-многообразия, построенного в (a).
Тогда наличие тройки полей следует из \ref{vesehint}.b.

\smallskip
{\bf \ref{vesethi}.} (a) Аналогично \ref{vesekle}.b, используя \ref{veseimm}.b.

\smallskip
{\bf \ref{veseimm}.} (a) http://en.wikipedia.org/wiki/Boy's\_surface

(b) Следует из (a) и теоремы классификации 2-многообразий.

\smallskip
{\bf \ref{veseobcy}.} (a) Рассмотрим окрестность объединения препятствующего цикла и множества вершин двойственной  триангуляции.
Дополнение до этой окрестности является окрестностью объединения двумерных многоугольников двойственного разбиения,
не являющихся красными.
Пару полей на этом дополнении можно продолжить на дополнение до препятствующего цикла по задаче \ref{vesehint}.a.

(b) Для данной вершины двойственной триангуляции рассмотрим граничную сферу соответствующего многогранника исходной триангуляции.
Четность количества ее двумерных граней, протыкаемых красными ребрами, равна сумме гомотопических классов отображений из границ граней в $SO_3$ и потому равна нулю.

(c) При изменении пары полей на одном ребре двойственной триангуляции на элемент  $x\in\pi_1(SO_3)\cong\Z_2$ на каждой грани двойственной триангуляции, соседней с этим ребром,
добавляется $x$.

\smallskip
{\bf \ref{veseobcyex}.} (a) Пусть $v$ --- поле на объединении ребер некоторого разбиения 2-многообразия $F$, для которого ненулевые элементы препятствующей расстановки стоят в вершинах $p_1,\dots,p_n$ двойственного разбиения и равны $\sgn\chi(F)$ (т.о. $n=|\chi(F)|$).
Пусть $v'$ --- единичное векторное поле на $S^1$. Тогда для `призматического' разбиения произведения $F\times S^1$ и пары $(v,v')$ препятствующий цикл является объединением окружностей
$p_i\times S^1$ по $i=1,\dots,n$.
(Для доказательства этого используйте то, что гомоморфизм включения $\Z\cong\pi_1(SO_2)\to\pi_1(SO_3)\cong\Z_2$ есть приведение по модулю 2.)

Замечание. Набор из $k$ касательных векторных полей на $X\times Y$ является набором общего положения, если обе его `проекции' на сомножители  --- общего положения.
(Поэтому если $e$ --- касательное векторное поле общего положения на $F$, а $e'$ --- единичное векторное поле на $S^1$, то пара $(e,e')$ не общего положения.)

\smallskip
{\bf \ref{vesesti}.} (b) Следует из (a) и задачи \ref{vesehint}.a.

(c) На 3-многообразии $N$ существует пара линейно независимых
касательных векторных полей тогда и только тогда, когда второй класс
Штифеля-Уитни $w_2(N)\in H_1(N,\partial N)$ нулевой.

\newpage
\section{Гомологии трехмерных многообразий}\label{03hom}

\subsection{Гомологии трехмерных многообразий}\label{03homhom}

\begin{pr}\label{3-home}
(a) $S^3$, $S^1\times S^2$, $S^1\times S^1\times S^1$ и $\R P^3$ попарно не гомеоморфны.

(b)  Если произведение сферы с $g$ ручками на окружность гомеоморфно
произведению сферы с $h$ ручками на окружность, то $g=h$.

(c) Если  $P$ и $Q$  2-многообразия (возможно, неориентируемые и имеющие край) и
$P\times S^1\cong Q\times S^1$, то $P\cong Q$.

(d) $S^3$ не гомеоморфно произведению 2-многообразия на окружность.

(e) $\R P^3$ не гомеоморфно произведению 2-многообразия на окружность.
\end{pr}

Для решения этой задачи \ref{3-home} нужны следующие понятия.
Они естественно возникли в \S\ref{0vesy3} и \S\ref{0vehima} при изучении векторных полей на 3-многообразиях.

{\it Одномерной группой гомологий $H_1(T)$ с коэффициентами $\Z_2$ клеточного разбиения $T$ 3-комплекса}\ называется
одномерная группа гомологий с коэффициентами $\Z_2$ двумерного остова
разбиения $T$.
В отличие от двумерного остова, эта группа (и аналогичные группы, определенные далее) инвариантна относительно гомеоморфности (т.е., относительно операции подразделения ребра) для 3-комплекса $T$ (задача \ref{3-homolindep}.a).

\begin{pr}\label{3-homolexa} (1)-(5) Найдите $H_1(T)$ для (построенных Вами) клеточных разбиений
$T$ из примеров (1)-(5) в \S\ref{03man}.

(6) $H_1(U)\cong H_1(P)$ для регулярной окрестности $U$ 2-комплекса $P$ в 3-комплексе $U$.

(a) $H_1(P\times I)\cong H_1(P)$, $H_1(P\times S^1)\cong H_1(P)\oplus\Z_2$ и $H_1(P\times G)\cong H_1(P)\oplus H_1(G)$ для любых связного 2-комплекса $P$, графа $G$ и соответствующих клеточных разбиений произведений $P\times I$,
$P\times S^1$ и $P\times G$.
\end{pr}

\begin{pr}\label{3-homolindep}
(a) Одномерные группы гомологий с коэффициентами $\Z_2$ гомеоморфных 3-комплексов изоморфны.

(b) Если $T'$ --- клеточное разбиение 3-комплекса $T$, то $H_1(T)\cong H_1(T')$.

(с) $H_1(M\#N)\cong H_1(M)\oplus H_1(N)$ для любых 3-многообразий $M$ и $N$.
\end{pr}

{\it Определение группы $H_1(T;\Z)$ для 2-комплекса $T$ с ориентированными ребрами.}
Расстановка целых чисел на ориентированных ребрах 2-комплекса $T$ называется {\it циклом},
если для каждой вершины сумма чисел на входящих ребрах равна сумме чисел на выходящих.
{\it Границей $\partial a$ грани $a=\{i,j,k\}$} называется расстановка {\it плюс} единиц на
ориентированных ребрах $(ij),(jk),(ki)$ и нулей на остальных ребрах.
Это означает, что если ребро $\{i,j\}$ ориентировано от $j$ к $i$, то на нем ставится {\it минус} единица.
Граница грани определена с точностью до умножения на $-1$.
Назовем {\it границей} линейную комбинацию границ нескольких граней с целыми коэффициентами.
Два цикла называются {\it гомологичными}, если их разность есть граница.
{\it Одномерной группой гомологий $H_1(T;\Z)$ с коэффициентами $\Z$ 2-комплекса $T$ с ориентированными
ребрами} называется группа циклов с точностью до гомологичности.

{\it Одномерной группой гомологий $H_1(T;\Z)$ с коэффициентами $\Z$ 3-комплекса $T$ с ориентированными
ребрами}\ называется одномерная группа гомологий с коэффициентами $\Z$ его двумерного остова.

\begin{pr}\label{3-homolzcx}
(1)-(3) Найдите $H_1(T;\Z)$ для (построенных Вами) 3-комплексов $T$ из примеров (1)-(3) выше.

(a) Группы $H_1(T;\Z)$ для разных наборов ориентаций ребер и граней изоморфны.

(b) Группы $H_1(T;\Z)$ гомеоморфных 3-комплексов изоморфны.
\end{pr}

{\it Определение приклеивающего слова грани $a$ клеточного разбиения $T$ 2-комплекса $T$ с ориентированными ребрами.}
Рассмотрим {\it дополнение 2-комплекса до открытой регулярной окрестности одномерного остова $T^{(1)}$ разбиения $T$}, т.е. объединение граней второго барицентрического подразбиения 2-комплекса, не пересекающих $T^{(1)}$.
Ввиду клеточности это дополнение является несвязным объединением дисков, соответствующих граням разбиения $T$.
Будем идти по краевой окружности диска, соответствующего грани $a$.
Для каждого ребра $e$ в геометрической границе грани $a$ будем дописывать к слову $e$, если проходим около этого ребра
в направлении, совпадающем с направлением (=ориентацией)
ребра, и $e^{-1}$, если проходим около этого ребра в направлении, противоположном направлению ребра.
Полученное слово $e_1^{\varepsilon_1}\dots e_k^{\varepsilon_k}$ называется {\it приклеивающим словом грани $a$}.
Здесь $e_1,\dots,e_k$ --- ориентированные ребра разбиения и $\varepsilon_1,\dots,\varepsilon_k\in\{+1,-1\}$.
Приклеивающее слово определено с точностью до умножения всех $\varepsilon_s$ на $-1$.

\begin{pr}\label{3-homoword}
Напишите приклеивающие слова
граней клеточных разбиений из примеров (1)-(5) в \ref{03hom}.
\end{pr}

{\it Определение группы $H_1(T;\Z)$ клеточного разбиения $T$ 2-комплекса с ориентированными ребрами.}
Расстановка целых чисел на ориентированных ребрах разбиения $T$ называется {\it циклом}, если для каждой вершины
сумма чисел на входящих в нее ребрах равна сумме чисел на выходящих (при этом одно и то же ребро может быть и входящим, и выходящим).
Для грани $a$ возьмем ее приклеивающее слово $e_1^{\varepsilon_1}\dots e_k^{\varepsilon_k}$ в разбиении $T$.
{\it Границей $\partial a$ грани $a$} называется расстановка $\varepsilon_1e_1+\dots+\varepsilon_k e_k$.
Назовем {\it границей} линейную комбинацию границ некоторых граней с целыми коэффициентами.
Два цикла называются {\it гомологичными}, если их разность есть граница.
{\it Одномерной группой гомологий $H_1(T;\Z)$ с коэффициентами $\Z$ клеточного разбиения $T$
с ориентированными ребрами} называется группа циклов с точностью до гомологичности.

{\it Одномерной группой гомологий $H_1(T;\Z)$ с коэффициентами $\Z$ клеточного разбиения $T$ 3-комплекса}\ называется
одномерная группа гомологий с коэффициентами $\Z$ двумерного остова
разбиения $T$.

\begin{pr}\label{3-homolzcxcell}
(1)-(5)
Найдите $H_1(T;\Z)$ для (построенных Вами) клеточных разбиений $T$ из примеров (1)-(5) в \ref{03man}.

(a) Группы $H_1(T;\Z)$ для разных наборов ориентаций ребер клеточного разбиения изоморфны.
\end{pr}

\begin{pr}\label{3-homolexaz} (6,a) Сформулируйте и докажите аналоги задач \ref{3-homolexa}.6,a для $H_1(T;\Z)$.
\end{pr}

\begin{pr}\label{3-homolindepz} (a,b,c) Сформулируйте и докажите аналоги задач \ref{3-homolindep}.a,b,c для $H_1(T;\Z)$.

(d) Одномерная группа гомологий с коэффициентами $\Z$ клеточного разбиения с ориентированными ребрами, имеющего одну
вершину, изоморфна абелевой группе, образующими которой являются ориентированные ребра, а соотношениями ---
приклеивающие слова граней.
\end{pr}

\begin{pr}\label{3-homol} Существует ли замкнутое 3-многообразие $N$, для которого $H_1(N)\cong$

(a) $\Z$ ? \qquad (b) $\Z^2$ ? \qquad (c) $\Z^3$ ? \qquad (c) $\Z\oplus\Z_2$ ? \qquad (d) $\Z_2\oplus\Z_3$ ?
\qquad (e) $\Z_4$ ?
\end{pr}

\begin{pr}\label{3-lens}
Любая конечно порожденная абелева группа изоморфна одномерной группе гомологий с коэффициентами в $\Z$ некоторого
замкнутого 3-многообразия.
\end{pr}

\begin{pr}\label{3-homreal}
(a) Для связного 3-комплекса $X$ любой элемент группы $H_1(X;\Z)$ {\it представляется}
некоторым отображением $f:S^1\to X$.

(b) Отображение $f:S^1\to X$ представляет нулевой элемент группы $H_1(X;\Z_2)$ тогда и только тогда, когда оно
продолжается на некоторое 2-многообразие $N$
с краем $\partial N=S^1$.

(c) Отображение $f:S^1\to X$ представляет нулевой элемент группы
$H_1(X;\Z)$ тогда и только тогда, когда оно продолжается на некоторое
ориентируемое 2-многообразие $N$ с краем $\partial N=S^1$.
\end{pr}


\begin{pr}\label{3-lensho}* (a) $L(5,1)\not \cong L(5,2)$. \qquad (b) $L(7,1)\not \cong L(7,2)$.
\end{pr}

{\bf Теорема Александера-Райдемайстера.}
{\it $L(p_1,q_1)\cong L(p_2,q_2)$ тогда и только тогда, когда $p_1=p_2$ и $q_1\equiv \pm q_1^{\pm1}\mod p$
(знаки $\pm$ не обязательно согласованы).}


\begin{pr}
Определите и классифицируйте двулистные накрытия над данным 3-многообразием.
\end{pr}

\begin{pr}*\label{vesew2w1} Используя задачу \ref{veseobcyex}.b, докажите, что
  $w_2(N)=w_1^2(N)$ для любого 3-многообразия $N$.
\end{pr}

\subsection{Фундаментальная группа и накрытия (набросок)}

Для полиэдра (т.е. тела симплициального комплекса) $X$ с отмеченной точкой $x_0$ {\it петлей} называется
отображение $f:[0,1]\to X$, для которого $f(0)=f(1)=x_0$.

\begin{pr}\label{hopf-equ}
Гомотопия в классе петель является отношением эквивалентности на множестве петель $[0,1]\to X$.
\end{pr}

\begin{figure}[h]\centering
\includegraphics{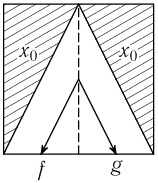}\qquad \includegraphics{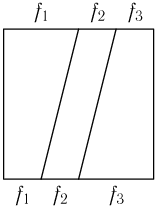}
\caption{Обратный элемент. Ассоциативность умножения.}\label{pr04-23}
\end{figure}

{\it Фундаментальной группой} $\pi_1(X,x_0)$ полиэдра $X$ с отмеченной точкой $x_0$ называется группа классов
гомотопности (в классе петель) петель.
Отмеченную точку часто не указывают в обозначениях.
Групповая структура на множестве классов гомотопности петель вводится следующим образом.
Определим {\it произведение} $[f_1][f_2]$ как гомотопический класс отображения
$$g(t):=\begin{cases} f_1(2t) & 0\le t\le1/2, \\
f_2(2t-1) & 1/2\le t\le1.\end{cases}$$
{\it Единичным элементом} называется гомотопический класс постоянного отображения.
Обратным элементом к $[f]$ называется гомотопический класс отображения
$\overline f(t):=f(1-t)$.
Равенства $[\bar f][f]=[f][\bar f]=e$ и $([f_1][f_2])[f_3]=[f_1]([f_2][f_3])$ ясны из рисунков~\ref{pr04-23}.

\begin{pr}\label{hopf-bp}
Если полиэдр $X$ 
связен, то $\pi_1(X,x_0)\cong \pi_1(X,x_1)$ для любых двух точек $x_0,x_1\in X$.
\end{pr}

\begin{pr}\label{hopf-exa}
Найдите фундаментальные группы следующих полиэдров

(a) кольца $S^1\times I$; \quad (b) листа Мебиуса; \quad (c) тора $S^1\times S^1$;

(d) проективной плоскости $\R P^2$; \quad
(e) `восьмерки' $S^1\vee S^1$; \quad

(f) бутылки Кляйна с дыркой; \quad
(g) бутылки Кляйна;  \quad (h) $\R P^3$; \quad (i) $L(p,q)$.

(j) the complement in $\R^3$ of the trivial knot;
\quad
(k) the complement in $\R^3$ of the trefoil knot.
\end{pr}

Для решения некоторых пунктов этой задачи нужны следующие понятия и результаты.

\begin{pr}\label{pi1-prop} (a) $\pi_1(X\times Y)\cong \pi_1(X)\oplus \pi_1(Y)$.

(b) Если полиэдр $X$ есть {\it деформационный ретракт} полиэдра $Y$ (ср. с задачей \ref{vfann})
или $X$ {\it сдавливается} на $Y$, то $\pi_1(X)\cong \pi_1(Y)$.
(Ср. со вторым способом решения задачи \ref{orihocl}.b.)
\end{pr}

\begin{pr}\label{hopf-s1s2}
(a) Любое отображение $S^1\to S^2$, образ которого не совпадает с $S^2$, гомотопно отображению в точку.

(b) Любое отображение $S^1\to S^2$ гомотопно кусочно-линейному (определите, что это такое) или гладкому
(определите, что это такое).

(с) Любое отображение $S^1\to S^2$ гомотопно отображению в точку.

(d) Любое отображение $S^1\to S^3$ гомотопно отображению в точку.
\end{pr}


Отображение $p:\t X\to X$ между полиэдрами (или, более общо, между подмножествами в $\R^n$) называется
{\it накрытием}, если для любой точки $x\in X$ существуют
ее окрестность $Ox$, число $n$, подмножество $F$ множества $\Z^n\subset\R^n$,
и гомеоморфизм $h:p^{-1}Ox\to Ox\times F$, для которого $\pro_F\circ h=p$.

\begin{pr}\label{hopf-cov} Ср. с задачей \ref{vfcov}.
Пусть $p:\t X\to X$ --- накрытие и $\t X,X$ связны.

(a) Для любых пути $s:[0,1]\to X$ и точки $\t x\in\t X$ с условием $p(\t x)=s(0)$ существует путь $\t s:[0,1]\to\t X$
({\it поднятие} пути $s$) такой, что $\t s(0)=\t x$ и $p\circ\t s=s$.

(b) Такое поднятие единственно, т.е. если $\t f_1,\t f_2:[0,1]\to\t X$ --- два поднятия одного и того же отображения
$f:[0,1]\to X$, причем $\t f_1(0)=\t f_2(0)$, то $\t f_1(t)=\t f_2(t)$ для любого $t$.

(c) {\it Лемма о поднятии гомотопии.}
Для любых гомотопии $\{f_t:[0,1]\to X\}_{t\in [0,1]}$ и поднятия $\t f_0:[0,1]\to\t X$
существует единственное поднятие $\{\t f_t:[0,1]\to\t X\}_{t\in [0,1]}$ гомотопии $f_t$.

(d) Если $\pi_1(\t X)=0$, то $|\pi_1(X)|=|p^{-1}(x_0)|$ для любой точки $x_0\in X$.

(e) $|\pi_1(X)|=|p^{-1}(x_0)|\cdot|\pi_1(\t X)|$ для любой точки $x_0\in X$.
\end{pr}

\begin{pr}\label{pi1-covsim} Пусть $p:\t X\to X$ --- накрытие, $\t X$ и $X$ связны, $x_0\in X$,
$\t x_0\in p^{-1}(x_0)$ и $\pi_1(\t X)=0$.

(a) Следующее определение задает групповую структуру в $p^{-1}(x_0)$.
Положим единичный элемент равным $x_0$.
Для точек $a,b\in p^{-1}(x_0)$ выберем пути $s_a,s_b:[0,1]\to\t X$, соединяющие $\t x_0$ с $a$ и $b$.
Положим $ab$ равным концу того поднятия пути $(s_a\circ p)(s_b\circ p)$, которое начинается в $\t x_0$.
Положим $a^{-1}$ равным концу того поднятия пути $(s_a\circ p)^{-1}$, которое начинается в $\t x_0$.

(b) $\pi_1(X)\cong p^{-1}(x_0)$.
\end{pr}

\begin{pr}\label{pi1-poin}
(a) {\it Теорема.} Фундаментальная группа тела клеточного разбиения с ориентированными ребрами, имеющего одну вершину,
изоморфна группе, образующими которой являются ориентированные ребра, а соотношениями --- приклеивающие слова граней.

(b) {\it Теорема Зейферта--Ван Кампена.} Пусть $U_1$ и $U_2$ --- тела подкомплексов комплекса с телом $X$.
Предположим, что группы $\pi_1(U_1)$, $\pi_1(U_2)$, $\pi_1(U_1\cap U_2)$ заданы образующими $S_1, S_2,S_0$ и
соотношениями $R_1,R_2,R_0$, соответственно, причем $S_k\cap S_l=\emptyset$ для любых $k,l$.
Тогда группа $\pi_1(U_1\cap U_2)$ задана образующими $S_1\cup S_2$ и соотношениями
$R_1,R_2,\{i_{1*}x=i_{2*}x\}_{x\in S_0}$, где $i_k:U_1\cap U_2\to U_k$ --- включение.
(Знать $R_0$ для нахождения $\pi_1(U_1\cap U_2)$ не обязательно.)

(c) {\it Теорема Пуанкаре.}
Для связного полиэдра $X$ выполнено $H_1(X;\Z)\cong\pi_1(X)/[\pi_1(X),\pi_1(X)].$
\end{pr}

\begin{pr}\label{pi1-tref}
Дополнения в $S^3$ до стандартной окружности и до трилистника не гомеоморфны.

(Из этого следует неизотопность тривиального узла и трилистника.)
\end{pr}

{\bf Теорема Адяна-Рабина.} {\it Не существует алгоритма распознавания, изоморфны ли две группы, заданные своими
конечными  представлениями.

Более формально, не существует алгоритма вычисления функции, которая по паре конечных представлений групп выясняет,
изоморфны ли эти группы.}



\begin{pr}\label{3-poin} {\it Сфера Пуанкаре.}
Если у трехмерного додекаэдра склеить противоположные грани с поворотом на $\pi/5$,
то получится 3-многообразие $N$, для которого $H_1(N;\Z)=0$, но $\pi_1(N)\ne0$.
(Имеется много других конструкций сферы Пуанкаре, см., например, задачу \ref{3-dehn}.e.)
\end{pr}


{\bf Теорема Перельмана.} {\it Если $\pi_1(N)=0$ для замкнутого 3-многообразия $N$, то $N\cong S^3$.}

\begin{pr}\label{3-pi1} Сушествует ли замкнутое 3-многообразие $N$, для которого $\pi_1(N)\cong$

(a) $\Z$ ?  \qquad (b) $\Z^3$ ?  \qquad (c) $\Z\oplus\Z_2$ ? \qquad (d) $\Z^2$ ? \qquad (e) $\Z^4$ ?
\end{pr}

\begin{pr}\label{3-bala} Любое связное 3-многообразие можно разбить на один многогранник так, чтобы была
ровно одна вершина и число одномерных граней было равно числу двумерных.
\end{pr}

\begin{pr}\label{3-fige}* Среди конечно порожденных абелевых групп только $\Z,\ \Z^3,\ \Z_n,\ \Z\oplus\Z_2$
изоморфны фундаментальной группе некоторого 3-многообразия.
\end{pr}

{\bf Теорема Столлингса.} \cite{Ma03} {\it Не существует алгоритма распознавания, является ли группа,
заданная конечным  представлением, фундаментальной группой некоторого замкнутого 3-многообразия.}


\smallskip
{\bf Теорема.} {\it Cуществует алгоритм распознавания тривиальности группы, заданной конечным  представлением
и являющейся фундаментальной группой некоторого (не заданного) замкнутого 3-многообразия.}

\smallskip
Это вытекает из \cite{So04} и теоремы Перельмана.

\comment

(2) Существование алгоритма, распознающего изоморфны ли заданные копредставлениями
фундаментальные группы 3-многообразия вытекает из положительного решения Перельманом
гипотезы Пуанкаре и моей работы
в которой доказано, что из алгоритмической неразрешимости этой проблемы следует, что гипотеза
Пуанкаре неверна. Прочитав работу (доказательство довольно короткое и -- мне кажется -- красивое)
можно понять, что некоторый алгоритм существует, но он безнадежно переборный.

По поводу вопроса о распознавании изоморфности групп, если известно, что они являются фундаментальными  группами
трехмерных многообразий. Этот вопрос можно понимать двояко: многообразия указаны или только сказано, что они есть.
Второй случай сводится к первому по Сосинскому, тривиальность группы здесь не причем.
А вот справедливость гипотезы Пуанкаре помогает определять только тривиальность группы (с помощью алгоритма
распознавания сферы по Рубинштейну-Томпсон).
Распознать изоморфность нетривиальных групп это не помогает. Так что этот вопрос скорее всего открыт.

Ваш Сергей Матвеев

АС: Задано 3-многообразие.

\endcomment

\subsection{Регулярные накрытия (набросок)}

Для отображения $\varphi:X\to Y$ между полиэдрами определим отображение $\varphi_*:\pi_1(X,x_0)\to\pi_1(Y,\varphi(x_0))$
формулой
$$\varphi_*[f:S^1\to X]:=[(\varphi\circ f):S^1\to Y].$$

\begin{pr}\label{hopf-inj}
(a) Придумайте инъективное отображение $\varphi:X\to Y$ между полиэдрами, для которого $\varphi_*$ не  инъективно.

(b) Если $p:\t X\to X$ --- накрытие, то $p_*:\pi_1(\t X)\to\pi_1(X)$ инъективно.
\end{pr}

Накрытие $p:\t X\to X$ называется {\it регулярным}, если $\t X,X$ связны и для любой петли $s:[0,1]\to X$ ее
поднятия, начинающиеся разных точках множества $p^{-1}(s(0))$, одновременно замкнуты или
незамкнуты.


\begin{pr}\label{hopf-reg} (a) Для связных $\t X,X$ регулярность накрытия $p:\t X\to X$ равносильна нормальности
подгруппы $p_*\pi_1(\t X)$ в $\pi_1(X)$.

(b) Приведите пример нерегулярного накрытия между связными графами.
\end{pr}

\begin{pr}\label{hopf-cov2} Пусть $x_0\in X$ и $p:\t X\to X$ --- регулярное накрытие.

(a) Для любой точки $\t x_0\in p^{-1}(x_0)$ определение из задачи \ref{pi1-covsim} задает групповую структуру в
$p^{-1}(x_0)$.

(b) {\it Теорема о накрытии.}
$p_*:\pi_1(\t X)\to\pi_1(X)$ мономорфизм и $\pi_1(X)/p_*\pi_1(\t X)\cong p^{-1}(x_0).$
\end{pr}

\subsection{Конструкции трехмерных многообразий}


\begin{pr}\label{3-dehn} Пусть $g:S^1\to S^3$ --- вложение.

(a) $g$ продолжается до такого вложения $G:S^1\times D^2\to S^3$, что $G(S^1\times0)=g(S^1)$ и $G(S^1\times a)$,
$a\in S^1$ имеют нулевой коэффициент зацепления.

(b) {\it Хирургия Дена.} Для $n\in\Z$ обозначим через $f_n$ автогомеоморфизм тора $S^1\times S^1$,
переводящий $S^1\times*$ в себя и $*\times S^1$ в кривую, гомологичную $n(S^1\times*)+(*\times S^1)$.
Положим
$$\chi_{g,n}(S^3):=(S^3-\Int G(S^1\times D^2))\cup_{f_n} D^2\times S^1.$$
(Это 3-многообразие не зависит от $G$ при данном $g$, чего мы не доказываем.)

Докажите, что $H_1(\chi_{g,n}(S^3))\cong \Z/n\Z$.

(c)  $\chi_{g,n}(S^3)\cong L(n,1)$ для стандартного узла $g$.

(d) Для трилистника $g:S^1\to S^3$ имеем $H_1(\chi_{g,1}(S^3))=0$, но $\chi_{g,1}(S^3)\not\cong S^3$.

(e) Определите аналогично хирургию Дена по зацеплению с набором чисел на компонентах.
Докажите, что любое 3-многообразия можно получить
при помощи такой хирургии.
\end{pr}

Аналогично определяется $\chi_{G,n}(M)$ для любых 3-многообразия $M$, вложения $G:S^1\times D^2\to M$ и $n\in\Z$.

Произведения 2-многообразий на окружность обобщаются до $S^1$-расслоений над 2-многообразиями (\S\ref{0vebun}) и до
расслоений над $S^1$ со слоем 2-многообразие.

(7) $S^1$-расслоение над 2-многообразием с краем является прямым произведением на $S^1$.
Тотальное пространство $S^1$-расслоения над замкнутым 2-многообразием $N$ гомеоморфно
\linebreak
$\chi_{S^1\times*,n}(S^1\times N)$ для некоторого $n$.

\begin{pr}\label{3-s1bun} Классифицируйте $S^1$-расслоения над замкнутыми 2-многообразиями с точностью

(a) до гомеоморфизма их тотальных пространств;
\qquad
(b) до изоморфизма расслоений.
\end{pr}

(8) Расслоение над $S^1$ со слоем 2-многообразие определяется как
$$S^1\widetilde\times_\sigma F:=F\times I/(x,0)\sim(\sigma(x),1)_{x\in F}$$
для некоторого автогомеоморфизма $\sigma:F\to F$.

\begin{pr}\label{3-skew}
(a) Существует 3-многообразие, не гомеоморфное $(S^1)^3$, которое является одновременно
$S^1$-расслоением над  $S^1\times S^1$ и $S^1\times S^1$-расслоением над $S^1$.

(b)* Пусть $N$ и $X$ --- сферы с ручками.
Если некоторое 3-многообразие $M$

$\bullet$ является $S^1$-расслоением над $N$,

$\bullet$ является $X$-расслоением над $S^1$,

$\bullet$ не гомеоморфно ни $N\times S^1$, ни $X\times S^1$,

то $N\cong X\cong S^1\times S^1$.
\end{pr}

\begin{pr}\label{3-heegard}
(a) Для любого линейного автогомеоморфизма $f$ тора $S^1\times S^1$ пространство
$D^2\times S^1\cup_fD^2\times S^1$ гомеоморфно либо $S^2\times S^1$, либо линзовому.

(b) Найдите $H_1(D^2\times S^1\cup_{f_A}D^2\times S^1)$,
где $f_A$ --- линейный автогомеоморфизм тора $S^1\times S^1$, заданный матрицей  $A$.

(c) {\it Диаграммы Хегора}. Для любого 3-многообразия $N$ найдется автогомеоморфизм $f$ сферы с ручками $S$,
ограничивающей полный крендель $X$, для которого $N\cong X\cup_f X$.
\end{pr}

\subsection{Указания и решения к некоторым задачам}

{\bf \ref{3-home}.} Используйте $H_1$.
Для (e) и для различения $S^1\times S^2$ и $\R P^3$ потребуется $H_1(\cdot;\Z)$.
Для (с) потребуются также ориентируемость и край.

\smallskip
{\bf \ref{3-homolexa}.} (1)-(5) $H_1(S^3)=H_1(D^3)=0$, \ $H_1(\R P^3)\cong H_1(D^2\t\times S^1)\cong\Z_2$,
\ $H_1(L(p,q))\cong\Z_{GCD(p,2)}$.

(6) Аналогично задаче \ref{pi1-prop}.b.

\smallskip
{\bf \ref{3-homolindep}.} Аналогично двумерному случаю.

\smallskip
{\bf \ref{3-homolexaz}.} $H_1(S^3;\Z)=H_1(D^3;\Z)=0$ и $H_1(\R P^3;\Z)\cong \Z_2$.

\smallskip
{\bf \ref{3-homolzcxcell}.} (4),(5) $H_1(L(p,q);\Z)\cong\Z_p$ и $H_1(D^2\t\times S^1;\Z)\cong\Z$.

\smallskip
{\bf \ref{3-homol}.} Ответы: да. Используйте задачу \ref{3-homolindepz}.c.

\smallskip
{\bf \ref{3-lens}.} Аналогично задаче \ref{3-homol}.

\smallskip
{\bf \ref{hopf-exa}.}
(a,b,c,d,e,g,h,i) Используйте накрытия, т.е. задачу \ref{pi1-covsim}.

(a,b,f,j) Используйте задачу \ref{pi1-prop}.b.

(a,c) Используйте задачу \ref{pi1-prop}.a.

(k) Используя теорему Зейферта--Ван Кампена о фундаментальной группе объединения,
для трилистника $f:S^1\to S^3$ нетрудно получить, что
$\pi_1(S^3-fS^1)=\left<x,y\ |\ xyx=yxy\right>=\left<a,b\ |\ a^2=b^3\right>$.

\smallskip
{\bf \ref{hopf-cov}.} (a--c) Аналогично задаче \ref{vfcov}.

(d) Определим отображение $\alpha:\pi_1(X)\to p^{-1}(x_0)$ так.
Возьмем $\t x_0\in p^{-1}(x_0)$.
Для петли $s:[0,1]\to X$ рассмотрим такое ее поднятие $\t s:[0,1]\to\t X$, для которого $\t s(0)=\t x_0$.
Положим $\alpha(s):=\t s(1)$.
Докажите, что это отображение биективно.

(e) Докажите, что у каждой точки ровно $|\pi_1(\t X)|$ прообразов при отображении $\alpha$ из (d).

\smallskip
{\bf \ref{pi1-covsim}.} (b) Аналогично задаче \ref{hopf-cov}.d.

\smallskip
Указания и решения к задачам \ref{hopf-equ}-\ref{pi1-poin} можно найти, например, в~\cite[\S2, \S11]{Pr04}.

\smallskip
{\bf \ref{pi1-tref}.}
Для различения трилистника и тривиального узла достаточно построить нетривиальный гомоморфизм
$\pi_1(S^3-fS^1)=\left<x,y\ |\ xyx=yxy\right>\to S_3$.
Его можно определить формулами $x\to(12)$, $y\to(23)$.
Значит, $\pi_1(S^3-fS^1)$ не абелева и не изоморфна группе $\Z$.
(Для тривиального узла $f_0:S^1\to S^3$ имеем $\pi_1(S^3-f_0S^1)=\Z$.)

\smallskip
{\bf \ref{3-poin}.} Для доказательства того, что
$\pi_1(N)\ne0$, постройте нетождественное накрытие $S^3\to N$.

\smallskip
{\bf \ref{3-pi1}.} Ответы: (a,b,c) да, (d,e) нет.

(d) Пусть существует 3-многообразие $N$, для которого $\pi_1(N)\cong\Z^4$.
Можно считать, что $N$ связно.
Возьмем клеточное разбиение, данное задачей \ref{3-bala}.
Приклеивая к $N$ клетки размерностей 3 и выше, получим клеточный комплекс $X\supset N$, для которого
$\pi_1(X)\cong\pi_1(N)\cong\Z^4$ и $\pi_k(X)=0$ для $k>0$.
Так как $(S^1)^4$ имеет то же свойство, то $X\simeq (S^1)^4$.
Поэтому $H_2(X;\Z_2)\cong H_2((S^1)^4;\Z_2)\cong\Z_2^6$.
Клеточный комплекс для $X$ имеет вид $\dots\to C_2\overset\partial\to C_1\to0$.
Так как
$$H_1(X;\Z_2)\cong \Z_2^4,\quad\text{то}\quad \dim\im\partial=\dim C_1-4=\dim C_2-4\quad\text{поэтому}$$
$$6=\dim H_2(X;\Z_2)\le \dim\ker\partial\le\dim C_2-\dim\im\partial=4.$$
Противоречие.

\smallskip
{\bf \ref{3-bala}.}
Возьмем произвольную триангуляцию 3-многообразия $N$.
Объединением ее тетраэдров вдоль дерева можно получить клеточный комплекс,
представляющий многообразие $N$, в котором ровно одна трехмерная клетка.
Стягиванием его ребер можно получить клеточный комплекс, представляющий многообразие $N$,
в котором ровно одна вершина и ровно одна трехмерная клетка.
Поскольку и при объединении тетраэдров, и при стягивании ребер эйлерова характеристика сохраняется,
в полученном комплексе число одномерных клеток равно числу двумерных.

\smallskip
{\bf \ref{hopf-cov2}.} (a) Аналогично задаче \ref{pi1-covsim}.a.

(b) Аналогично задачам \ref{hopf-cov}.e и \ref{hopf-inj}.b.

\smallskip
{\bf \ref{3-dehn}.} (d) Равенство $H_1(\chi_{g,1}(S^3))=0$ доказывается при помощи последовательности Майера-Виеториса.
Для доказательства негомеоморфности используйте фундаментальную группу.

\smallskip
{\bf \ref{3-skew}.} (a) Возьмем $e\in\Z-\{0\}$ и
$S^1\widetilde\times_\sigma(S^1\times S^1)$, где  $\sigma$ ---  автогомеоморфизм тора, полученный из автоморфизма
плоскости, заданного формулой $(x,y)\mapsto(x+ey,y)$.

(b) Вычислите гомологии 3-многообразия через
$S^1$-расслоение над $N$ и через $X$-расслоение над $S^1$.

\smallskip
{\bf \ref{3-lensho}.} (a) Для $N=L(5,1)$ или $N=L(5,2)$ определите {\it форму зацеплений}
$\lk:H_1(N)\times H_1(N)\to\Z_5\quad\text{формулой}\quad$
$$\lk([a],[b]):=a\cap b'\mod5,\quad\text{где}\quad b'\in C_2(N)\quad\text{и}\quad\partial b'=5b.$$
Докажите, что форма зацеплений корректно определена.
Проверьте, что формы зацеплений многообразий $L(5,1)$ и $L(5,2)$ не изоморфны.

\newpage
\section{Наборы векторных полей}\label{0vesy}

\subsection{О существовании наборов касательных полей}

{\bf Теорема о векторных полях.} {\it Если $n+1=2^rm$, где $m$ нечетно, то
на $\R P^n$ не существует набора из $2^r$ линейно независимых касательных векторных полей.}

\smallskip
{\bf Теорема об алгебрах с делением.}
{\it Если на $\R^n$ имеется структура алгебры с делением,
то $n$ есть степень двойки.}

\smallskip
Доказательство Хопфа, не использующее систем векторных полей, было получено одновременно с
доказательством Штифеля, использующим их \cite{Hi95}.

Более точно, {\it алгебры с делением на $\R^n$ имеются только при $n
=1,2,4,8$}.
Эта знаменитая теорема Ботта-Милнора-Кервера (см. ссылки в~\cite[\S4]{MS74})
доказывается также с использованием топологии (но гораздо более продвинутой) \cite{Hi95}.

В этом параграфе приводится построение характеристических классов и наброски доказательств всех этих теорем.
Другие исторически первые эффектные применения характеристических классов --- теоремы Уитни о невложимости
и Понтрягина-Тома о некобордантности многообразий --- описаны в \S\ref{0venemb} и \S\ref{0cobord}, ср.~\cite[\S4]{MS74}.

Для доказательства следующих элементарно формулируемых фактов (и задачи \ref{vesenpro} ниже)
также необходимы характеристические классы (кроме \ref{veseprs1}.a, части `тогда' в задаче \ref{veseprs1}.b и
части `только тогда' в задаче \ref{veseprs1}.c).
Далее в соответствующих местах написано, что можно вернуться к решению этих задач.

Далее в этом параграфе пара (тройка, четверка, ...) ортонормированных касательных векторных полей называется просто
парой (тройкой, четверкой, ...) полей.

Для $n$-многообразия $N$ обозначим через $N_0$ дополнение до внутренности некоторого $n$-мерного шара в $N$.

\begin{pr}\label{vesepro4}
Пусть $F$ и $F'$ --- замкнутые 2-многообразия.

(a) На $(F\times F')_0$ есть тройка полей тогда и только тогда, когда
одно из них ориентируемо, а у другого эйлерова характеристика четна.

(b) Если на $(F\times F')_0$ есть пара полей, то
одно из них ориентируемо или у обоих эйлерова характеристика четна.

(c)* Верно ли обратное к (b)?
\end{pr}

\begin{pr}\label{veseprs1}
  Пусть $M$ --- замкнутое 3-многообразие.

(a) На $M\times S^1$ существует пара полей.

(b) На $M\times S^1$ существует тройка полей тогда и только тогда, когда на $M$ существует пара полей.

(c) На $M\times S^1$ существует четверка полей тогда и только тогда, когда $M$ ориентируемо.

Указание. Выведите часть `тогда' из теоремы Штифеля.
\end{pr}

\begin{pr}\label{vesepr4}
На дополнении замкнутого ориентируемого 4-многообразия до шара есть пара полей.
\end{pr}

 \subsection{Характеристические классы для 4-многообразий}\label{0vesy4}

Формально, содержание этого пункта не используется в дальнейшем.

Для решения следующих задач \ref{ves42}--\ref{ves42bd} необходимо пространство $SO_4$ положительных реперов в $\R^4$,
равенства $\pi_1(SO_4)\cong\Z_2$, $\pi_2(SO_4)=0$, и изоморфность гомоморфизма $\pi_1(SO_3)\to\pi_1(SO_4)$,
индуцированного включением.
См. \S\ref{0hopfbun}.


\begin{pr}\label{ves42}
Пусть $N$ --- замкнутое связное 4-многообразие.

(a) Определите $H_2(N)$ и $w_2(N)\in H_2(N)$ как препятствие к построению {\it тройки} полей.

(b) Вне окрестности любого непустого цикла, представляющего $w_2(N)$, существует тройка полей.

(c) Препятствие $w_2(N)$ не полно.

(d) $w_2(M\times S^1)=w_2(M)\times S^1$ для любого замкнутого 3-многообразия $M$.
(Определите сами, что такое $\times$.)

(e') $w_2(F\times S^1\times S^1)=\rho_2\chi(F) p\times S^1\times S^1$ для любого связного замкнутого
2-многообразия $F$ и $p\in F$.

(e) Найдите $w_2(F\times F')$ для замкнутых 2-многообразий $F$ и $F'$.
Ответ: $w_2(F\times F')=w_2(F)p\times F'+w_1(F)\times w_1(F')+w_2(F')F\times p'$, где
$p\in F$, $w_2(F):=\rho_2\chi(F)$ и то же для $F'$.
\end{pr}

При помощи пункта (d) доказывается часть `только тогда' задачи \ref{veseprs1}.b, а пункта (e) --- задача
\ref{vesepro4}.a.

\begin{pr}\label{vesw2}
 Следующие три условия на замкнутое связное 4-многообразие $N$ равносильны.

$w_2(N)=0$;

существует тройка полей на дополнении в $N$ до некоторого графа;

существует тройка полей на $N_0$.
\end{pr}

Для замкнутого связного ориентируемого 4-многообразия $N$ препятствие к продолжению тройки полей с $N_0$ на $N$  лежит в
$\pi_3(SO_4)\cong\Z\oplus\Z$, т.е. является парой чисел (см. \S\ref{0hopfbun}).
Эти числа --- эйлерова характеристика $\chi(N)$ и (с точностью до множителя) сигнатура $\sigma(N)$ формы пересечений
$\cap:H_2(N;\Z)\times H_2(N;\Z)\to\Z$ (ее определение аналогично приведенному в \S\ref{02homint}, см. \S\ref{0vesyint}).

\smallskip
{\bf Теорема.} {\it На замкнутом связном ориентируемом
4-многообразии $N$ существует четверка ортонормированных касательных векторных
полей тогда и только тогда, когда $\chi(N)=\sigma(N)=0\in\Z$ и $w_2(N)=0\in H_2(N)$.}

\smallskip
(Доказательство выходит за рамки этой книги.)

\begin{pr}\label{ves42bd}
(a--f) Сформулируйте и докажите аналоги задачи \ref{ves42} для
связных ориентируемых 4-многообразий с непустым краем.
Аналог пункта (c) следующий: на связном ориентируемом 4-многообразии $N$ с непустым краем существует
четверка линейно независимых касательных векторных полей тогда и только тогда,
когда $w_2(N)=0\in H_2(N,\partial N)$.
\end{pr}

\begin{pr}\label{ves43}
Пусть $N$ --- замкнутое связное  4-многообразие (не обязательно ориентируемое).
Для решения этой задачи необходимо многообразие Штифеля $V_{4,2}$ ортонормированных пар векторов в $\R^4$,
равенства $\pi_1(V_{4,2})=0$, $\pi_2(V_{4,2})\cong\Z$, и изоморфность гомоморфизма
$\pi_2(S^2)\cong\pi_2(V_{3,1})\to\pi_2(V_{4,2})$, индуцированного включением.
См. \S\ref{0hopfbun}.

(a) Определите $H_1(N;\Z)$ и $W_3(N)\in H_1(N;\Z)$ как препятствие к построению {\it пары} полей для ориентируемого $N$.

(b) Вне окрестности любого непустого графа, представляющего $W_3(N)$, можно построить пару полей.

(c) Препятствие $W_3(N)$ не полно.

(d) $W_3(N)=0$ тогда и только тогда, когда на $N_0$ существует пара полей.

(e) Найдите $\rho_2W_3(F\times F')$ для замкнутых 2-многообразий $F$ и $F'$.
Ответ:
$\rho_2W_3(F\times F')
=w_2(F)p\times w_1(F')+w_2(F')w_1(F)\times p'$, где
$p\in F$, $w_2(F)=\rho_2\chi(F)$ и то же для $F'$.
\end{pr}

При помощи пункта (e) решается задача \ref{vesepro4}.b.

\begin{pr}\label{ves43bd}
(a-d)Решите аналог предыдущей задачи для многообразий с непустым краем.
\end{pr}

\begin{pr}\label{ves4bock} Пусть $N$ --- замкнутое связное ориентируемое 4-многообразие.

(a) Если $w_2(N)=0$, то $W_3(N)=0$.

(b)* Если $w_2(N)$ представляется ориентируемым 2-многообразием, то $W_3(N)=0$.


(с)* Класс $w_2(N)$ представляется ориентируемым 2-многообразием.
(Значит, $W_3(N)=0$; отсюда вытекает утверждение задачи \ref{vesepr4}.)
\end{pr}

 \begin{pr}\label{ves4int}* Пусть $N$ --- замкнутое связное ориентируемое 4-многообразие.

(a) Докажите, что $a\cap a=w_2(N)\cap a$ для любого $a\in H_2(N)$.
(Ср. с задачей \ref{intstwh}.a.)

(b) $w_2(N)\cap w_2(N)=\chi(N)\mod2$. (Ср. с задачей \ref{intstwh}.b.)

(c) Определите {\it число Понтрягина} $p_1(N)\in\Z$ как полное препятствие к существованию
{\it четверки касательных векторных полей, имеющей ранг не менее трех в каждой точке}.
(См. задачу \ref{cobhirz}.a.)
\end{pr}

\begin{pr}\label{ves4imm} $\R P^2\times \R P^2$
\quad
(a) не погружается в $\R^5$.
\quad
(b) не вкладывается в $\R^6$.
\end{pr}

 \subsection{Общее определение групп гомологий}\label{0vesyhom}


Дадим {\it определение групп гомологий (с коэффициентами в $\Z_2$)},
независимое от рассуждений, в которых оно появилось.
Пусть $T$ --- симплициальный комплекс, имеющей $c_k$ симплексов размерности $k$, $k=0,1,2,\dots$.
(Необходимые изменения для разбиения $T$ читатель легко сделает самостоятельно.)
Обозначим через $C_k=C_k(T)$ группу расстановок нулей и единиц на $k$-мерных симплексах
(относительно операции покомпонентного сложения).
Ясно, что $C_k\cong\Z_2^{c_k}$.

Для произвольного ребра $a$ обозначим через $\partial_0a$
расстановку единиц в вершинах этого ребра и нулей в остальных вершинах.
Продолжим $\partial_0$ до линейного отображения $\partial_0:C_1\to C_0$.

Аналогично, для произвольной $k$-мерной грани $a$ обозначим через $\partial_{k-1}a$ расстановку
единиц на ребрах границы этой грани.
Продолжим $\partial_{k-1}$ до линейного отображения $\partial_{k-1}:C_k\to C_{k-1}$.

Группы $\partial_{k-1}^{-1}(0)$ и $\partial_kC_{k+1}$ называются группами
$k$-циклов и $k$-границ, соответственно.

\begin{pr}\label{}
$\partial_k\partial_{k+1}=0$.
Иными словами, любая граница является циклом.
\end{pr}

Положим
$$H_0(N):=C_0/\partial_0C_1,
\quad\text{и}\quad \quad H_k(N):=\partial_{k-1}^{-1}(0)/\partial_kC_{k+1}.$$
В некоторых местах этого параграфа мы {\it указываем} коэффициенты $\Z_2$ в обозначениях.
При этом пропуск коэффициентов по-прежнему означает, что предполагаются коэффициенты $\Z_2$.

Чтобы определить гомологии с коэффициентами в $\Z$, фиксируем дополнительно {\it ориентацию}
(т.е. порядок вершин с точностью до четной перестановки) каждого симплекса триангуляции $T$
(никакой согласованности ориентаций разных симплексов не предполагается).
Обозначим через $C_k(\Z)=C_k(N;\Z)$ группу расстановок
целых чисел на ориентированных $k$-мерных симплексах триангуляции $T$
(относительно операции покомпонентного сложения).
Ясно, что $C_k(\Z)\cong\Z^{c_k}$.
Для $k$-симплекса $a$ с вершинами $a_0\dots a_k$ обозначим через $\widehat{a_i}$
симплекс размерности $k-1$ с вершинами $a_0,\dots,a_{i-1},a_{i+1},\dots,a_k$.
Положим $\varepsilon_i=-1$, если ориентация симплекса $\widehat{a_i}$ задается порядком
$(a_0,\dots,a_{i-1},a_{i+1},\dots,a_k)$, и $\varepsilon_i=1$, иначе.
Положим
$\partial_{k-1}a=\sum\limits_{i=0}^k(-1)^i\varepsilon_i\widehat{a_i}$.

\begin{pr}\label{}
$\partial_k\partial_{k+1}=0$.
Иными словами, любая граница является циклом.
\end{pr}

После этого гомологии с коэффициентами $\Z$ определяются аналогично гомологиям с коэффициентами $\Z_2$.
Аналогично определяются гомологии с коэффициентами $\Z_p$ и $\Q$.


\begin{pr}\label{}
{\bf Теорема инвариантности гомологий.}
{\it Гомологии кусочно-линейно гомеоморфных симплициальных комплексов изоморфны}.
\end{pr}

 \subsection{Общее определение формы пересечений}\label{0vesyint}

 {\it Определение формы пересечений
$$\cap:H_i(N)\times H_{n-i}(N)\to\Z_2$$
для замкнутого $n$-многообразия $N$.}
(Аналогично \S\ref{02homint}.)
Для $i$-мерной грани $x$ разбиения $T$ и $(n-i)$-мерной грани $y$ двойственного разбиения $T^*$ положим $x\cap y=1$,
если $x$ и $y$ двойственны, и $x\cap y=0$ иначе.
Для набора $X$ из $i$-мерных граней разбиения $T$ и набора $Y$ из $(n-i)$-мерных граней разбиения $T^*$ положим
$X\cap Y$ равным количеству их точек пересечения по модулю 2 (т.е. продолжим умножение $\cap$ на ребрах до билинейной
симметричной формы $\cap:C_i(T)\times C_{n-i}(T^*)\to\Z_2$).
Тогда формула $[a]\cap[b]:=a\cap b$ корректно задает форму пересечений
ввиду теоремы инвариантности групп гомологий и следующей задачи.

\begin{pr}\label{vesint}
(a) Пересечение цикла и границы всегда нулевое.

(b) Определение формы пересечений корректно.

(с) $(\alpha+\beta)\cap\gamma=\alpha\cap\gamma+\beta\cap\gamma$ и
$\alpha\cap(\beta+\gamma)=\alpha\cap\beta+\alpha\cap\gamma$, т.е.
пересечение на гомологиях билинейно.

(d)* Если $N$ --- замкнутое ориентируемое 4-многообразие, то $x^4=0$ для любого $x\in H_3(N)$.
\end{pr}

\begin{pr}\label{vesingen} Дано $n$-многообразие $N$ с краем.

(a) Форма пересечений $\cap:H_i(N)\times H_{n-i}(N)\to\Z_2$ определена корректно формулой
$[a]\cap[b]:=a\cap b$.

(b) Она может быть вырожденной.

(c) Определите целочисленную форму пересечений $\cap:H_i(N;\Z)\times H_{n-i}(N;\Z)\to\Z$.

(d) Определите билинейное умножение $\cap:H_i(N,\partial N)\times H_{n-i}(N)\to\Z_2.$

(e) Для $n$-многообразия $N$ с краем (возможно, пустым) формула
$[a]\cap[b]:=a\cap b$ корректно определяет билинейные умножения
$$H_i(N)\times H_j(N)\to H_{i+j-n}(N),
\quad H_i(N)\times H_j(N,\partial N)\to H_{i+j-n}(N)$$
$$\text{и}
\quad H_i(N,\partial N)\times H_j(N,\partial N)\to H_{i+j-n}(N,\partial N).$$
\end{pr}

\subsection{Характеристические классы для $n$-многообразий}\label{0vesych}

Следующая теорема обобщает теоремы ориентируемости (\S\ref{02homcy}, \S\ref{03manor}),
Эйлера-Пуанкаре (\S\ref{0ve2ma}) и Хопфа (\S\ref{0vehima}),
лемму Штифеля о препятствии (\S\ref{0vesy3ch}) и некоторые задачи из \S\ref{0vesy4}.

Через $\Z_{(i)}$ обозначается группа $\Z$ для четного $i$ и $\Z_2$ для нечетного~$i$.

\smallskip
{\bf Теорема Штифеля о препятствии.}
{\it Если на замкнутом $n$-многообразии $N$ существует набор из $k$
линейно независимых касательных векторных полей
($1<k<n$), то} $(n-k+1)$-й класс Штифеля-Уитни
$W_{n-k+1}(N)\in H_{k-1}(N,\Z_{(n-k)})$ \ {\it нулевой.}
\footnote{Это еще одно преимущество `гомологического' подхода над `когомологическим': некоторые `когомологические'
классы Штифеля-Уитни для неориентируемых многообразию принимают значения в {\it скрученных} когомологиях, которые
изоморфны обычным гомологиям ввиду {\it скрученного} изоморфизма Пуанкаре.}

\smallskip
{\it Набросок определения класса $W_{n-k+1}(N)$, использующего общее положение.}
Рассмотрим набор общего положения из $k$ касательных векторных полей на $n$-многообразии $N$.
Ввиду общности положения множество точек многообразия, в которых векторы этого набора
линейно зависимы, является объединением $(k-1)$-мерных подмногообразий.
Если $n-k$ четно, то на них можно ввести ориентацию.
Эти подмногообразия несут цикл с коэффициентами в $\Z_{(n-k)}$.
Классом $W_{n-k+1}(N)$ называется класс гомологичности этого цикла (\S\ref{0vesyhom}).

\smallskip
{\it Набросок определения класса $W_{n-k+1}(N)$.}
Обозначим через $V_{n,k}$ многообразие Штифеля ортонормированных $k$-реперов в $\R^n$.
Тогда (\S\ref{0hopfbun})
$$\pi_i(V_{n,k})=0\quad\text{для}\quad i<n-k\quad\text{и}\quad\pi_{n-k}(V_{n,k})
\cong\Z_{(n-k)}\quad\text{для}\quad 1<k<n.$$
Ана\-ло\-гич\-но трехмерному случаю (\S\ref{0vesy3ch}) первое из этих равенств влечет, что набор
из $k$ линейно независимых касательных векторных полей всегда строится на $(n-k)$-остове.
Ввиду второго из этих равенств продолжению полученного набора полей на $(n-k+1)$-остов препятствует
некоторая расстановка $\omega_{n-k+1}$ на $(k-1)$-мерных многогранниках двойственного разбиения
целых чисел при четном $n-k$, или нулей и единиц при нечетном $n-k$.

Аналогично предыдущему (см. ссылки перед формулировкой теоремы) среди всех таких расстановок выделяются {\it циклы} и
{\it границы}, а также определяется {\it гомологичность} циклов.
Группа циклов с точностью до гомологичности называется группой гомологий $H_{k-1}(N,\Z_{(n-k)})$.
Ср. \S\ref{0vesyhom}.
Классом Штифеля-Уитни называется класс гомологичности препятствующей расстановки
$$W_{n-k+1}(N)=[\omega_{n-k+1}]\in H_{k-1}(N,\Z_{(n-k)}).$$

\begin{pr}\label{vesnw}
(a) Препятствующая расстановка является циклом.

(b) Докажите теорему Штифеля о препятствии.

(c) $W_{n-k+1}(N)=0$ тогда и только тогда, когда набор из $k$ полей существует на дополнении в
$N$ до некоторого $(k-1)$-мерного подкомплекса.
\end{pr}

Обратная теорема (к теореме Штифеля о препятствии) неверна, т.е. препятствие
$W_{n-k+1}$ уже не всегда является полным (см. задачи \ref{ves42}.c и \ref{ves43}.c).

Для замкнутого $n$-многообразия $N$ аналогично определяются классы $W_1(N)\in H_{n-1}(N,\Z_2)$
и $W_n(N)
\in H_0(N,\Z)$.
Эти определения равносильны данным в \S\S\ref{02hom},\ref{0ve2ma},\ref{0vehima}: $W_1(N)=w_1(N)$ и $W_n(N)=\chi(N)$.

Через $\rho_2:H_s(N;\Z)\to H_s(N;\Z_2)$ обозначается приведение по модулю два.
Обозначим $w_i(N):=\rho_2W_i(N)\in H_{n-i}(N;\Z_2)$.
Этот класс легче вычисляется.
Удобно считать, что $w_k(N)=0$ при $k>n$ и что $w_0(N)=[N]\in H_n(N)$ представляется объединением $n$-мерных
многогранников произвольного разбиения многообразия $N$.

Как выразить классы Штифеля-Уитни произведения многообразий через через классы Штифеля-Уитни сомножителей?

\begin{pr}\label{vesnww}
(a) {\it Формула Уитни-Ву}. $w_s(P\times P')=\sum\limits_{k=0}^s w_k(P)\times w_{s-k}(P')$.
Эту формулу удобно переформулировать следующим образом.
{\it Полный класс Штифеля-Уитни} замкнутого многообразия $N$ определяется как
$$w(N):=1+w_1(N)+w_2(N)+\dots\in H_n(N)\oplus H_{n-1}(N)\oplus H_{n-2}(N)\dots$$
Здесь единица в сумме (но не в индексе!) обозначает $[N]$; это удобно, поскольку $[N]\cap x=x$
для любых $s$ и $x\in H_s(N)$.
В этих обозначениях формула Уитни-Ву записывается как $w(P\times P')=w(P)\times w(P')$.

(b) Вычислите $w((\R P^2)^k)$.

(c) Вычислите классы Штифеля-Уитни произведения любого числа 2-многообразий.

(d)* Если $w_1(N)=w_2(N)=\dots=w_{[n/2]}(N)=0$ для замкнутого $n$-многообразия $N$, то $w_s(N)=0$ для любого $s$.
\end{pr}

\begin{pr}\label{vesncalc}
(a) $w_1(\R P^{2k+1})=0$ и $w_1(\R P^{2k})=[\R P^{2k-1}]$.

(b) $w_2(\R P^n)=0$ для $n\equiv0,3\mod4$ и $w_2(\R P^n)=[\R P^{n-2}]$ для $n\equiv1,2\mod4$, $n\ge2$.

(c) $w_i(\R P^n)=\displaystyle{n+1\choose i}[\R P^{n-i}]$ для $0\le i\le n$.
\end{pr}

\begin{pr}\label{vesntr}* Для любой триангуляции замкнутого $n$-многообразия $N$ объединение $k$-мерных симплексов ее
барицентрического подразделения является циклом, представляющим класс $w_{n-k}(N)$.
\end{pr}

\begin{pr}\label{vespont}*
(a) Для замкнутого ориентируемого $n$-многообразия $N$ определите {\it $j$-й класс Понтрягина} $p_j(N)\in H_{n-4j}(N;\Z)$
как препятствие к построению набора $n-2j+2$ векторов, имеющей ранг не менее $n-2j+1$ в каждой точке.

(b) $2(p_s(P\times P')-\sum\limits_{k=0}^s p_k(P)\times p_{s-k}(P'))=0$.

(c) $p_i(\C P^n)=\displaystyle{n+1\choose i}[\C P^{n-2i}]$ для $1\le i\le n/2$.
\end{pr}

\subsection{Указания и решения к некоторым задачам}

{\bf \ref{vesepro4}.} (a)  По задаче \ref{ves42}.ae
$$0=w_2(F\times F')=w_2(F)p\times F'+w_1(F)\times w_1(F')+w_2(F')F\times p',\text{ где }p\in F,\ w_2(F):=\rho_2\chi(F)$$
и то же для $F'$.
Пересекая правую часть (в смысле формы пересечений) с $p\times F'$, получаем $w_2(F')=0$.
Аналогично, пересекая правую часть с $F\times p'$, получаем $w_2(F)=0$.
Значит, $w_1(F)\times w_1(F')=0$.
Если ни $F$, ни $F'$ не ориентируемо, то $w_1(F)\ne0$ и $w_1(F')\ne0$.
Тогда существуют $\alpha\in H_1(F)$ и $\alpha'\in H_1(F')$, для которых $\alpha\cap w_1(F)=1$ и
$\alpha'\cap w_1(F')=1$.
Поэтому $w_1(F)\times w_1(F')\cap \alpha\times\alpha'=1\ne0$.
Полученное противоречие доказывает утверждение задачи.

(Для решения этой задачи группу  $H_2(F\times F')$ вычислять не нужно!)


\smallskip
{\bf \ref{ves42}.} (a) Аналогично лемме Штифеля.

(b) Аналогично задаче \ref{veseobcy}.a.

(c) В качестве контрпримера годится $N=S^4$, на которой нет даже одного поля.

(d) Аналогично задаче \ref{veseobcyex}.a.
Пусть $u,v$ --- пара полей на объединении ребер некоторого разбиения 3-многообразия $M$, для которой препятствующий
цикл является объединением ребер $e_1,\dots,e_n$ двойственного разбиения.
Пусть $v'$ --- единичное касательное векторное поле на $S^1$.
Тогда для `призматического' разбиения произведения $M\times S^1$ и тройки $(u,v,v')$ препятствующий цикл является
объединением колец $e_i\times S^1$ по $i=1,\dots,n$.

(e) Аналогично задаче \ref{ves43}.e.
Нужно рассмотреть тройку полей $(u,v',v+u')$.

\smallskip
{\bf \ref{vesw2}.} Аналогично лемме Штифеля с использованием равенства $\pi_2(SO_4)=0$.

\smallskip
{\bf \ref{ves42bd}.} Аналогично задачам \ref{ves42} и \ref{vesw2}.

\smallskip
{\bf \ref{ves43}.} (a)-(d) Аналогично задачам \ref{ves42} и \ref{vesw2}.

(e) Возьмем разбиение 2-многообразия $F$, для которого $Оp$ есть внутренность многоугольника.
Рассмотрим такую пару полей $u,v$ (не общего положения!) на $F$, что $u\ne0$ вне $p$,
$u\perp v$, $v=0$ на объединении ребер некоторого представителя $\omega$ класса $w_1(F)$, лежащего в двойственной
триангуляции, и $p\not\in\omega$.
Рассмотрим пару полей $v',u'$ общего положения на $F'$ с аналогичными свойствами.
Рассмотрим пару полей $(u+u',v+v')$ (не общего положения!) на $F\times F'$.
Множество точек линейной зависимости этой пары --- множество точек, где пары $(u,v)$ и $(u',v')$
одинаково линейно зависимы.

В точках множества $Op\times F'_0$ пара $(u',v')$ линейно зависима только при $u'=0$, т.е. на
$Op\times\omega'$.
Соответствующая линейная зависимость между $u$ и $v$ имеет вид $u=0$ и имеется только в $p$.
В точках множества $F_0\times Op'$ пара $(u,v)$ линейно зависима только при $v=0$, т.е. на
$\omega\times Op'$.
Соответствующая линейная зависимость между $u'$ и $v'$ имеет вид $v'=0$ и имеется только в $p'$.

Поэтому пара $(u',v')$ невырождена вне $p\times\omega'\cup\omega\times p'$.
Рассмотрим разбиение 4-многообразия $F\times F'$ на многоугольники, являющееся произведенем
разбиений 2-многообразий $F$ и $F'$.
Красными (по модулю 2!) могут быть только ребра из набора $p\times\omega\cup\omega\times p'$ двойственного разбиения.
Ребра набора $p\times\omega'$ ($\omega\times p'$) будут красными тогда и только тогда, когда
$\chi(F)$ ($\chi(F')$) четно.
Это доказывает искомую формулу.

\smallskip
{\bf \ref{ves43bd}.} 
Для многообразий с краем (ориентируемых или нет) получится $W_3(N)\in H_1(N,\partial N;\Z)$.

\smallskip
{\bf \ref{ves4bock}.} (a) Следует из \ref{vesw2} и \ref{ves43}.d.

(с) Пусть класс в $H_2(N)$ представлен замкнутым 2-многообразием $A$
(на самом деле, любой класс в $H_2(N)$ может быть так представлен).
Рассмотрим 1-цикл на $A$, препятствующий ориентируемости 2-многообразия $A$.
На нем можно ввести ориентацию (подумайте, как!).
Гомологический класс в $N$ полученного ориентируемого 1-цикла и есть $\beta[A]$.

{\it Гомоморфизм Бокштейна} $\beta$ определяется так.
Если $a=a_1+\dots+a_s$ --- гомологический 2-мерный цикл, где $a_1,\dots,a_s$ являются многоугольниками
некоторого разбиения, то возьмем расстановку $\t a=a_1+\dots+a_s$ целых чисел на многоугольниках, положим
$\beta[a]:=[\frac12\partial\t a]$ и проверим
корректность (т.е. проверим, что $\beta(\partial b)=0$).
Ср. \S\ref{0homolex}.

\smallskip
{\bf \ref{ves4int}.}
(b) Вычет $w_2(N)\cap a$ есть препятствие к тому, чтобы на цикле,
представляющем $a$, можно построить тройку линейно независимых касательных (к $N$) векторных полей.
Вычет $a\cap a$ есть препятствие к тому, чтобы на 2-многообразии
$A\subset N$, представляющем $a$, можно было бы построить {\it стабильное
нормальное поле}, т.е. пару полей в $\R^{n+1}$, нормальных к $A$, где $A\subset N\subset\R^n$.

(d) Вот набросок определения, использующего общее положение.
Рассмотрим четверку общего положения касательных векторных полей на $N$.
Ввиду общности положения данная четверка векторов имеет ранг менее трех лишь в конечном количестве точек
(разберитесь, почему!).
У этих точек можно определить знак (разберитесь, как!).
Количество этих точек с учетом знака и называется $p_1(N)$.

\smallskip
{\bf \ref{ves4imm}.}
Обозначим $N:=\R P^2\times\R P^2$ и $a:=[\R P^1]\in H_1(\R P^2)$.
Пусть $f:N\to\R^5$ --- погружение.
Рассмотрим препятствие $x$ к наличию четверки линейно независимых векторных полей на $f(N)$ (имеются в виду
векторные поля в $\R^5$, на них не накладывается ни условие касательности, ни условие нормальности к $f(N)$).
Тогда цепочка равенств в $H_2(N)$
$$0=x=w_2(N)+w_1^2(N)=(a^2\times1+a\times a+1\times a^2)+(a\times1+1\times a)^2=a\times a\ne0$$
дает противоречие.
Здесь

$\bullet$ первое равенство выполнено, поскольку такая четверка существует.

$\bullet$ второе равенство доказывается аналогично теореме Уитни-Ву и второму доказательству теоремы Штифеля
(см. детали ниже).

$\bullet$ третье равенство верно по задаче \ref{ves42}.e (или поскольку $w(N)=(1+a+a^2)\times(1+a+a^2)$
по теореме Уитни-Ву, задача \ref{vesnww}.a).

$\bullet$ четвертое равенство очевидно.

$\bullet$ пятое следует из $(a \times a)\cap( a\times a)=1\in\Z_2\ne0$.

Для доказательства второго равенства возьмем четверку касательных векторных полей $v_1,v_2,v_3,v_4$ на $N$,

которые линейно зависимы на некотором трехмерном подкомплексе $\omega_1$ некоторой триангуляции $T$
многообразия $N$, представляющем класс $w_1(N)$, и

для которых тройка $v_2,v_3,v_4$ линейно зависима на некотором двумерном подкомплексе $\omega_2$
триангуляции $T$, представляющем класс $w_2(N)$.

Возьмем также четверку касательных векторных полей $u_1,u_2,u_3,u_4$ на $N$, которые линейно зависимы на
некотором трехмерном подкомплексе $\omega_1'$ триангуляции $T^*$, представляющем класс $w_1(N)$.
Возьмем нормальное касательное векторное полей $u$ на $N$, модуль которого в любой точке равен
четырехмерному объему параллелепипеда, натянутого на  $u_1,u_2,u_3,u_4$, а направление (при ненулевом модуле)
определяется условием $u\perp u_i$ и положительности репера $u,u_1,u_2,u_3,u_4$.

Теперь рассмотрим множество линейной зависимости четверки $u+v_1,v_2,v_3,v_4$ векторных полей на
$f(N)$. Это множество есть $\omega_2\cup(\omega_1\cap\omega_1')$.
Из этого следует второе равенство (аналогично другому простому доказательству
теоремы Штифеля и задачам \ref{ves42}.e, \ref{ves43}.e).

\bigskip
{\bf \ref{vesnww}.} (a) Аналогично задачам \ref{ves42}.e и \ref{ves43}.e.

Обозначим $k=\dim P$ и $k'=\dim P'$.
Возьмем $m\ge k+k'$.
Рассмотрим следующие $m-n+1$ (касательных векторных) полей на $P$.
Первое поле $v_1$ обращается в ноль только на конечном количестве точек, представляющих класс $w_k(P)$.
Второе поле $v_2$ перпендикулярно первому и обращается в ноль только на некотором одномерном подкомплексе некоторой
триангуляции многообразия $P$, представляющем класс $w_{k-1}(P)$.
Третье поле $v_3$ перпендикулярно первым двум и обращается в ноль только на некотором двумерном подкомплексе некоторой
триангуляции многообразия $P$, представляющем класс $w_{k-2}(P)$.
И т.д.
$k$-е поле $v_k$ перпендикулярно первым $k-1$ и обращается в ноль только на некотором $(k-1)$-мерном подкомплексе
некоторой триангуляции многообразия $P$, представляющем класс $w_1(P)$.
Поля $v_{k+1},\dots,v_{m-n+1}$ нулевые.
Аналогично строятся поля $v_1',\dots,v_{m-n+1}'$ на $P'$.
Теперь, рассматривая множество линейной зависимости набора полей $v_1+v_{m-n+1}',v_2+v_{m-n}',\dots,v_{m-n+1}+v_1'$,
получаем искомую формулу.

(d) Обозначим $w_i:=w_i(N)$.

Пусть $n=4$. Докажем, что $w_3=0$ (потом аналогично получится, что $w_4=0$).
Достаточно доказать, что $w_3\cap[F]=0$ для любого 3-подмногообразия $F$.
Имеем
$$w_3\cap[F]=w_3(\tau_N|_F)=w_3(\tau_F\oplus\nu_{F\subset N})=w_3(F)+w_2(F)w_1(\nu_{F\subset N})=0,$$
ибо $w_3(F)=0$ и $w_2(F)=0$.

Рассмотрим общий случай. Обозначим $k=[n/2]$. Желательно доказать, что $w_{k+1}=0$
(потом аналогично получится, что $w_{k+2}=0,\dots,w_n=0$).
Докажем, что $w_{k+1}\cap[F]=0$ для любого $(k+1)$-подмногообразия $F$ (это более слабое свойство).
Обозначим $\nu:=\nu_{F\subset N}$.
Имеем $w(N)\cap[F]=w(\tau_N|_F)=w(\tau_F\oplus\nu)=w(\tau_F)w(\nu)$.
Из этого и условия задачи следует, что $w_s(\nu)=\overline w_s(F)$ для $s\le k$.
Далее,
$$w_{k+1}\cap[F]=w_{k+1}(F)+w_k(F)w_1(\nu)+\dots+w_1(F)w_k(\nu)+w_{k+1}(\nu).$$
Так как $w_{k+1}(\nu)=0$, то это равенство влечет $w_{k+1}\cap[F]=\overline w_{k+1}(F)$.
Последнее число нулевое, ибо $F$ вложимо в $\R^{2k+2}$.

Чтобы приведенное доказательство равенства $w_{k+1}\cap x=0$ заработало для класса $x$, не реализуемого
подмногообразием, нужно определить классы Штифеля-Уитни для `нормального
расслоения' такого класса.
Попытки это сделать приводят к определению {\it стинродовых квадратов}.
(Возможно, именно так они и были придуманы. Соотвествующая работа Стинрода посвящена другой задаче, но там
приводится формальное определение стинродовых квадратов без мотивировки.)
Было бы интересно вычленить из всей теории стинродовых квадратов то, что нужно для решения данной задачи и то,
что нужно для $\overline w_n(N)=0$.

\smallskip
{\bf \ref{vespont}.} Вот набросок определения класса $p_j(N)$, использующего общее положение.
Рассмотрим набор общего положения из $n-2j+2$ касательных векторных полей на $N$.
Ввиду общности положения множество точек многообразия, в которых данный набор
$n-2j+2$ векторов имеет ранг менее $n-2j+1$, является объединением
$(n-4j)$-мерных подмногообразий.
На них можно ввести ориентацию.
Объединение этих подмногообразий с точностью до гомологичности и называется классом $p_j(N)$.

\newpage
\section{Непогружаемость и невложимость}\label{0venemb}

\subsection{Основные результаты о невложимости}

В 1935 году Хопф рассказал о результатах Штифеля о наборах касательных
векторных полей на Международной топологической конференции в Москве.
Там выяснилось, что Хасслер Уитни около 1934 г. тоже естественно пришел к
определению характеристических классов, изучая классическую проблему топологии
о вложимости многообразий.

Мы работаем в гладкой категории, т.е. все многообразия,
векторные поля и отображения предполагаются гладкими.
Отображение $f:N\to\R^m$ гладкого многообразия $N$ называется {\it гладким
погружением}, если $df(x)$ невырожден для любой точки $x\in N$.
Гладкое погружение $f:N\to\R^m$ называется {\it гладким вложением}, если оно инъективно.

Уитни начал с доказательства вложимости $n$-многообразий в $\R^{2n}$ (\S\ref{0emb}).
Его построение препятствий к вложимости ($n$-многообразий в $\R^m$ при $m<2n$)
основано на рассмотрении {\it трубчатой окрестности} подмногообразия.
Он создал теорию сферических расслоений и ввел так называемые классы Штифеля--Уитни
и нормальные (двойственные) классы Штифеля--Уитни
гладкого $n$-многообразия $N$.
С тех пор они играют большую роль в топологии и дифференциальной геометрии.
Используя их, Уитни доказал невложимость некоторых проективных пространств.

\smallskip
{\bf Теорема Уитни.}
{\it Если $\R P^n$ вложимо в $\R^m$ или погружаемо в $\R^{m-1}$, то $\displaystyle{m\choose n}$ четно.}

\smallskip
{\bf Следствие.}
{\it Если $n$ есть степень двойки, то $\R P^n$ не погружаемо в $\R^{2n-2}$ и не вложимо в $\R^{2n-1}$.}

\begin{pr}\label{veserpn}
(a) Выведите Следствие из теоремы Уитни.

(b) Если $\R P^n$ вложимо в $\R^m$ или погружаемо в $\R^{m-1}$, то $\displaystyle{i\choose n}$ четно для любого
$m\le i\le 2n$
(и даже $\displaystyle{i\choose s}$ четно для любых $1\le s\le n$ и $m\le i\le 2s$, что вытекает из предыдущего и
тождества Паскаля, а поэтому не дает новой информации).
\end{pr}

{\bf Теорема Уитни для коразмерности 1.}
{\it Если $\R P^n$ погружаемо в $\R^{n+1}$, то либо $n+1$, либо $n+2$ являются степенями двойки.}

\smallskip
Далее в этом параграфе приводятся наброски доказательств этих результатов.

Позже было проделано много других вычислений, дающих интересные следствия.

\begin{pr}\label{vesenpro}
(a) $(\R P^2)^k$ не погружаемо в $\R^{3k-1}$ (и не вложимо в $\R^{3k}$).

(b) Произведение любых $k$ 2-многообразий погружаемо в $\R^{3k}$ и вложимо в $\R^{3k+1}$.

(c) Произведение любых $k$ замкнутых неориентируемых многообразий размерностей $n_1,\dots,n_k$
не погружаемо в $\R^{n_1+\dots+n_k+k-1}$
(и не вложимо в $\R^{n_1+\dots+n_k+k}$).

(d) Найдите наименьшую размерность евклидова пространства, в которое вложимо данное произведение 2-многообразий.
\end{pr}

{\bf Теорема Хефлигера-Хирша.}
{\it Пусть $N$ --- многообразие размерности $n\ge3$.
Если либо $n$ не есть степень двойки, либо $N$ ориентируемо, либо $N$
незамкнуто, то $N$ вложимо в $\R^{2n-1}$ и погружаемо
в $\R^{2n-2}$~\cite{RS99},~\cite{Sk08}.}

\smallskip
Для $n=3$ и $n=4$ эта теорема была доказана позже и другими авторами, см. ссылки в~\cite{Sk08}.
Из доказательства теорема Хефлигера-Хирша мы сможем наметить там только более
простую часть (т.е. доказательство теоремы Масси).
Более сложная часть --- теорема, обратная к теореме Уитни о препятствии.

Теорему Хефлигера-Хирша обобщают следующая гипотеза.

\smallskip
{\bf Гипотеза Масси.} {\it Любое $n$-мерное многообразие вложимо в $\R^{2n+1-\alpha(n)}$.

(b) Любое $n$-мерное многообразие погружаемо в~$\R^{2n-\alpha(n)}$.
(Некоторые считают, что в доказательстве~\cite{Co85} этой гипотезы имеются пробелы.)}



\smallskip
{\it Нормальным векторным полем} для погружения $f:N\to\R^m$ называется семейство нормальных к
$f(N)$ векторов $v(x)$ в точках $x\in f(N)$, непрерывно зависящих от точки $x\in N$.

\begin{pr}\label{vesenotr}*
(a) При $m\ge3n/2+1$ или $m=n+3$ для любого вложения $S^n$ в $\R^m$ существует
набор из $m-n$ линейно независимых нормальных векторных полей~\cite{Ke59}, ~\cite{Ma59}.

(b) Существует вложение $S^7$ в $\R^{11}$, для которого нет четверки линейно независимых нормальных векторных
полей~\cite{Ha67}.

(с) Для любого гладкого вложения $f:N\to\R^6$ замкнутого ориентируемого
3-многообразия $N$ существует тройка линейно независимых нормальных векторных полей,  ср.~\cite{DW59};

(d) Существует гладкое вложение $\C P^2\to\R^8$, для которого нет четверки линейно
независимых нормальных полей.


(e) {\it Проблема Хирша.} Для каких $m$ и многообразий $N$ для любого вложения $N\to\R^m$
существует набор из $m-n$ линейно независимых нормальных векторных полей?
\end{pr}

Из этих утверждений читатель сможет доказать (с) для $N=S^3$ сейчас, а (с) для общего случая и (d) --- далее
(в соответствующих местах написано, что можно вернуться к решению этих задач).
Доказательство утверждений (a,b) выходит за рамки этой книги, а полное решение проблемы Хирша неизвестно.

\subsection{Доказательства непогружаемости и невложимости}

Попробуйте решить задачу \ref{ves4imm} и прочитать ее решение.
\footnote{Возможно, читателю вместо знакомства с этим пунктом будет интересно довести до конца доказательства
результатов, не требующие нормальных классов. Тогда нужно перейти к следующему параграфу.
Впрочем, поработать с наборами нормальных полей полезно для мотивировки понятия расслоения из следующего параграфа.}


\begin{pr}\label{vesnoex} Если замкнутое $n$-многообразие $N$

(a) погружаемо в $\R^{n+1}$, то $w_2(N)+w_1^2(N)=0$.

(b) вложимо в $\R^{n+2}$, то $w_2(N)+w_1^2(N)=0$.
\end{pr}

{\it Доказательство непогружаемости $(\R P^2)^3$ в $\R^8$.}
Обозначим $N:=(\R P^2)^3$ и $a:=[\R P^1]\in H_1(\R P^2)$.
Пусть $f:N\to\R^8$ --- погружение.
Рассмотрим препятствие $x$ к наличию шестерки линейно независимых векторных полей на
$f(N)$
(имеются в виду векторные поля в $\R^8$, на них не накладывается ни условие
касательности, ни условие нормальности).
 Тогда цепочка равенств в $H_3(N)$
$$0=y=w_3(N)+w_1^3(N)=
([a^2\times a\times 1]+a\times a\times a)+[a\times1\times1]^3=a\times a\times a\ne0$$
дает противоречие.
Здесь

$\bullet$ квадратные скобки обозначают `симметризацию' выражения.

$\bullet$ первое равенство выполнено, поскольку такая шестерка существует.

$\bullet$ второе равенство доказывается аналогично теореме Уитни-Ву и второму доказательству теоремы Штифеля
(см. детали ниже).

$\bullet$ третье равенство верно, поскольку $w(N)=(1+a+a^2)^3$ по теореме Уитни-Ву (задача \ref{vesnww}.a).

$\bullet$ четвертое равенство очевидно.

$\bullet$ пятое следует из $(a \times a \times a)\cap( a\times a\times a)=1\in\Z_2\ne0$.

Для доказательства второго равенства возьмем шестерку касательных векторных полей $v_1,\dots,v_6$ на $N$, для которых

$v_1\perp v_i$ и $v_1=0$ на некотором 5-мерном подкомплексе $\omega_1$ некоторой триангуляции $T$, представляющем класс
$w_1(N)$, и

пятерка $v_2,\dots,v_6$ линейно зависима на некотором 4-мерном подкомплексе $\omega_2$ триангуляции $T$, представляющем
класс $w_2(N)$.

четверка $v_3,\dots,v_6$ линейно зависима на некотором 3-мерном подкомплексе $\omega_3$ триангуляции $T$,
представляющем класс $w_3(N)$.

Возьмем также пару $u_1,u_2$ нормальных к $f(N)$ векторных полей на $f(N)$, для которых

$u_1=0$ некотором 4-мерном подкомплексе $\overline\omega_2$ триангуляции $T^*$, представляющем класс $\overline w_2(f)$
(определите, что это такое, аналогично пункту `нормальные векторные поля' предыдущего параграфа).

$u_2\perp u_1$ и $u_2=0$ на некотором 5-мерном подкомплексе $\overline\omega_1$ триангуляции $T^*$, представляющем класс
$\overline w_1(f)$ (определите, что это такое, аналогично пункту `нормальные векторные поля' \S\ref{0ve2ma}).

Теперь рассмотрим шестерку $u_1+v_1,u_2+v_2,v_3,v_4,v_5,v_6$ векторных полей на $f(N)$.
Эта шестерка линейно зависима в одном из следующих трех случаев:

$v_3,v_4,v_5,v_6$ линейно зависимы;

$u_2=0$ и $v_2,v_3,v_4,v_5,v_6$ линейно зависимы;

$u_1=0$ и $v_1=0$.

Поэтому множество линейной зависимости такой шестерки есть
$\omega_3\cup(\omega_2\cap\overline\omega_1)\cup(\omega_1\cap\overline\omega_2)$.
Из этого и пунктов (a,b) следующей задачи вытекает второе равенство (аналогично другому простому доказательству
теоремы Штифеля и задачам \ref{ves42}.e, \ref{ves43}.e).
QED

\begin{pr}\label{vesnoex4}
(a) Если замкнутое $n$-многообразие $N$ погружаемо в $\R^{n+2}$, то  $w_2(N)+w_1^3(N)=0$.

(b) Если замкнутое $n$-многообразие $N$ вложимо в $\R^{n+3}$, то  $w_2(N)+w_1^3(N)=0$.

(с) $(\R P^2)^4$ не погружаемо в $\R^{11}$ и не вложимо в $\R^{12}$.

(d) Если замкнутое $n$-многообразие $N$ погружаемо в $\R^{n+3}$ или вложимо в $\R^{n+4}$, то
$(w_4+w_2^2+w_2w_1^2+w_1^4)(N)=0$.
\end{pr}

\subsection{Нормальные классы Уитни}

Как получить аналогичное необходимое условие для погружаемости замкнутого $n$-многообразия в $\R^{n+k}$?

\begin{pr}\label{vesnowh}
(a) Для любого замкнутого $n$-многообразия $N$ существуют и единственны элементы
$\overline w_s(N)\in H_{n-s}(N)$, $s=0,1,2,\dots,n$, для суммы $\overline w(N)$ которых
выполнено $\overline w(N)\cap w(N)=1$.

Эти элементы называются {\it нормальными классами Уитни.}

(b) $\overline w_1(N)=w_1(N)$, $\overline w_2(N)=w_2(N)+w_1^2(N)$, $\overline w_3(N)=w_2(N)+w_1^3(N)$,
$\overline w_1(N)=(w_4+w_2^2+w_2w_1^2+w_1^4)(N)$.

(c) Докажите утверждение о погружаемости в следующей теореме.
\end{pr}

 {\bf Теорема Уитни о препятствии.}
{\it Пусть $N$ --- замкнутое $n$-многообразие.

Если $N$ погружаемо в $\R^m$, то $\overline w_s(N)=0\in H_{n-s}(N)$ при $s>m-n$.

Если $N$ вложимо в $\R^m$, то $\overline w_s(N)=0\in H_{n-s}(N)$ при $s\ge m-n$.}

\begin{pr}\label{vesnoemb} Пусть $N$ --- замкнутое ориентируемое $n$-многообразие,
$f:N\to\R^m$ --- вложение и $f':N\to\R^m$ --- близкое к $f$ отображение.

(a) Постройте препятствие Уитни
$$\overline W_{m-n}(N):=[\Sigma(f)]\in H_{2n-m}(N;\Z_{(m-n)})$$
к вложимости $n$-многообразия $N$ в $\R^m$ как гомологический класс множества самопересечений
для отображения общего положения~(аналогично \S\ref{0ve2ma} и \cite[\S2]{Sk08}).

(b) Этот класс не зависит от $f$.

(c) Если $N$ вложимо в $\R^m$, то $\overline W_{m-n}(N)=0\in H_{2n-m}(N;\Z_{(m-n)})$.
\end{pr}


 Чтобы использовать эти теоремы для получения конкретных результатов
(например, невложимости $\R P^4$ в $\R^7$), нужны вычисления.

\begin{pr}\label{vesnorp}
(a) Вычислите $\overline w((\R P^2)^k)$.

(b) Вычислите нормальные классы Уитни произведения любого числа 2-многообразий.
\end{pr}

 Теперь Вы уже можете решить задачу \ref{vesenpro}.

Приведем в качестве задачи прямое геометрическое определение классов $\overline w_s(N)$.

\begin{pr}\label{vesnodir}
Пусть $f:N\to\R^m$ --- погружение замкнутого $n$-многообразия $N$.
Аналогично параграфу 'векторные поля' можно определить {\it нормальный класс
Эйлера} $\overline W_{m-n}(f)\in H_{2n-m}(N;\Z)$ как препятствие к построению
ненулевого нормального векторного поля (ср. пункт 'нормальные векторные поля').
Аналогично многомерной теореме Штифеля о препятствии можно определить
{\it нормальный класс Штифеля-Уитни}
$\overline W_k(f)\in H_{n-k}(N;\Z_{(k-1)})$ как препятствие к построению
набора из $m-n+1-k$ линейно независимых нормальных векторных полей
($1<k<m-n$).
Напомним, что через $\Z_{(i)}$ обозначается группа $\Z$ для четного $i$ и
$\Z_2$ для нечетного~$i$.

(a) $\overline W_k(f)$ не зависит от $f$ при $k<m-n$.
 \footnote{А вот класс $\overline W_{m-n}(f)$ {\it может зависеть} от выбора
погружения $f$: например, существуют вложения замкнутых неориентируемых
2-многообразий в $\R^4$ с разными классами Эйлера.
\newline
Другое доказательство независимости $\rho_2\overline W_k(f)$ от $f$ при $k<m-n$,
не использующее регулярной гомотопии, следует из пункта (f).}

(b) Отображение
$\pi_l(V_{u+l,u})\to\pi_l(V_{u+l+1,u+1})$, индуцированное
добавлением одного вектора, является изоморфизмом при $u>1$ или четном $l$,
и является эпиморфизмом $\Z\to\Z_2$ при $u=1$ и нечетном $l$.

Указание. Рассмотрите точную последовательность расслоения
$V_{u+l,u}\to V_{u+l+1,u+1}\to S^{u+l}$.
Воспользуйтесь \cite[Теорема 11.4]{Pr06}.

(с) При фиксированном $N$ обозначим  $\overline W_k(N):=\overline W_k(f)$ при $m>k+n$.
Если $N$ погружаемо в $\R^m$, то $\overline W_s(N)=0\in H_{n-s}(N;\Z_{(m-n)})$ при $s>m-n$.

(d) Если $N$ вложимо в $\R^m$, то $\overline W_{m-n}(N)=0\in H_{2n-m}(N;\Z_{(m-n)})$.

(e) Это определение препятствия Уитни $\overline W_{m-n}(N)$ эквивалентно приведенному в задаче \ref{vesnoemb}.
\end{pr}

\begin{pr}\label{vesnoncalc}
(a) $\overline W_1(N)=w_1(N)$

(b) $\rho_2\overline W_2(N)=w_2(N)+w_1(N)^2$

(c) $\rho_2\overline W_k(N)=\overline w_k(N)$ для $1\le k\le m-n$.
\end{pr}


Аналогично нормальным классам Штифеля-Уитни можно определить {\it нормальные классы Понтрягина}
$\overline p_k(N)\in H_{n-4k}(N;\Z)$ для замкнутого ориентируемого $n$-многообразия $N$.
Они препятствуют погружаемости и вложимости:

{\it если $N$ погружаемо в $\R^m$, то $\overline p_k(N)=0$ для $2k>m-n$.}

{\it если $N$  вложимо в $\R^{n+2k}$, то $\overline p_k(N)=0$.}

\begin{pr}\label{vesno64}*
Пусть $N$ --- замкнутое ориентируемое связное 4-подмногообразие замкнутого ориентируемого 6-многообразия $M$.

(a) Определите группу $H_4(M;\Z)$ и форму пересечений $\cap:H_4(M;\Z)^3\to\Z$.

(b) Определите {\it число Понтрягина} $\overline p_1\in\Z$ как
полное препятствие к существованию {\it пары нормальных (к $N$) векторных
полей (касательных к $M$), имеющей ранг не менее одного в каждой точке}.

(c) Обозначим через $[N]\in H_4(M;\Z)$ гомологический класс подмногообразия $N$ и
через $\overline e\in H_2(N;\Z)$ и класс Эйлера нормального расслоения $N$ в $M$.
Тогда $[N]^3=\overline e\cap\overline e=\overline p_1$.

(d) $p_1(M)\cap[N]=p_1(N)+\overline p_1.$
\end{pr}



\subsection{Указания к некоторым задачам}

{\bf \ref{vesenotr}.} (с) для $N=S^3$. Докажите, что препятствие к существованию единичного векторного поля,
нормального и к $f(S^3)$ и к уже построенному нормальному полю, нулевое.

\smallskip
{\bf \ref{vesnoex4}.} (d) Следует из (e).

(e) Пусть $f:N\to\R^{11}$ --- погружение.
Рассмотрите препятствие $z$ к наличию набора из $n$ линейно независимых векторных полей на
$f(N)$. Докажите, что
$$0=z=w_4(N)+w_3(N)\overline w_1(N)+w_2(N)\overline w_2(N)+w_1(N)\overline w_3(N)$$
и воспользуйтесь пунктами (a,b,c).

\smallskip
{\bf \ref{vesnowh}.} (a) Индукция по $s$.

(b) Аналогично рассуждению в начале пункта и \ref{vesnoex}.b.


\smallskip
{\bf \ref{vesnodir}.} Обозначим $n=\dim N$ и возьмем вложение $g:N\to \R^{n+k}$.

(a) Возьмем два погружения $f:N\to\R^m$ и $g:N\to\R^p$.
Обозначим через $i_m$ включение $i:\R^m
\subset\R^{\max\{m,p,2n+1\}}$.
Аналогично определим $i_p$.
Тогда
$$\overline  W_k(f)=\overline  W_k(i_m\circ f)=\overline W_k(i_p\circ g)=
\overline  W_k(g).$$
Второе равенство следует из того, что по общему положению погружения $i_m\circ f$ и $i_p\circ g$ регулярно гомотопны
(определение напомнено в пункте `векторные расслоения').

Докажем первое равенство (третье доказывается аналогично).
Мы используем обозначения из доказательства теоремы Штифеля о препятствии (см. предыдущий пункт).
Рассмотрим отображение $V_{m-n,m-n-k+1}\to V_{s+m-n,s+m-n-k+1}$, определенное добавлением $s$
векторов.
Тогда индуцированное отображение групп $\pi_{k-1}$ является изоморфизмом (пункт (b)).
Поэтому препятствие $\overline W_k(f)$ равно препятствию
$\overline W_k(i_m\circ f)$ к построению семейства из $s+m-n+1-k$ линейно
независимых нормальных векторных полей для погружения $i_m\circ f$.

{\it Другое доказательство независимости $\overline w_k(f)$ от $f$.}
Возьмем два погружения $f,g:N\to\R^m$ и включение $i:\R^m\to\R^{2m}$.
Используем понятия расслоения из следующего параграфа.
Так как
$$\tau_N\oplus\nu_f\cong\tau_N\oplus\nu_g\cong m\varepsilon,\quad\text{то}
\nu_f\oplus m\varepsilon\cong\nu_g\oplus m\varepsilon.$$
То есть нормальные расслоения погружений $i\circ f$ и $i\circ g$ изоморфны.
Значит, $\overline  w_k(f)=\overline  w_k(i\circ f)=\overline w_k(i\circ g)=
\overline  w_k(g)$.

(c) Для композиции
$N\xrightarrow{f}\R^m\xrightarrow{i}\R^{m+n}$ погружения и включения
существует семейство из $n$ линейно независимых нормальных векторных полей,
поэтому $\overline W_{m-n+1}(N)=\overline W_{m-n+1}(i\circ f)=0$.

(d) Возьмем гомотопию $F:N\times I\to\R^m$ общего положения между $f'$ и отображением, образ которого
не пересекает $f(N)$.
Тогда $\Cl f^{-1}FN$ есть ориентированное погруженное подмногообразие (многообразия $N$) с краем $\Cl f^{-1}f'N$.
Поэтому $[\Cl f^{-1}f'N]=0$.

(e) Аналогично  случаю $m=2n$, рассмотренному в \S\ref{0ve2ma}.
Гомологический класс ориентированного погруженного подмногообразия
$\Cl f^{-1}f'N$ (многообразия $N$) равен классу $\overline W_{m-n}(N)$.

\smallskip
{\bf \ref{vesnoncalc}.}
(a) Одновременность обнуления обоих классов следует из того, что
обнуление касательного (нормального) класса равносильно возможности ввести
согласованную ориентацию в касательных (нормальных) пространствах.
Равенство этих классов --- более сильное условие.
(Поэтому, например,~\cite[доказательство Леммы 1 на стр. 178]{FF89} неполно.)

Построим препятствие к наличию $n+k$ линейно независимых векторных полей на $N\subset\R^{n+k}$
(на векторные поля не накладывается ни условие касательности, ни условие нормальности).
Это препятствие будет нулевым элементом группы  $H_{n-1}(N)$.

По наборам $e_1,\dots,e_n$ касательных векторных полей общего положения и $f_1,\dots,f_k$
нормальных векторных полей общего положения построим набор $e_1,\dots,e_n,f_1,\dots,f_k$
векторных полей на $N$. Гомологический класс множества вырождения последнего набора есть
$w_1(N)+\overline w_1(N)=0$.

(b) Построим препятствие к наличию $n+k-1$ линейно независимых векторных полей на $N\subset\R^{n+k}$
(на векторные поля не накладывается ни условие касательности, ни условие нормальности).
Это препятствие будет нулевым элементом группы  $H_{n-2}(N)$.
Рассуждение аналогично пункту (a) и задаче \ref{ves42}.e.

(c) Аналогично пунктам (a,b) и формуле Уитни-Ву (задача \ref{vesnww}.a).

\smallskip
{\bf \ref{vesno64}.}
(b) Число $\overline p_1$ --- сумма (со знаком) тех точек, в которых ранг пары
нормальных векторов общего положения меньше 1, т.е. в которых {\it оба} вектора
нулевые.

(c) Пусть $N_1,N_2\subset M$ --- близкие к $N$ подмногообразия, находящиеся в
общем положении относительно друг друга и $N$. Тогда
$$[N]^3=[N_1]\cap[N_2]\cap[N]=\#(N_1\cap N_2\cap N)=
\#[(N_1\cap N)\cap(N_2\cap N)]=\overline e\cap \overline e=\overline p_1.$$
Здесь $\#$ означает алгебраическую сумму точек пересечения в $M$, а
$N_1\cap N$, $N_2\cap N$ --- {\it ориентированные} пересечения в $N$.
Для доказательства последнего равенства возьмем два сечения нормального
расслоения $N$ в $M$, находящиеся в общем положении.
Тогда $\overline e$ --- гомологический класс надлежащим образом
ориентированного нулевого многобразия каждого сечения.
Поэтому $\overline e\cap\overline  e=\overline p_1$.

(d) Аналогично второму доказательству теоремы Штифеля (или ввиду $\tau_M|_N\cong\tau_N\oplus\nu(N\subset M)$, см. ниже).


\bigskip
\newpage
\section{Расслоения и их применения}\label{0vebun}

\subsection{Простейшие расслоения}

В этом пункте мы не докажем новых теорем, но
введем некоторые объекты, важные для дальнейшего и вообще в математике.
(Эти объекты близки к {\it утолщениям}, см. \cite{Sk}, параграфы `утолщения графов' и
`трехмерные утолщения двумерных полиэдров'.)

\begin{figure}[h]\centering
\includegraphics{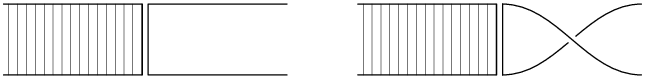}
\caption{Склейка концов ленточек}\label{7}
\end{figure}

\smallskip
{\it Определение $I$-расслоения над графом.}
Для графа $N$ с $V$ вершинами и $E$ ребрами возьмем $E$ ленточек, занумерованных
числами $1,2,\dots,E$.
Для каждого $i$ зафиксируем соответствие между концами $i$-й ленточки и
вершинами $i$-го ребра.
Склеим концы ленточек, соответствующие одной и той же вершине.
Заметим, что для каждой вершины это можно сделать несколькими способами
(рис.~\ref{7}).
Полученный занумерованный раскрашенный двумерный объект называется
{\it $I$-расслоением (или расслоением со слоем отрезок)} над графом $N$.

Примеры: цилиндр или лента Мебиуса над окружностью.

Два $I$-расслоения над графами называются {\it эквивалентными},
если можно так расставить стрелочки на отрезках склейки, что для
любых двух ленточек одного цвета в этих двух расслоениях стрелки
на концах обеих ленточек одновременно
сонаправлены или одновременно противонаправлены.

\smallskip
{\it Эквивалентное определение $I$-расслоения над графом для знакомых с двулистными накрытиями
(см. конец параграфа 'инволюции').}
Для двулистного накрытия $p:N'\to N$  и для
каждой точки $x\in N$ соединим отрезком две точки из $p^{-1}(x)$.
Полученный двумерный объект называется {\it расслоением со слоем отрезок (или
$I$-расслоением)} над $N$.
$I$-расслоения {\it эквивалентны}, если соответствующие двулистные накрытия
эквивалентны.

\smallskip
{\it Эквивалентное определение $I$-расслоения над графом для знакомых с основами топологии.}
Отображение $p:M\to N$ называется {\it $I$-расслоением} над $N$,
если любая точка $x\in N$ имеет такую окрестность $Ox$, что
$p^{-1}x\cong Ox\times I$ и композиция этого гомеоморфизма с проекцией на
$Ox$ совпадает с $p$.
Расслоения $p_1:M_1\to N$ и $p_2:M_2\to N$ называются
{\it эквивалентными}, если существует такой гомеоморфизм
$h:M_1\to M_2$, что $p_1=p_2\circ h$.)

\begin{figure}[h]\centering
\begin{tabular}{c@{\hskip3em}c}
\hbox{\includegraphics{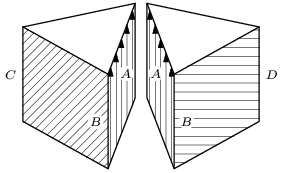}}
&
\hbox{\includegraphics{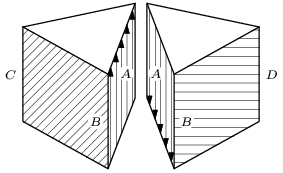}}\\
\hbox{\includegraphics{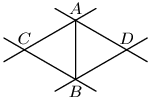}}
&
\hbox{\includegraphics{resko02.172.eps}}
\end{tabular}
\caption{Склейка призмочек}\label{8}
\end{figure}

\smallskip
{\it Определение $I$-расслоения над 2-многообразием $N$.}
Пусть дана триангуляция 2-многообразия $N$.
Раскрасим ее грани, ребра и вершины в разные цвета. Возьмем
призмы тех же цветов, что и грани триангуляции.
Раскрасим боковые грани и боковые ребра каждой призмы в
те же цвета, что соответствующие им ребра и вершины
соответствующей грани. Склеим боковые грани призм,
окрашенные в один и тот же цвет, так, чтобы склеились
боковые ребра, окрашенные в один и тот же цвет.
Для каждой боковой грани призмы это можно сделать
двумя способами (рис.~\ref{8}). Полученный раскрашенный
трехмерный объект называется {\it $I$-расслоением} над $N$.


Например, регулярная окрестность (т.е. окрестность, являющаяся утолщением)
2-многообразия в 3-многообразии является пространством $I$-расслоения (т. е.
является $I$-расслоением при некоторой раскраске).
Ср. \S\ref{0vesy3}.

Два $I$-расслоения над одним и тем же 2-многообразием называются {\it эквивалентными},
если на призмах можно так расставить стрелки, сонаправленные друг с другом и
параллельные боковым ребрам, что стрелки в одноцветных боковых гранях
призм этих двух расслоений одновременно сонаправлены или одновременно
противонаправлены.

\smallskip
Ясно, что {\it существуют взаимно однозначные соответствия между классами
эквивалентности инволюций над $N$, двулистных накрытий над $N$ и
$I$-расслоений над $N$.}

{\it $S^1$-расслоения над графами} и их эквивалентность определяются
аналогично $I$-расслоениям с заменой прямоугольников на кольца, являющиеся
прямыми произведениями окружности $S^1$ на ребра графа.

{\it $S^1$-расслоения над 2-многообразиями} и их эквивалентность определяются
аналогично $I$-расслоениям с заменой призмочек на полнотория, являющиеся
прямыми произведениями окружности $S^1$ на грани триангуляции.

\begin{pr}\label{}
Пространствами $S^1$-расслоений являются:
\quad

(a) $X\times S^1$ над $X$.
\quad
(b) Бутылка Кляйна над окружностью.
\quad
(c) Сфера $S^3$ над $S^2$.
\quad

(d)* Граница {\it регулярной окрестности} 2-многообразия $N$, вложенного в $\R^4$, над $N$.


(e) Пространство, полученное из несвязного объединения двух эквивалентных
$I$-расслоений над $N$ отождествлением их граничных двулистных накрытий по некоторой эквивалентности.
\end{pr}

\subsection{Векторные расслоения}

В топологии играет огромную роль следующее обобщение.
Мы используем его для вывода

$\bullet$ теоремы Штифеля о векторных полях и об алгебрах с делением из теоремы Штифеля о препятствии,

$\bullet$ теоремы Уитни из теоремы Уитни о препятствии,

$\bullet$ нижеследующей теоремы Масси, используемой при доказательстве теоремы Хефлигера-Хирша.

Мы следуем~\cite[\S2]{MS74}, ср.~\cite[\S19]{FF89}, \cite[Глава 5]{Pr04}.

{\it Векторным расслоением} называется отображение $p:E\to B$ вместе с заданием
структуры векторного пространства над $\R$ в множестве $p^{-1}b$ для каждой
точки $b\in B$, если выполнено следующее условие {\it локальной тривиальности}:

для любой точки $b\in B$ найдутся такие ее окрестность $Ob\subset B$, целое число $n\ge0$ и гомеоморфизм
$h_b:Ob\times\R^n\to p^{-1}Ob$, что для любой $a\in Ob$ сужение $h_b|_{a\times R^n}$ является изоморфизмом
векторных пространств $\R^n$ и $p^{-1}a$.

При этом пространство $B$ называется {\it базой} векторного расслоения, а отображение $p$ называется {\it проекцией}.

\begin{pr}\label{}
Следующие конструкции действительно дают векторные расслоения.
\end{pr}

{\it Примеры векторных расслоений (мы строим только отображения $p$, а
структуры векторных пространств в $p$-прообразах очевидны).}

($\varepsilon$) {\it Тривиальным} векторным расслоением называется проекция
$p:B\times\R^n\to B$.
Далее мы обозначаем через $\varepsilon$ одномерное тривиальное расслоение
(база которого ясна из контекста).

($\tau$) Пусть многообразие $N$ гладко вложено в $\R^m$.
Обозначим через $T_xN$ касательную ($n$-мерную) плоскость к многообразию $N$
в точке $x\in N$.
Определим
$$TN:=\{(x,v)\in N\times\R^m\ |\ v\in T_xN\}.$$
{\it Касательным расслоением} $\tau_N$ многообразия $N$ называется
отображение $\tau_N:TN\to N$, определенное формулой $\tau_N(x,v):=x$.

($\nu$) Пусть $f:N\to\R^m$ --- гладкое погружение.
Определим
$$E(f):=\{(x,v)\in N\times\R^m\ |\ v\perp T_{f(x)}f(N)\}.$$
{\it Нормальным расслоением} $\nu_f$ погружения $f$ называется
отображение $\nu_f:E(f)\to N$, определенное формулой $\nu_f(x,v):=x$.

\smallskip
Через $Dp$ и $Sp$ обозначим $D^n$- и $S^{n-1}$-подраслоения расслоения $p$.

\smallskip
{\bf Теорема о трубчатой окрестности.}
{\it  Если $N\subset\R^m$ --- замкнутое гладкое компактное подмногообразие и $\nu:E\to N$ --- нормальное расслоение,
то существует гладкое вложение $D\nu\to\R^m$, при котором нулевое сечение переходит в $N$.}

\smallskip
($\zeta_n$) Напомним, что $\R P^n$ есть пространство прямых в $\R^{n+1}$, проходящих через начало координат.
Определим
$$E(\zeta_n):=\{(y,v)\in \R P^n\times\R^{n+1}\ |\ v\in y\}.$$
{\it Тавтологическим расслоением} называется отображение
$$\zeta_n:E(\zeta_n)\to\R P^n,\quad\text{определенное формулой}\quad\zeta_n(y,v):=y.$$
\quad
(*) Пусть $B=\cup_i U_i$ --- открытое покрытие и
$\varphi_{ij}:U_i\cap U_j\to O_n$ непрерывные отображения, для которых
$$\varphi_{ii}=\id,\quad \varphi_{ij}=\varphi_{ji}^{-1}\quad\text{и}\quad
\varphi_{ik}=\varphi_{ij}\varphi_{jk}.$$
Положим
$$E:=\bigcup\limits_{\{(x,s)=(x,\varphi_{ij}s)\}_{x\in U_i\cap U_j}}
U_i\times\R^n\quad\text{и}\quad p[x,s]=x.$$

 Векторные расслоения $p_1:E_1\to B$ и $p_2:E_2\to B$ с одинаковой базой называются {\it изоморфными} (обозначение:
$p_1\cong p_2$), если существует гомеоморфизм $f:E_1\to E_2$, для которого $p_2f=p_1$ и сужение
$f|_{p_1^{-1}b}:p_1^{-1}b\to p_2^{-1}b$ есть изоморфизм векторных пространств для любой $b\in B$.

\begin{pr}\label{bun-exa}
(a) $\tau_{S^1}$, $\tau_{S^3}$ и $\tau_{S^7}$ тривиальны.

(b) $\tau_N$ тривиально тогда и только тогда, когда имеется $n$ касательных
ортонормированных векторных полей на $N$.

(c) $\tau_{S^{2k}}$ нетривиально.

(d) $\tau_{(S^1)^n}$ тривиально.

(e) $E(\zeta_1)$ есть лента Мебиуса, а $\zeta_1$ --- проекция на его среднюю линию.

(f) $\zeta_n\cong\nu(\R P^n\subset\R P^{n+1})$.

(g) Любое векторное расслоение над шаром $B^n$ изоморфно тривиальному.

(h)* Любое векторное расслоение (над конечным полиэдром $B$) может быть получено конструкцией (*).
(Поэтому конструкция (*) доставляет эквивалентное определение векторного расслоения.)
\end{pr}

Пусть задано векторное расслоение $p:E\to B$.
Отображение $s:B\to E$ называется {\it сечением расслоения $p$}, если $p\circ s=\id_B$.
Сечения называются {\it послойно гомотопными}, если они гомотопны в классе сечений.

Ясно, что касательное (нормальное) векторное поле --- то же, что сечение
касательного (нормального) расслоения.

\begin{pr}\label{}
(a) У проекции $p$ листа Мебиуса на его среднюю линию существует
ровно одно сечение (с точностью до послойной гомотопности).

(b) Обозначим через $S(p)$ множество классов послойной гомотопности сечений расслоения $p$.
Опишите $S(p)$ для стандартной проекции $p:K\to S^1$ бутылки Кляйна на окружность.

(с) $\zeta_n$ не является тривиальным ни при каком $n$
(указание: докажите, что $\zeta_n$ не имеет ненулевого сечения).
\end{pr}

Характеристические классы
Штифеля-Уитни $w_i(\xi)$ векторного расслоения $\xi$ определяются аналогично вышеизложенному как
препятствия к построению семейств линейно-независимых сечений.
Ясно, что $w_i(\tau_N)=w_i(N)$ и $w_i(\nu_f)=\overline w_i(N)$.

Пусть $p_1:E_1\to B$ и $p_2:E_2\to B$ --- векторные расслоения с одинаковой базой.
Определим
$$E(p_1\oplus p_2):=\{(x,y)\in E_1\times E_2\ |\ p_1(x)=p_2(y)\}.$$
{\it Суммой Уитни} векторных расслоений $p_1$ и $p_2$ называется отображение
$$p_1\oplus p_2:E(p_1\oplus p_2)\to B,\quad\text{определенное формулой}\quad p_1\oplus p_2(x,y):=p_1(x)=p_2(y).$$

\begin{pr}\label{bun-sum}
(a) Это действительно векторное расслоение.

(b) $\zeta_1\oplus\zeta_1\cong2\varepsilon$
\quad
(c) $\tau_{S^2}\oplus\varepsilon\cong3\varepsilon$. \quad
\quad
(d) $\tau\oplus\nu_f\cong m\varepsilon$.
\end{pr}

\begin{pr}\label{}
(a) $w_i(\xi\oplus n\varepsilon)=w_i(\xi)$.

(b) $w_1(\xi\oplus\eta)=w_1(\xi)\oplus w_1(\eta)$.

(c) Если $\dim\xi=m$ и $\dim\eta=n$, то $w_{m+n}(\xi\oplus\eta)=w_m(\xi)w_n(\eta)$.

(d) {\it Формула Уитни-Ву.}
{\it Полный класс Штифеля-Уитни} расслоения $\xi$ определяется как $w(\xi):=1+w_1(\xi)+w_2(\xi)+\dots$.
Тогда $w(\xi\oplus\eta)=w(\xi)w(\eta)$.

(e) {\it Формула Уитни-Ву для классов Понтрягина.}
{\it Полный класс Понтрягина} расслоения $\xi$ определяется как $p(\xi):=1+p_1(\xi)+p_2(\xi)+\dots$.
Если $H_i(B;\Z)$ не имеют 2-кручения, то $p(\xi\oplus\eta)=p(\xi)p(\eta)$.
\end{pr}

\begin{pr}\label{bun-rpt}
(a) $\tau_{\R P^n}\oplus\varepsilon\cong(n+1)\zeta_n$.

(b) $w(\R P^n)=(1+[\R P^{n-1}])^{n+1}$.

(c) $w_i(\R P^n)=\displaystyle{n+1\choose i}[\R P^{n-i}]$ для $0\le i\le n$.

(d) Если $\displaystyle{m\choose i}$ четно для любого $i=1,2,\dots,m-1$, то $m$ есть степень двойки.

(e) Докажите теорему Хопфа-Штифеля об алгебрах с делением.

(f) Докажите теорему Штифеля о векторных полях.

Указание к (e), (f): $[\R P^k]\ne0$ для $1\le k\le n$, поскольку $[\R P^k]\cap[\R P^{n-k}]=1$.
\end{pr}

\begin{pr}\label{bun-rpn}
 (a) $\overline w_i(\R P^n)=\displaystyle{-n-1\choose i}[\R P^i]={n+i\choose i}[\R P^i]$ для $0\le i\le n$.

(b) Докажите теорему Уитни.


(c) Если $n$-многообразие $N$ погружаемо в $\R^{n+1}$, то $w_i(N)=(w_1(N))^i$.

(d) Докажите теорему Уитни для коразмерности 1.
\end{pr}

\begin{pr}\label{bun-cpn}
(a) $\tau_{\C P^n}\oplus\varepsilon_\C\cong(n+1)\zeta_{n,\C}$ как {\it комплексные} векторные расслоения.

(b) $p(\C P^n)=(1+[\C P^{n-2}])^{n+1}$.

(c) $p_i(\C P^n)=\displaystyle{n+1\choose i}[\C P^{n-2i}]$ для $1\le i\le n/2$.

(d) $\overline p_i(\C P^n)=\displaystyle(-1)^i{ n+i\choose i}[\C P^{n-2i}]$ для $1\le i\le n/2$.


(e) Если $\C P^n$ вложимо в $\R^m$ или погружаемо в $\R^{m-1}$, то $\displaystyle{n+i\choose n}=0$ при $m-n\le 2i\le n$.
\end{pr}

{\it Произведением} векторных расслоений $p_1:E_1\to B_1$ и $p_2:E_2\to B_2$ называется отображение
$$p_1\oplus p_2:E_1\times E_2\to B_1\times B_2,\quad\text{определенное формулой}
\quad (p_1\times p_2)(x,y):=(p_1(x),p_2(y))$$
(с очевидно определяемыми структурами векторных пространств в слоях).

\begin{pr}\label{}
  (a) $\tau_{N_1\times N_2}=\tau_{N_1}\times\tau_{N_2}$.

(b) Если $\dim\xi=m$ и $\dim\eta=n$, то
$w_{m+n}(\xi\times\eta)=w_m(\xi)\times w_n(\eta)$.

(c) $w(\xi\times\eta)=w(\xi)\times w(\eta).$

(d)* $p(\xi\times\eta)=p(\xi)\times p(\eta).$
 \end{pr}



\subsection{Степени двойки и классы Штифеля-Уитни (набросок)}

Для нормальных классов Штифеля-Уитни произвольного замкнутого многообразия выполняются интересные
соотношения.
(На этом основаны теорема Хефлигера-Хирша и многие другие результаты.)
Например, старший нормальный класс всегда тривиален: $\overline W_n(N)=0$.
Это следует из вложимости любого $n$-многообразия в $\R^{2n}$ и теоремы Уитни о препятствии.
(Это также является частным случаем нижеследующей формулы Ву \ref{bun-wufor}.b.)

\smallskip
{\bf Теорема Масси.}
{\it Для замкнутого гладкого $n$-многообразия $N$

(1) Если либо $n$ не есть степень двойки, либо $N$ ориентируемо, то
$\overline W_{n-1}(N)=0$~\cite{Ma60},~\cite{Ma62}.

(2) Если $q$ меньше количества единиц в двоичной записи
числа $n$, то $\overline w_{n-q}(N)=0$~\cite{Ma60}.}

\smallskip
Доказательство теоремы Масси интересно тем, что является примером применения
следующего важного понятия.
Пусть класс $[X]\in H_{n-k}(N)$ реализуется вложением $i:X\subset N$.
Положим
$$\Sq_s[X]:=i_*w_s(\nu_N(X))\in H_{n-k-s}(N)$$
равным $i_*$-образу $s$-го класса Штифеля-Уитни нормального расслоения
вложения (это определение эквивалентно общепринятому~\cite[следствие на странице 263]{FF89}).

Далее замкнутое многообразие $N$ фиксировано и пропускается из обозначений
характеристических классов.
Задачи достаточно решить для классов, реализуемых подмногообразиями.
(См. пункт 'реализация циклов подмногообразиями' в \S\ref{0homol}.)
А справедливы они и для общего случая, для которого $\Sq_s$ определяются более сложно.
Для алгебраических вычислений этого пункта удобно {\it обозначить}
$H^s:=H_{n-s}$.
(Здесь никакого другого смысла, кроме переобозначений, не подразумевается, но вообще-то он имеется.)
Тогда умножение (т.е. билинейное отображение пересечения) в гомологиях
$n$-многообразия устроено так: $H^k\times H^l\to H^{k+l}$.


\begin{pr}\label{bun-sq}
(a) Проверьте корректность определения операции $\Sq_s$ относительно
гомологичности, реализуемой подмногообразием.

(b) {\it Общее положение.} Если $x\in H_{n-k}(N)$, то $\Sq_sx=0$ для $s>k$.

(c) {\it Случай класса Эйлера.} Если $x\in H_{n-k}(N)$, то $\Sq_kx=x\cap x$.

(d) $\Sq_1=\rho_2\beta$, где $\rho_2$ --- приведение по модулю 2 и $\beta$ --- гомоморфизм Бокштейна (\S\ref{0homol}).
Здесь многообразие $N$ ориентируемо (или коэффициенты скручены).

(e) {\it Формула Картана для $\Sq_1$.} $\Sq_1(x\cap y)=x\cap\Sq_1y+(\Sq_1x)\cap y$.

(f)* {\it Формула Картана.} $\Sq_k(x\cap y)=\sum_{s=0}^k(\Sq_sx)\cap(\Sq_{k-s}y)$.

(g) $\Sq_kx^{2^s}$ равно 0 при $k\ne0,2^s$.
 \end{pr}

\begin{pr}\label{}
Пусть $N$ ориентируемо.

(a) $\Sq_1w_{n-1}=0$.
\quad
(b)* $\Sq_1w_{2j+1}=0$.
\quad
(c)* $\Sq_1w_{2j}=w_{2j+1}$.
\quad
(d)* $\beta w_{2j}=W_{2j+1}$.
 \end{pr}

\begin{pr}\label{bun-wucla}
  (a) Для любого $k$ существует и единственен такой $v_k\in H_{n-k}(N)$, что $\Sq_kx=x\cap v_k$ при любом $x\in H_k(N)$.

(b) $v_1=w_1$ (значит, $\Sq_1w_{n-1}=w_1\cap w_{n-1}$).

(c) Если $N$ ориентируемо, то $v_2=w_2=\overline w_2$.
 \end{pr}

\begin{pr}\label{bun-wufor}*
{\it Формулы Ву.}
(a) $w_k=\sum\limits_{s=0}^k \Sq_sv_{k-s}$~\cite[Theorem 11.14]{MS74}.

(b) Если $x\in H_k(N)$ и $k>0$, то
$\overline w_k\cap x=\sum\limits_{i=1}^k\overline w_{k-i}\cap\Sq_ix$~\cite{Ma60}.

(c) $\Sq_kw_m=\sum\limits_{i=0}^k\displaystyle{ m-i-1\choose k-i}w_i\cap w_{m+k-i}$~\cite[30.2.A]{FF89}.
\end{pr}



 Поскольку операции $\Sq_s$ могут быть определены гомотопически, из формулы Ву
\ref{bun-wufor}.a вытекает, что классы Штифеля-Уитни многообразия гомотопически инвариантны.
См. другие формулы Ву в~\cite[Problem 11.E]{MS74},~\cite[Problem 11.F]{MS74}.

Чтобы проиллюстрировать {\it идею доказательства пункта (1) теоремы Масси}, покажем,
что {\it $\overline w_3=0$ для ориентируемого 4-многообразия $N$.}
Для любого $x\in H_3(N)$ имеем
$$
\overline w_3\cap x =
\overline w_2\cap\Sq_1x+\overline  w_1\cap\Sq_2x+\Sq_3x =
\overline w_2\cap x^2 =x^4=\Sq_1 x^3= 0.
$$
Здесь первое равенство верно по формуле Ву (\ref{bun-wufor}.b), второе по общему положению и
поскольку $\Sq_1x=x^2$ (случай класса Эйлера \ref{bun-sq}.c), третье поскольку
$\overline w_2=v_2$, четвертое по формуле Картана для $\Sq_1$ (\ref{bun-sq}.e), а последнее
поскольку и $N$, и окружность, представляющая класс $x^3$, ориентируемы.
 QED

\begin{pr}\label{bun-part1}
  (a) $\overline w_{n-1}\cap x=\overline w_{n-2^k}\cap x^{2^k}$
для $x\in H_{n-1}(N)$.

(b) Докажите пункт (1) теоремы Масси по модулю 2.

(c) Докажите пункт (1) теоремы Масси.
 \end{pr}

\begin{pr}\label{bun-adem}*
  (a) $\Sq_1\Sq_1=0$.

(b) {\it Соотношения Адема.} Если $0<a<2b$, то
$$\Sq_a\Sq_b=\sum\limits_{j=0}^{[a/2]}{b-1-j\choose a-2j}\Sq_{a+b-j}\Sq_j\mod2.$$
\quad
(с) Докажите пункт (2) теоремы Масси.
 \end{pr}

\subsection{Классификация расслоений}

\begin{pr}\label{bun-cl1}
  Множество классов эквивалентности одномерных векторных расслоений над
замкнутым ориентируемым $n$-многообразием $N$ находится во взаимно-однозначном
соответствии с $H_{n-1}(N)$.
 \end{pr}

 Векторное расслоение $p:E\to B$ называется {\it ориентируемым}, если в
векторных пространтвах $p^{-1}b$ заданы ориентации, причем изоморфизм
$x\mapsto h_b(b,x)$ сохраняет ориентацию.

Ориентированные векторные расслоения $p_1:E_1\to B$ и $p_2:E_2\to B$ с
одинаковой базой называются {\it ориентированно изоморфными}, если существует
гомеоморфизм $f:E_1\to E_2$, для которого $p_2f=p_1$ и сужение
$f|_{p_1^{-1}b}:p_1^{-1}b\to p_2^{-1}b$ есть сохраняющий ориентацию изоморфизм
векторных пространств.

\begin{pr}\label{}
(a) Какие расслоения из предыдущего пункта являются ориентируемыми?

(b) Любые два одномерных ориентируемых векторных расслоения над одной базой
ориентируемо эквивалентны.

(c) Определите класс Эйлера $e(\xi)$ как препятствие к построению ненулевого
сечения для ориентируемого векторного расслоения $\xi$.

(d) $e(\xi\oplus\eta)=e(\xi)\cap e(\eta)$ и
$e(\xi\times\eta)=e(\xi)\times e(\eta)$.

(e) Любые два ориентируемых двумерных векторных расслоения над
связным 2-многообразием с непустой границей ориентируемо эквивалентны.

(f) Множество ориентируемых двумерных векторных расслоений над замкнутым
ориентируемым связным 2-многообразием  с точностью до ориентируемой
эквивалентности находится во взаимно однозначном соответствии с $\Z$.


(g)* Как изменится результаты пунктов (e) и (f) для произвольных
2-многообразий?

(h)* Как изменится результат пункта (f) для векторных двумерных
расслоений?
Ср.~\cite{MF90}.
 \end{pr}

Для непрерывного отображения $f:B'\to B$ и векторного расслоения $p:E\to B$ определим
$E(f^*p):=\{(x,y)\in B'\times E\ |\ fx=py\}.$
{\it Индуцированным расслоением} называется отображение
$f^*p:E(f^*p)\to B'$, заданное формулой $(f^*p)(x,y):=x$
(с очевидно задаваемыми структурами векторных пространств в слоях).

Определение {\it гомотопности} отображений напомнено в начале \S\ref{0hopf}.
В следующих двух задачах $p:E\to B$ --- одномерное векторное расслоение.

\begin{pr}\label{}
 (a) Гомотопные отображения индуцируют эквивалентные расслоения.
Иными словами, формула $[f]\mapsto f^*\zeta_m$ корректно задает отображение
множества $[B,\R P^m]$ гомотопических классов отображений $B\to \R P^m$ в
множество классов эквивалентности $n$-мерных векторных расслоений над
полиэдром $B$.

(b) Отображение $f:B\to\R P^m$, для которого $p\simeq f^*\zeta_m$,
существует тогда и только тогда, когда существует {\it линейный мономорфизм}
$F:E\to\R^{m+1}$ (т.е. отображение, которое на каждом слое расслоения
$p:E\to B$ является линейным мономорфизмом).

(c) Для каких $m$ существует линейный мономорфизм $E(\zeta_1)\to\R^m$ (относительно расслоения $\zeta_1$)?

(d) Для любой точки $x\in B$ существует ее окрестность $Ox$ в $B$ и линейный мономорфизм $F:p^{-1}Ox\to\R$.

(e) Отображение $F:E\to\R^{m+1}$ является линейным мономорфизмом тогда и только тогда, когда для любой точки
$x\in B$ существует такое $i$, что отображение $\pro_i\circ F|_{p^{-1}x}:p^{-1}x\to\R$ является линейным мономорфизмом.
 \end{pr}

{\it Линейным функционалом} $\varphi:E\to\R$ называется отображение, линейное (но не обязательно мономорфное)
на слоях расслоения $p$.

\begin{pr}\label{bun-cla}
(a) Для любого подполиэдра $B_1\subset B$ любой линейный функционал $\varphi_1:p^{-1}B_1\to\R$ продолжается до
линейного функционала $\varphi:E\to\R$.

(b) Для некоторого $m$ существует линейный мономорфизм $\varphi:E\to\R^{m+1}$.
Иными словами, отображение $[\varphi]\mapsto \varphi^*\zeta_m$ сюръективно.

(c) Отображение $[\varphi]\mapsto \varphi^*\zeta_m$ инъективно при $m\ge2\dim B+2$.

(d) Для некоторого $m$ отображение $[\varphi]\mapsto \varphi^*\zeta_m$ биективно.

(e) Докажите следующую теорему.
\end{pr}

 {\it Многообразием Грассмана} $G(m,n)$ называется пространство всех $n$-мерных
линейных подпространств в $\R^m$.
Например, $G(m,1)\cong\R P^{m-1}$.

\smallskip
{\bf Теорема классификации векторных расслоений.}
{\it Множество классов эквивалентности $n$-мерных векторных расслоений над
полиэдром $B$ находится во
взаимно-однозначном соответствии с множеством $[B,G(m,n)]$ гомотопических
классов отображений $B\to G(m,n)$ при достаточно большом $m$.}

\smallskip
Взаимно-однозначное соответствие из сформулированной теоремы легко описать.
Определим
$$E(\zeta_{m,n}):=\{(\alpha,v)\in G(m,n)\times\R^m\ |\ v\in\alpha\}.$$
Определим {\it тавтологическим расслоение}
$$\zeta_{m,n}:E(\zeta_{m,n})\to G(m,n)\quad\text{формулой}\quad
\zeta_{m,n}(\alpha,v):=\alpha.$$
Соответствие из теоремы задается формулой $[f]\mapsto f^*\zeta_{m,n}$.

{\it Ориентируемым многообразием Грассмана} $G_+(m,n)$ называется пространство всех
ориентируемых $n$-мерных линейных подпространств в $\R^m$.
Например, $G_+(m,1)\cong S^{m-1}$.

\smallskip
{\bf Теорема классификации ориентируемых векторных расслоений.}
{\it Множество классов ориентируемой эквивалентности $n$-мерных ориентируемых
векторных расслоений над полиэдром $B$
при достаточно большом $m$ находится во взаимно-однозначном соответствии с
множеством $[B,G_+(m,n)]$ гомотопических классов отображений $B\to G_+(m,n)$.}

\begin{pr}\label{bun-claor}
  Докажите этот результат и явно постройте взаимно-однозначное соответствие.
 \end{pr}

 При помощи приведенных теорем характеристические классы можно получать из
{\it когомологий} пространств $G(m,n)$ или $G_+(m,n)$.
Иногда характеристические классы так и определяются.
Такое абстрактное алгебраическое определение неудобно для геометрических применений
(но, возможно, удобно для других целей).

\newpage
\subsection{Указания к некоторым задачам}

 {\bf \ref{bun-sum}.} (d) Нужный изоморфизм расслоений задается формулой $(x,v_1)\oplus(x,v_2)\mapsto(x,v_1+v_2)$.

\smallskip
{\bf \ref{bun-rpt}.} (a) Используем соглашения из примера ($\zeta_n$) перед задачей \ref{bun-exa}.
Обозначим через $l\subset\R^{n+1}$ прямую, проходящую через начало координат.
Касательный вектор в точке $l\in\R P^n$ можно канонически отождествить с линейным отображением
$l\to l^\perp\subset\R^{n+1}$.
(Можно и с $l^\perp$, но это отождествление не будет каноническим.)
Тогда пару из касательного вектора в точке $l\in\R P^n$ и числа можно канонически отождествить с парой линейных
отображений $l\to l$ и $l\to l^\perp$.
А последнюю пару --- с линейным отображением $l\to\R^{n+1}$, т.е. с упорядоченным набором длины $n+1$ линейных
фунционалов $l\to\R$.
Ввиду наличия скалярного произведения в $\R^{n+1}$ линейный функционал  $l\to\R$  можно канонически отождествить
с элементом из  $l$.

\smallskip
{\bf \ref{bun-wucla}.} (a) Докажите, что $\Sq_k:H_k(N)\to\Z_2$ --- гомоморфизм и используйте
двойственность Пуанкаре.

\smallskip
{\bf \ref{bun-wufor}.} (b) В случае $k=n=\dim N$ нужно доказать, что $\overline w_n(N)=0$, а это вытекает из теоремы Уитни.
Покажем, как свести утверждение для класса $x$, реализованного $k$-подмногообразием  $X$, к случаю $k=n=\dim N$.
Обозначим через $\nu_X$ нормальное расслоение $X$ в $N$,
через $\nu_N$ нормальное расслоение $N$ в $\R^m$ и через $i:X\to N$ включение.
Можно считать, что $\overline w_s=w_s(\nu_N)$.
Тогда
$$\sum\limits_{s=0}^k\overline w_{k-s}\cap_N\Sq_sx =
i\sum\limits_{s=0}^k(\overline w_{k-s}\cap x)\cap_Xw_s(\nu_X) =
i\sum\limits_{s=0}^k(\overline w_{k-s}(\nu_N|_X))\cap_Xw_a(\nu_X) =$$
$$=w_k(\nu_N|_X\oplus\nu_X) = w_k(\nu_{X\subset\R^m})=\overline w_k(X)=0.$$

{\bf \ref{bun-part1}.} (a) Индукция по $k$ с использованием формулы Ву, общего положения и формулы Картана
(задачи \ref{bun-wufor}.b и \ref{bun-sq}.bf).

(c) Используйте соотношение $\beta w_{2j}=W_{2j+1}$.

\smallskip
{\bf \ref{bun-adem}.} (c)
 Пересечение с классом $\overline w_{n-1}$ дает отображение
$H_{n-1}(N)\to H_0(N)$.
По формуле Ву его можно представить как сумму {\it итерированных стинродовых
квадратов} $\Sq_I:=\Sq_{i_1}\Sq_{i_2}\dots\Sq_{i_k}$ для разных
$I=(i_1,i_2,
\dots,i_k)$.
Из соотношений Адема следует, что любой  итерированный стинродов квадрат есть
сумма {\it допустимых}, т.е. удовлетворяющих условию $i_1\ge2i_2,\dots,
i_2\ge2i_3\dots$.
Поэтому если $\overline w_{n-q}\ne0$, то для некоторых $q\in H_{n-1}(N)$ и
допустимой последовательности $I$ имеем $\Sq_Ix\ne0$.


\smallskip
{\bf \ref{bun-cl1}.} Начните с $n=2$.
Можно делать аналогично классификации инволюций (\S\ref{0inv}) или свести к ней.

\smallskip
{\bf \ref{bun-cla}.} (b) Используйте компактность полиэдра $B$.

(e) Аналогично только что доказанному одномерному случаю.


\newpage
\section{Общие свойства гомологий (набросок)}\label{0homol}

\hfill{\it `The wise men followed the star and found the house.}

\hfill{\it But if I followed the star, should I find the house?'}

\hfill{\it `It depends, perhaps,' I said smiling, `on whether you are a wise man'.}

\hfill{\it G. K. Chesterton, Manalive.
\footnote{\it `Мудрецы следовали за звездой и нашли дом.
Но если бы я последовал за звездой, нашел бы я дом?' ---
`Возможно, это зависит' сказал я, улыбаясь, `от того, мудрец ли ты.'
Г.К. Честертон, Живой, пер. автора.}
}

\subsection{Простейшие свойства}\label{0homolsi}

Мы опускаем обозначения коэффициентов из записи групп гомологий,
если формула верна для произвольных коэффициентов (которые мы обозначаем $G$).
Советуем сначала проработать материал для коэффициентов $\Z_2$, затем для $\Z$.
Класс гомологий цикла $a$ обозначается через $[a]$.
Напомним, что $D^0$ --- точка.

\begin{pr}\label{prod0}
(a) $H_0(D^0;G)\cong G$ и $H_k(D^0;G)=0$ при $k>0$.

(b) $H_0(X;G)\cong G^{c(X)}$, где $c(X)$ --- количество компонент связности полиэдра $X$.
 \end{pr}

Положим $\t H_k(X):=H_k(X)$ при $k>0$ и $\t H_0(X):=\varepsilon^{-1}(0)/\partial_0C_1$, где $\varepsilon:C_0\to G$
--- сумма чисел, стоящих в вершинах.

\begin{pr}\label{prosn}
(a) $\t H_k(D^n)=0$.
\quad

(b) $\t H_k(\Con X)=0$, где $\Con$ --- конус.

(c) $\t H_k(D^0)=0$.

(d) Для $n>0$ группа $H_k(S^n;G)$ есть $0$ при $k\ne0,n$ и $G$ при $k=0,n$.

(e) $\t H_k(X)\cong\t H_{k+1}(\Sigma X)$, где $\Sigma$ --- надстройка.

(f) Группа $\t H_k(S^n_1\vee S^n_2\vee\dots\vee S^n_l;G)$ есть $0$ при $k\ne n$ и $G^l$ при $k=n$.

(g) Обозначим $B^n_l=S^n_1\vee S^n_2\vee\dots\vee S^n_l$.
Для симплициального
отображения $g:B^n_l\to B^n_m$ обозначим через $d_{pq}(g)$ степень композиции
$S^n_p\overset{\subset}\to B^n_l\overset{g}\to B^n_m \overset{proj}\to S^n_q$.
Тогда $g_*:H_n(B^n_l)\to H_n(B^n_m)$ --- линейное отображение с матрицей $d_{pq}(g)$ (в стандартных базисах).

(h) $H_k(X\vee Y)\cong H_k(X)\oplus H_k(Y)$ для $k>0$.

(i) Найдите $H_k(S^p\times S^q)$.
 \end{pr}

\begin{pr}\label{prof}
(a) Придумайте, как по отображению $f:X\to Y$ построить ({\it индуцировать}) отображение $f_*:H_k(X)\to H_k(Y)$.

(b) Отображение $f_*$ является гомоморфизмом.

(c)* Гомотопные отображения индуцируют одинаковые изоморфизмы.

(d) {\bf Теорема гомотопической инвариантности гомологий.}
{\it Если $X$ гомотопически эквивалентно $Y$, то $H_k(X)\cong H_k(Y)$.}
\end{pr}

\subsection{Гомологии пары, вырезание и точная последовательность}\label{0homolpa}

Пусть $A$ и $B$ --- подкомплексы комплекса $X$.
Тогда $C_k(A)\subset C_k(X)$.
Положим
$$Z_k(X,A):=\{a\in C_k(X)\ |\ \partial_{k-1}a\in C_{k-1}(A)\}\quad\text{и}\quad H_k(X,A):=\frac{Z_k(X,A)}{\partial_kC_{k+1}(X)\oplus C_k(A)}.$$
При $k=0$ условие $\partial_{k-1}a\in C_{k-1}(A)$ пусто.
(Такие группы появлялись при построении препятствий на многообразиях с краем.)

\begin{pr}\label{proexc}
(a) $H_0(X,A;G)\cong G^{c(X)}$, где $c(X)$ --- количество компонент связности полиэдра $X$, не содержащих подполиэдра $A$.

(b) $H_k(X)\cong H_k(X,\emptyset)$. \qquad

(c) {\bf Связь абсолютных и относительных гомологий.} $H_k(X,A)\cong\t H_k(X/A)$.

(d) {\bf Изоморфизм вырезания.} $H_k(A\cup B,A)\cong H_k(B,B\cap A)$.

{\it Более точно, отображение $H_k(B,B\cap A)\to H_k(A\cup B,A)$,
индуцированное включением $(B,B\cap A)\to(A\cup B,A)$, является изоморфизмом.}

(e) Если $A$ есть деформационный ретракт в $X$ (т.е. если существует гомотопия $f_t:X\to X$, для которой $f_0=\id X$, $f_1(X)\subset A$ и
$f_1(a)=a$ для любого $a\in A$), то $H_k(X,A)=0$ при $k\ge1$.
 \end{pr}

{\it Набросок другого решения задачи \ref{prosn}.d об $H_k(S^n)$.}
Индукция по $n\ge0$.
Докажем шаг индукции.
Ввиду изоморфизма вырезания $H_k(S^n)\cong H_k(D^n,S^{n-1})$.
Определим отображение $\partial_*:H_k(D^n,S^{n-1})\to H_{k-1}(S^{n-1})$ формулой $\partial_*[a]=[\partial a]$.
В задаче \ref{prod} Вы проверите корректность этого определения.
Достаточно доказать, что отображение $\partial_*$ является изоморфизмом.
Для доказательства эпиморфности отображения $\partial_*$ возьмем цикл $a$ и $[a]\in H_{k-1}(S^{n-1})$.
По задаче \ref{prosn}.a $H_{k-1}(D^n)=0$.
Поэтому существует такая цепь $b$, что $a=\partial b$.
Тогда $b$ есть относительный цикл пары $(D^n,S^{n-1})$ и $[a]=\partial_*[b]$.
QED

\begin{pr}\label{prod} (a) Отображение $\partial_*:H_k(D^n,S^{n-1})\to H_{k-1}(S^{n-1})$ корректно определено формулой $\partial_*[a]=[\partial a]$.

(b) Отображение $\partial$ мономорфно.
\end{pr}

{\it Гомологической последовательностью пары $(X,A)$}
$$\dots\to H_k(A)\overset{i_*}\to H_k(X)\overset{j_*}\to H_k(X,A)
\overset{\partial_*}\to H_{k-1}(A)\to\dots$$
называется последовательность групп и гомоморфизмов, определенных так.
Гомоморфизм $i_*$ индуцирован включением.
При $k\ge1$ отображения $\partial_*$ и $j_*$ определены формулами
$\partial_*[a]=[\partial a]$ и $j_*[a]=[a]$, соответственно.

\begin{pr}\label{procopa}
(a) Отображение $\partial_*$ корректно определено и является гомоморфизмом.

(b) Если $H_k(X)=H_{k-1}(X)=0$ (например, если $X$ стягиваемо), то $\partial_*$ изоморфизм.

(c) При $k\ge1$ отображение $j_*$ корректно определено формулой $j_*[a]=[a]$ и является гомоморфизмом.

(d) Если $H_k(A)=H_{k-1}(A)=0$ (например, если $A$ стягиваемо), то $j_*$ изоморфизм.

(e) Если $A$ есть ретракт пространства $X$ (т.е. если существует отображение $f:X\to A$, тождественное на $A$), то $i_*$ --- мономорфизм,
$j_*$ --- эпиморфизм и $\partial_*=0$.

(f) Если $A$ стягиваемо по $X$ в точку (т.е. если существует отображение $\con A\to X$, сужение которого на основание есть включение), то
$i_*=0$, $j_*$ --- мономорфизм и $\partial_*$ --- эпиморфизм.

(g) Если существует гомотопия $f_t:X\to X$, для которой $f_0=\id$ и $f_1(X)\subset A$, то $i_*$ --- эпиморфизм, $j_*=0$ и $\partial_*$
--- мономорфизм.
 \end{pr}

Оказывается, все факты, сформулированные в предыдущей задаче, являются частными случаями {\it одной} теоремы.

Последовательность групп и гомоморфизмов называется {\it точной}, если ядро
каждого ее гомоморфизма совпадает с образом предшествующего гомоморфизма.

\begin{pr}\label{proesp}
(a) {\bf Теорема о точности последовательности пары.}
{\it Гомологическая последовательность пары $(X,A)$ точна.}

(b) Вычислите гомологическую последовательность пар $(S^n,D^n)$ и $(S^n,S^{n-1})$.
 \end{pr}

\begin{pr}\label{geho-fiex}
(a) Последовательность $0\to G\overset\varphi\to H$ точна
тогда и только тогда, когда $\varphi$ есть мономорфизм.
Последовательность $G\overset\varphi\to H\to0$ точна
тогда и только тогда, когда $\varphi$ есть эпиморфизм.
Последовательность $0\to G\overset\varphi\to H\to0$ точна
тогда и только тогда, когда $\varphi$ есть изоморфизм.

(b) {\it Короткая точная последовательность}
$0\to G\overset\varphi\to H\overset\psi\to K\to0$ точна тогда и только
тогда, когда $\varphi$ есть мономорфизм, $K\cong H/\varphi(G)$ и $\psi$ есть
естественная проекция.
\end{pr}

\subsection{Другие точные последовательности}\label{0homolex}

\begin{pr}\label{procomv} (a) Если $A\cup B$ стягиваемо, то при $k\ge1$ разность гомоморфизмов включений
$\alpha=i_A-i_B:H_k(A\cap B)\to H_k(A)\oplus H_k(B)$ есть изоморфизм.

(b) Если $A\cap B$ стягиваемо, то при $k\ge1$ сумма гомоморфизмов включений
$\beta=I_A\oplus I_B:H_k(A)\oplus H_k(B)\to H_k(A\cup B)$ есть изоморфизм.

(c) Композиции
$$H_k(A\cup B)\overset{j_*}\to H_k(A\cup B,A)\cong H_k(B,A\cap B)
\overset{\partial_*}\to H_{k-1}(A\cap B)\quad\text{и}\quad$$
$$H_k(A\cup B)\overset{j_*}\to H_k(A\cup B,B)\cong H_k(A,A\cap B)
\overset{\partial_*}\to H_{k-1}(A\cap B)$$
(в которых $\cong$ есть изоморфизм вырезания) совпадают с точностью до знака.

(d) Если $A$ и $B$ стягиваемы, то при $k\ge1$ каждый гомоморфизм из (c) есть
изоморфизм $H_k(A\cup B)\to H_{k-1}(A\cap B)$.

(e) Докажите следующий результат для случая, когда $A$ и $B$ --- многообразия
размерности $n$ с краем, пересекающиеся по $(n-1)$-мерному подмногообразию края.
 \end{pr}

{\bf Теорема о последовательности Майера-Виеториса.}
{\it Последовательность
$$\to H_k(A\cap B)\overset{i_A-i_B}\to H_k(A)\oplus H_k(B)
\overset{I_A\oplus I_B}\to H_k(A\cup B)\overset{\gamma}\to H_{k-1}(A\cap B)\to$$
точна.
Здесь гомоморфизмы $i_A-i_B$, $I_A\oplus I_B$ и $\gamma$ определены в предыдущей задаче.}

\smallskip
Эта последовательность является аналогом формулы включений-исключений.

\begin{pr}\label{procom}
Вычислите гомологии

(a) Дополнения к узлу в трехмерном пространстве.

(b) Дополнения к вложенной сфере с ручками в четырехмерном пространстве.
Указание: найдется окрестность сферы с ручками, граница которой диффеоморфна $N\times S^1$.
\end{pr}

 \begin{pr}\label{pro}
Пусть $f_A$ --- линейный автоморфизм тора $T^2$, заданный матрицей $A$.
Вычислите гомологии

(a) Пространства $N\cup_{f_A}D^2\times S^1$, где $N$ --- данное 3-многообразие
с краем $\partial N\supset T^2$.

(b) Пространства $T^2\times I/(x,0)\sim(f_A(x),1)$.
 \end{pr}

Пусть $X$ --- полиэдр.
Если задан гомеоморфизм $f:X\to X$  на себя, то {\it пространством $X$-расслоения над окружностью} называется пространство,
полученное из $X\times I$ отождествлением точек $(x,0)$ и $(f(x),1)$.
Эквивалентность $X$-расслоений определяется аналогично пункту `простейшие расслоения' в \S\ref{0vebun}.

\begin{pr}\label{probuci}*
(a) Пусть $N$ и $X$ --- замкнутые ориентируемые 2-многообразия.
Если некоторое 3-многообразие $M$ является одновременно $S^1$-расслоением
над $N$, не эквивалентным проекции $N\times S^1\to N$, и $X$-расслоением над
окружностью, не эквивалентным проекции $N\times S^1\to N$, то $N=X=T^2$.


(b) Существует 3-многообразие, которое является одновременно непрямым $S^1$-расслоением над $T^2$ и непрямым
$T^2$-расслоением над окружностью.

(c) `Опишите' все расслоения из (b).
 \end{pr}

Как выразить гомологии с коэффициентами $\Z_2$ через гомологии с коэффициентами $\Z$?

\begin{pr}\label{proun}
{\bf Формулы универсальных коэффициентов.}
{\it Для простого $p$
$$H_k(X;\Z_p)\cong(H_k(X;\Z)\otimes\Z_p)\oplus(H_{k-1}(X;\Z)*\Z_p),$$
где $G\otimes\Z_p:=G/pG$ и $G*\Z_p$
--- совокупность элементов порядка 1 или $p$ в $G$.

$\text{Если}\quad H_k(X;\Z)\cong\Z^b\oplus\Tors H_k(X;\Z),\quad\text{то}\quad H_k(X;\Q)\cong\Q^b.$}
 \end{pr}

{\it Замечание.}
Аналогичные простые формулировки и доказательства имеют остальные {\it формулы универсальных коэффициентов},
выражающие (ко)гомологии с коэффициентами в $\Z$, $\Z_p$ и $\Q$ через гомологии  с коэффициентами в $\Z$.
Чтобы сделать все эти простые и полезнейшие результаты менее доступными, их
обычно формулируют и доказывают на языке абстрактных тензорных, периодических и $Ext$-произведений.
Впрочем, этот язык (как и более экзотические коэффициенты) полезен, например, для изучения обобщений примера Понтрягина
двумерных компактов $A$ и $B$ с трехмерным произведением $A\times B$ (Бокштейн, Кузьминов, Дранишников).

\begin{pr}\label{probock}
  (a) {\it Последовательность Бокштейна}
$$\dots\to H_{k+1}(X;\Z_2)\overset{\beta_k}\to H_k(X;\Z)\overset2\to H_k(X;\Z)
\overset{\rho_2}\to H_k(X;\Z_2)\to\dots$$
точна.
Здесь $2$ --- умножение на 2, $\rho_2$ --- приведение по модулю 2, а
{\it гомоморфизм Бокштейна} $\beta_k$ определяется так.
Если $a=a_1+\dots+a_s$ --- цикл по модулю 2, где $a_1,\dots,a_s$ являются
$(k+1)$-симплексами некоторой триангуляции, то возьмем целочисленную цепь
$\t a=a_1+\dots+a_s$, положим $\beta_k[a]:=[\frac12\partial\t a]$ и проверим
корректность.

(b) Докажите формулы универсальных коэффициентов для $p=2$.

(с) То же для произвольного $p$.

(d) То же для коэффициентов $\Q$.

(e)* Выразите форму пересечений в гомологиях с коэффициентами в $\Z_2$ через форму пересечение и форму зацеплений
в гомологиях с коэффициентами в $\Z$ (форма зацеплений определена в \S\ref{0dual}).
 \end{pr}

\begin{pr}\label{proku}
{\bf Формула Кюннета.}
{\it Для простого $p$}
$$H_k(X\times Y;\Z_p)\cong\oplus_{i=0}^k H_i(X;\Z_p)\otimes H_{k-i}(Y;\Z_p).$$
 \end{pr}

\subsection{Указания к некоторым задачам}


{\bf \ref{prosn}.} (a) Аналогично (b) или следует из (b).

(b) $\partial\Con x=x$.

(e) Любой $k$-цикл в $\Sigma X$ гомологичен надстройке над некоторым циклом в $X$.

\smallskip
{\bf \ref{prof}.} См. детали в ~\cite{Pr06}.

\smallskip
{\bf \ref{proexc}.} (d) Сведите к случаю, когда $A-B$ есть внутренность одного симплекса старшей размерности.

\smallskip
{\bf \ref{prod}.} (b) Аналогично эпиморфности.

\smallskip
{\bf \ref{proesp}.} Аналогично приведенному наброску доказательства теоремы о гомологиях сфер и задаче \ref{procopa}.

\smallskip
{\bf \ref{procomv}.}
(a) Эпиморфность отображения $i_A$ доказывается так.
Пусть дан произвольный цикл $a$ в $A$.
Так как $A\cup B$ стягиваемо, то существует цепь $c$ в $A\cup B$, для которой $\partial c=a$.
Тогда $c$ есть относительный цикл в $(A\cup B,A)$.
При изоморфизме вырезания ему соответствует относительный цикл $c'=c\cap B$ в $(B,A\cap B)$.
Тогда $\partial c'$ есть цикл в $A\cap B$, гомологичный в $A$ циклу $a$.

(d) Рассмотрите случай $A\cup B=D^n_+\cup_{S^{n-1}}D^n_-$.

\smallskip
{\bf \ref{probuci}.} Вычислите гомологии 3-многообразия $M$ через $S^1$-расслоение над $N$ и через $X$-расслоение над $S^1$.

\smallskip
{\bf \ref{probock}.} (a)
Докажем, что $\rho_2\circ\beta_k=0$.
Рассмотрим цикл $\alpha$, представляющий некоторый элемент группы $H_{k+1}(X;\Z)$.
Значит,  $\partial\alpha=0$.
Возьмем $\widetilde{\rho_2[\alpha]}:=\alpha$.
Тогда $\beta_k(\rho_2[\alpha])=[\frac12\partial\widetilde{\rho_2[\alpha]}]=[\frac12\partial\alpha]=0$.

Теперь докажем, что если $\beta_k(a)=0$, то $a=\rho_2(\alpha)$ для некоторого цикла $\alpha$.
Если $\beta_k(a)=0$, то по определению $[\frac12\partial\widetilde a]=0$. \
Значит, $\frac12\partial\widetilde a=\partial\alpha'$ для некоторой цепи $\alpha'$.
Отсюда $\partial\widetilde a=2\partial\alpha'$ и $\partial(\widetilde a-2\alpha')=0$.
Поэтому $\alpha:=\widetilde a-2\alpha'$ является циклом.
При этом $a=\rho_2(\widetilde a)=\rho_2(\widetilde a-2\alpha')$.


\newpage
\section{Двойственности Пуанкаре и Александера-Понтрягина}\label{0dual}

\subsection{Простая часть двойственности Пуанкаре}

Возьмем разбиение $T$ \ $n$-многообразия $N$ на многогранники и двойственное разбиение $T^*$.
Будем использовать обозначения, введенные в начале \S\ref{0homol}.
Тогда количество $i$-мерных клеток в $T^*$ равно $c_{i*}:=c_{n-i}$.
Используя $T^*$ вместо $T$, определим аналогично \S\ref{0homol} группы $C_i(T^*)$ и отображения
$\partial_{i*}:C_{i+1}(T^*)\to C_i(T^*)$.

\smallskip
{\it Доказательство изоморфизма $H_1(N)\cong H_2(N)$ для замкнутого 3-мно\-го\-об\-ра\-зия $N$.}
Этот изоморфизм вытекает из теоремы инвариантности гомологий и
$$\dim H_1(T)=\dim\partial_0^{-1}(0)-\dim\partial_1(C_2)=
c_1-\rk\partial_0-\rk\partial_1\overset{(1)}=
 c_{2,*}-\rk\partial_{2,*}-\rk\partial_{1,*}\overset{(2)}=\dim H_2(T^*).$$
Равенство (2) доказывается аналогично первым двум.
Докажем равенство (1).
Ясно, что $\alpha\subset\beta$ тогда и только тогда, когда $\beta^*\subset\alpha^*$.
Значит, матрицы операторов $\partial_{0*}$, $\partial_{1*}$ и $\partial_{2*}$ в стандартных базисах групп
$C_i(T^*)$ получаются {\it транспонированием} из матриц операторов $\partial_2$, $\partial_1$ и $\partial_0$ в стандартных базисах групп $C_i(T)$.
Это влечет равенство (1).
QED

\begin{pr}\label{pdsimz2}
{\bf Теорема двойственности Пуанкаре по модулю 2 (простая часть).}
{\it Для замкнутого $n$-мно\-го\-об\-ра\-зия $N$ выполнено $H_i(N)\cong H_{n-i}(N)$. }
\end{pr}

\begin{pr}\label{pdsim3z}
Для замкнутого ориентируемого 3-мно\-го\-об\-ра\-зия $N$

(a) кручение группы $H_2(N;\Z)$ нулевое.

(b) $H_2(N;\Z)\cong H_1(N;\Z)/Tors$.
\end{pr}

\begin{pr}\label{pdsimz}
(a) {\bf Теорема двойственности Пуанкаре (простая часть).}
{\it Для замкнутого ориентируемого $n$-мно\-го\-об\-ра\-зия $N$

$\bullet$ свободные части групп $H_i(N;\Z)$ и $H_{n-i}(N;\Z)$ изоморфны;

$\bullet$ кручения групп $H_i(N;\Z)$ и $H_{n-i-1}(N;\Z)$ изоморфны.}

(b) {\bf Теорема двойственности Лефшеца  по модулю 2 (простая часть).}
{\it Для (компактного) $n$-многообразия $N$ с краем $H_i(N)\cong H_{n-i}(N,\partial N)$.}

(c) {\bf Теорема двойственности Лефшеца (простая часть).}
{\it Для ориентируемого $n$-многообразия $N$ с краем

$\bullet$ свободные части групп $H_i(N;\Z)$ и $H_{n-i}(N,\partial N;\Z)$ изоморфны;

$\bullet$ кручения групп $H_i(N;\Z)$ и $H_{n-i-1}(N,\partial N;\Z)$ изоморфны.}
\end{pr}

Заметим, что иногда (с целью усложнения доказательства) простую часть теоремы
двойственности Пуанкаре доказывают с использованием индекса пересечения.

 \subsection{Сложная часть двойственности Пуанкаре}

Впрочем, индекс пересечения действительно крайне полезен:
например, для `сложной части' теоремы двойственности Пуанкаре.
Определение формы пересечений приведено перед задачей \ref{vesint}.
Следующие задачи и теоремы интересны даже для $n=3$.

\smallskip
{\bf Теорема двойственности Пуанкаре по модулю 2 (сложная часть).}
{\it Для замкнутого $n$-мно\-го\-об\-ра\-зия $N$ форма пересечений по модулю 2
невырождена, т.е. для любого $\alpha\in H_i(N)-\{0\}$ существует такой
$\beta\in H_{n-i}(N)$, что $\alpha\cap\beta=1\in\Z_2$.}

\begin{pr}\label{pdcom1}
(a) Докажите теорему для $i=n-1$.

(b) Выведите случай $i=1$ из случая $i=n-1$.
\end{pr}

Простая идея, работающая для $i=1,n-1$, не годится для других случаев.
Приведем в виде цикла задач план другого (чуть более сложного) доказательства случая $n=3$, $i=1$,
которое годится для общего случая (см. детали в~\cite[\S71]{ST04}).

\begin{pr}\label{}
(a) Пусть базис $x_1,\dots,x_{c_1}$ группы $C_1$ получен из
стандартного базиса этой группы оператором $A$.
Базис $y_1,\dots,y_{c_1}$ группы $C_{2*}$ получен из стандартного базиса этой
группы оператором $(A^{-1})^*$ тогда и только тогда, когда $x_i\cap y_j=\delta_{ij}$.

Указание. При преобразованиях базисов в $C_1$ и $C_{2*}$, заданных матрицами
$A$ и $B$, матрица $M$ билинейной формы переходит в $A^TMB$.

(b) Пусть $x_1,\dots,x_{c_1}$ и $y_1,\dots,y_{c_1}$ --- такие базисы групп
$C_1$ и $C_{2*}$, что $x_i\cap y_j=\delta_{ij}$.
Возьмем стандартные базисы в группах $C_0$ и $C_{3*}$.
Тогда матрицы операторов $\partial_0$ и $\partial_2^*$ получаютя друг из друга транспонированием.

(c) Пусть $x_1,\dots,x_{c_1}$ --- такой базис группы $C_1$, что

$x_1,\dots,x_{b_1}$ --- базис группы 1-границ и $x_1,\dots,x_{z_1}$ --- базис группы 1-циклов.

Возьмем такой базис $y_1,\dots,y_{c_1}$ группы $C_{2*}$, что $x_i\cap y_j=\delta_{ij}$.
Тогда

$y_{b_1+1},\dots,y_{c_1}$ --- базис группы 2-циклов и $y_{z_1+1},\dots,y_{c_1}$ --- базис группы 2-границ.

(d) Для любого базиса $\alpha_1,\dots,\alpha_k$ группы $H_1(N)$ существует базис
$\beta_1,\dots,\beta_k$ группы $H_2(N)$, для которого $\alpha_i\cap\beta_j=\delta_{ij}$.

(e) Докажите сложную часть теоремы двойственности Пуанкаре по модулю 2.

(f) {\bf Теорема двойственности Лефшеца по модулю 2.}
{\it Для $n$-многообразия $N$ с краем форма пересечений по модулю 2
$H_i(N)\times H_{n-i}(N,\partial N)\to\Z_2$ невырождена.}
\end{pr}

\begin{pr}\label{pdcomlk}
(a) Для замкнутого ориентируемого $n$-мно\-го\-об\-ра\-зия $N$ целочисленная форма пересечений $\cap$ унимодулярна, т.е. для любого примитивного
(т.е. не делящегося на целое число, большее 1) элемента $\alpha\in H_i(N;\Z)$ существует такой $\beta\in H_{n-i}(N;\Z)$, что $\alpha\cap\beta=1\in\Z$.

(b) Для замкнутого ориентируемого $n$-мно\-го\-об\-ра\-зия $N$  определим форму зацеплений
$$\lk:\Tors H_i(N;\Z)\times\Tors H_{n-1-i}(N;\Z)\to\Q/\Z\quad\text{формулой}\quad
\lk([a],[b]):=\{\dfrac{a\cap B}n\},$$
$\quad\text{где}\quad\partial B=nb.$
Это определение корректно.

(c) Найдите форму зацеплений для $\R P^3$.

(d) Форма зацеплений невырождена.
\end{pr}

{\bf Теорема двойственности Пуанкаре.}
{\it Для замкнутого ориентируемого $n$-мно\-го\-об\-ра\-зия

$\bullet$ целочисленная форма пересечений $\cap$  унимодулярна;

$\bullet$ форма зацеплений $\lk$ невырождена.}

\begin{pr}\label{pdcomz}
{\bf Теорема двойственности Лефшеца.}
{\it Для ориентируемого $n$-многообразия $N$ с краем

$\bullet$  целочисленная форма пересечений $\cap:H_i(N;\Z)\times H_{n-i}(N,\partial N;\Z)\to\Z$ унимодулярна,

$\bullet$ форма зацеплений $\lk:\Tors H_i(N;\Z)\times\Tors H_{n-1-i}(N,\partial N;\Z)\to\Q/\Z$ невырождена.}
 \end{pr}

Надеюсь, следующее шутливое замечание будет полезно начинающему.
`Сложную часть' теоремы двойственности Пуанкаре часто либо не доказывают (в~\cite{FF89} теорема 6 в \S17 на стр. 148 названа очевидной),
либо доказывают с использованием когомологий и формулы универсальных коэффициентов (с целью усложнения доказательства).
Когомологии действительно полезны при работе с дифференциальными формами, при изучении алгебраической геометрии
или гомотопической топологии многообразий с краем или произвольных комплексов.
А во многих учебниках и лекционных курсах {\it когомологии многообразий}
вводятся намного раньше тех проблем, для решения которых они нужны.
В итоге когомогогии применяются для усложнения доказательств,
а также чтобы дать возможность студентам почаще допускать ошибки (что легче
сделать, если обозначать канонически изоморфные группы по-разному).


\subsection{Двойственность Александера и ее применения}

Препятствия к {\it вложимости} полиэдра в $S^m$ можно получить, основываясь на следующей идее, восходящей к Александеру.
Рассматривая дополнение $S^m-N$ полиэдра $N\subset S^m$, можно вывести необходимые условия на сам полиэдр $N$.
Эта идея видна уже на доказательстве невложимости графов в плоскость с помощью формулы Эйлера.
В общем случае формула Эйлера заменяется на ее обобщение --- двойственность Александера (см. ниже).

Ясно, что $H_1(S^3-f_0S^1)=\Z$ для стандартного вложения (т.е. тривиального узла) $f_0:S^1\to S^3$.
Для другого узла $f:S^1\to S^3$ группу $H_1(S^3-fS^1)$ можно вычислить, располагая весь узел, кроме `переходов', в плоскости проекции, а `переходы' над плоскостью, и применяя последовательность Майера-Виеториса.
По-видимому, Александер, проделывая такие вычисления, заметил, что $H_1(S^3-fS^1)\cong\Z$ для любого узла $f$.
Это привело его к открытию двойственности, носящей его имя.
(Если вместо
группы $H_1$
использовать {\it фундаментальную группу}, то узлы уже можно различать, см. \S\ref{03man} или ~\cite[VI]{Pr04}.)

Полезно также использовать двойственность с учетом пересечений и зацеплений (т.е. двойственность между {\it кольцами гомологий} многообразия $N\subset S^m$ и его дополнения $S^m-N$).
Так Х. Хопф доказал в 1940, что $\R P^n$ не вложимо в $S^{n+1}$.
См. задачи \ref{alex-appl}, ср.~\cite[III, 2.1]{SE62},~\cite{ARS01}.
Р. Том получил условия на кольцо (ко)гомологий замкнутого $(m-1)$-многообразия, необходимые для вложимости в $S^m$~\cite{Th51}.
Петерсон изучал двойственность между {\it (ко)гомологическими операциями} пространств $N$ и $S^m-N$
и получил некоторые интересные теоремы невложимости.
Применения двойственности Александера приведены также в \S\ref{0emb}.

\begin{pr}\label{alex-appl}
(a) $\C P^2$ не вложимо в $S^5$.

(b) $\C P^2\#\C P^2$ не вложимо в $S^5$.

(c) $\R P^3$ (и трехмерное линзовое пространство $L(p,q)$) не вложимо в $S^4$.

(d) $\R P^3\# \R P^3$ не вложимо в $S^4$.

(e) $\R P^n$ не вложимо в $S^{n+1}$ для $n\ne4k+1$.

(f)* $\R P^5$ не вложимо в $S^6$.
\qquad
(f')* $\R P^n$ не вложимо в $S^{n+1}$.

(g)* $\R P^3\times\R P^3$ не вложимо в $S^7$.
\qquad
(g')* $(\R P^3)^r$ не вложимо в $S^{3r+1}$.
\end{pr}

\begin{pr}\label{alex-mvs} Пусть замкнутое связное $n$-многообразие $N$ лежит в $S^{n+1}$.
Обозначим через $A$ и $B$ замыкания компонент дополнения $S^{n+1}-N$.
Тогда $H_n(A)=H_n(B)=0$ и $i_A\oplus i_B:H_s(N)\to H_s(A)\oplus H_s(B)$ изоморфизм для любого $0<s<n$.
\end{pr}

\begin{pr}\label{alex-even} Пусть замкнутое связное $4k$-многообразие $N$ вложимо в $S^{4k+1}$.

(a) Эйлерова характеристика многообразия $N$ четна (т.е. $\rk H_{2k}(N;\Z)$ четен).

(b) Свободная часть группы $H_{2k}(N;\Z)$ разлагается в сумму двух прямых слагаемых, на каждом из которых
форма пересечений тривиальна (т.е. сигнатура $\sigma(N)=0$).
\end{pr}


Аналогичное необходимое условие можно получить и для вложимости ориентируемых $(4k+2)$-многообразий в $S^{4k+3}$.
Но оно автоматически выполнено ввиду кососимметричности формы пересечений.
Аналогичное необходимое условие можно получить и с коэффициентами $\Z_2$.
Но часть (a) этого условия равносильна целочисленной
версии
по формуле универсальных коэффициентов.

В задаче \ref{alex-even} условие вложимости $N$ в $ S^{4k+1}$ может быть заменено на более слабое:
достаточно того, чтобы $N$ было краем некоторого ориентируемого $(4k+1)$-многообразия (см. \S\ref{0cobord}).
Пример $N=\R P^3$ показывает, что следующее утверждение \ref{alex-odd} таким образом ослабить нельзя.

\begin{pr}\label{alex-odd} Пусть замкнутое связное $(2l+1)$-многообразие $N$ вложимо в $S^{2l+2}$.

(a) Подгруппа кручения группы $H_l(N;\Z)$ изоморфна $G\oplus G$ для некоторой группы $G$.

(b) Для некоторого разложения $H_l(N;\Z)\cong G\oplus G$ форма зацепления обращается в ноль
на каждом из двух прямых слагаемых $G$.
\end{pr}

В этом и следующем пунктах коэффициенты групп гомологий, если они не указаны, могут быть произвольными.
Полагаем $H_l(X)=0$ при $l<0$.

\smallskip
{\bf Теорема двойственности Александера.}
{\it Для любых $s$, $m>n+1$ и замкнутого гладкого ориентируемого $n$-подмногообразия $N\subset S^m$
определенное  ниже отображение
$$AD:H_{s+n-m+1}(N)\to H_s(S^m-N)$$
является  изоморфизмом. }

\smallskip
{\it Определение отображения $AD$ для $m>n+1$.}
Рассмотрим $(s+n-m+1)$-цикл $x\subset N$ (симплициальный цикл в некоторой гладкой триангуляции).
Концы нормальных к $N$ векторов длины $\varepsilon$, начала которых пробегают $x$, образуют подмножество в $S^m-N$, обозначаемое $AD_\varepsilon(x)$.
При достаточно малом $\varepsilon$ множество $AD_\varepsilon(x)$ является носителем $s$-цикла в $S^m-N$.
Значит, можно определить $AD(x):=[AD_\varepsilon(x)]$ и проверить корректность.

\smallskip
В теорему двойственности Александера обычно включают также случай $m=n+1$, равносильный
теореме Жордана вместе с утверждением задачи \ref{alex-mvs}.

 Обозначим через $ON$ множество концов нормальных векторов длины не более $\varepsilon$.
Положим $C_N:=S^m-\Int ON$.
(В этом тексте не определялись гомологии некомпакотных пространств, и если читатель не знает их определения, то может просто всюду заменить $S^m-N$ на $C_N$.)
Эти пространства $ON$ и $C_N$ являются $m$-мерными многообразиями с общим краем
(край образован концами векторов длины ровно $\varepsilon$).
Обозначим через $p:ON\to N$ отображение, ставящее в соответствие точке нормального вектора его начало
(этим отображение корректно определено при достаточно малых $\varepsilon$).

Отображение $AD$ является композицией
$H_{s+n-m+1}(N)\overset {p|_{\partial C_N}^!}\to H_s(\partial C_N)\overset{i_*}\to H_s(C_N)$
гомоморфизмов включения и прообраза.
Вообще говоря, ни один из них не является изоморфизмом.
\footnote{Изоморфизм $AD$ совпадает с композицией `обычной' двойственности Александера
\cite{Pr06} и изоморфизма Пуанкаре.
Действительно, $AD(x)$ является границей $(s+1)$-цепи $p^{-1}(x)$, пересечение которой в $S^m$
с любым $(n-1-(s+n-m+1))$-циклом из $N$ совпадает с пересечением в $N$ этого цикла и $x$.}

 Два доказательства теорема двойственности Александера намечены в двух следующих задачах.

\begin{pr} \label{}
Рассмотрим композицию гомоморфизмов `прообраза' (Тома), вырезания и граничного:
$$H_{s+n-m+1}(N)\overset{p^!}\to H_{s+1}(ON,\partial ON)\overset e\to H_{s+1}(S^m,C_N)\overset\partial\to H_s(C_N).$$
\ \quad
(a) Дайте определение отображения $p^!$.

(b) $AD=\partial ep^!$.

(c) $p^!$ и $\partial$ --- изоморфизмы.
\end{pr}

\begin{pr} \label{}
Для $s$-цикла $x$ в $C_N$ возьмем $(s+1)$-цепь $x'$  в $S^m$ общего положения относительно $N$, для которой $\partial x'=x$.
(Общность положения можно заменить на то, что $x'$ --- цепь в клеточном разбиении сферы $S^m$, двойственном к тому,
в котором $N$ является  подкомплексом.)
Тогда пересечение $x'\cap N$ будет $(s+1+n-m)$-циклом в $N$.
Положим $AD^{-1}[x]:=[x'\cap N]$.

(a) Это определение осмысленно, т.е. $\partial(x'\cap N)=0$.

(b) Это определение корректно, т.е. не зависит от выбора $a$ и $a'$.

(c) Построенное отображение $AD^{-1}$ действительно обратно к $AD$ и действительно является изоморфизмом.
\end{pr}

\begin{pr}\label{alex-gene}
(a) Для любого гладкого $n$-подмногообразия $N\subset S^m$ с краем
выполнено $H_s(N,\partial N)\cong H_{s+m-n-1}(S^m-N)$.

(b) Аналог теоремы двойственности Александера верен и для неориентируемых многообразий с коэффициентами $\Z_2$.
\end{pr}

\subsection{Другая двойственность Александера и двойственность Понтрягина}

Для подполиэдра $N\subset S^m$ обозначим через $ON$ его {\it регулярную}
окрестность в $S^m$.
(Регулярная окрестность определяется аналогично \S\ref{03man}.)
Положим $C_N:=S^m-\Int ON$.
Обозначим через $p:ON\to N$ ретракцию.
Эти определения обобщают определения, данные в предыдущем пункте.


\begin{pr}\label{alex-thm} Пусть  $N\subset S^m$ ---  подполиэдр.

(a) {\bf Теорема двойственности Александера.}
$H_s(N)\cong H_{s+1}(C_N,\partial C_N)$.
Более точно, изоморфизмом является композиция
$H_{s+1}(C_N,\partial C_N)\overset\partial\to H_s(\partial C_N)
\overset{p_*}\to H_s(N).$
\quad

(b) {\bf Следствие.} $H_s(N;\Z_2)\cong H_{m-s-1}(S^m-N;\Z_2)$;

свободные части групп $H_s(N;\Z)$ и $H_{m-s-1}(S^m-N;\Z)$ изоморфны;

кручения групп $H_s(N;\Z)$ и $H_{m-s-2}(S^m-N;\Z)$ изоморфны.

(c) Для многообразий $N$ изоморфизмы Александера
$p\partial$ и $ip|_{\partial C_N}^!$ 'не коммутируют': для
вложения $\xi:N\to\partial C_N$ вообще говоря,
$ip|_{\partial C_N}^!\xi^!\partial\ne i\xi p\partial$.

В этом смысле для многообразий имеется {\it две разные}
двойственности Александера $p\partial$ и $ip|_{\partial C_N}^!$.

(d) $H_s(\partial C_N)\cong H_s(N)\oplus H_s(S^m-N)$ при $0<s<m$.
\end{pr}

Теорема двойственности Александера из предыдущего пункта, в которой $N$ заменено на $X$, вытекает из теоремы двойственности Александера \ref{alex-thm}.a, если взять $N=C_X$ и применить изоморфизм Тома
$p^!:H_{s+1}(OX,\partial OX)\to H_{s+1+n-m}(X)$.

Задача \ref{alex-thm}.b является `простой частью' двойственности Александера-Понтрягина
(ср. с `простой частью' двойственности Пуанкаре).
`Сложная часть' была открыта Л. С. Понтрягиным около 1930 г. (это было одним из первых его научных достижений).

\smallskip
{\bf Теорема двойственности Понтрягина.}
{\it Для подполиэдра $N\subset S^m$

$\bullet$ форма зацеплений
$$\lk:H_s(N;\Z)\times H_{m-s-1}(S^m-N;\Z)\to\Z$$
невырождена и унимодулярна.

$\bullet$ форма вторичных зацеплений
$$\lk:Tors H_s(N;\Z)\times Tors H_{m-s-2}(S^m-N;\Z)\to\Q/\Z$$
невырождена.}

\smallskip
{\it Определение формы зацеплений.}
Для $(m-s-1)$-цикла $b$ в $S^m-N$ возьмем $(m-s)$-цепь $b'$ общего
положения в $S^m$, для которой $\partial b'=b$.
(Формально, пусть $N$ --- подкомплекс сферы $S^m$ в некоторой триангуляции,
тогда цикл $b'$ берем в двойственной триангуляции.)
Для $s$-цикла $a$ в $N$ положим $\lk([a],[b]):=a\cap b'$.

\smallskip
{\it Определение формы вторичных зацеплений.}
Для $(m-s-2)$-цикла $b$ в $S^m-N$ и некоторого $B\in\Z$ возьмем $(m-s)$-цепь $b'$ общего
положения в $S^m$, для которой $\partial b'=Bb$.
Аналогично, для $s$-цикла $a$ в $N$ и некоторого $A\in\Z$ возьмем $(s+1)$-цепь $a'$ общего
положения в $S^m$, для которой $\partial a'=Aa$.
Положим $\lk([a],[b]):=\{\dfrac{a'\cap b'}{AB}\}$.

\begin{pr} \label{}
(a) Эти определения корректны.
\quad

(b) Докажите
теорему двойственности Понтрягина.
\end{pr}

\subsection{Указания к некоторым задачам}

{\bf \ref{pdsim3z}.} (a) Возьмем такие 3-цепь $y$ и 2-цикл $z$, что $\partial y=kz$ для некоторого целого $k$.
Кратность (в цепи $y$) 3-симплекса, не входящего в $y$, равна нулю и поэтому делится на $k$.
Если кратность некоторого 3-симплекса делится на $k$, то кратность любого соседнего с ним 3-симплекса делится на $k$.
Поэтому $y=ky_1$.

\smallskip
{\bf \ref{pdsimz}.} Аналогично задачам \ref{pdsimz2} и \ref{pdsim3z}, ср.~\cite[\S69]{ST04}.

\bigskip
{\bf \ref{pdcom1}.} (a) 
По задаче \ref{pdrealc1m2} класс $\alpha$ представляется подмногообразием $A$.
Оно не разбивает многообразие $N$.
Действительно, иначе $A$ было бы границей каждой компоненты связности
дополнения $N-A$, а значит, было бы гомологично нулю.

Возьмем маленький отрезок, трансверсально пересекающий подмногообразие $A$ ровно в одной точке.
Так как $A$ не разбивает $N$, то концы этого отрезка можно соединить ломаной вне $A$.
Объединение этих ломаной и отрезка является 1-циклом, трансверсально
пересекающим $A$ ровно в одной точке.
Гомологический класс этого 1-цикла является искомым.

(b)  Следует из (a) и простой части двойственности Пуанкаре.


\smallskip
{\bf \ref{pdcomlk}.} (b) {\it Независимость от выбора цепи $B$.} Пусть $\partial B'=\partial B=nb$.
Тогда $a\cap B'-a\cap B=a\cap(B'-B)=0$, поскольку $\partial (B'-B)=0$ и $a$ имеет конечный порядок.

{\it Независимость от выбора цепи $b$} вытекает из независимости от выбора цепи $B$, поскольку
$\partial(B+nc)=n(b+\partial c)$.

{\it Независимость от выбора цикла  $a$.} Имеем
$$(a+\partial A)\cap B-a\cap B=\partial A\cap B=\pm A\cap\partial B=\pm n A\cap b.$$
\ \quad
(c) Ответ: $\lk(a,a)=1/2$ для образующей $a\in H_1(\R P^3)$.

(d) (При написании этого решения использованы~\cite[\S69, \S71, Предложение 2]{ST04} и текст С. Аввакумова.)
 Выберем на многообразии некоторую триангуляцию вместе с двойственной к ней. Для каждого $k$ выберем естественный базис в группе $k$-мерных цепей триангуляции и двойственный ему естественный базис в группе $(n-k)$-мерных цепей двойственной триангуляции. Через $d_k$ обозначим граничный оператор, действующий из группы $(k+1)$-мерных цепей в группу $k$-мерных цепей. Через $d^*_k$ обозначим граничный оператор, действующий из группы $(k+1)$-мерных цепей двойственной триангуляции в группу $k$-мерных цепей двойственной триангуляции.

Тогда $d_k d_{k+1}=0$ и $d^*_k d^*_{k+1}=0$.
Кроме того, выбранные базисы {\it дуальны}, т.е. $d^*_{n-k} = d^T_{k-1}$ и $a\cap b^* = \delta_{ab}$ для любых элемента $a$ базиса $k$-мерных цепей и элемента $b^*$ базиса $(n-k)$-мерных двойственных цепей.

Пусть $a$ и $a^*$, \ $b$ и $b^*$ --- двойственные элементы дуальных базисов.
При каждом из следующих двух преобразований дуальность базисов сохраняется:

$a':=-a$, \ $(a^*)':=-a^*$.

$a':=a+b$, $(b^*)'=b^*-a^*$.

С помощью этих двух преобразований можно добиться того, чтобы для каждого $k$ матрицы операторов $d_{k-1}$ имели нормальную форму (рис 1).
???Здесь через $C_k$ обозначен базис цепей, через $Z_k$ циклов и через $B_k$ элементы конечного (возможно, равного 1!) порядка.
Ввиду $d^*_{n-k} = d^T_{k-1}$ матрицы операторов $d^*_{n-k}$ также примут нормальную форму
(рис 2; ???$\dim C_k=\dim B^*_{n-k}$).
Ввиду $a\cap b = \delta_{ab}$ матрица формы пересечений примет следующий вид (рис 3).

Пусть теперь $b\in B_{k-1}$ имеет порядок $m>1$ (? $m>1$ не используется?).
Тогда $mb=d_{k-1}c$ для некоторого $c\in C_k$.
Тогда $c\cap c^*=\delta_{cc}=1$.
Значит, форма зацепления невырождена.

\smallskip
{\bf \ref{alex-appl}.}
(a) Примените \ref{alex-even}.a.

(b) Примените \ref{alex-even}.b.

Две формы на $H_2(\C P^2\#\C P^2)$ со значениями в $\Z$, имеющие матрицы
$\left(\begin{matrix}1&0\\0&1\end{matrix}\right)$ и
$\left(\begin{matrix}0&1\\1&0\end{matrix}\right)$ в некоторых базисах, не изоморфны,
поскольку для первой формы имеются нечетные квадраты, а для второй --- нет.

(c,e) Примените \ref{alex-odd}.a.

(d) Примените \ref{alex-odd}.b.
Две формы на $H_1(\R P^3\#\R P^3)$ со значениями в $\Z_2$, имеющие матрицы
$\left(\begin{matrix}1&0\\0&1\end{matrix}\right)$ и
$\left(\begin{matrix}0&1\\1&0\end{matrix}\right)$ в некоторых базисах, не изоморфны,
поскольку для первой формы имеются ненулевые квадраты, а для второй --- нет.

(f,f') Применим форму пересечений.
 Suppose to the contrary that $\R P^n\subset S^{n+1}$ is an embedding.
Let $A$ and $B$ be the closures of the connected components of $S^{n+1}-\R P^n$
and let $i_A$, $i_B$ be the inclusion-induced homomorphisms.
We omit the coefficients $\Z_2$.
Using the Mayer-Vietoris sequence for $S^{n+1}=A\cup B$, one sees that
$\partial_A+\partial_B:H_n(A,\partial)\oplus H_n(B,\partial)\to H_{n-1}(\R P^n)$ is an isomorphism.
The\-re\-fo\-re by relabelling, if necessary, we can assume that there is an element
$a\in H_n(A,\partial)$ such that $\partial_Aa=[\R P^{n-1}]$.
So, $\partial_Aa^n=[\R P^{n-1}]^n\ne0$.
But from the above Mayer-Vietoris sequence it follows that $H_1(A,\partial)=0$, which is a contradiction.

(g,g') Применим форму пересечений. См. \cite{ARS01}.
Let $N=(\R P^3)^r$.
Suppose to the contrary that $N\subset S^{3r+1}$ is an embedding.
Let $A$ and $B$ be the closures of the connected components of $S^{3r+1}-N$
and let $i_A$, $i_B$ be the inclusion-induced homomorphisms.
We omit the coefficients $\Z_2$.
Using the Mayer-Vietoris sequence for $S^{3r+1}=A\cup B$, one sees that
$\partial_A+\partial_B:H_{2r+1}(A,\partial)\oplus H_{2r+1}(B,\partial)\to H_{2r}(N)$ is an isomorphism.
We have $H_*(N)=\left<x_1,\dots,x_r\ |\ x_i^4=0\right>$, where $x_i$ is represented by the Cartesian product of $\R P^2$ at the $i$-th place and $\R P^3$ at other places.
The\-re\-fo\-re by relabelling, if necessary, we can assume that there is an element
$a\in H_{2r+1}(A,\partial)$ such that $\partial_Aa=x_1\dots x_r+\dots$, where dots
denote summands containing $x_l^2$ for some $l$.
So, $\partial_Aa^2=(x_1\dots x_r)^2$ and
$\partial_Aa^3=(x_1\dots x_r)^3\ne0$.
But from the above Mayer-Vietoris sequence it follows that $H_1(A,\partial)=0$, which is a contradiction.

\smallskip
{\bf \ref{alex-mvs}.} Примените последовательность Майера-Виеториса.


\smallskip
{\bf \ref{alex-even}.} Считаем, что $N\subset S^{4k+1}$.
Используем обозначения из задачи \ref{alex-mvs} и ее результат для $s=l:=2k$.
Теперь утверждение (a) вытекает из
$$Free H_l(A)\cong Free H_{l+1}(A,\partial)\quad\text{и}
\quad H_{l+1}(A,\partial)\cong H_{l+1}(S^{2l+1},B)\cong H_l(B).$$
Докажем (b).
Если $x,y\in Free H_l(N)$ и $i_Ax=i_Ay=0$, то $x=\partial_AX$ и $y=\partial_AY$
для некоторых $X,Y\in H_{l+1}(A,\partial)$.
Тогда
$$x\cap_N y=\partial_AX\cap_N\partial_AY=X\cap_A i_A\partial_AY=X\cap0=0.$$
(Или $x\cap_N y=\partial_A(X\cap_AY)=0$, поскольку $X\cap_AY\in H_1(A,\partial)=0$.)

\smallskip
{\bf \ref{alex-odd}.} Считаем, что $N\subset S^{2l+2}$.
Используем обозначения из задачи \ref{alex-mvs} и ее результат для $s=l$.
Теперь утверждение (a) вытекает из
$$\Tors H_l(A)\cong \Tors H_{l+1}(A,\partial)\quad\text{и}
\quad H_{l+1}(A,\partial)\cong H_{l+1}(S^{2l+2},B)\cong H_l(B).$$
Докажем (b).
Если $x,y\in\Tors H_l(N)$ и $i_Ax=i_Ay=0$, то

$\bullet$ $x=\partial X$ для некоторой $(l+1)$-цепи $X$ в $A$,

$\bullet$ $y=\partial Y$ для некоторой $(l+1)$-цепи $Y$ в $A$, и

$\bullet$ $nx=\partial X'$ для некоторой $(l+1)$-цепи $X'$ в $N$.

Тогда
$$\lk(x,y)=\frac1n X'\cap_N y=\frac1n X''\cap_A Y=X\cap_A Y=0\in\Q/\Z,$$
где $X''$ получено из $X'$ сдвигом внутренности внутрь $A$ и последнее равенство верно,
поскольку $\partial(nX)=nx=\partial X''$.
Ср. \cite[Exercise 4.5.12.d]{GS99}.

\smallskip
{\bf \ref{alex-thm}.}
(a) Теорема следует из того, что композиция $p_*\partial$ обратна следующей:
$$H_s(N)\cong H_s(ON)\cong H_{s+1}(S^m,ON)\cong H_{s+1}(C_N,\partial C_N).$$
Здесь первый изоморфизм получается ввиду гомотопической инвариантности
гомологий, второй из точной последовательности пары $(S^m,ON)$, а
третий из теоремы вырезания.

(b) Примените двойственность Лефшеца к $m$-многообразию $C_N$.

(c) Рассмотрите, например, стандартную окружность $N\subset\R^3$.

(d) Примените последовательность Майера-Виеториса (или задачу \ref{alex-mvs}).

 \newpage
\section{Препятствия к кобордантности}\label{0cobord}

\hfill{\it He looked at the snail. Can it see me? he wondered. }

\hfill{\it Then he felt, how little I know, and how little it is possible to know;}

\hfill{\it  and with this thought he experienced a moment of joy.}

\hfill{\it I. Murdoch, The Flight From The Enchanter.
\footnote{\it Он посмотрел на улитку. 'Может ли она видеть меня?' задумался он.
Тогда он почувствовал, как мало он знает, и как мало возможно знать;
и эта мысль принесла ему радость. А. Мердок,
Бегство от волшебника, пер. автора.}
}

\subsection{Введение}\label{0cobint}

При помощи характеристических классов были получены блестящие достижения теории классификации многообразий.
Классификация многообразий с точностью до {\it кобордизма} (определение см. в \S\ref{0cobeul}) была начата Львом
Семеновичем Понтрягиным в 1930-х в качестве первого шага к важным и трудным проблемам классификации
{\it многообразий с точностью до гомеоморфизма}, а также классификации {\it оснащенных многообразий с точностью до
оснащенного кобордизма} (определение см. в \S\ref{0hopfs3s2}).
При помощи кобордизмов ученик Понтрягина Владимир Абрамович Рохлин доказал следующий результат.

\smallskip
 {\bf Теорема Рохлина о сигнатуре.} {\it Сигнатура гладкого замкнутого почти
параллелизуемого 4-мерного многообразия делится на 16.}

\smallskip
В 1942 Понтрягин придумал новые характеристи\-ческие клас\-сы.
Рохлин показал, что эти классы естественно появляются при обобщении задачи Штифеля:
при построении набора касательных векторных полей, имеющей в каждой точке ранг, не меньший заданного.
При дальнейшем обобщении появляются другие характеристи\-ческие клас\-сы, но они уже `выражаются' через
классы Штифеля-Уитни и Понтрягина.
Поэтому среди этих классов выделены {\it классы Понтрягина}, образующие `максимальную независимую систему'.

Последующие замечательные результаты Рене Тома о классификации многообразий с точностью до кобордизма
были удостоены Филдсовской медали 1952 г.
С помощью них была получена знаменитая формула Фридриха Хирцебруха, связывающая сигнатуру формы пересечений с классами
Понтрягина.

\smallskip
{\bf Следствие формулы Хирцебруха о сигнатуре.} {\it Сигнатура гладкого замкнутого почти параллелизуемого 8-мерного
многообразия делится на 224.}

\smallskip
Эта формула сделала возможной следующее открытие.

\smallskip
{\bf Пример сферы Милнора.} {\it Существует замкнутое 7-многообразие,
гомотопически эквивалентное (и даже гомеоморфное) сфере $S^7$, но не диффеоморфное ей.}

\smallskip
Многообразия $M$ и $N$ называются {\it гомотопически эквивалентными}, если
существуют непрерывные отображения $f:M\to N$, $g:N\to M$ такие, что
$f\circ g$ гомотопно $\id_N$  и $g\circ f$ гомотопно $\id_M$.
См. построение в \S\ref{0hosp-con}.

Классификация гомотопических сфер стала основой теории хирургии, полезной для
классификации многообразий с точностью до {\it диффеоморфизма}.
Все эти теории применяются для изучения {\it гладких} многообразий и их отличий от {\it топологических}
многообразий~\cite[\S28]{DNF84}.

Мы пропускаем {\it целые }  коэффициенты из обозначений групп гомологий.
Мы работаем в гладкой категории, если не оговорено противное.
Слово `гладкое' пропускается.
Напомним, что все многообразия в книге предполагаются компактными.
Под {\it многообразием} можно понимать подмногообразие евклидова пространства
Определения $n$-мерных {\it гладких} подмногообразий, касательных и нормальных векторных полей и полей направлений
на них, а также гомотопности векторных полей, аналогичны случаю 2-многообразий (\S\ref{0ve2ma}).

\subsection{Эйлерова характеристика}\label{0cobeul}


Для
решения большей части задач этого пункта достаточно начальных сведений об эйлеровой характеристике
(задача \ref{veheul}, см. также \S\S \ref{02maneul}, \ref{03maneul}).

\begin{pr} \label{cob12}
(a) Если $A$ и $B$ --- замкнутые многообразия и $A=\partial M$, то $A\times B=\partial(M\times B)$.

(b) Любое одномерное или двумерное замкнутое ориентируемое многообразие является краем некоторого многообразия.

(c) Любое двумерное замкнутое неориентируемое многообразие четной эйлеровой
характеристики является краем некоторого многообразия.
\end{pr}

{\bf Теорема.}
{\it Любое замкнутое 3-многообразие является краем некоторого многообразия.}
(Для ориентируемых~\cite[VII, Theorem 2]{Ki89}, для произвольных~\cite[Chapter 6]{St68}.)

\begin{pr} \label{cobeul}
(a) $\R P^2$ не является краем многообразия.

(b) Замкнутое 2-многообразие является краем многообразия тогда и только тогда,
когда его эйлерова характеристика четна.

(c) {\bf Теорема.} {\it Если замкнутое многообразие является краем многообразия, то его эйлерова характеристика четна.}

(Эта теорема интересна только для четномерных многообразий.)

(d) Если замкнутое $2k$-многообразие $N$ является краем многообразия, то $\rk H_k(N)$ четен.

(e) {\it Аддитивность.} $\chi(M\cup N)=\chi(M)+\chi(N)-\chi(M\cap N)$.

(f) {\it Мультипликативность.}  $\chi(M\times N)=\chi(M)\chi(N)$.
\end{pr}

\begin{pr} \label{cobpro}
(a) $\R P^{2k+1}$ является краем некоторого многообразия.

(b) $\R P^n$ является краем многообразия тогда и только тогда,
когда $n$ нечетно (\cite[Упражнение 29]{FF89} можно усилить).

(c) $\C P^{2k}$ не является краем многообразия.

(d) $\C P^{2k+1}$ является краем  многообразия (даже ориентируемого).

(e) $\C P^n$ является краем
многообразия тогда и только тогда, когда $n$ нечетно
(\cite[Упражнение 29]{FF89} можно усилить).

(f) $\R P^2\times\R P^2$ не является краем многообразия.

(g) При каком условии $\R P^{n_1}\times\dots\times\R P^{n_k}$ является краем многообразия?
\end{pr}

Замкнутые многообразия $N_1$ и $N_2$ называются {\it кобордантными}, если существует многообразие ({\it кобордизм})
с границей $N_1\sqcup N_2$.
При этом одно из многообразий может быть пустым.

\subsection{Сигнатура}\label{0cobsig}

Для
решения большей части задач этого пункта достаточно владения основами теории гомологий и пересечений
(задача \ref{veheul}, \S\S \ref{0vesyhom}, \ref{0vesyint}, \ref{0homolpa}, см. также \S\S \ref{02maneul},
\ref{02homhom}, \ref{02homint}, \ref{03maneul}, \ref{03homhom}).

{\it Ориентированным} многообразием называется ориентируемое многообразие с фиксированной ориентацией.
Для ориентированного многообразия $N$ через $-N$ обозначается ориентированное многообразие, полученное из $N$
изменением ориентации.

\begin{pr} \label{cobori12}
(a) Край ориентируемого многообразия замкнут и ориентируем.

(b) Любое одномерное или двумерное замкнутое ориентируемое
многообразие является краем некоторого ориентируемого многообразия.

(с) Любое одномерное или двумерное замкнутое ориентированное многообразие
является ориентированным краем краем некоторого ориентированного многообразия.

(d) Для любого замкнутого ориентированного многообразия $N$ многообразие $N\sqcup(-N)$
 является ориентированным краем некоторого ориентированного многообразия.
\end{pr}

{\bf Теорема.}
{\it Любое замкнутое ориентированное 3-многообразие является ориентированным краем некоторого ориентированного
многообразия.}

\begin{pr} \label{cobcp2}
Ориентированное (произвольно) многообразие $\C P^2\sqcup\C P^2$
не является ориентированным краем ориентированного многообразия.
(Иными словами, многообразия $\C P^2$ и $-\C P^2$ не являются ориентированно кобордантными.)

Указание: если не получается, то см. следующие задачи.
\end{pr}

\begin{pr} \label{cobinbdev}
Для ориентированного $2k$-мерного многообразия $M$ с краем, гомоморфизма включения $i:H_k(\partial M)\to H_k(M)$ и
формы пересечений $\cap:H_k(M)\times H_k(M)\to\Z$ (\S\S\ref{02homint},\ref{0vesyint})

(a) $ix\cap ix=0$ при любом $x\in H_k(\partial M)$;

(b) $\im i=H_k(M)^\perp$.
\end{pr}

 \begin{pr} \label{cobinbdod} (ср. \cite[теорема 8.16]{Pr06})
Для ориентированного $(2k+1)$-многообразия $M$ с краем, гомоморфизмов
$H_{k+1}(M,\partial)\overset\partial\to H_k(\partial M)\overset i\to H_k(M)$ и
формы пересечений $\cap:H_k(\partial M)\times H_k(\partial M)\to\Z$ (\S\S\ref{02homint},\ref{0vesyint})

(a) $x\cap\partial y=ix\cap y$ для любых $x\in H_k(\partial M)$ и $y\in H_{k+1}(M,\partial)$,
где в правой части подразумевается билинейное пересечение $\cap:H_k(M)\times H_{k+1}(M,\partial)\to\Z$

(b) $\partial y\cap\partial y=0$ при любом $y\in H_{k+1}(M,\partial)$.

(с) $\im\partial=(\im\partial)^\perp$.

(d) $2\rk\im\partial=\rk H_k(\partial M)$.
 \end{pr}

\begin{pr} \label{cobsign}
(a) {\bf Теорема Понтрягина.}
{\it Если замкнутое ориентированное $4k$-многообразие $N$ является ориентированным краем ориентируемого
многообразия, то сигнатура $\sigma(N)$ (формы пересечений на $H_{2k}(N;\Z)$) равна нулю.}
(Или, эквивалентно, если замкнутые ориентированные $4k$-многообразия $N_1$ и
$N_2$ ориентированно кобордантны, то $\sigma(N_1)=\sigma(N_2)$.)

(b) Если для замкнутого ориентированного $4k$-многообразия $N$ существует сохраняющий ориентацию
диффеоморфизм $N\to-N$, то $\sigma(N)=0$.
 \end{pr}

\begin{pr} \label{cobsignpro}
(a) {\it Аддитивность.} $\sigma(M\sqcup N)=\sigma(M)+\sigma(N)$.

(b) {\it Мультипликативность.}  $\sigma(M\times N)=\sigma(M)\sigma(N)$.

(c)* {\it Аддитивность Новикова-Рохлина.}
$\sigma(M\bigcup\limits_{\partial M=\partial N}N)=\sigma(M)+\sigma(N)$ \cite{Pr06}.
\end{pr}

Замкнутые ориентированные многообразия $N_1$ и $N_2$ {\it ориентированно кобордантны}, если существует ориентированное
многообразие ({\it ориентированный кобордизм}) с ориентированным краем $N_1\sqcup(-N_2)$.
При этом одно из многообразий может быть пустым.

\begin{pr} \label{cobord} (a) Если многообразия гомеоморфны, то они кобордантны.
Обратное неверно.

(b) Если непустые многообразия кобордантны, то их размерности равны.
Размерность кобордизма на 1 больше размерности того из многообразий $M_1,M_2$, которое непусто.

(c) Многообразия кобордантны тогда и только тогда, когда их несвязное
объединение является краем некоторого многообразия.

(d) Существуют непустые некобордантные многообразия одинаковой размерности.
 \end{pr}

\begin{pr}\label{cobconn}
(a) $M\sqcup N$ кобордантно $M\#N$.

(b) Любое многообразие кобордантно связному (пустое считается связным).
\end{pr}

 \subsection{Числа Штифеля-Уитни}

Необходимые здесь сведения о классах Штифеля-Уитни содержатся в \S\ref{0vesy3ch} и \S\ref{0vesych},
см. также \S\ref{02homcy}-\S\ref{02homhom} и \S\ref{0vesy4}.

\begin{pr} \label{cobpont}
(a) Не существует многообразия с границей $N:=\R P^2\times\R P^2\sqcup\R P^4$.
Иными словами, многообразия $\R P^2\times\R P^2$ и $\R P^4$ не кобордантны.

(b) Не существует многообразия, границей которого является непустое объединение
одного, двух или трех различных многообразий из $N_1:=\R P^6$,
$N_2=\R P^4\times\R P^2$ и $N_3:=\R P^2\times\R P^2\times\R P^2$.

(c) Ни для какого $n$ не существует многообразия, границей которого является непустое объединение
некоторого числа различных многообразий из множества
$$X_n:\ =\ \{\R P^{2i_1}\times\dots\times\R P^{2i_k}\ :\ i_1,\dots,i_k\in\Z_{>0},\ n=i_1+\dots+i_k\}.$$
(Указание: докажите и используйте следующую задачу.)
\end{pr}

{\it Числом Штифеля-Уитни} замкнутого $n$-многообразия $N$, отвечающим разбиению $n=i_1+\dots+i_s$ на целые
положительные слагаемые, называется число
$$w_{i_1,\dots,i_s}(N):=\Sigma w_{i_1}(N)\cap\dots\cap w_{i_s}(N).$$
Здесь $\Sigma$ означает суммирование компонент вектора $w_{i_1}(N)\cap\dots\cap w_{i_s}(N)$, которые
нумеруются компонентами связности многообразия $N$.

\begin{pr} \label{cobponth}
(a) $w_k(\partial M)=\partial w_k(M)$.

(b) {\bf `Неориентируемая' теорема Понтрягина.} {\it Если замкнутое $n$-многообразие $N$ является краем многообразия,
то $w_{i_1,\dots,i_s}(N)=0$ для любого разбиения $n=i_1+\dots+i_s$.

Иными словами, если замкнутые $n$-многообразия $N_1$ и $N_2$ кобордантны, то
$w_{i_1,\dots,i_s}(N_1)=w_{i_1,\dots,i_s}(N_2)$ для любого разбиения $n=i_1+\dots+i_s$.}
\end{pr}

Обратное утверждение тоже верно.
Оно принадлежит Тому и доказывается гораздо труднее~\cite{St68}.

На множестве $\Omega^O$ классов кобордизма замкнутых многообразий (всевозможных размерностей) имеются операции
несвязного объединения и декартова
произведения.

\begin{pr} \label{cobomo}
 (a) Эти операции действительно корректно определены.

(b) Множестве $\Omega^O$ с этими операциями является градуированным
коммутативным кольцом с единицей и соотношением $a+a=0$.

(с)* Для целых положительных чисел рассмотрим гиперповерхность
$H\subset\R P^4\times \R P^2$ степени $(1,1)$.
Тогда $H$ является 5-многообразием и не является краем никакого многообразия.

(d)* Для целых положительных чисел рассмотрим гиперповерхность
$H_{p,q}\subset\R P^{2^{p+1}\times q}\times \R P^{2^p}$ степени $(1,1)$.
Тогда $H_{p,q}$ является многообразием размерности $2^p(2q+1)-1$ и не является
краем никакого многообразия.
 \end{pr}

 {\bf Классификационная теорема Тома.}
{\it Градуированное кольцо $\Omega^O$ изоморфно
градуированному кольцу полиномов над $\Z_2$ от образующих $x_i$ размерности
$i$, где $i$ пробегает все положительные целые числа, не равные $2^s-1$. }

\smallskip
В качестве $x_{2i}$ можно взять класс кобордизма многообразия $\R P^{2i}$.
Нечетномерные образующие были построены, например, Дольдом (ориентируемые) и Милнором, см. задачу \ref{cobomo}.cd.

\subsection{Числа Понтрягина и формула Хирцебруха}

Необходимые здесь сведения о классах Понтрягина содержатся в задачах \ref{vespont}, см. также задачу \ref{ves4int}.с.

\begin{pr} \label{cobcp4}
$\C P^2\times\C P^2\sqcup(-\C P^4)$ не является ориентированным краем ориентированного многообразия.
Иными словами, многообразия $\C P^2\times\C P^2$ и $\C P^4$ не являются ориентированно кобордантными.
(Указание: докажите и используйте следующую задачу.)
 \end{pr}

{\it Числом Понтрягина} замкнутого ориентированного $4k$-многообразия $N$, отвечающим разбиению $4k=i_1+\dots+i_s$
на целые положительные слагаемые, называется число $p_{i_1,\dots,i_s}(N):=\Sigma p_{i_1}(N)\cap\dots\cap p_{i_s}(N)$.
Здесь $\Sigma$ означает суммирование компонент вектора $p_{i_1}(N)\cap\dots\cap p_{i_s}(N)$
(которые нумеруются компонентами связности многообразия $N$).

\begin{pr} \label{cobpontor}
(a) $p_k(\partial M)=\partial p_k(M)$.

(b) {\bf `Ориентированная' теорема Понтрягина.}
{\it Если замкнутое ориентированное $4k$-многообразие $N$ является ориентированным краем ориентируемого многообразия,
то $p_{i_1,\dots,i_s}(N)=0$ для любого разбиения $4k=i_1+\dots+i_s$.

Или, эквивалентно, если замкнутые ориентированные $4k$-многообразия $N_1$ и $N_2$ ориентированно кобордантны, то
$p_{i_1,\dots,i_s}(N_1)=p_{i_1,\dots,i_s}(N_2)$ для любого разбиения $4k=i_1+\dots+i_s$.}

(c) Если для замкнутого ориентированного $4k$-многообразия $N$ существует сохраняющий ориентацию диффеоморфизм $N\to-N$,
то $p_{i_1,\dots,i_s}(N)=0$ для любого разбиения $4k=i_1+\dots+i_s$.
 \end{pr}


{\bf Теорема Тома.} {\it Если все числа Понтрягина замкнутого ориентированного многообразия нулевые
(или если его размерность не делится на 4), то объединение нескольких экземпляров многообразия $X$
(которые все берутся с одинаковой ориентацией) является ориентированным
краем некоторого замкнутого ориентированного многообразия.}

\smallskip
По поводу описания кольца $\Omega^{SO}$ классов ориентированной кобордантности
ориентированных многообразий см., например,~\cite[\S19.6.D, \S43.2]{FF89}.

При доказательстве следующих результатов можно использовать теорему Тома.

\begin{pr} \label{cobhirz}
(a) Для любого замкнутого ориентируемого 4-многообразия выполнено $p_1=3\sigma$.

(b) Для любого замкнутого ориентируемого 8-многообразия выполнено $7p_2-p_{1,1}=45\sigma$.

(c) {\bf Теорема Хирцебруха о сигнатуре (упрощенная).}
{\it Для любого $k$ существуют такие рациональные числа $a_{i_1,\dots,i_r}$, отвечающие разбиениям
$4k=i_1+\dots+i_s$ на целые положительные слагаемые, что для любого замкнутого ориентируемого $4k$-многообразия
$M$ выполнено $\sigma(M)=\sum\limits_{i_1+\dots+i_r=4k}a_{i_1,\dots,i_r}p_{i_1,\dots,i_r}(M)$.}
\end{pr}

\subsection{Указания и решения к некоторым задачам}

 {\bf \ref{cobeul}.}
(a) {\it Первое решение.}
Пусть, напротив, $M$ --- 3-многообразие и $\partial M\cong\R P^2$.
Обозначим через $M'$ копию многообразия $M$.
Тогда $0=\chi(M\cup_{\R P^2} M')=\chi(M)+\chi(M')-\chi(\R P^2)$, откуда
$\chi(\R P^2)$ четно.
Противоречие.

{\it Второе решение.}
Пусть, напротив, $M$ --- 3-многообразие и $\partial M\cong\R P^2$.
Тогда
$$\tau(M)|_{\partial M}\cong\tau(\partial M)\oplus\nu(\partial M\subset M)
\cong\tau(\partial M)\oplus\varepsilon.$$
Поэтому $0=e(M)\cap\partial M=e(\partial M)\ne0$, где $e$ --- класс Эйлера по
модулю 2 и равенства означают сравнения по модулю 2.

{\it Третье решение.}
Пусть, напротив, $M$ --- 3-многообразие и $\partial M\cong\R P^2$.
Тогда $w_1(\partial M)=\partial w_1(M)$, где
$\partial:H_2(M,\partial M)\to H_1(M)$ --- граничный гомоморфизм.
Пусть $(\omega,\partial\omega)\subset(M,\partial M)$ --- 2-подмногообразие c
краем, реализующее класс $w_1(M)$.
Тогда $\partial\omega\subset\partial M$ реализует класс $w_1(\partial M)$.
Значит, $\omega\cap_M\omega$ есть 1-подмногообразие с краем
$\partial\omega\cap_{\partial M}\partial\omega$.
Поэтому $w_1(\partial M)\cap w_1(\partial M)=0$, что неверно для $\partial M=\R P^2$.

\smallskip
{\bf \ref{cobpro}.}  (d) Постройте и используйте расслоение $\C P^{2k+1}\to H P^k$  со слоем $S^2$.
Или постройте и используйте инволюцию на $\C P^{2k+1}$ без неподвижных точек.

(f) $\chi(\R P^2\times\R P^2)=\chi(\R P^2)^2\equiv1\mod2$.

\smallskip
{\bf \ref{cobinbdev}.}  (b) Используйте теорему двойственности Пуанкаре (сложную часть).

\smallskip
{\bf \ref{cobinbdod}.}
(c) Используйте (a,b) и двойственность Пуанкаре.

Если $x\in H_k(\partial M)$ и $x\cap\im\partial=0$, то
$ix\cap y=x\cap\partial y=0$ для любого $y\in H_{k+1}(M,\partial)$.
Значит, по двойственности Пуанкаре $ix=0$, т.е. $x\in\im\partial$.

(d) Используйте (c) и двойственность Пуанкаре.

\smallskip
{\bf \ref{cobsign}.} (a) Следует из задачи \ref{cobinbdod}.bd.

(b) Если $N\cong-N$, то $N\sqcup N$ кобордантно пустому.

\smallskip
{\bf \ref{cobpont}.} (a) Аналогично третьему доказательству того, что $\R P^2$ не ограничивает.
По задаче \ref{vesnww}.ab $w_1(N)=[\R P^2\times\R P^1+\R P^1\times\R P^2\sqcup \R P^3]$.
Тогда
$$w_1(N)^2=[\R P^1\times\R P^1+\R P^1\times\R P^1\sqcup \R P^2]=[\emptyset\sqcup \R P^2]
\quad\text{и}\quad w_1(N)^4=[\emptyset\sqcup \R P^0]\ne0.$$
\quad
(b) Аналогично (a). Используйте $w_2(\R P^n)={n+1\choose 2}[\R P^{n-2}]$ (задача \ref{vesncalc}).

\smallskip
{\bf \ref{cobponth}.} (b) Следует из (a).

\smallskip
{\bf \ref{cobcp4}.} Используйте задачи \ref{vespont}.bc и \ref{cobpontor}.b.

\smallskip
{\bf \ref{cobpontor}.} (b) Следует из (a).

(c) Если $N\cong-N$, то $N\sqcup N$ кобордантно пустому.

\smallskip
{\bf \ref{cobhirz}.}
(a,b) По теореме Тома достаточно проверить (a) для $M=\C P^2$, а (b) для $M=\C P^4$ и $M=\C P^2\times \C P^2$.
Используйте задачи \ref{vespont}.bc.


\newpage
 \section{Гомотопическая классификация отображений }\label{0hopf}

\hfill{\it Увы! Что б ни сказал потомок просвещенный,}

\hfill{\it все так же на ветру, в одежде оживленной,}

\hfill{\it к своим же Истина склоняется перстам,}

\hfill{\it с улыбкой женскою и детскою заботой,}

\hfill{\it как будто в пригоршне рассматривая что-то,}

\hfill{\it из-за ее плеча невидимое нам.}

\hfill{\it В. Набоков, Дар.
\footnote{\it Alas! In vain historians pry and probe: The same wind blows, and
in the same live robe Truth bends her head to fingers curved cupwise; And with
a woman's smile and a child's care Examines something she is holding there
Concealed by her own shoulders from our eyes. V. Nabokov, Gift.}
}

 \subsection{Определения и исторические замечания}

Непрерывные отображения $X\to Y$ определены в начале \S\ref{0vepl}.
Их гомотопность определяется дословно так же.

Проблема классификации непрерывных отображений с точностью до
гомотопности  --- одна из важнейших в  топологии, поскольку к ней сводятся многие другие
задачи топологии и ее приложений.
Эта проблема изучается в этом и следующем параграфах.
Основная теорема топологии о гомотопической классификации отображений $S^n\to S^n$ была
доказана Хайнцем Хопфом в 1926 г.
Он использовал инвариант 'степень', введенный Лейтзеном Эгбертом Яном Брауэром в 1911 г.
В 1932 г. Хопф обобщил этот результат до классификации отображений $n$-полиэдра в $S^n$
('по заказу' Павла Сергеевича Александрова).
Приводимые здесь формулировка и доказательство теоремы Хопфа принадлежит Хасслеру Уитни (1937).
Дальнейшее развитие теорема Хопфа-Уитни получила в работах Сэмюэля Эйленберга и Сондерса Маклейна (1940),
Льва Семеновича Понтрягина (1941), Нормана Стинрода (1947), Джона Генри Константина Уайтхеда
(1949) и Михаила Михайловича Постникова (1950).

Напомним, что все отображения считаются непрерывными и прилагательное `непрерывное' опускается.
Обозначим через $[X,Y]$ множество непрерывных отображений $X\to Y$, переводящих отмеченную точку в отмеченную точку,
с точностью до гомотопности в классе таких отображений.
Слово `непрерывный' опускается.

 \subsection{Групповая структура}

Обозначим
$$S^n:=\{(x_1,\dots,x_{n+1})\in\R^{n+1}\ :\ x_1^2+\dots+x_{n+1}^2=1\},\quad \pi_n(X):=[S^n,X]
\quad\text{и}\quad \pi^n(X):=[X;S^n].$$

\begin{figure}[h]\centering
\includegraphics{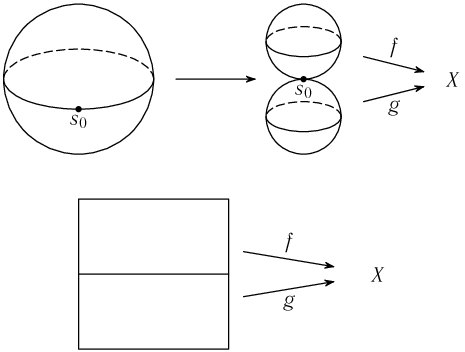}
\caption{Сумма сфероидов }\label{summa}
\end{figure}

Пусть $X$ --- произвольное многообразие (или даже тело симплициального комплекса).

{\it Суммой} гомотопических классов отображений $f,g:S^2\to X$ называется гомотопический класс композиции
$S^2\overset c\to S^2\vee S^2\overset{f\vee g}\to X$ стягивания $c$ экватора в точку и отображения $f\vee g$
(рис.~\ref{summa} вверху).
Если под отображением $S^2\to X$ понимать отображение $I^2\to X$, переводящее границу куба $I^2$ в точку, то сумма будет выглядеть как на рис.~\ref{summa} внизу.

{\it Нулевым элементом} множества $\pi_2(X)$ называется гомотопический класс отображения в точку.
{\it Обратным элементом} к гомотопическому классу отображения $f:S^2\to X$ называется гомотопический класс
композиции $s\circ f$ отображения $f$ с зеркальной симметрией $s:S^2\to S^2$.
Нетрудно проверить, что эти определения превращают множество $\pi_2(X)$ в абелеву группу.

\begin{pr}\label{hopf-cosu}
(a) {\it Суммой} гомотопических классов отображений $f:S^1\times S^1\to X$ и $g:S^2\to X$ называется
гомотопический класс композиции $S^1\times S^1\overset c\to S^1\times S^1\vee S^2\overset{f\vee g}\to X$
cтягивания $c$ малой окружности в точку и отображения $f\vee g$.
Докажите, что сумма определяет действие группы $\pi_2(X)$ на множестве $[S^1\times S^1;X]$.

(b) Cумма с отображением в точку определяет взаимно-однозначное соответствие
$\pi_2(S^2)=\pi^2(S^2)\to\pi^2(S^1\times S^1)$.

(c) Введите структуру группы на $\pi_n(S^n)$.

(d)  $\deg:\pi_n(S^n)\to\Z$ является изоморфизмом групп.

(e) Определите действие группы $\pi_n(S^n)$ на $\pi^n(N)$ для компактного $n$-многообразия $N$.
\end{pr}

\begin{pr}\label{hopfpi2} Найдите $\pi_2(N)$ для
\quad
(a) компактного 2-многообразия $N$;
\quad
(b) $N=S^2\times S^1$;
 \end{pr}

\begin{pr}\label{geho-fipr}
$\pi_k(X\times Y)\cong\pi_k(X)\oplus\pi_k(Y)$.
\end{pr}

 \subsection{Теорема Фрейденталя о надстройке}

\begin{pr}\label{hopf-s1r4}
Любое гладкое (или кусочно-линейное) вложение $S^1\sqcup S^1\to\R^4$ гладко (или кусочно-линейно) изотопно
такому, у которого образы компонент лежат по разные стороны от некоторой гиперплоскости.
\end{pr}

\begin{pr}\label{hopf-sus}
Обозначим через $N$ и $S$ северный и южный полюс сферы $S^2$, соответственно.
Любое отображение $S^2\to S^2$ гомотопно такому отображению $f$, для которого

(a) $f^{-1}(N)$ и $f^{-1}(S)$ лежат внутри северной и южной полусфер $D^2_+$ и $D^2_-$, соответственно.


(b) $f(D^2_+)\subset D^2_+$ и $f(D^2_-)\subset D^2_-$.

(c) $f(\cos\alpha\cos\theta,\cos\alpha\sin\theta,\sin\alpha)=(g(\cos\theta,\sin\theta)\cos\alpha,\sin\alpha)$
для некоторого отображения $g:S^1\to S^1$.
\end{pr}

Отображение из задачи \ref{hopf-sus}.c называется {\it надстройкой} $\Sigma g$ над отображением $g:S^1\to S^1$.

\begin{pr}\label{hopf-sus2}
Надстройка $\Sigma:\pi_1(S^1)\to\pi_2(S^2)$ \quad (a) корректно определена; \quad (b) сюръективна.
\end{pr}

\begin{pr}\label{hopf-susg}
(a) Надстройка $\Sigma:\pi_q(S^n)\to\pi_{q+1}(S^{n+1})$ корректно определена.

(b) Надстройка сюръективна для $q\le 2n-1$.

(с) Любая гомотопия между надстройками над двумя отображениями $S^2\to S^2$ гомотопна, неподвижно на концах,
гомотопии в классе надстроек.

(d) Надстройка является взаимно-однозначным соответствием для $q\le 2n-2$.

(e) {\it Теорема Александера.}
Любой сохраняющий ориентацию кусочно-линейный гомеоморфизм $D^n\to D^n$ изотопен тождественному.
\end{pr}

\begin{pr}\label{hopf-sus1} (a) Определите инвариант Хопфа $H(F)$ для гомотопии $F:S^2\times I\to S^2$ между двумя
отображениями $S^2\to S^2$, являющимися надстройками.

(b) Если этот инвариант равен нулю, то $F$ гомотопна надстроечной гомотопии неподвижно на концах.

(c) Существует гомотопия $W:S^2\times I\to S^2$ между тождественными отображениями, для которой $H(W)=1$.

(d) Если надстройки над двумя петлями $S^1\to S^1$ гомотопны, то одна петля гомотопна другой.
(Значит, надстройка $\Sigma:\pi_1(S^1)\to\pi_2(S^2)$ является изоморфизмом.)

(e) $\ker(\Sigma:\pi_3(S^2)\to\pi_4(S^3))$ порождается гомотопическим классом отображения $pw$ из задачи \ref{hopf-wh}. \end{pr}

\begin{pr}\label{hopf-retw}
(a) $\Sigma(S^1\times S^1)$ ретрагируется на $\Sigma(S^1\vee S^1)$.

(b) $S^1\times S^1$ не ретрагируется на $S^1\vee S^1$.

(a',b') То же с заменой $S^1$ на $S^2$.
\end{pr}

\begin{pr}\label{hopf-wh}
Разложим $S^3=S^1\times D^2\bigcup\limits_{S^1\times S^1}D^2\times S^1$.
Рассмотрим композицию $S^3\overset w\to S^2\vee S^2\overset p\to S^2$
объединения проекций на $S^2$ и `схлопывания'.

(a) $S^2\times S^2\cong D^4/\sim$, где $x\sim y\ \Leftrightarrow \ x,y\in S^3$ и $w(x)=w(y)$.

(b) Найдите коэффициент зацепления $pw$-прообразов северного и южного полюса
(полюса отличны от отмеченной точки).

(c) $\Sigma w:S^4\to S^3\vee S^3$ гомотопно отображению в точку.

(d) $\Sigma(pw)$ гомотопно отображению в точку.
\end{pr}

\begin{pr}\label{hopf-nret}
(a) $pw$ не гомотопно отображению в точку.

(b) $w$ не гомотопно отображению в точку.

(с) $\pi_3(S^2\vee S^2)\cong\Z\oplus\Z\oplus\Z$.
\end{pr}

\begin{pr}\label{hopf-al}
(a) Коэффициент зацепления равен $\alpha$-инварианту для зацеплений $S^1\sqcup S^1\to\R^3$,
у которых обе компоненты незаузлены.

(b) Коэффициент зацепления равен $\alpha$-инварианту для зацеплений $S^1\sqcup S^1\to\R^3$.

(c)* $\alpha=\pm\Sigma^\infty\lambda$ для зацеплений $S^p\sqcup S^q\to\R^m$.
 \end{pr}

\subsection{Отображения трехмерной сферы в двумерную}\label{0hopfs3s2}

\begin{pr}\label{geho-grp}
Введите структуру группы на $\pi_3(S^2)$ (и на $\pi_k(S^n)$)
аналогично вышеописанной (\S\ref{0hopf}) структуре группы на $\pi_2(S^2)$.
\end{pr}

Представим $S^3=\{(z_1,z_2)\in\C^2\ :\ |z_1|^2+|z_2|^2=1\}$ и отождествим
$S^2$ с $\C P^1$.
Определим {\it отображение Хопфа}
$
\eta:S^3\to S^2\quad\text{формулой}\quad \eta(z_1,z_2)=(z_1:z_2).
$
Доказательство того, что это отображение не гомотопно отображению в точку (задачи \ref{geho-inho}, \ref{geho-hoet}.a) основано на следующей конструкции.

{\it Инвариантом Хопфа} отображения $S^3\to S^2$ называется коэффициент зацепления прообразов двух точек общего положения при гладкой (или кусочно-линейной) аппроксимации данного отображения.
(Коэффициент зацепления $\lk$ определен, например, в \cite{Sk}, параграф `реализуемость двумерных комплексов',
пункт `коэффициент зацепления'.)
Приведем аккуратное определение.
Любое отображение $f:S^3\to S^2$ гомотопно гладкому отображению $g$~\cite[18.4]{Pr04}.
Точка $y\in S^2$ называется {\it регулярным значением} гладкого отображения $g$, если $\rk dg(x)=2$ для любой точки $x\in g^{-1}y$.
Для гладкого отображения $g:S^3\to S^2$ имеются регулярные значения $y_1,y_2\in S^2$~\cite[18.4]{Pr04}.
Тогда $g^{-1}y_i=S^1_{i1}\sqcup S^1_{i2}\sqcup\dots\sqcup S^1_{ik_i}$ есть несвязное объединение окружностей~\cite[18.4]{Pr04}.
Определим {\it инвариант Хопфа}
$H(f):=\sum_{i=1,j=1}^{k_1,k_2}\lk(S^1_{1i},S^1_{2j}).$
Корректность этого определения вытекает из следующей задачи.

\begin{pr}\label{geho-inho} Сумма в правой части не зависит от

(a) гладкой гомотопии отображения $g$.

(b) выбора отображения $g$.

(c) выбора точки $y_1$ (и точки $y_2$).

(d) гомотопии отображения $f$.
\end{pr}

\begin{pr}\label{geho-hoet}
(a) $H(\eta)=1$.

(b) Инвариант Хопфа определяет гомоморфизм $\pi_3(S^2)\to\Z$ (тоже называемый инвариантом Хопфа).

(c) Существует бесконечно много попарно негомотопных отображений $S^3\to S^2$ (Хопф, 1931).

(d) $\R P^2$ не ретрагируется на $\R P^1$.

(e) $\C P^2\cong D^4/\sim$, где $x\sim y\ \Leftrightarrow \ x,y\in S^3$ и $\eta(x)=\eta(y)$.

(f) $\C P^2$ не ретрагируется на $\C P^1$. \quad

(f') $\pi_3(\C P^2)=0$.

(g)* Любое отображение $S^3\to S^2$ с тривиальным инвариантом Хопфа гомотопно надстройке
над некоторым отображением $S^2\to S^1$.

(i) {\bf Теорема Хопфа-Понтрягина-Фрейденталя} (1938).
{\it Инвариант Хопфа является изоморфизмом $\pi_3(S^2)\to\Z$.}
\end{pr}

Другой способ доказательства изоморфизма $\pi_3(S^2)\cong\Z$ намечен в следующей задаче \ref{geho-fr3}.


\begin{pr}\label{geho-fr0}
{\it Оснащенным набором точек} в $S^2$ называется набор точек в $S^2$ вместе с парой неколлинеарных
касательных к $S^2$ векторов в каждой точке набора.
{\it Оснащенным кобордизмом} называется компактное одномерное подмногообразие в $S^2\times I$, ортогонально подходящее к $S^2\times\{0,1\}$, вместе с парой неколлинеарных нормальных полей на этом подмногообразии.
{\it Границей} оснащенного кобордизма $(L,\xi)$ называется пара
$(L\cap S^2\times\{0,1\},\xi|_{L\cap S^2\times\{0,1\}})$ оснащенных наборов точек в $S^2$.
Два оснащенных набора точек в $S^2$ называются {\it оснащенно кобордантными}, если существует оснащенный кобордизм, границей которого является эта пара оснащенных наборов точек.

(a) Множество $\pi_2(S^2)$ находится во взаимно-однозначном соответствии с множеством
оснащенных наборов точек в $S^2$ с точностью до оснащенной кобордантности.

Это соответствие (и его обобщения) называется {\it соответствием Понтрягина}.

(b) Во что при соответствии Понтрягина переходит операция суммы?

(c) $\pi_2(S^2)\cong\Z$. (Докажите с помощью оснащенных наборов точек!)
\end{pr}

\begin{pr}\label{geho-fr3}
{\it Оснащенным зацеплением в $S^3$} называется замкнутое одномерное подмногообразие в $S^3$ вместе с парой неколлинеарных нормальных полей на этом подмногообразии.
{\it Оснащенным кобордизмом} называется компактное двумерное подмногообразие в $S^3\times I$, ортогонально подходящее к $S^3\times\{0,1\}$, вместе с парой неколлинеарных нормальных полей на этом подмногообразии.
{\it Границей} оснащенного кобордизма $(L,\xi)$ называется пара
$(L\cap S^3\times\{0,1\},\xi|_{L\cap S^3\times\{0,1\}})$ оснащенных зацеплений в $S^3$.
Два оснащенных зацепления в $S^2$ называются {\it оснащенно кобордантными}, если существует оснащенный кобордизм, границей которого является эта пара оснащенных зацеплений.

(a) Множество $\pi_3(S^2)$ находится во взаимно-однозначном соответствии с множеством
оснащенных зацеплений в $S^3$ с точностью до оснащенного кобордизма.

(b) $\pi_3(S^2)\cong\pi_1(SO_2)\cong\Z$. (Докажите с помощью оснащенных зацеплений!)
\end{pr}

\begin{pr}\label{geho-fr4}
(a) Множество $\pi_4(S^3)$ находится во взаимно-однозначном соответствии с множеством
{\it оснащенных зацеплений} в $S^4$ с точностью до  {\it оснащенного кобордизма}.
(Определите самостоятельно, что это такое.)

(b) $\pi_4(S^3)\cong\pi_1(SO_3)\cong\Z_2$.

(c) {\bf Теорема Понтрягина.} (1938). $\pi_{n+1}(S^n)\cong\Z_2$ для $n\ge3$.

(d) Во что при соответствии Понтрягина переходит гомоморфизм надстройки?

(e) Гомоморфизм $\Sigma:\pi_3(S^2)\to\pi_4(S^3)$ изоморфен приведению по модулю 2 (и, в частности, не инъективен).

(f) $\alpha$-инвариант вложения $f\sqcup g:S^3\sqcup S^3\to\R^6$ представляется (с точностью до $\Sigma^2$)
оснащенным подмногообразием $\overline f^{-1}g(S^3)\subset D^4$, где $\overline f:D^4\to\R^6$ отображение общего положения, продолжающее $f$, и на $g(S^3)$ взято произвольное нормальное оснащение.
\end{pr}

\begin{pr}\label{geho-frhi}
(a) Определите инвариант Хопфа $H:\pi_5(S^3)\to\Z$.

(b) Этот инвариант нулевой.

(c) Множество $\pi_{n+k}(S^n)$ находится во взаимно-однозначном соответствии с множеством
{\it оснащенных $k$-подмногообразий} в $S^{n+k}$ с точностью до  {\it оснащенного кобордизма}.
(Определите самостоятельно, что это такое.)

(d) Гомоморфизм $\Sigma:\pi_4(S^2)\to\pi_5(S^3)$ сюръективен.

(e)* $\pi_4(S^2)\cong\pi_4(S^3)$.

(f)* $|\pi_{n+2}(S^n)|\le2$ для любого $n$.
\end{pr}

{\bf Теорема Понтрягина.} (1950) $\pi_{n+2}(S^n)\cong\Z_2$ при $n\ge 4$.

\smallskip
{\bf Теорема Рохлина.} (1951) $\pi_{n+3}(S^n)\cong\Z_{24}$ при $n\ge 5$.

 \subsection{Отображение Хопфа как обобщение накрытия}

Еще один способ доказательства изоморфизма $\pi_3(S^2)\cong\Z$ намечен в следующих задачах \ref{geho-file}.bcd'e'.

Пусть $X$ --- произвольное многообразие (или даже полиэдр, т.е. тело симплициального комплекса).
Отображение $\t f:X\to S^3$ называется {\it поднятием} отображения $f:X\to S^2$, если $f=\eta\circ\t f$.

\begin{pr}\label{geho-file}
(a) Для любого $a\in S^2$ выполнено $\eta^{-1}a\cong S^1$.

(b) 
{\bf Лемма о локальной тривиальности.}
Для любой точки $x\in S^2$ существует гомеоморфизм $h:\eta^{-1}(S^2-\{x\})\to(S^2-\{x\})\times S^1$, для
которого $pr_1\circ h=\eta$.

(c) {\bf Лемма о поднятии пути.} Любой путь $s:[0,1]\to S^2$  имеет поднятие $\t s:[0,1]\to S^3$.

(d) Любое отображение $D^3\to S^2$ имеет поднятие $D^3\to S^3$.

(e) Любое отображение $S^3\to S^2$ гомотопно такому, которое имеет поднятие.

(f) Отображение $\eta_*:\pi_3(S^3)\to\pi_3(S^2)$, определенное композицией с отображением Хопфа,
корректно определено является гомоморфизмом.
(Даже эпиморфизмом ввиду (e).)

(d')   
{\bf Лемма о поднятии гомотопии.}
Для любых отображения $F_0:S^3\to S^3$ и гомотопии $f_t:S^3\to S^2$ отображения $f_0=\eta\circ F_0$ существует
гомотопия $F_t:S^3\to S^3$ отображения $F_0$, для которой $f_t=\eta\circ F_t$.

(e') Отображение $\eta_*$ инъективно.
\end{pr}

Лемма о поднятии гомотопии называется также {\it леммой о накрывающей гомотопии.}
Она справедлива с заменой $S^3$ на любой полиэдр $X$.

Определим отображение $\partial:\pi_2(S^2)\to\pi_1(S^1)$.
Представим произвольный элемент $\varphi\in\pi_2(S^2)$ отображением $(D^2,1)\to (S^2,1)$.
Поднятие этого отображения переводит $S^1=\partial D^2$ в $\eta^{-1}(1)=S^1$.
Обозначим через $\partial\varphi\in\pi_1(S^1)$ гомотопический класс полученного отображения $S^1\to S^1$.


\begin{pr}\label{geho-pi2s2}
(a) $\partial$ корректно определено.
\qquad
(b) $\partial$ гомоморфизм.

(c) $\partial$ мономорфизм.
\qquad
(d) $\partial$ эпиморфизм.
\end{pr}

\begin{pr}\label{geho-n3}
(a) Отображение сужения $\deg:[S^1\times S^2;S^2]\to [1\times S^2;S^2]=\Z$ сюръективно.

(b) Любое отображение $S^1\times S^2\to S^2$, сужение которого на $1\times S^2\to S^2$ гомотопно постоянному,
гомотопно такому, которое имеет поднятие $S^1\times S^2\to S^3$.

(с) Постройте
биекцию $\deg^{-1}(0)\to\Z$.

(d)*  $|\deg^{-1}(k)|=2k$ при $k\ne0$.

(c') Постройте сюръекцию $\deg:[S^1\times S^1\times S^1;S^2]\to\Z^3$ и биекцию $\deg^{-1}(0)\to\Z$.

(d')*  $|\deg^{-1}(p,q,r)|=2GCD(p,q,r)$ при $(p,q,r)\ne(0,0,0)$.
\end{pr}


\begin{pr}\label{geho-n3g} {\bf Теорема Понтрягина.}
(a) Для ориентируемого 3-многообразия $N$ постройте сюръекцию $\deg:[N;S^2]\to H_1(N;\Z)$ и биекцию $\deg^{-1}(0)\to\Z$.

(b)*  $|\deg^{-1}(\gamma)|$ есть наибольший делитель класса $[2\gamma]\in H_1(N;\Z)/Tors$.

(a') Для 3-полиэдра $N$ постройте сюръекцию $\deg:[N;S^2]\to H^2(N;\Z)$ и биекцию $\deg^{-1}(0)\to H^3(N;\Z)$.

(b')* 
Для любого $\gamma\in H^2(N;\Z)$ имеется биекция $\deg^{-1}(\gamma)\hm\to\dfrac{H^3(N;\Z)}{2\gamma\cup H^1(N;\Z)}$.

`Определение' произведения $\cup:H^1(N;\Z)\times H^2(N;\Z)\to H^3(N;\Z)$: число на симплексе 1234 равно произведению
чисел на симплексе 12 и на симплексе 234.
Впрочем, это определение естественно появляется при изучении множества $[N,S^2]$, поэтому его можно придумать, и не зная определения.

(c)* Как по $\gamma$ быстро описать $2\gamma\cup H^1(N;\Z)$?
\end{pr}

Теорема Понтрягина для ориентируемых 3-многообразий равносильна теореме классификации векторных полей
на ориентируемых  3-многообразиях из \S\ref{0vehima} (ввиду теоремы Штифеля о параллелизуемости).

\begin{pr}\label{geho-fihh}
(a) Придумайте кватернионный аналог $S^7\to S^4$ отображения Хопфа.

(b) Отображение $\partial:\pi_6(S^4)\to\pi_5(S^3)$, определенное аналогично отображению Хопфа $S^3\to S^2$,
корректно определено является гомоморфизмом.

(c) Сформулируйте и докажите аналог леммы о локальной тривиальности (задача \ref{geho-file}.b) для отображения Хопфа $S^7\to S^4$.

(d) Сформулируйте и докажите аналог леммы о поднятии гомотопии (задача \ref{geho-file}.c) для отображения Хопфа
$S^7\to S^4$.

(e) $\partial$ изоморфизм.

(f) $\Sigma:\pi_5(S^3)\to\pi_6(S^4)$ изоморфизм.
\end{pr}

\subsection{Точная последовательность расслоения }\label{0hopfbun}\nopagebreak

\begin{pr}\label{geho-rpn}
$\pi_1(\R P^n)\cong\Z_2$ для $n\ge2$ и $\pi_k(\R P^n)=0$ для $k=2,3,\dots,n-1$.
\end{pr}

Отображение $p:E\to B$ между многообразиями (или даже телами симпилициальных комплексов) называется
{\it расслоением (в смысле) Серра}, если для любых $n$, отображения $F_0:S^n\to E$ и гомотопии $f_t:S^n\to B$ отображения $f_0=p\circ F_0$ существует гомотопия $F_t:S^n\to E$ отображения $F_0$, для которой $f_t=p\circ F_t$.
(Ср. с леммой о поднятии гомотопии.)

Примеры: проекция на сомножитель произведения, накрытие, отображения Хопфа $S^3\to S^2$ и $S^7\to S^4$,
{\it отображение Хопфа} $S^{2n+1}\to\C P^n$, заданное формулой $\eta(z_0,z_1,\dots,z_n)=(z_0:z_1:\dots:z_n)$.

\begin{pr}\label{geho-ficx}
(a) Отображение Хопфа $S^{2n+1}\to\C P^n$ действительно является расслоением Серра.

(b) $\pi_2(\C P^n)\cong\Z$ и $\pi_k(\C P^n)=0$ для $k=1,3,4,5,\dots,2n$.
\end{pr}

Многообразия $E$ и $B$ называются {\it тотальным пространством} и {\it базой} расслоения Серра.
Множество $F=p^{-1}(x)$ для некоторого $x\in B$ называется его {\it слоем}.
(Формально, это определение некорректно.)
Расслоение Серра обозначается как $F\to E\to B$.

\begin{pr}\label{geho-fiexb} {\bf Точная последовательность расслоения.}
{\it Если $F\overset \subset\to E\overset p\to B$ расслоение Серра, то имеется точная последовательность}
$$\dots\to\pi_{n+1}(B)\overset{\partial}\to\pi_n(F)\overset{i}\to\pi_n(E)\overset{p}\to\pi_n(B)
\overset{\partial}\to\pi_{n-1}(F)\dots\to$$
Определите ее гомоморфизмы!
\end{pr}

Точность последовательности определена в \S\ref{0homolpa}.
Следующие задачи \ref{geho-fiso}.b для $n=3$  и \ref{geho-fist}.сde использовались при построении и изучении  характеристических классов (\S\ref{0vesy}).

\begin{pr}\label{geho-fiso}
(a) Придумайте расслоение Серра $SO_n\to SO_{n+1}\to S^n$.

(b) Гомоморфизм включения $\pi_k(SO_n)\to \pi_k(SO_{n+1})$ является изоморфизмом при $k\le n-2$ и
эпиморфизмом при $k=n-1$.
\end{pr}

{\bf Теорема Ботта о вещественной периодичности.}
{\it Группы $\pi_n(SO):=\pi_n(SO_{n+2})$ задаются следующей таблицей:} \qquad
\begin{tabular}{c|cccccccc}
$n\mod8$ &~0&~1&~2~&~3~&4    &5    &6   &7    \\
\hline
$\pi_n(SO)$ &$\Z_2$ &$\Z_2$  &0  &$\Z$  &0  &0 &0 &$\Z$
\end{tabular}

\begin{pr}\label{geho-fist} Обозначим через $V_{m,n}$ многообразие Штифеля ортонормированных наборов $n$ векторов в $\R^m$.

(a) $V_{m,1}\cong S^{m-1}$ и $V_{m,m-1}\cong V_{m,m}\cong SO_m$.

(b) Придумайте расслоение Серра $V_{m,n}\to V_{m+1,n+1}\to S^m$.

(c) Гомоморфизм включения $\pi_k(V_{m,n})\to \pi_k(V_{m+1,n+1})$ является изоморфизмом при $k\le m-2$ и
эпиморфизмом при $k=m-1$.

(d) $\pi_k(V_{m,n})=0$ для $k<m-n$.

(e) Граничный гомоморфизм $\partial:\pi_m(S^m)\to\pi_{m-1}(S^{m-1})$ расслоения Серра из (b) для $n=1$ является умножением на $\chi(S^{m-1})=1+(-1)^{m-1}$.

(f) $\pi_{m-n}(V_{m,n})=\begin{cases}\Z &n=1\text{ или $m-n$ четно}\\ \Z_2 &\text{ иначе}\end{cases}$.

(g) $\pi_4(V_{5,2})$ конечна.
\end{pr}

\begin{pr}\label{geho-pi3so}
(a) Придумайте расслоение Серра $SO_3\to SO_5\to V_{5,2}$.

(b) $\pi_3(SO_3)\cong\pi_3(SO_5)\cong\Z$, причем гомоморфизм включения $\pi_3(SO_3)\to\pi_3(SO_5)$ является умножением на 2.
\end{pr}


\subsection{Точная последовательность вложения (или Баррата-Пуппе)}\nopagebreak

\begin{pr}\label{geho-act}
(a) Определите действие $\#$ группы $\pi_q(S^m)$ на $[S^p\times S^q;S^m]$ аналогично \ref{hopf-cosu}.a.

(b) Последовательность множеств $\pi_{p+q}(S^m)\overset\#\to [S^p\times S^q;S^m]\overset r\to [S^p\vee S^q;S^m]$ точна.
(Отмеченные точки --- гомотопические классы отображений в точку; $r$ --- сужение.)

(c) Отображение $\#$ инъективно.
\end{pr}

\begin{pr}\label{geho-pup} Пусть $X,Y,Z$ комплексы и $A\subset X$ подкомплекс.

(a) $X/A\simeq X\cup_A\con A$.

(b) $\Sigma A\simeq \con X\cup_{\con A} \Sigma A$.

(c) {\bf Точная последовательность вложения (или Баррата-Пуппе).}
{\it
Последовательность множеств
$$\dots\to[\Sigma X;Z]\overset r\to[\Sigma A;Z]\overset \varphi\to[X/A;Z]\overset\#\to [X;Z]\overset r\to [A;Z]$$
точна.
Здесь  $r$ --- сужение; \ $\#$ --- композиция с проекцией $X\to X/A$; \ $\varphi$ --- композиция с композицией
$X/A\to X\cup_A\con A\to\Sigma A$ гомотопической эквивалентности и сжатия пространства $X$ в точку; отмеченные точки --- гомотопические классы отображений в точку.}

(d) $\Sigma (X\times Y)$ ретрагируется на $\Sigma(X\vee Y)$.

(e) $[X\times Y/X\vee Y;Z]\overset\#\to [X\times Y;Z]$ инъективно.
\end{pr}

 \subsection{Реализация циклов подмногообразиями (набросок)}

Приведем план наглядного изложения простейших результатов о гомотопической классификации отображений многообразий
и их следствий о реализуемости гомологических классов подмногообразиями.

\begin{pr}\label{} Пусть $N$ --- замкнутое $n$-многообразие.

(a) При $n\ge3$ любой гомологический класс по модулю 2 размерности 1 реализуется
вложенной окружностью, и гомологичность по модулю 2 между окружностями реализуется
2-подмногообразием многообразия $N\times I$.

(b) При $n\ge5$ любой гомологический класс по модулю 2 размерности 2 реализуется
2-подмногообразием (естественно, замкнутым), и гомологичность по модулю 2 между 2-подмногообразиями реализуется 3-подмногообразием (естественно, с краем) многообразия $N\times I$.
\end{pr}

\begin{pr}\label{}
Пусть $N$ --- замкнутое ориентируемое $n$-многообразие.

(a) Сформулируйте и докажите аналоги предыдущей задачи.

(b) Отображение
$$[N;S^1]\to H_{n-1}(N;\Z),\quad [f]\mapsto[f^{-1}(x)],$$
ставящее в соответствие гомотопическому классу отображения $f:N\to S^1$
гомологический класс прообраза $f^{-1}(x)$ регулярного значения, корректно
определено.

(c) Для любого $y\in H_{n-1}(N;\Z)$ существует и единственно (с точностью до
гомотопии) отображение $f:N\to S^1$, для которого гомоморфизм
$f_*:H_1(N)\to H_1(S^1)\cong\Z$ задается формулой $f_*(a)=a\cap y$.

(d) Отображения из предыдущих пунктов являются взаимно обратными биекциями
$[N;S^1]\leftrightarrow H_{n-1}(N;\Z)$.

(e) Любой целочисленный гомологический класс коразмерности 1 реализуется ориентируемым
подмногообразием, и целочисленная гомологичность между ориентируемыми подмногообразиями
коразмерности 1 реализуется ориентируемым подмногообразием
многообразия $N\times I$.
\end{pr}

\begin{pr}\label{}
  Пусть $N$ --- замкнутое ориентируемое $n$-многообразие.

(a) Отображение
$$[N;\C P^n]\to H_{n-2}(N;\Z),\quad [f]\mapsto[f^{-1}(\C P^{n-1})]$$
ставящее в соответствие гомотопическому классу отображения $f:N\to \C P^n$,
трансверсального вдоль $\C P^{n-1}\subset\C P^n$, гомологический класс
прообраза $f^{-1}(\C P^{n-1})$, корректно определено и является биекцией.

(b) Любой целочисленный гомологический класс коразмерности 2 реализуется ориентируемым
подмногообразием, и целочисленная гомологичность между ориентируемыми подмногообразиями
коразмерности 2 реализуется реализуется ориентируемым подмногообразием
многообразия $N\times I$.

(c) Если $n\le5$, то любой целочисленный класс гомологий с целыми коэффициентами реализуется
ориентированным подмногообразием.
\end{pr}

 {\bf Теорема Тома.} {\it Некоторое кратное любого $k$-мерного целочисленного
гомологического класса замкнутого ориентируемого $n$-многообразия реализуется
отображением ориентируемого многообразия; при $k<n/2$ это отображение можно
считать вложением}~\cite[\S27, Следствие 3]{DNF84}.

\begin{pr}\label{pdrealc1m2}
 Пусть $N$ --- замкнутое $n$-многообразие.

(a) Отображение
$$[N;\R P^{n+1}]\to H_{n-1}(N),\quad[f]\mapsto[f^{-1}(\R P^n)]$$
ставящее в соответствие гомотопическому классу отображения $f:N\to\R P^{n+1}$,
трансверсального вдоль $\R P^n\subset\R P^{n+1}$, гомологический класс
прообраза $f^{-1}(\R P^n)$, корректно определено и является биекцией.

(b) Любой гомологический класс по модулю 2 коразмерности 1 реализуется
подмногообразием, и гомологичность по модулю 2 между подмногообразиями
коразмерности 1 реализуется подмногообразием многообразия $N\times I$.

(c)* При $n\le4$ любой гомологический класс по модулю 2 реализуется
подмногообразием, но не любая гомологичность по модулю 2 между
2-подмногообразиями реализуется 3-подмногообразием в $N\times I$~\cite[II.1]{Ki89}.
\end{pr}

\begin{pr}\label{pdrealhopf}
  Пусть $N$ --- замкнутое ориентируемое $n$-многообразие.

(a) Любое ориентируемое 2-многообразие, кроме сферы, ретрагируется на некоторую окружность.

(b) {\it Теорема Хопфа.} $N$ ретрагируется на окружность $\Sigma\subset N$
тогда и только тогда, когда $[\Sigma]\subset H_1(N;\Z)$ {\it примитивен}, т.е.
не делится ни на одно число, отличное от $\pm1$.

(c) $N$ ретрагируется на некоторую окружность тогда и только тогда, когда
группа $H_1(N;\Z)$ бесконечна (А. Т. Фоменко и И. Н. Шнурников).
\end{pr}


\subsection{Указания и решения к некоторым задачам}

\smallskip
{\bf \ref{hopf-sus}, \ref{hopf-susg}.} \cite[\S10]{FF89}.

\smallskip
{\bf \ref{hopf-sus2}.} (b) Следует из \ref{hopf-sus}.c.

\smallskip
{\bf \ref{hopf-sus1}.} (b) Примените оснащенные зацепления.

(c) Тождественное отображение гомотопно коммутатору $aba^{-1}b^{-1}$ петель.

(d) Следует из (b,c).

(e) Аналогично (a,b,c,d).

\smallskip
{\bf \ref{hopf-retw}.}
(a) Рассмотрим композицию
$$\Sigma(S^1\times S^1)\overset c\to \Sigma(S^1\times S^1)\vee\Sigma(S^1\times S^1)
\overset{\Sigma pr_1\vee \Sigma pr_2}\to \Sigma(S^1\times*)\vee\Sigma(*\times S^1)=\Sigma(S^1\vee S^1),$$
где $c$ --- стягивание экватора в точку.
Докажите, что ее сужение на $\Sigma(S^1\vee S^1)$ гомотопно тождественному отображению.
Примените теорему Хопфа о продолжении гомотопии.

(a') Аналогично (a).

(b) Примените $\pi_1$. \qquad

(b') Примените \ref{hopf-wh}.a и \ref{hopf-nret}.b.

\smallskip
{\bf \ref{hopf-wh}.}
(b) Ответ: 2.
\qquad
(с) Ввиду (a) это равносильно \ref{hopf-retw}.a'.

\smallskip
{\bf \ref{hopf-nret}.}
(a) Примените инвариант Хопфа или оснащенные зацепления.

(b,c) Примените обобщение инварианта Хопфа или двуцветные оснащенные зацепления.

(d) Следует из (c).

\smallskip
{\bf \ref{geho-inho}.} Аналогично степени \ref{geho-inho}.

(b,d) Следует из (a): гомотопируйте произвольную гомотопию к гладкой.

(c) Следует из (a): покрутите сферу $S^2$.

\smallskip
{\bf \ref{geho-hoet}.} (c) Следует из (a) и (b).

(d) Примените $\pi_1$.

(f) Ввиду (e) это равносильно негомотопности нулю отображения $\eta$.

(g) Удобнее всего делать через оснащенные зацепления.

(i) Следует из (a,b,g) и гомотопности любого отображения $S^2\to S^1$ постоянному.

\smallskip
{\bf \ref{geho-fr0}.} Аналогично  \cite[18.5]{Pr04}.

\smallskip
{\bf \ref{geho-fr3}.} (a) \cite[18.5]{Pr04}.

(b) Сначала докажите, что любое зацепление {\it кобордантно} тривиальному.
Для этого используйте ориентированную поверхность в $\R^3$, ориентированной границей которой является
ориентированной зацепление Хопфа.
Потом позаботьтесь об оснащении.

\smallskip
{\bf \ref{geho-fr4}.} (a) Аналогично предыдущей задаче.

(b) Докажем, что гомоморфизм $J:\pi_1(SO_3)\to\pi_4(S^3)$, определенный оснащением стандартной окружности, является изоморфизмом.
Ввиду общности положения любое оснащенное зацепление в $S^4$ оснащенно кобордантно такому, носителем которого является стандартная окружность.
Значит, $J$ сюръективен.
Для доказательства его инъективности рассмотрим двумерный оснащенный кобордизм в $S^4\times I$ между как-то
оснащенными стандартными окружностями.
Ввиду общности положения {\it тривиальные}
оснащения на его краях продолжаются до оснащения на всем кобордизме.
Поэтому {\it заданное} оснащение кобордизма можно рассматривать как отображение из сферы с ручками и дырками в $SO_3$. Сумма `степеней' сужений такого отображения на его граничные окружности равна нулю.
Поэтому стандартные окружности были оснащены одинаково.

{\it Замечание.} Гомоморфизм $J:\pi_n(SO_k)\to\pi_{n+k}(S^k)$ определяется оснащением стандартной сферы.
Вообще говоря, он не является ни мнонморфизмом, ни эпиморфизмом.
Какие его свойства можно доказать аналогично предыдущему?

(с) Аналогично (b). Или следует из (b) и теоремы Фрейденталя о настройке.

(d) В гомоморфизм, индуцированный включением экватора в сферу.

(e) Гомоморфизм $\pi_1(SO_2)\to \pi_1(SO_3)$, индуцированный включением $SO_2\to SO_3$,
изоморфен приведению по модулю 2.

\smallskip
{\bf \ref{geho-frhi}.}
(b) Покрутите сферу $S^2$, чтобы северный полюс совместился с южным.

(c) \cite[18.5]{Pr04}; ср. ~\cite{Po76}, \cite{CRS}.

(d) Выведите из (b) и (c).

(e) Аналогично задаче \ref{geho-file}.

(f) Следует из (d,e), \ref{geho-fr4}.b, теоремы Фрейденталя о надстройке и
гомотопности любого отображения $S^2\to S^1$ постоянному.

\smallskip
{\bf \ref{geho-file}.}
(c,d,d',f) 
Аналогично доказательству для накрытий.

(e) Представим данное отображение $\varphi:S^3\to S^2$ отображением $D^3_+\to S^2$, переводящим
$\partial D^3_+$ в точку 1.
По (d) существует поднятие $\t\varphi_+:D^3_+\to S^3$ последнего отображения.
Так как $\t\varphi_+(\partial D^3_+)\subset\eta^{-1}(1)=S^1$, то
$\t\varphi_+|_{\partial D^3_+}$ продолжается до отображения $\t\varphi_-:D^3_-\to \eta^{-1}(1)$.
Отображения $\t\varphi_+$ и $\t\varphi_-$ образуют отображение $\t\varphi:S^3\to S^3$.
Отображение $\eta\circ\t\varphi$ имеет поднятие и гомотопно $\varphi$.

(e') Аналогично (e).
Пусть для отображения $F:S^3\to S^3$ композиция $\eta\circ F$ гомотопна отображению в точку $1\in S^2$.
По (d) существует накрывающая гомотопия $F_t$.
Тогда $F_1(S^3)\subset\eta^{-1}(1)=S^1$.
Значит, $F_1$ пропускается через $S^1$ и потому гомотопно отображению в точку.
Следовательно, и $F$ гомотопно отображению в точку.

\smallskip
{\bf \ref{geho-pi2s2}.} Аналогично задаче \ref{geho-file}.

\smallskip
{\bf \ref{geho-n3}.} (b) Аналогично задаче \ref{geho-file}.e.
\qquad
(c) Аналогично задаче \ref{geho-file}.fe'.

(c') Аналогично (a,b,c).
\qquad
(d,d') \cite{CRS}.

\smallskip
{\bf \ref{geho-n3g}.}
(a,a') \cite[3.1.4]{Pr06}. \qquad (b,b') \cite{CRS}.

\smallskip
{\bf \ref{geho-rpn}.} Используйте накрытие $S^n\to\R P^n$.

\smallskip
{\bf \ref{geho-fist}.}
(c) Из точной последовательности расслоения (b).

(d) Индукция по $n$.

(f),(g) Из точной последовательности расслоения (b) для $n=1$ и (e).

\smallskip
{\bf \ref{geho-pi3so}.} (b) Из точной последовательности расслоения (a) и \ref{geho-fist}.eg.

\smallskip
{\bf \ref{geho-act}.} (c) См. \ref{geho-pup}.e.

\smallskip
{\bf \ref{geho-pup}.} (c) Точность в $[X,Z]$ проверяется несложно (аналогично \ref{geho-act}.b).
Точность в остальных членах аналогична, ибо последовательность можно начать с $[X,Z]$,
положив $A':=X$ и $X':=X\cup\con A$.
Не забудьте доказать для этой редукции коммутативность квадратов и треугольников, например, что сужение на $X$
композиции гомотопической эквивалентности $X\cup\con A\to X/A$ и отображения $f:X/A\to Z$
гомотопно композиции стягивания $X\to X/A$ и отображения $f$.

\smallskip
{\bf \ref{pdrealc1m2}.} (b)
Приведем рассуждения для $n=3$.
Возьмем некоторые триангуляцию $T$ многообразия $N$ и 2-цикл $a$, представляющий данный класс из $H_2(T)$.
К каждому ребру триангуляции $T$ примыкает четное число граней цикла $a$.
`Растащим' эти грани на пары, получим 2-цикл в некотором измельчении $T'$ триангуляции $T$,
гомологичный 2-циклу $a$ и представленный 2-комплексом $K$, каждая точка которого имеет в $K$ окрестность,
изоморфную конусу над несвязным объединением окружностей.
Тех точек, для которых окружностей больше одной, конечное число.
Для каждой из таких точек `растащим' конусы, отвечающие разным окружностям.
Получим 2-цикл в некотором измельчении $T''$ триангуляции $T'$, гомологичный 2-циклу $a$ и представленный несамопересекающимся связным замкнутым локально евклидовым 2-подкомплексом (не обязательно ориентируемым).

\smallskip
{\bf \ref{pdrealhopf}.} (b) Докажите, что оба этих условия равносильны наличию
$(n-1)$-подмногообразия в $N$, трансверсально пересекающего $\Sigma$ ровно в одной точке.

\newpage
\section{Классификация погружений}\label{0imm}

\subsection{Выворачивание сфер наизнанку и классификация погружений}\nopagebreak

Пусть $N$ --- гладкое подмногообразие в $\R^M$ (см. определение в \S\ref{0ve2ma}).
Мы работаем в гладкой категории, т.е., опускаем прилагательное `гладкий' при словах `вложение', `погружение' и т.д.

Отображение $f:N\to\R^m$ называется {\it погружением}, если $df(x)\ne0$ для любой точки $x\in N$.

Например, на рис. \ref{tordyr} справа,
\ref{td} справа и \ref{kleinb} изображены погружения тора с дыркой в $\R^2$, тора с тремя дырками в $\R^2$ и
бутылки Клейна в $\R^3$.

Например, любое вложение является погружением.

Погружения $f_0,f_1:N\to\R^m$ называются {\it регулярно гомотопными}, если существует такое погружение
$H:N\times I\to S^m\times I$, что

$\bullet$ $H(x,0)=(f_0(x),0)$ и $H(x,1)=(f_1(x),1)$ для любого $x\in N$,

$\bullet$ $H(N\times\{t\})\subset \R^m\times\{t\}$  для любого $t\in I.$

(Это не то же самое, что гомотопность в классе погружений.)

Например,

$\bullet$ стандартное вложение $S^1\to\R^2$ не регулярно гомотопно композиции стандартного вложения и отражения относительно прямой.

$\bullet$ вложения цилиндра в $\R^3$ на рис. \ref{5-1} в середине не регулярно гомотопны;
если нижнюю ленточку еще раз перекрутить, то полученная ленточка будет регулярно гомотопна верхней;

$\bullet$ вложения листа Мебиуса в $\R^3$ на рис. \ref{glu} в предпоследней колонке регулярно гомотопны;
они не регулярно гомотопны вложению на рис. \ref{glu} в последней колонке.

Обозначим через $I^m(S^n)$ множество погружений $S^n\to\R^m$ c точностью до регулярной гомотопности.
В решении задач этого пункта (кроме теоремы общего положения) можно использовать без доказательства то,
что определенные Вами отображения, аналогичные расслоениям Серра из доказательства теоремы Смейла, являются расслоениями Серра.

\begin{pr}\label{em-wh}
(a) {\bf Теорема Уитни для плоскости.}
{\it Существует биекция
$I^2(S^1)\to\Z$.
Эта биекция определяется количеством оборотов касательного вектора. }

(b) $|I^3(S^1\times[0,1])|=2$.
\qquad

(c) $|I^3(S^1\times S^1-\Int D^2)|=4$.

(d) {\bf Теорема общего положения.} {\it Любые два погружения $S^n\to\R^m$ регулярно гомотопны при $m\ge 2n+1$.}
\end{pr}

Знаменитый результат Смейла состоит в том, что двумерная сфера выворачивается наизнанку в трехмерном евклидовом  пространстве (т.е. стандартное вложение сферы $S^2$ в $\R^3$ регулярно гомотопно композиции стандартного вложения и отражения относительно гиперплоскости).

\smallskip
{\bf Трехмерная теорема Смейла.} {\it Любые два погружения $S^2\to\R^3$ регулярно гомотопны.}

\smallskip
{\bf Теорема Кайзера.}
{\it Любые два погружения $S^n\to\R^{n+1}$ регулярно гомотопны тогда и только тогда, когда $n\in\{0,2,6\}$.}

\smallskip
Более того, при $n\not\in\{0,2,6\}$ сфера $S^n$ не выворачивается наизнанку в $\R^{n+1}$.
Этот результат \cite{Ka84}
связан с параллелизуемостью
сфер $S^1$, $S^3$ и $S^7$.

\smallskip
{\it Набросок доказательства трехмерной теоремы Смейла.}
Обозначим через $\Sigma_2$ пространство погружений $S^2\to\R^3$, для которых в окрестности точки
$1\in\partial D^2$ погружение стандартно.
Достаточно доказать, что $\Sigma_2$ связно.
Пространство $\Delta_2$ определяется аналогично пространству $\Sigma_2$ с заменой $S^2$ на $D^2$.
Очевидно, что $\Delta_2$ стягиваемо.
Обозначим через $\Sigma_1$ пространство погружений $S^1\to\R^3$ с единичным нормальным векторным полем, для которых в окрестности точки $1\in S^1$  и погружение, и поле стандартны.
Рассмотрим отображение $\Delta_2\to\Sigma_1$, которое ставит в соответствие погружению диска его сужение на
границу вместе с векторным полем, смотрящим из точек границы диска внутрь диска.
Слои этого отображения гомотопически эквивалентны $\Sigma_2$.
Наиболее содержательная и трудная часть доказательства теоремы Смейла состоит в том, что описанное отображение
является расслоением Серра.
Если это доказано, то из точной последовательности расслоения (учитывая, что $\Delta_2$ стягиваемо) получаем
$\pi_0(\Sigma_2)\cong\pi_1(\Sigma_1)$.

Пространство $\Delta_1$ определяется аналогично пространству $\Sigma_1$ с заменой $S^1$ на $D^1$.
Очевидно, что $\Delta_1$ стягиваемо.
Рассмотрим отображение $\Delta_1\to SO_3$, которое ставит в соответствие погружению отрезка с нормальным векторным полем репер из нормального вектора в точке $-1$,  касательного вектора к погружению в точке $-1$ и третьего вектора, составляющего с ними ортонормированный репер.
Слои этого отображения гомотопически эквивалентны $\Sigma_1$.
Аналогично предыдущему это отображение является расслоением Серра.
Поэтому из точной последовательности расслоения (учитывая, что $\Delta_1$ стягиваемо) получаем
$\pi_1(\Sigma_1)\cong\pi_2(SO_3)=0$.
QED

\begin{pr}\label{em-sm}
(a) $|I^3(S^1\times S^1)|=4$.

(a')** Сколько классов регулярной гомотопности вложений $S^1\times S^1\to\R^3$: 3 или 4?

(b) Используя $\pi_6(SO_7)=0$, докажите теорему Кайзера для $n=6$.

(c)* Докажите теорему Кайзера для $n=3$.

(d) {\bf Теорема Смейла.} {\it Существует взаимно-однозначное соответствие $I^m(S^n)\to \pi_n(V_{m,n})$. }


(e) {\bf Теорема Смейла-Хирша.} {\it
Если $n$-многообразие $N$ параллелизуемо, то существует взаимно-однозначное соответствие $I^m(N)\to [N,V_{m,n}]$.
Эта биекция определяется формулой $f\mapsto (df:\R^n\to\R^m)$.}

(f) {\bf Теорема Уитни.}
{\it Существует биекция  $I^{2n}(S^n)\to\Z$ для $n=1$ или $n$ четного и $I^{2n}(S^n)\to\Z_2$ для остальных $n$.
Эта биекция определяется количеством точек самопересечения аппроксимации общего положения, со знаком для $n=1$ или $n$ четного. }
\end{pr}

\smallskip
Заметим, что биекции из теоремы Уитни для плоскости и теоремы Уитни для $n=1$ различны.

Биекция из  теоремы Смейла определяется так.
Для погружений $f,g:S^n\to\R^m$ совместим $df,dg:TS^n\to\R^m$ гомотопией на $TD^n_-$ в классе линейных мономорфизмов.
Тогда два полученных отображения $TD^n_+\to\R^m$ совпадают на границе.
Значит, они дают отображение $S^n\to V_{m,n}$.
Обозначим через $\Omega(f,g)\in\pi_n(V_{m,n})$ его гомотопический класс.
(Это {\it препятствие} к гомотопности дифференциалов $df,dg:TS^n\to\R^m$ в классе линейных мономорфизмов.)
Тогда для любого погружения $g:S^n\to\R^m$ отображение $f\mapsto\Omega(f,g)$ определяет биекцию
$I^m(S^n)\to \pi_n(V_{m,n})$.

\smallskip
{\bf Теорема Кервера.} {\it Любые два вложения $S^n\to\R^m$ регулярно гомотопны при $2m\ge 3n+2$.}

\smallskip
Набросок доказательства приведен в следующем пункте.

Из теоремы Кервера следует, что {\it при $2m\ge 3n+2$ любое вложение $S^n\to\R^m$ имеет тривиальное нормальное расслоение.}
Это не так для любых $n=4l-1$ и каждого $m=4l+2,4l+3,\dots,6l-1$ \cite[6.8]{Ha66}.
Значит, размерностное ограничение $2m\ge 3n+2$ в теореме Кервера ослабить нельзя.

\begin{pr}\label{em-hhg}* Для некоторых $n\not\equiv-1\mod4$ попробуйте ослабить размерностное ограничение
в теореме Кервера, используя усиления теоремы Фрейденталя о надстройке (`трудную часть' Уайтхеда,
теорему Джеймса о двойной надстройке,  EHP-последовательность,...).
\end{pr}

\begin{pr}\label{em-lin}
Множество линейных мономорфизмов $TS^n\to\R^m$ с точностью до гомотопии (в классе линейных мономорфизмов)
находится во взаимно-однозначном соответствии с группой $\pi_n(V_{m,n})$.
\end{pr}

Классифицировать погружения произвольного многообразия можно на языке векторных расслоений (\S\ref{0vebun})~\cite{Hi60}, см. также~\cite{HP64}.

Ясно, что если $n$-многообразие $N$ погружаемо в $\R^{n+k}$, то существует векторное $k$-расслоение $\nu$ над $N$ такое, что $\nu\oplus\tau_N\cong(n+k)\varepsilon$.
(Действительно, так как $\nu(f)\oplus\tau(N)=(n+k)\varepsilon$ для погружения $f:N\to\R^{n+k}$, то можно взять $\nu=\nu(f)$ .
Необходимое условие из теоремы Уитни о препятствии фактически является
необходимым для существования такого расслоения $\nu$ --- ввиду формулы Уитни-Ву).
Это необходимое условие оказывается {\it достаточным}.

\smallskip
{\bf Теорема Смейла-Хирша.}
{\it Пусть $N$ --- гладкое $n$-многообразие.

(a) Если существует векторное $k$-расслоение $\nu$ над $N$ такое, что $\nu\oplus\tau_N\cong(n+k)\varepsilon$,
то $N$ погружаемо в $\R^{n+k}$.

(b) Если $f,g:N\to\R^{n+k}$ --- погружения такие, что $df,dg:TN\to\R^{n+k}$ гомотопны
в классе линейных мономорфизмов, то $f$ и $g$ регулярно гомотопны.}

\subsection{Набросок доказательства теоремы Кервера}\nopagebreak

Для достаточно малой окрестности $O\Delta$ диагонали $\Delta$ в $S^n\times S^n$
обозначим
$$SS^n=O\Delta-\Delta.$$
Для погружения $h:S^n\to\R^m$ определим эквивариантное отображение
$$\t h:SS^n\to S^{m-1}\quad\text{формулой}\quad\t h(x,y)=\frac{hx-hy}{|hx-hy|}.$$

\begin{pr}\label{em-hh}
(a) Если погружения $h_0$ и $h_1$ регулярно гомотопны, то $\t{h_0}\simeq_{eq}\t{h_1}$.

(b) Если $h_0,h_1:S^n\to\R^m$ --- вложения, то $\t{h_0}\simeq_{eq}\t{h_1}$.

(c) {\bf Теорема Хефлигера-Хирша.}
{\it Если $h_0,h_1:S^n\to\R^m$ --- два погружения, $\t{h_0}\simeq_{eq}\t{h_1}$ и $2m\ge 3n+2$,
то $h_0$ и $h_1$ регулярно гомотопны.}
\end{pr}

Теорема Хефлигера-Хирша справедлива с заменой $S^n$ на произвольное многообразие; доказательство получается из нижеприведенного применением теории препятствий.

Теорема Кервера вытекает из наблюдения \ref{em-hh}.b и теоремы Хефлигера-Хирша \ref{em-hh}.c.

Для доказательства теоремы Хефлигера-Хирша введем следующее понятие.
{\it Эквивариантное многообразие Штифеля} $V^{eq}_{m,n}$ --- пространство эквивариантных
(относительно антиподальных инволюций) отображений $S^{n-1}\to S^{m-1}$.

\begin{pr}\label{em-lhh}
(a) $V_{m,1}^{eq}\cong S^{m-1}$.

(b) Множество эквивариантных отображений $SS^n\to S^{m-1}$ с точностью до эквивариантной гомотопии
находится во взаимно-однозначном соответствии с группой $\pi_n(V^{eq}_{m,n})$.
\end{pr}

Теорема Хефлигера-Хирша вытекает из построения биекции в теореме Смейла (\ref{em-sm}.d), наблюдений \ref{em-lin} и \ref{em-lhh}.b, а также следующей леммы для $k=n$.

\smallskip
{\bf Лемма.} \cite[1.1]{HH62} {\it Гомоморфизм включения $\rho_n:\pi_k(V_{m,n})\to\pi_k(V^{eq}_{m,n})$ есть изоморфизм
для $0\le k\le2(m-n)-2$ и эпиморфизм для $k=2(m-n)-1$.}

\begin{pr}\label{em-5l}
(a) {\it 5-лемма.} Если в коммутативной диаграмме
$$
\begin{array}{ccccc}
A_1&\to A_2&\to A_3&\to A_4&\to A_5\\
\varphi_1\downarrow&\varphi_2\downarrow&\varphi_3\downarrow&\varphi_4\downarrow
&\varphi_5\downarrow\\
B_1&\to B_2&\to B_3&\to B_4&\to B_5\\
\end{array}
$$
с точными строками отображения $\varphi_1,\varphi_2,\varphi_4$ и $\varphi_5$
являются изоморфизмами, то и $\varphi_3$ --- изоморфизм.

(b) Из восьми предположений 5-леммы ($\varphi_1,\varphi_2,\varphi_4$ и
$\varphi_5$ --- мономорфизмы и $\varphi_1,\varphi_2,\varphi_4$ и
$\varphi_5$ --- эпиморфизмы) для доказательства мономорфности отображения
$\varphi_3$ требуются только три (какие?), а для доказательства его
эпиморфности --- другие три (какие?). Два оставшихся предположения являются,
таким образом, вообще лишними.

(c) Если, сохранив все предположения 5-леммы, исключить из ее диаграммы
стрелку $\varphi_3$, будет ли верно, что $A_3\cong B_3$?
 \end{pr}

\begin{pr}\label{em-res}
(a) Слой отображения сужения $r^{eq}:V^{eq}_{m,2}\to V^{eq}_{m,1}$ гомеоморфен пространству
$\Omega S^{m-1}$ всех отображений $(S^1,1)\to (S^{m-1},1)$.

(b) Слой отображения сужения $r^{eq}:V^{eq}_{m,n}\to V^{eq}_{m,n-1}$ гомеоморфен $\Omega^{n-1}S^{m-1}$.

(с)* Отображение сужения $r^{eq}:V^{eq}_{m,n}\to V^{eq}_{m,n-1}$ является расслоением Серра.
\end{pr}

{\it Доказательство леммы.}  Индукция по $n$.
Для $n=1$ лемма верна, поскольку $V_{m1}\cong V_{m1}^{eq}\cong S^{m-1}$.
Теперь предположим, что $n\ge2$.
Возьмем расслоения Серра
$$S^{m-n}\to V_{mn}\overset r\to V_{m,n-1}\quad\text{и}\quad
\Omega^{n-1}S^{m-1}\to V^{eq}_{mn}\overset{r^{eq}}\to V^{eq}_{m,n-1},$$
где $r$ и $r^{eq}$ --- сужения.
Ясно, что $r^{eq}\rho_n=\rho_{n-1}r$.
Включение $\rho:V_{m,n}\to V^{eq}_{m,n}$ на слоях сужений индуцирует композицию
$$\pi_i(S^{m-n})\overset{\Sigma^{n-1}}\to\pi_{i+n-1}(S^{m-1})\overset\cong\to \pi_i(\Omega^{n-1}S^{m-1}).$$
Поэтому $\rho_n$ индуцирует отображение точных гомотопических последовательностей раслоений:
$$\xymatrix{\pi_{i+1}(V_{m,n-1}) \ar[r] \ar[d]^{\rho_{n-1}} & \pi_i(S^{m-n}) \ar[r] \ar[d]^{\Sigma^{n-1}} &
\pi_i(V_{m,n}) \ar[r] \ar[d]^{\rho_n} & \pi_i(V_{m,n-1}) \ar[r] \ar[d]^{\rho_{n-1}}
& \pi_{i-1}(S^{m-n}) \ar[d]^{\Sigma^{n-1}} \\
\pi_{i+1}(V^{eq}_{m,n-1}) \ar[r] & \pi_{i+n-1}(S^{m-1}) \ar[r] & \pi_i(V^{eq}_{mn}) \ar[r] &
\pi_i(V^{eq}_{m,n-1}) \ar[r] & \pi_{i+n-2}(S^{m-1})
}.$$
Предположим сначала, что $i\ge1$.
По предположению индукции $\rho_{n-1}$ --- изоморфизмы.
По теореме Фрейденталя о надстройке $\Sigma^{n-1}$ --- изоморфизмы.
По 5-лемме, $\rho_n$ --- изоморфизм для $i\le2(m-n)-2$.

Аналогично, $\rho_n$ --- эпиморфизм для $i=2(m-n)-1$.

Доказательство для $i=0$ аналогично, только правый столбец рассматриваемой
диаграммы должен быть заменен нулями, поскольку $r$ и $r^{eq}$ --- сюръекции.
\footnote{Заметим, что то, что мы обозначаем через $m$ и $n$, в \cite{HH62}
обозначалось через $n$ и $m$, соответственно.
Заметим также, что в \cite{HH62} имеются следующие опечатки: в формулировке условия $(1.1_{n-1})$
и в доказательстве утверждения (1.2) $\pi_i(V^{eq}_{m,n},V_{m,n})=0$
должно быть для $0\le i\le2(m-n)-1$, а не только для $0<i<2(m-n)-1$).}
QED


\newpage
\section{Вложения и заузливания}\label{0emb}

\hfill{\it It is a riddle which shares with the universe}

\hfill{\it the merit of having no answer.}

\hfill{\it W.\?Maugham, The Moon and Sixpence.
\footnote{\it Это загадка, имеющая общее с мирозданием достоинство:
отсутствие ответа. У. С. Моэм, Луна и грош, пер. автора.}
}

\subsection{Введение: проблемы вложимости и заузливания}

Как писал Е. К. Зиман, тремя классическими проблемами топологии являются

\smallskip
(1) {\it Проблема гомеоморфизма.}
Когда данные два пространства $N$ и $M$ гомеоморфны?
Как описать множество гомеоморфических классов мно\-го\-об\-ра\-зий из
заданного класса, например, заданной размерности $n$?

(2) {\it Проблема вложимости.}
Какие пространства $N$ вложимы в $S^m$ для данного $m$?

(3) {\it Проблема заузливания.} Какие вложения $f,g:N\to S^m$ изотопны?
Как описать множество изотопических класов вложений $N\to S^m$?

\smallskip
Идеи и методы, применяемые для изучения проблем вложимости и заузливания,
применяются и для проблемы гомеоморфизма (и для других проблем топологии и ее
приложений).


Определение гладкого и кусочно-линейного вложения приведены
в пунктах 'нормальные классы Уитни' и 'препятствие Уитни к вложимости'.
Два вложения $f,g:N\to S^m$ называются {\it (объемлемо) изотопными}, если
существует такой гомеоморфизм $F: S^m\times I\to S^m\times I$, что

(i) $F(y,0)=(y,0)$ для любого $y\in S^m,$

(ii) $F(f(x),1)=(g(x),1)$ для любого $x\in N$, и

(iii) $F(S^m\times\{t\})= S^m\times\{t\}$  для любого $t\in I.$

Этот гомеоморфизм $F$ называется (объемлющей) изотопией.
(Объемлющей) изотопией также называют гомотопию $S^m\times I\to S^m$ или
семейство отображений $F_t: S^m\to S^m$, очевидным образом порожденные
отображением $F$.

Через CAT будем обозначать гладкую (DIFF) или ку\-соч\-но-ли\-ней\-ную (PL)
категорию.
Если определение или утверждение имеет силу в обеих категориях, то CAT
опускается.
Если любые два вложения $N\to S^m$ (объемлемо) изотопны, то $N$ называется
{\it незаузленным в $ S^m$}.


\begin{pr}\label{em-namb}
 Два PL вложения $f,g:N\to S^m$ называются {\it PL не объемлемо
изотопными}, если существует такое PL вложение $F:N\times I\to S^m\times I$,
что

\quad $F(x,0)=(f(x),0)$ для любого $x\in N$,

\quad $F(x,1)=(g(x),1)$ для любого $x\in N$, и

\quad $F(N\times\{t\})\subset S^m\times\{t\}$  для любого $t\in I.$

(a) Любые два PL вложения $S^1\to S^3$ (т.е. узла) являются PL не объемемо
изотопными.

(b) То же для $S^n\to S^m$.

(c)* Определение {\it гладкой не объемлемой изотопности} гладких вложений
получается из предыдущего определения заменой PL на DIFF.
Существуют два гладких вложения $S^1\to S^3$ (т.е. узла), не являющиеся гладко не объемемо изотопными.
\end{pr}

\subsection{Общее положение}

 {\bf Теорема.} {\it Любой $n$-мерный полиэдр вложим в $\R^{2n+1}$.}

\smallskip
{\it Набросок доказательства.}
Это доказательство достаточно разобрать для $n=1$ или $n=2$.
Поставим в соответствии вершинам данного $n$-полиэдра $N$ точки пространства
$\R^{2n+1}$ {\it общего положения}, т.е. такие, что никакие $2n+2$ из них не
лежат в одной $2n$-мерной плоскости.
(Это возможно, т.к. если точки $A_1,\dots,A_k$ выбраны, то можно взять
точку $A_{k+1}$ вне объединения $2n$-мерных плоскостей, натянутых на наборы
$2n+1$ точек.)
Отобразим теперь каждый симплекс полиэдра $P^n$ линейно на симплекс,
натянутый на соответсвующие точки в $\R^{2n+1}$.
Получим отображение $f$.

Проверим, что $f$ --- вложение.
Нужно доказать, что образы любых двух симплексов не имеют 'лишних' (т.е.
отличных от образа пересечения) общих точек.
Пусть, напротив образы какого-то $k$-мерного и $l$-мерного симплекса
пересекаются не по образу их общей $s$-мерной грани ($s=-1$ для
непересекающихся симплексов).
Тогда $k+l+1-s$ образов их вершин лежат в некотором $(k+l-s-1)$-мерном
пространстве.
Добавим к этим образам еще $2n+2-(k+l+1-s)$ образов других вершин.
Тогда полученные $2n+2$ точки лежат в $2n$-мерной плоскости, что противоречит
построению отображения $f$.
 QED

\begin{pr}\label{em-ges1}
(a) Любые два вложения окружности в $\R^4$ (объемлемо) изотопны.

(b) Любые две окружности в $\R^4$ не зацеплены, т.е. их можно заключить в
непересекающиеся четырехмерные шары.
\end{pr}

 {\bf Теорема.} {\it Любое $n$-многообразие вложимо в $\R^{2n}$.}

\smallskip
{\it Набросок доказательства для PL категории.}
Будем считать, что $n\ge3$ (случай $n\le2$ разберите самостоятельно) и
что данное $n$-многообразие $N$ связно.
Существует кусочно линейное отображение $f:N\to\R^{2n}$ {\it общего положения},
т.е. имеющее конечное число попарно различных двойных точек
$x_1,y_1,\dots,x_k,y_k$, где $f(x_i)=f(y_i)$ (a).
Тогда на $N-\{x_1,y_1,\dots,x_k,y_k\}$ отображение $f$ является вложением.
Так как $n\ge2$, то можно соединить дугой $l_1$ точку $x_1$ с $y_1$ так, чтобы
дуга $l_1$ не пересекала других точек $x_i,y_i$.
Тогда $f(l_1)$ --- замкнутая кривая в $\R^{2n}$.
Так как $n\ge2$, то эта кривая незаузлена в $\R^{2n}$.
Поэтому cуществует диск $D^2\subset\R^{2n}$ такой, что $\partial D^2=f(l_1)$.
Поскольку $n+2<2n$, то мы можем взять диск $D^2$ {\it общего положения}, т.е.
$D^2\cap f(N^n)=f(l_1)$ (b).
Существует окрестность $B^{2n}$ диска $D^2$ в $\R^{2n}$, гомеоморфная
$2n$-шару, для которой

(1) $B^n=f^{-1}(B^{2n})$ гомеоморфно $n$-шару, причем

(2) граница и внутренность шара $B^n$ переходят при отображении $f$ в границу и
внутренность шара $B^{2n}$, соответственно (с).

Определим отображение $g:B^n\to B^{2n}$ формулой $g(s,r):=(rf(s),r)$, где
$s\in\partial B^n$.
Заменим отображение $f$ на $B^n$ отображением $g$.
Таким образом мы избавимся от двойных точек $x_1,y_1$.
Аналогично для остальных $i=2,\dots,k$ избавимся
от двойных точек $x_i,y_i$ и получим вложение $N$ в $\R^{2n}$.
 QED

\begin{pr}\label{em-ge2n} (a)--(с) Докажите утверждения, отмеченные в наброске соответствующими буквами.
Указание к (c): нужно брать {\it регулярную} окрестность, т.е. 'тесную'
окрестность без 'лишних' дыр.
\end{pr}

 Заметим, что доказательство Уитни {\it гладкой} вложимости любого
$n$-многообразия в $\R^{2n}$ более сложно~\cite{Ad93},~\cite{Pr06}.

Не существует одной теоремы или одного определения, формализующего идею
{\it общего положения}.
Формально, для каждой конкретной задачи определение общего положения свое.
Подробнее см.~\cite[1.6.D]{Ru73},~\cite{RS72}.

Дальнейшие задачи, сформулированны для произвольного $n$, интересно решить
даже для наименьшего $n$, удовлетворяющего условию.

\begin{pr}\label{em-zee}
В этой задаче подразумевается PL категория.

(a) Вложение $S^n\to\R^m$ объемлемо изотопно стандартному тогда и
только тогда, когда оно продолжается до вложения $D^{n+1}\to\R^m$.

(b) Любые два вложения $S^n\to\R^{2n+1}$ объемлемо изотопны при $n\ge2$.

(c) Любые два вложения $S^n\to\R^{2n}$ объемлемо изотопны при $n\ge4$.

(d)* Любые два вложения $S^3\to\R^6$ объемлемо изотопны.

(e) Любые два вложения $S^n\to\R^m$ объемлемо изотопны при $2m\ge3n+4$.

(f)* {\it Теорема Зимана.}
Любые два вложения $S^n\to\R^m$ объемлемо изотопны при
$m\ge n+3$~\cite[\S2]{Sk08}.
\end{pr}

Заметим, что в гладкой категории утверждение \ref{em-zee}.e верно, а утверждения \ref{em-zee}.d,f --- нет!

Для доказательства {\it изотопических аналогов} приведенных теорем (т.е. задач
\ref{em-isot}.ab ниже) полезен следующий факт (доказательство которого непросто).
Определение {\it конкордантности} вложений получается из определения
изотопности отбрасыванием условия сохранения уровней.

\smallskip
{\bf Теорема.}
{\it В коразмерности больше двух конкордантность влечет объемлемую изотопию}.

\smallskip
Доказательство этой теоремы является обобщением доказательства теоремы Зимана о
незаузленности сфер (поэтому очевидное решение задачи \ref{em-zee} с его помощью неразумно).

\begin{pr}\label{em-isot}
(a) Любые два вложения $n$-мерного полиэдра в $\R^{2n+2}$
изотопны при $n\ge2$.
(Эту и следующую теоремы интересно доказать даже для $n=2$.)

(b) Любые два вложения связного $n$-мерного
многообразия в $\R^{2n+1}$ изотопны при $n\ge2$.

(c) Связная сумма трилистника со своим зеркальным образом
конкордантна тривиальному узлу (но не изотопна ему).
\end{pr}

 \subsection{Идея дополнения}

О применении идеи дополнения к {\it вложимости} полиэдра в $S^m$ см. \S\ref{0dual}, пункт `двойственность Александера'. Изотопические нварианты вложений можно получить, рассматривая дополнение $S^m-fN$ до образа многообразия $N$ при вложении в $S^m$.
В самом деле, если $f,g:N\to S^m$ --- два изотопных вложения, то $S^m-fN\cong S^m-gN$.
Поэтому любой топологический инвариант пространства $S^m-fN$ является инвариантом изотопии вложения $f$.
Впервые эту идею  применял Дж. Александер около 1910 г. к изучению узлов в трехмерной сфере (\S\ref{03man}).
Изложение развития идеи дополнения в теории узлов не входит в нашу задачу.
Ограничимся лишь формулировками достаточных условий полноты инварианта дополнения.

\smallskip
{\bf Теорема Папакириякопулоса.} {\it Вложение $f:S^1\to S^3$ изотопно
стандартному тогда и только тогда, когда $\pi_1(S^3-fS^1)\cong\Z$}~\cite{Pa57}.

\smallskip
{\bf Теорема Левина.} {\it Для любого $n\ne2$ гладкое вложение
$f:S^n\to S^{n+2}$ гладко изотопно стандартному тогда и только
тогда, когда $S^{n+2}-fS^n$ гомотопически эквивалентно окружности}~\cite{Le65}.

\smallskip
Для $n=1$ теорема Левина равносильна теореме Папакириякопулоса, поскольку условие $(S^3-fS^1)\simeq S^1$ равносильно условию $\pi_1(S^3-fS^1)\cong\Z$ ввиду асферичности пространства $S^3-fS^1$ \cite{Pa57}.

Аналог теоремы Левина справедлив в кусочно-линейной (PL) и топологической (TOP) локально плоских категориях
при $n\ne2,3$~\cite{Pa57}, \cite{St63}, \cite{Le65}, см. также~\cite{Gl62}.
Напомним, что (PL или TOP) вложение $S^n\subset S^m$ называется (PL или TOP)
{\it локально-плоским}, если каждая точка из $S^n$ обладает такой окрестностью
$U$ в $S^m$, что пара $(U,U\cap S^n)$ (PL или TOP) гомеоморфна паре
$(B^n\times B^{m-n},B^n\times0)$.

Краткая история изучения узлов {\it в коразмерности 1}, т.е. вложений $S^n\to S^{n+1}$, такова.
Известная теорема Жордана, впервые доказанная Брауэром, утверждает, что сфера $S^n$, содержащаяся в $S^{n+1}$, разбивает сферу $S^{n+1}$ на две компоненты.
Легко доказывается, что {\it вложение $f:S^n\to S^{n+1}$ незаузлено тогда и только тогда, когда замыкания этих компонент дополнения являются шарами}.
(Теорема Левина в некотором смысле является аналогом этого утверждения.)
 В 1912 г. Шенфлис доказал, что любая окружность $S^1\subset S^2$ незаузлена.

Утвержение о незаузленности $S^n$ в $S^{n+1}$ получило название {\it гипотезы Шенфлиса.}
В 1921 г.  Александер объявил, что он доказал гипотезу Шенфлиса для произвольного $n$.
Однако, в 1923 г. он нашел контрпример --- знаменитую рогатую сферу Александера
(с помощью той же идеи дополнения~\cite{Ru73}).
Тем не менее, он доказал, что в кусочно-линейной категории гипотеза Шенфлиса верна для $n=2$.


После появления рогатой сферы Александера гипотезой Шенфлиса стали также
называть утверждение о незаузленности {\it локально плоского} вложения
$f:S^n\to S^{n+1}$.
Она была доказана только в 1960-1973 годах:

$\bullet$ для гладкого случая при $n\ne3$ в~\cite{Sm61},~\cite{Ba65}.

$\bullet$ для кусочно-линейного локально плоского случая при $n\ne3$ (следует из TOP
локально плоского случая и результатов Кирби-Зибенмана~\cite{Ho90}).

$\bullet$ для топологического локально плоского случая в~\cite{Br60}, \cite{Ma59}, \cite{Mo60}.

Заметим, что элегантное и короткое доказательство Брауна топологического локально плоского аналога теоремы Смейла-Бардена послужило началом теории 'клеточных множеств', которая стала важной частью геометрической топологии.

Кусочно-линейный случай гипотезы Шенфлиса для $n\ge3$ (или, эквивалентно,
локально плоский кусочно-линейный случай гипотезы Шенфлиса для $n=3$)
остается известной и трудной нерешенной {\it проблемой Шенфлиса}~\cite{RS72}.

\subsection{Общий инвариант дополнения}

Для $n$-мерного полиэдра (т.е., тела симплициального $n$-мерного комплекса) $N$ и вложения $f:N\to S^{n+k}$ положим
$$C_f:=S^{n+k}-Of(N)\simeq S^{n+k}-f(N),$$
где $Of(N)$ --- регулярная (трубчатая) окрестность образа $f(N)$ в $S^{n+k}$.
Топологический и гомотопический тип пространства $C_f$ являются инвариантами вложения $f$.
По двойственности Александера гомологические группы пространства $C_f$ не зависят от $f$ (см. также задачи
\ref{com-simp}.cd).
Поэтому гомотопический тип пространства $C_f$ --- относительно слабый инвариант для $k\ge3$.
Однако он полезен для построения других инвариантов.

Обозначим через $i:S^m\to S^{m+1}$ стандартное включение.

\begin{pr}\label{com-simp}
(a) Для стандартного вложения $f:S^q\to S^m$ имеем $C_f\simeq S^{m-q-1}$.

(b) При $m>p+q$ для стандартного вложения $f:S^p\times S^q\to S^m$ имеем
$C_f\simeq S^{m-p-q-1}\vee S^{m-q-1}\vee S^{m-p-1}$.

(c) При $m-n\ge3$ пространство $C_f$ односвязно.

(d) $C_{i\circ f}\simeq\Sigma C_f$.

(e) $C_f$ ретрагируется на (`меридиональную') окружность для $m=n+2$.

(f) Если $N=S^{n_1}\sqcup\dots\sqcup S^{n_r}$, $m-n_i\ge3$ и
$f:N\to S^m$ --- произвольное вложение, то $C_f\simeq S^{m-n_1-1}\vee\dots\vee S^{m-n_r-1}$.

(g)* {\it Пример Ламбрехтса.} Гомотопический тип дополнения к образу тора $N=S^p\times S^q$ может
зависеть от вложения в $S^m$ при $m-p-q\ge3$.

(h)* Существуют гладкие вложения $f:\R P^2\to S^4$ и $g:\C P^2\to S^7$,
для которых $C_f\simeq \R P^2$ и $C_g\simeq \C P^2$.
(Это `стандартные' вложения; для $\C P^2$ `других' и не существует~\cite{Sk10}.)
\end{pr}

Для решения пункта (f) нужен следующий результат \ref{com-whi}.f.

\begin{pr}\label{com-whi}
(a) $S^1\times S^1$ и $S^1\vee S^1\vee S^2$ имеют изоморфные гомологические группы, но не гомотопически эквивалентны.

(b) $S^2\times \R P^3$ и $S^3\times \R P^2$ имеют изоморфные гомотопические группы, но не гомотопически эквивалентны.

(d) Если все гомотопические группы полиэдра тривиальны, то он стягиваем.

(c) {\it Теорема Уайтхеда.} Отображение (конечномерных) полиэдров, индуцирующее изоморфизм всех гомотопических групп, является гомотопической эквивалентностью.

(e) Односвязный полиэдр, группы гомологий $H_q$ которого тривиальны при $q\ge2$ (т.е. изоморфны группам гомологий точки), стягиваем.

(f) Односвязный полиэдр, группы гомологий  которого изоморфны группам гомологий сферы (или букета сфер),
гомотопически эквивалентен этой сфере (или этому букету сфер).

(g) {\it Теорема Уайтхеда.} Отображение односвязных полиэдров, индуцирующее эпиморфизм на $\pi_2$ и изоморфизм всех гомологических групп, является гомотопической эквивалентностью.
\end{pr}

Для решения пунктов (e,f,g) нужны следующие результаты \ref{com-hur}.de.

Для полиэдра $X$ определим {\it гомоморфизм Гуревича}
$$h=h_{X,k}:\pi_k(X,*)\to H_k(X;\Z)\quad\text{формулой}\quad h(f:S^k\to X):=f_*[S^k].$$
Вот эквивалентное определение.
Заменим отображение $f:S^k\to X$ на гомотопное симплициальное $f':S^k\to X$ для некоторой триангуляции $T$ сферы $S^k$.
Поставим на ориентированный $k$-симплекс $\sigma\subset X$ число симплексов $\tau\in T$, для которых $f'\tau=\sigma$, со знаком.
(Знак симплекса $\tau$ положителен, если $f'|_\tau:\tau\to\sigma$ сохраняет ориентацию, и отрицателен, если меняет.)
Докажем, что эта расстановка является циклом.
Положим $h(f)$ равным ее гомологическому классу.

\begin{pr}\label{com-hur} (a) Это определение корректно.

(b) Даже для $k\ge2$ гомоморфизм Гуревича может не быть мономорфизмом и может не быть эпиморфизмом.

(c) {\it Теорема Гуревича.} Если  $X$ --- полиэдр, $n\ge 1$ и $\pi_k(X,*)=0$ для $0\le k\le n$, то
$H_k(X;\Z)=0$ для $0\le k\le n$ и $h:\pi_{n+1}(X,*)\to H_{n+1}(X;\Z)$ изоморфизм.

(Ср. с теоремой Пуанкаре \ref{0homol}.4e.)

(d) {\it Обратная теорема Гуревича.} Если  $X$ --- связный односвязный полиэдр, $n\ge 1$ и $H_k(X;\Z)=0$ для
$2\le k\le n$, то $\pi_k(X)=0$ для $2\le k\le n$ и $h:\pi_{n+1}(X,*)\to H_{n+1}(X;\Z)$ изоморфизм.

(d') В ситуации из пунктов (c,d) $h_{n+2,X}$ эпиморфизм при $n\ge2$.

(f) Определите относительный гомоморфизм Гуревича $h_{k,X,A}:\pi_k(X,A)\to H_k(X,A;\Z)$.

(e) {\it Относительная теорема Гуревича.} Пусть  $X,A$ --- связные односвязные полиэдры, $A$ подполиэдр в $X$ и $\pi_2(X,A)=0$.

Если $n\ge2$ и $\pi_k(X,A)=0$ для $3\le k\le n$, то $H_k(X,A;\Z)=0$ для $3\le k\le n$ и
$h:\pi_{n+1}(X,A)\to H_{n+1}(X,A;\Z)$ изоморфизм.

Если $n\ge2$ и $H_k(X,A;\Z)=0$ для $3\le k\le n$, то $\pi_k(X,A)=0$ для $3\le k\le n$ и
$h:\pi_{n+1}(X,A)\to H_{n+1}(X,A;\Z)$ изоморфизм.

(e') В ситуации из пункта (e) $h_{n+2,X,A}$ эпиморфизм при $n\ge2$.
\end{pr}

Для решения пунктов (c-e) полезен следующий результат \ref{com-cell}.
Ср. \cite[10.1]{Pr06}.

Пусть в 2-комплексе задан 1-подкомплекс, т.е. граф.
Рассмотрим {\it дополнение до  открытой регулярной окрестности графа}, т.е. объединение симплексов второго барицентрического подразбиения 2-комплекса, не пересекающих граф.
Вложение графа в 2-комплекс называется {\it клеточным}, если каждая связная компонента этого дополнения
гомеоморфна диску (или, эквивалентно, разбивается любой ломаной с концами на границе этой грани).
Например, вложение точки в сферу клеточно.

{\it Клеточным разбиением} 2-комплекса называется клеточное вложение некоторого графа этот 2-комплекс.
Возможно, начинающему будет более понятен термин {\it разбиение на многоугольники}.
Например, 2-комплекс является клеточным разбиением себя.

Пусть в 3-комплексе задан 2-подкомплекс.
Рассмотрим {\it дополнение до  открытой регулярной окрестности 2-подкомплекс}, т.е. объединение симплексов второго барицентрического подразбиения 2-комплекса, не пересекающих 2-подкомплекс.
Вложение 2-подкомплекса в 3-комплекс называется {\it клеточным}, если каждая связная компонента этого дополнения
гомеоморфна 3-диску.
Например, вложение точки в 3-сферу клеточно.

{\it Клеточным разбиением} 3-комплекса называется клеточное вложение некоторого 2-комплекса в этот 3-комплекс, вместе с клеточным разбиением этого 2-комплекса.
Возможно, начинающему будет более понятен термин {\it разбиение на многогранники}.
Например, 3-комплекс является клеточным разбиением себя.

Аналогично определяется {\it клеточное разбиение} $n$-комплекса.

\begin{pr}\label{com-cell}
(a) Если  $X$ --- полиэдр и $\pi_k(X,*)=0$ для $0\le k\le n$, то $X$ гомотопически эквивалентен
полиэдру, имеющему клеточное разбиение c единственной вершиной и без клеток размерностей $1,2,\dots,n$.

(b) Если  $(X,A)$ --- полиэдральная пара и $\pi_k(X,A)=0$ для $0\le k\le n$, то $(X,A)$ гомотопически эквивалентна
полиэдральной паре, имеющей клеточное разбиение  без клеток размерностей $1,2,\dots,n$ в $X-A$.
\end{pr}

\begin{pr}\label{com-bh} Пусть $f:N\to S^7$ --- вложение замкнутого связного 4-многообразия $N$.

(a) $\pi_2(C_f)\cong\Z$.

(b) Если $N$ односвязно, то  $\pi_3(C_f)\cong\Z_{d(\varkappa(f))}$.
\footnote{Можно даже доказать, что
$C_f\simeq C_{\varkappa (f)}:=S^2\cup_{\varkappa (f)}(D^4_1\sqcup\dots\sqcup D^4_{b(N)})$ \cite{Sk10}.
Here $b(N):=\rk H_2(N;\Z)$ and we identify $H_2(N;\Z)$ and $\Z^{b(N)}$ by any isomorphism, so $\varkappa (f)$ is identified with an ordered set of $b(N)$ integers, which set defines a homotopy class of maps
$\partial(D^4_1\sqcup\dots\sqcup D^4_{b_2(N)})\to S^2$.
(This ordered set depends on the identification of $H_2(N;\Z)$ and $\Z^{b(N)}$, but the homotopy type of
$C_{\varkappa (f)}$ does not.
The homotopy equivalence $C_f\simeq C_{\varkappa (f)}$ is not canonical.)}
Чтобы определить $d(\varkappa(f))$, обозначим

$\bullet$ через $A_{f,k}=A_k:H_k(N;\Z)\to H_{k+1}(C_f,\partial;\Z)$ двойственность Александера,

$\bullet$ $\varkappa(f):=A_2^{-1}(A_4[N]\cap A_4[N])\in H_2(N;\Z)$ {\it инвариант Беша-Хефлигера},

$\bullet$ $d(0):=0$ и $d(x):=\max\{k\in\Z\ |\ \text{there is }y\in H_2(N;\Z):\ x=ky\}\quad\text{для}\quad x\ne0.$
\end{pr}


По соображениям общего положения, $C_f$ не зависит от $f$ для $k\ge n+2$.
Это пространство обозначим через $C_k(N)$.
Так как
$$C_{i\circ f}\simeq\Sigma C_f,\quad\text{то}\quad
\Sigma^{K-k}C_f\simeq C_K(N)\quad\text{для}\quad K\ge n+2.$$
{\it Необходимым условием дополнения} для вложимости полиэдра $N$ в $S^{n+k}$ является $(K-k)$-де\-над\-стра\-и\-ва\-е\-мость пространства $C_K(N)$, то есть существование пространства $C$ такого, что $\Sigma^{K-k}C\simeq C_K(N)$.
Если это условие выполняется для $K=n+2$, то оно автоматически выполняется для любого $K\ge n+2$.

\subsection{Комбинация инвариантов дополнения и окрестности}

В этом и следующем пунктах мы работаем в гладкой категории.

По соображениям общего положения нормальное векторное расслоение $\nu_f$ погружения $f:N\to S^{n+k}$
замкнутого $n$-многообразия $N$ не зависит от $f$ для $k\ge n+1$.
Это {\it стабильное нормальное расслоение} обозначим через $\nu_k(N)$.
Так как
$$\nu_{i\circ f}=\nu_f\oplus1,\quad\text{то}\quad
\nu_f\oplus(K-k)\varepsilon=\nu_K(N)\quad\text{для}\quad K\ge n+1.$$
Теорему Смейла-Хирша можно переформулировать так:
{\it если существует $k$-расслоение $\nu$ над $N$ такое, что
$\nu\oplus(K-k)\varepsilon\cong\nu_K(N)$, то $N$ погружаемо в $\R^{n+k}$.}

\begin{pr}\label{com-norm}  $C_k\simeq D\nu_k/S\nu_k$ для $k\ge n+2$.
\end{pr}

  Нормальное расслоение $\nu_f$ погружения $f$ с точностью до эквивалентности является инвариантом регулярной
гомотопии погружения $f$ (и  инвариантом изотопии погружения $f$, если оно является вложением).
Этот изотопический инвариант не очень сильный, потому что, например, нормальные расслоения
различных вложений cтабильно эквивалентны.
Для вложения $f$ пространство расслоения $\nu_f$ и регулярная окрестность
$Of(N)$ образа $f(N)$ в $S^m$ могут быть отождествлены при помощи гомеоморфизма
$\kappa$, при котором нулевое сечение переходит в $f(N)$.
Заметим, что этот гомеоморфизм $\kappa$ определен неоднозначно
(более точно, различающие элементы лежат в группах $H^l(N;\pi_l(\SO_{k-1}))$).
О близком понятии утолщения и утолщаемости см. \cite{Sk}, параграфы `утолщения графов' и
`трехмерные утолщения двумерных полиэдров'.


{\it Совмещая} идею дополнения и идею окрестности, Дж. Левин, С. П. Новиков и
В. Браудер получили {\it редукцию} проблем вложимости и изотопии к
алгебраическим проблемам (впрочем, достаточно сложным).
Доказательство достаточности {\it комбинации} инвариантов дополнения и окрестности условий --- одно из наиболее
важных приложений {\it хирургии} к топологии многообразий.
По поводу развития этого подхода и преодоления возникающих трудностей см.~\cite{Sk08'},~\cite{Sk10}, \cite{CS11}.

Дадим определение пормальной системы и сформулируем теорему Браудера-Левина.
Композиция $S\nu_f\overset\kappa=\partial C_f\subset C_f$ называется
{\it приклеивающим отображением} $a(f,\kappa)$.
Тройка $S(f,\kappa)=(\nu_f,C_f,a(f,\kappa))$ называется {\it нормальной системой пары $(f,\kappa)$}.
Вообще, {\it нормальная система коразмерности $k$ на многообразии $N$} --- это тройка $S=(\nu,C,a)$, состоящая из

(i) векторного $k$-расслоения $\nu$;

(ii) пространства $C$; и

(iii) непрерывного отображения $p:S\nu\to C$.

Будем говорить, что нормальные системы $(\nu,C,p)$ и $(\nu_1,C_1,p_1)$
{\it эквивалентны}, если существуют изоморфизм расслоений $b:\nu\to\nu_1$
и гомотопическая эквивалентность $r:C\to C_1$, для которых
$r\circ p\simeq p_1\circ Sb$.
Очевидно, класс эквивалентности нормальной системы $(f,\kappa)$ не зависит от $\kappa$.
Этот класс называется {\it нормальной системой $S(f)$ вложения $f$}.
Нормальная система вложения является изотопическим инвариантом.
Так как любые два вложения $N\to S^{n+K}$ изотопны при $K>n+1$, то $S(f)$ не зависит от вложения $f$ при $K>n+1$.
Эта нормальная система называется {\it стабильной нормальной системой многообразия $N$} и обозначается $S_K(N)$.

{\it Надстройка} $\Sigma S$ над нормальной системой $(\nu,C,p)$ --- это нормальная система $(\nu\oplus1,\Sigma C,p')$, в которой $p'$ является надстройкой над $p$ на каждом слое.
{\it Необходимое условие Браудера-Левина} для вложимости многообразия $N$ в $S^{n+k}$ состоит в существовании нормальной системы $S$ на $N$, такой что $\Sigma^{K-k}S\sim S_K(N)$ (ясно, что это условие не зависит от $K$).

\smallskip
{\bf Теорема  Браудера-Левина.}
{\it Пусть $N$ --- замкнутое $n$-многообразие, $K>n+1$, $k\ge3$ и существует такая нормальная система $S=(\nu,C,p)$
на $N$, что $\Sigma^{K-k}S\sim S_K(N)$ и $\pi_1(C)=0$.
Тогда существует вложение $f:N\to S^{n+k}$ такое, что $S(f)\sim S$.}

\smallskip
Для гомотопической сферы $N$ теорема Браудера-Левина доказана в~\cite{Le65'}, а в общем случае --- в~\cite{Br68}.
См. также~\cite{Gi71}.
Теорема Браудера-Левина вытекает из следующей леммы~\cite{CRS04}.

\smallskip
{\bf Лемма о компрессии и денадстройке.}
{\it Пусть $F:N\to S^{n+k+1}$ --- вложение замкнутого гладкого $n$-многообразия,
$k\ge3$, $n+k\ge5$ и существует такая нормальная система $S=(\nu,C,p)$
на $N$, что $\Sigma S\sim S(F)$ и $\pi_1(C)=0$.
Тогда существует такое вложение $f:N\to S^{n+k}$, что $S(f)\sim S$.}

Доказательство леммы  о компрессии и денадстройке (принадлежащее Браудеру)
использует
понятие нормального инварианта и приводится в~\cite{CRS04}.

\subsection{Указания и решения к некоторым задачам}

 {\bf \ref{com-simp}.} (b) Используйте (d).

(c) Используйте общее положение.

(e) Постройте ретракцию по остовам некоторого клеточного (или симплициального) разбиения дополнения $C_f$.

(f) Следует из \ref{com-whi}.f ввиду (c) и двойственности Александера (\S\ref{0dual}).

(g) Рассмотрим расслоение Хопфа $S^3\to S^7\overset h\to S^4$.
Возьмем стандартное вложение $S^2\subset S^4$.
Его дополнение гомотопически эквивалентно $S^1$.
Имеем $h^{-1}(S^1)\cong S^1\times S^3\subset S^7$.
Дополнение до этого вложения $f:S^1\times S^3\to S^7$ есть
$$C_f\simeq h^{-1}(S^2)\cong S^2\times S^3\not\simeq S^2\vee S^3\vee S^5\simeq C_{f_0},$$
где $f_0:S^1\times S^3\to S^7$ --- стандартное вложение.

Аналогично строятся два вложения $S^3\times S^7\subset S^{15}$, дополнения
которых не гомотопически эквивалентны.

\smallskip
{\bf \ref{com-whi}.} (a) Гомотопические группы не изоморфны.

(b) Гомологические группы не изоморфны.

(d) Возьмем некоторую триангуляцию полиэдра.
Упорядочим симплексы в соответствии с возрастанием размерности.
Докажем индукцией по $n$, что объединение включение в полиэдр объединения первых $n$ симплексов гомотопно постоянному отображению.

(c) Постройте обратную гомотопическую эквивалентность аналогично (d).

(e) Следует из (d) и обратной теоремы Гуревича \ref{com-hur}.d.

(f) Приведем план решения для сферы, для букетов сфер решение аналогично.
Примените обратную теорему Гуревича \ref{com-hur}.d.
Она даст отображение $f:S^q\to X$, индуцирующее изоморфизм в гомологических группах.
Примените относительную теорему Гуревича \ref{com-hur}.e к цилиндру $\cyl f$ этого отображения и (гомологическую и гомотопическую) точную последовательность пары $(\cyl f, X)$.
Получится, что $f$ индуцирует изоморфизм и в гомотопических группах.
Значит, по (c) $f$ --- гомотопическая эквивалентность.

(g) Аналогично (f).

\smallskip
{\bf \ref{com-hur}.} (c) Следует из теорем о первой `нетривиальной' гомотопической и гомологической группах
полиэдра, представленного клеточным разбиением из задачи \ref{com-cell}.a.

(d) Следует из (c).

(e) Первое утверждение аналогично (c) следует из теорем о первой `нетривиальной' гомотопической и гомологической группах полиэдральной пары,  представленной клеточным разбиением из задачи \ref{com-cell}.b.
Второе утверждение следует из первого.

\smallskip
{\bf \ref{com-bh}.} (a) Следует из односвязности (задача \ref{com-simp}.c), двойственности Александера и обратной теоремы Гуревича (задача \ref{com-hur}.d).

(b) Возьмем отображение $h=h_f:C_f\to\C P^\infty$, соответствующее элементу $A_4[N]\in H_5(C_f,\partial)$.
Так как $\cyl h\simeq\C P^\infty$, получаем
$$\pi_3(C_f)\cong\pi_4(\cyl h,C_f)\cong H_4(\cyl h,C_f)\cong H_4(\cyl h)/h_*H_4(C_f)\cong\Z_{d(\varkappa(f))}.$$
Здесь

$\bullet$ первый и третий изоморфизмы получаются из точных последовательностей пар;

$\bullet$ второй изоморфизм выполнен по относительной теореме Гуревича \ref{com-hur}.e;

$\bullet$ четвертый изоморфизм выполнен, поскольку для двойственного отображения $h^*:H^4(\C P^\infty)\to H^4(C_f)$
и образующей $a\in H^2(\C P^\infty)$ имеем $h^*(a\cup a)=h^*a\cup h^*a=PDA_4[N]\cup PDA_4[N]=PDA_2\varkappa(f)$.

\newpage
\section{Нестандартные гомотопические сферы}\label{0hosp}

\subsection{Пример нестандартной гомотопической сферы }\label{0hosp-con}

В этом пункте приводится построение знаменитого примера Милнора нестандартной гомотопической сферы (\S\ref{0cobint}).
Для
решения большей части задач достаточно владения основами теории гомологий (\S\S \ref{0vesyhom}, \ref{0vesyint}, \ref{0homolpa}, см. также \S\S \ref{02homhom}, \ref{02homint}, \ref{03homhom}).

Обозначим $T_n:=\{(x,y)\in S^n\times S^n\ :\ |x,y|\le1\}$ (трубчатую окрестность диагонали).

\begin{pr} \label{cobtn}
(a) Чему гомеоморфно $T_1$?
\qquad
(b) Чему гомеоморфны $T_3$ и $T_7$?

(с) $T_2$ не гомеоморфно $S^2\times D^2$.
\qquad
(d) $T_4$ не гомеоморфно $S^4\times D^4$.
 \end{pr}

Обозначим $T:=T_4$.

\begin{pr} \label{cobt4con}
Любое ли отображение

(a) $S^1\to\partial T$; \qquad (b) $S^2\to\partial T$; \qquad (c) $S^3\to\partial T$

гомотопно отображению в точку?
\end{pr}

Докажем, что $H_3(\partial T)\ne0$.
Из этого будет вытекать, что $\partial T$ не является гомотопически эквивалентным
сфере $S^7$ и необходимо усложнение конструкции.

Для ориентируемого $2k$-многообразия $M$ обозначим через
$\cap:H_k(M)\times H_k(M)\to\Z$  и $\cap:H_k(M,\partial)\times H_k(M)\to\Z$
его {\it форму пересечений} и {\it билинейное отображение пересечений} (\S\S\ref{02homint},\ref{0vesyint}).
Обозначим через $[\Delta]\in H_4(T)$ гомологический класс диагонали $\Delta\subset S^4\times S^4$.
Рассмотрим фрагмент точной последовательности пары
$$\xymatrix{H_4(T) \ar[r]^j & H_4(T,\partial)  \ar[r]^\partial \ar[d]^{PD} & H_3(\partial T) \ar[r]^i
& H_3(T)\ar[d]^{\cong}\\
& H_4(T)\cong H_4(\Delta)\cong\Z & & H_3(\Delta)=0
}$$
Здесь $i$ --- гомоморфизм включения, $j$ --- гомоморфизм `позволяющий границу', $\partial$ --- граничный
гомоморфизм, $PD$ --- неканонический изоморфизм, полученный из двойственности Лефшеца (задача \ref{pdsimz}.c).


\begin{pr} \label{cobtnhom}
(a') $T$ деформационно ретрагируется (см. определение в задаче \ref{proexc}.e) на $\Delta$.

(a) $[\Delta]$ порождает $H_4(T)$ и $[\Delta]\cap [\Delta]=2$.

(b) Для любых $x,y\in H_4(T)$ выполнено $x\cap jy=x\cap y$.

(c) Существует $x\in H_4(T,\partial)$, для которого $x\cap [\Delta]=1$.

(d) $H_3(\partial T)\ne0$.

(e) $H_3(\partial T)\cong\Z_2$.
 \end{pr}


Обозначим через $p:T\to S^4$ сужение на $T$ проекции $S^4\times S^4\to S^4$ на первый сомножитель.
Обозначим через $(T',p')$ копию пары $(T,p)$.

\begin{pr} \label{cobdiff}
Существует диффеоморфизм $f:(p')^{-1}D^4\to p^{-1}D^4$, переводящий
пересечение с диагональю в $p^{-1}z$, $z\in\Int D^4$.
 \end{pr}

Обозначим $V:=T\cup_fT'$ (в копиях $T$ и $T'$ склеиваются сужения проекций на четырехмерные диски, причем диск из
базы одного сужения отождествляется  со слоем другого.)
После сглаживания получится 8-многообразие $V$ с краем (которое называется
{\it водопроводным соединением} двух копий многообразия $T$).

\begin{pr} \label{cobvhom}
(a) Любое отображение $S^1\to\partial V$ или $S^2\to\partial V$ гомотопно отображению в точку.

(b) $V$ деформационно ретрагируется на $\Delta\vee\Delta'\cong S^4\vee S^4$.

(c) $H_4(V)\cong\Z\oplus\Z$ имеет базис $s_1:=[\Delta]$, $s_2:=[\Delta']$, в котором матрица формы пересечений $H_4(V)\times H_4(V)\to\Z$ имеет вид $\left(\begin{matrix} 2&1\\1&2\end{matrix}\right)$.

(d) Для этого базиса существует такой базис $x_1,x_2$ группы $H_4(T,\partial)$, что $s_i\cap x_j=\delta_{ij}$.

(e) $H_3(\partial V)\ne0$.

(f) $H_3(\partial V)\cong\Z_3$.
\end{pr}

Итак, $\partial V$ не является гомотопически эквивалентным сфере $S^7$ и необходимо усложнение конструкции.


\smallskip
{\it Построение сферы Милнора.}
Рассмотрим граф с вершинами
$1,\dots,8$ и ребрами $12,23,34$, $45,56,67$ и $58$.
Для каждой вершины $a$ графа

возьмем свой экземпляр $(T_a,p_a)$ пары $(T,p)$, и

выберем столько непересекающихся дисков $D^4_{ab}\subset S^4$, сколько
вершин $b$ соединено с $a$.

Для каждого ребра $ab$ склеим сужения расслоений на четырехмерные диски,
соответствующие концам этого ребра, отождествляя базу одного сужения расслоения со слоем другого.
Т.е. склеим $p_a^{-1}D^4_{ab}$   и $p_b^{-1}D^4_{ba}$, как при построении многообразия $V$.

 После сглаживания получится 8-многообразие $W$ с краем.
(Оно называется {\it водопроводным соединением} копий многообразия $T$
сообразно графу.)
Край $\partial W$ этого многообразия и есть 7-многообразие из примера Милнора.

 \smallskip
Пример сферы Милнора вытекает из следующих задач \ref{cobw}.e и \ref{hosp-connd}.a.

\begin{pr} \label{cobw}
(a) Любое отображение $S^1\to\partial W$ или $S^2\to\partial W$ гомотопно отображению в точку.

(c) Найдите матрицу формы пересечений многообразия $W$.

(b) Выберите базисы  $s_1,\dots,s_8$ и $x_1,\dots,x_8$ групп $H_4(W)$ и
$H_4(W,\partial)$, для которых $s_i\cap x_j=\delta_{ij}$.

(d) $H_3(\partial W)=0$.

(e) Край $\partial W$ гомотопически эквивалентен $S^7$.

(Заметим, что по теореме Смейла $\partial W$ топологически гомеоморфно $S^7$.)

(f) Сигнатура (формы пересечений) многообразия $W$ равна 8.
 \end{pr}

Читатель, не знакомый с понятием расслоения, может не решать следующей задачи.

\begin{pr}\label{cobss} (a) $S^n\times S^n$ вложимо в $\R^{2n+1}$.

(b) Сумма касательного расслоения к $S^n\times S^n$ и одномерного тривиального
расслоения над $S^n\times S^n$ тривиальна.

(c) Если $p:E\to B$ --- векторное расслоение со слоем $\R^n$, расслоение $p\oplus\varepsilon$ тривиально и
$b:=\dim B<n$, то расслоение $p$ тривиально.

(d) Если нормальное расслоение к многообразию в $\R^m$ тривиально, то сумма касательного и одномерного тривиального тривиальна.

(e) Для связного многообразия с непустым краем если сумма касательного и одномерного тривиального расслоений тривиальна, то касательное расслоение тривиально.

(f) Многообразие $T_n$ {\it параллелизуемо}, т.е. имеет семейство из $2n$
касательных векторных полей, линейно независимых в каждой точке.

(g) Многообразие $V$ параллелизуемо.

(h) Многообразие $W$ параллелизуемо.
\end{pr}

В следующей задаче можно использовать предыдущую и следствие формулы Хирцебруха о сигнатуре
(\S\ref{0cobint}; от него достаточно делимости на 7).

\begin{pr}\label{hosp-connd}
(a) $\partial W$ не диффеоморфно $S^7$.

(b) $\partial W$ не вложимо гладко в $S^8$.
\end{pr}

\begin{pr} \label{cob4w}
Построим аналогично 4-мерное многообразий $W$ с краем.

(a) $H_1(\partial W)=0$.

(b) $\partial W$ не диффеоморфно $S^3$.

(c) $\partial W$ не вложимо гладко в $S^4$.

В пунктах (b) и (c) можно использовать теорему Рохлина о сигнатуре, см. \S\ref{0cobint}.
(Впрочем, (b) делается и без этого, с использованием $\pi_1$.)
\end{pr}


\begin{pr}\label{hosp-con-2}
{\it Второе построение сферы Милнора}~\cite[28]{DNF84}.

(a) Постройте расслоение $S^7\overset{S^3}\to S^4$ (ср. \ref{geho-fihh}.a). Найдите его класс Эйлера.

(b) Пространство $S^3$-расслоения $\xi$ над $S^4$ гомотопически эквивалентно
$S^7$ тогда и только тогда, когда $e(\xi)=1$.

(c) Определим отображение $f_{hj}:S^3\to SO_4$ формулой $f_{hj}(u)v:=u^hvu^j$.
Для соответствующего расслоения $\xi_{hj}:S^7\overset{S^3}\to S^4$ имеем $e(\xi_{hj})=h+j$.

(d)* $p_1(\xi_{hj})=\pm2(h-j)$.

(e) Выведите пример сферы Милнора из предыдущих пунктов и следствия формулы Хирцебруха о сигнатуре.
\end{pr}

 \subsection{Конечность множества гомотопических сфер (набросок)}\label{hosp-cla}

В этом пункте приводится
набросок доказательства знаменитой теоремы Кервера-Милнора о конечности множества гомотопических сфер.

\smallskip
{\bf Теорема Кервера-Милнора.}
{\it Для любого $n\ge6$ множество $n$-многообразий, гомотопически эквивалентных $S^n$,
с точностью до диффеоморфизма, конечно.}

\begin{pr}\label{}
Этот результат равносилен следующему: для любого $n\ge6$ множество $\theta_n$ ориентированных $n$-многообразий, гомотопически эквивалентных $S^n$ (гомотопических сфер), с точностью до сохраняющего ориентацию диффеоморфизма, конечно.
\end{pr}

 {\bf Лемма.} {\it Нормальное расслоение вложения гомотопической сферы в $\R^m$ тривиально для большого $m$.}

\smallskip
Препятствие к тривиальности лежит в $\pi_{n-1}(SO)$.
Поэтому лемма верна для $n\equiv3,5,6,7\bmod8$ ввиду теоремы периодичности Ботта (\S\ref{0hopfbun}).
Для других $n$ доказательство более сложно \cite{KM63}.

О конструкции Понтрягина см., например, \S\ref{0hopf}, \cite[\S18]{Pr04}.

Для гомотопической сферы $N\subset \R^m$ c нормальным оснащением $\zeta$ обозначим через $p(N,\zeta)\in\pi_m(S^{m-n})\cong\pi_n^S$ класс оснащенного кобордизма.

Для стандартной сферы $S^n\subset \R^m$ и $x\in\pi_n(SO_{m-n})$ рассмотрим оснащение нормального
расслоения, полученное из стандартного оснащения при помощи $x$.
Обозначим через $J(x)\in\pi_m(S^{m-n})$ класс оснащенного кобордизма полученного оснащенного многообразия.
В \S\ref{0hopf} доказано, например, что $J:\pi_1(SO_2)\to\pi_3(S^2)$ есть изоморфизм.

\begin{pr}\label{hosp-j}
(a) $J:\pi_n(SO_{m-n})\to\pi_n(S^{m-n})$ гомоморфизм.

(b) $p(N,x\zeta)=p(N,\zeta)+J(x)$.
\end{pr}

 Для большого $m$ имеем $\pi_n(SO_{m-n})\cong\pi_n(SO)$.
Поэтому предыдущая конструкция дает отображение $J:\pi_n(SO)\to\pi_n^S$.
В \S\ref{0hopf} доказано, например, что $J:\pi_1(SO)\to\pi_1^S$ есть изоморфизм.
Кроме того, $J:\pi_3(SO)\to\pi_3^S$ есть приведение по модулю 24.

Ввиду задачи \ref{hosp-j} отображение $p:\theta_n\to \pi_n^S/\im J$ корректно определено формулой
$p(N):=p(N,\zeta)+\im J$.

\begin{pr}\label{hosp-clap}
(a) Операция связного суммирования превращает $\theta_n$ в группу.

(b) $p$ гомоморфизм.

(c) $p(N)=0$ тогда и только тогда, когда $N$ является границей параллелизуемого многообразия.
\end{pr}

Ввиду (a,b) и поскольку группа $\pi_n^S$ конечна для $n>0$, достаточно доказать, что $\ker p$ конечно.

См. задачу \ref{cobconn}.

\begin{pr}\label{}
(a) {\it Перестройка} дает многообразие, кобордантное исходному.

(b)* Обратно, если многообразия кобордантны, то одно можно получить
из другого перестройками.
\end{pr}

\begin{pr}\label{}
(a) Любое ориентируемое многообразие размерности $\ne3$ кобордантно односвязному.

(b) Любое спинорное (т.е. параллелизуемое в окрестности двумерного
остова некоторой триангуляции) многообразие размерности $\ge6$ кобордантно двусвязному.

(c) Любое спинорное многообразие размерности $\ge8$ параллелизуемо в окрестности трехмерного остова некоторой триангуляции и потому кобордантно трехсвязному.
\end{pr}

В следующих задачах $N=\partial W$ гомотопическая $n$-сфера и $W$ параллелизуемо.
По поводу задач, отмеченных звездочками, см. \cite{KM63}.

\begin{pr}\label{} В этой задаче $n\ge5$.

(a) Если $W$ стягиваемо, то $N\cong S^n$.

(b)* Можно так выбрать оснащение сферы $S^i\subset W$, чтобы результат
перестройки по этой сфере с этим оснащением был параллелизуемым.

(c)  Перестройками можно добиться того, чтобы $W$ стало $([n/2]-1)$-связным.
\end{pr}

\begin{pr}\label{}
Пусть $n=2k$ и $W'$ получено из $W$ перестройкой сферы $S^k\subset W$.

(a) Существует $x\in H_k(W')$, для которого $H_k(W')/x\cong H_k(W)/[S^k]$.

(b) Если $[S^k]\in H_k(W)$ примитивен, то $H_k(W')\cong H_k(W)/[S^k]$.

(c) Если $k\ge4$ четно, то $\rk H_k(W')\ne\rk H_k(W)$.

(d)* Если $k\ge4$ четно, то $N\cong S^n$.

(e)* Если $k\ge3$ нечетно, то $N\cong S^n$.
\end{pr}

\begin{pr}\label{}
(a)* Если $n=4l-1\ge7$ и $\sigma(W)=0$, то $N\cong S^n$.

(b)* Существует замкнутое почти параллелизуемое $4l$-многообразие $M$,
для которого $\sigma(M)\ne0$.

(c) Если $n=4l-1\ge7$ и $\sigma(W)$ делится на $\sigma(M)$, то $N\cong S^n$.
\end{pr}

\begin{pr}\label{}
Пусть $n=4l+1\ge17$.
Перестройками можно добиться того, чтобы $W$ стало $2l$-связным (и осталось параллелизуемым).
Реализуем элемент $x\in H_{2l+1}(W)$ вложением $x:S^{2l+1}\to W$.
Обозначим через
$$q(x)\in\ker[i_*:\pi_{2l}(SO_{2l+1})\to\pi_{2l}(SO)]\cong\Z_2$$
препятствие к тривиальности нормального расслоения вложения $x$.
Обозначим $Arf(N):=Arf(q):=\sum_iq(a_i)q(b_i)\in\Z_2$, where
$a_1,\dots,a_m,b_1,\dots,b_m$ is a symplectic basis of $H_{2l+1}(W)$.

(a) $q$ квадратичная форма над $\Z_2$.

(b) $Arf(N)$ действительно зависит только от $N$.

(c) Если $Arf(N)=0$, то $N\cong S^n$.
\end{pr}

Заметим, что $Arf(N)=0$ для $l\ne1,3,7,15,31$.
Это решение знаменитой проблемы Кервера, полученное около 2008 г.

\begin{pr}\label{}
 Выведите теорему Кервера-Милнора из предыдущего.
\end{pr}

 \subsection{Указания и решения к некоторым задачам}

{\bf \ref{cobtn}.} (a) $S^1\times D^1$.

(b) $S^3\times D^3$ и $S^7\times D^7$.

(c) $H_1(\partial T_2)\ne0$ аналогично задаче \ref{cobtnhom}.d.
Или $\pi_1(\partial T_2)\ne0$ из точной гомотопической последовательности расслоения
(задача \ref{geho-fiexb}).

(d) Аналогично (c) следует из задачи \ref{cobtnhom}.d или из $\pi_3(\partial T_4)\ne0$.

\smallskip
{\bf \ref{cobt4con}.} Ответы: (a), (b) да, (c) нет.

(a) По задаче \ref{cobtnhom}.a' $\pi_1(T)\cong\pi_1(\Delta)$.
(Ср. со вторым способом решения задачи \ref{orihocl}.b.)
Так как $\pi_1(\Delta)=\pi_1(S^4)=0$, то любое отображение $S^1\to\partial T$ продолжается до отображения
$D^2\to T$.
Так как $2+4<8$, то отображение $D^2\to T$ можно изменить малым шевелением, чтобы
оно перестало пересекать диагональ $\Delta$.
Отображение $D^2\to T$, не пересекающее $\Delta$, гомотопно отображению $D^2\to\partial T$.

Замечание. Эту задачу можно решать, используя точную гомотопическую последовательность расслоения.
Для $S^1$ и $S^2$ получится по сути то же решение, что и выше, но более сложно изложенное.

(b) Аналогично (a).

(c) По теореме Гуревича (задача \ref{com-hur}.c) следует из (a,b) и задачи \ref{cobtnhom}.d.
Или $\pi_3(\partial T)\ne0$ из точной гомотопической последовательности расслоения (задача \ref{geho-fiexb}).

\smallskip
 {\bf \ref{cobtnhom}.}
(a) По (a') $H_4(T)\cong H_4(\Delta)=\left<\Delta\right>$.
(Ср. со вторым способом решения задачи \ref{orihocl}.b.)
Число $[\Delta]\cap [\Delta]=2$ равно сумме индексов особых точек касательного векторного поля общего положения на $S^4$.

(c) $x:=[T\cap(*\times S^2)]$.

(d) Так как $x\cap[\Delta]=1$ нечетно, то по (a,b) $x\not\in\im j=\ker\partial$.
Поэтому $\partial\ne0$.
Значит, $H_3(\partial T)=0$.

\smallskip
{\bf \ref{cobvhom}.}
(a) Аналогично задаче \ref{cobt4con}.ab.

{\it Более сложный способ.}
Используя теорему Зейферт-Ван Кампена, докажите $\pi_1(\partial V)=0$.
Из точной последовательности пары (ср. диаграмму перед задачей \ref{cobtnhom}), (b) и двойственности Лефшеца
(задача \ref{pdsimz}.c) получаем $H_2(\partial V)=0$.
По теореме Гуревича (задача \ref{com-hur}.c) $\pi_2(\partial V)=0$.

(с) Аналогично задаче \ref{cobtnhom}.b.

(d) Аналогично задаче \ref{cobtnhom}.c.

(e) Аналогично задаче \ref{cobtnhom}.d.
Пусть, напротив, $H_3(\partial V)=0$.
Тогда в точной последовательности пары $(V,\partial V)$ имеем $\partial=0$.
Значит, $j$ эпиморфизм.
Поэтому найдутся такие $a,b\in\Z$, что $x_1=j(as_1+bs_2)$.
Из (c,d) получаем $2a+b=1$, $a+2b=0$.
Складывая, получаем, что 1 делится на 3.
Противоречие.

(f) Аналогично задаче \ref{cobtnhom}.e.

\smallskip
{\bf \ref{cobw}.}
(a) Аналогично задаче \ref{cobt4con}.ab.

(c) Ответ: $a_{ii}=2$, $a_{ij}=1$ или $a_{ij}=0$ сообразно тому, соединены ли в графе вершины $i,j$  ребром или нет.

(b) Аналогично задаче \ref{cobtnhom}.c.

(d) (Ср. [Br72, V.2.7].)
Достаточно доказать, что $j:H_4(W)\to H_4(W,\partial)$ эпиморфно.
Для произвольного $x'\in H_4(W,\partial)$ определим линейную функцию $f_{x'}:H_4(W)\to\Z$ формулой
$f_{x'}(y)=x'\cap y$.
По (c) форма пересечений $H_4(W)\times H_4(W)\to\Z$ унимодулярна.
Поэтому существует $x\in H_4(W)$, для которого $f_{x'}(y)=x\cap y$ для любого $y\in H_4(W)$.
По двойственности Пуанкаре $x'-jx=0$.
Итак, $j$ эпиморфно.

\smallskip
{\bf \ref{cobss}.}  (b) Следует из (a).

(c)  Препятствия к тривиальности обоих расслоений одинаковы (и, значит, нулевые), поскольку
отображение включения $\pi_b(SO_n)\to \pi_b(SO_{n+1})$ является изоморфизмом (задача \ref{geho-fiso}.b).

(d,e) Следует из (c).

(f) Следует из (b,e).

\smallskip
{\bf \ref{hosp-con-2}.}  (b) Гомоморфизм $\pi_4(S^4)\to\pi_3(S^3)$ из точной последовательности
расслоения является умножением на $e(\xi)$.

\smallskip
{\bf \ref{hosp-clap}.} Используйте задачу \ref{cobss}.de.

\newpage
\section{Приложение: Классификация сечений и интегрируемые системы}\label{0sec}

 \subsection{Классификация сечений и зейфертовых сечений}

Обозначим через $\pi_X:X\times Y\to X$ проекцию на первый сомножитель.
Отображение $s:X\to X\times Y$ называется {\it сечением}, если
$\pi_X\circ s=\id_X$.
Сечения называются {\it послойно гомотопными}, если они гомотопны в
классе сечений.
Обозначим через $S(\pi_X)$ множество классов послойной гомотопности сечений.
Равенство между множествами означает наличие взаимно-однозначного
соответствия между ними.

\begin{pr}\label{}
(a) $S(\pi_{S^1}:S^1\times I\to S^1)=\{0\}$.

(b) $S(\pi_{S^1}:S^1\times S^1\to S^1)=\Z$.

(c) $S(\pi_N:N\times S^1\to N)=H_1(N;\Z)$ для замкнутого
2-многообразия~$N$.

(d) $S(\pi_X)$ находится во взаимно-однозначном соответствии с множеством
$[X;Y]$ отображений $X\to Y$ с точностью до гомотопности (см. определение в начале \S\ref{0hopf}).
\end{pr}

 Понятия векторного поля и сечения (прямого произведения) являются важнейшими
частными случаями понятия {\it сечения расслоения}, см. \S\ref{0vebun}.

 Рассмотрим 2-многообразие $P^*$ с непустым краем и инволюцией $t$ на
$P^*$, имеющей лишь конечное число неподвижных точек.

\begin{pr}\label{}
  Существуют такие негомеоморфные 2-многообразия $P^*$ и $P^*_1$ с
непустыми краями и инволюциями $t$ и $t_1$, имеющими лишь конечное число
неподвижных точек, для которых $P^*/t\cong P^*_1/t_1$.
\end{pr}

 Сечение $s:P^*\to P^*\times S^1$ называется {\it симметричным}, если
$p_{S^1}\circ s=p_{S^1}\circ s\circ t$.
Множество симметричных сечений с точностью до гомотопии в классе
симметричных сечений обозначается кратко $S(P^*,t)$.
Положим $P=P^*/t$.

\begin{pr}\label{}
$S(P^*,t)=H_1^t(P^*,\partial P^*;\Z)\cong H_1(P,\partial P;\Z)\cong\Z^{1-\chi(P)}$,
если инволюция $t$ имеет лишь конечное число неподвижных точек.
Для пространства с инволюцией симметричные гомологические группы
обозначаются добавлением верхнего индекса $t$ к обычному обозначению.
Их определение естественно появляется при изучении множества $S(P^*,t)$.
Читатель может придумать его сам или подсмотреть в \cite{Sk}, параграфе `гомотопическая классификация отображений', пункте 'эквивариантные отображения'.
\end{pr}

В приложениях (см. следующий пункт) появляются аналоги сечений для расслоений с особенностями.
Приведем определение и результат для частного случая (мотивированного приложениями).

Положим $Q=Q(P^*):=P^*\times I/\{(a,0)\sim(t(a),1)\}$.
(Выражаясь научно, это расслоение над $S^1$ со слоем $P^*$ и сшивающим
отображением $t$~\cite[определение 2.2]{BFM90}.)

\begin{pr}\label{}
  $Q(P^*)$ зависит только от $P$, а не от $P^*$.
\end{pr}

 Обозначим через $p:P^*\to P$ проекцию. Определим
отображение $\pi:Q\to P$ формулой $\pi[(a,t)]=p(a)$. (Это
расслоение Зейферта, имеющее сингулярные слои только над
звездочками и только типа $(2,1)$.)

Вложение $\lessmskips{2mu}f:P^*\to Q$ называется {\it зейфертовым
сечением}, если $\lessmskips{2mu}\pi\circ f\hm=p$. (В гладкой
категории нужно дополнительно предполагать, что $f$ трансверсально
слоям отображения $\pi$. В~\cite{BF94} сечения Зейферта
назывались трансверсальными площадками.)

\smallskip
{\bf Теорема классификации зейфертовых сечений.} \cite{RS99'}
{\it При фиксированных $P^*$ и $t$ множество $X$ зейфертовых сечений с
точностью до изотопии над $\pi$ находится во взаимно-однозначном соответствии
с $H_1(P,\partial P;\Z)\cong\Z^{1-\chi(P)}$.}

\begin{figure}[h]\centering
\includegraphics{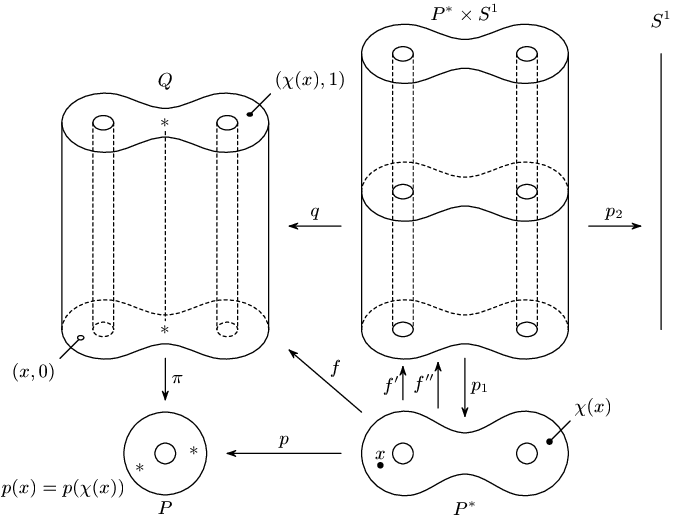}
\caption{Классификация зейфертовых сечений}\label{seifert}
\end{figure}

{\it Доказательство.} Определим отображение
$$
q:P^*\times S^1\cong\frac{P^*\times I}{\{(a,0)\sim(a,1)\}}
\to\frac{P^*\times I}{\{(a,0)\sim(t(a),1)\}}\cong Q
$$
формулой
$$
q[(a,u)]=\begin{cases} [(a,2u)]&0\le u\le\frac12 \\
[(t(a),2u-1)]& \frac12\le u\le1\end{cases}.
$$
Поскольку $t$ инволюция, $q$ корректно определено и непрерывно.

Пусть $f:P^*\to Q$ --- сечение Зейферта.
Для каждой точки $x\in P^*$, не являющейся неподвижной точкой инволюции
$t$, существует {\it единственная} такая точка
$$
f'(x)\in P^*\times S^1,\quad\text{что}\quad qf'(x)=f(x).
$$
Для каждой неподвижной точки $x$ инволюции $t$ существуют {\it две} такие
точки
$$
u,v\in S^1,\quad\text{что}\quad q(x,u)=q(x,v)=f(x).
$$
Поскольку маленькая проколотая дисковая окрестность в $P^*$
неподвижной точки $x$ связна, мы можем взять в качестве $f'(x)$
либо $(x,u)$, либо $(x,v)$, так что отображение $f':P^*\to
P^*\times S^1$ будет непрерывным. Построенное отображение $f'$ ---
(обычное) сечение тривиального расслоения $P^*\times S^1\to P^*$.

Поскольку $f$ --- вложение, $p_{S^1}f'(x)\ne-p_{S^1}f'(t(x))$
ни для какой точки $x\in P^*$.
Поэтому $f'$ послойно гомотопно симметричному сечению $f''$.

Обратно, для любого симметричного сечения $F:P^*\to P^*\times S^1$,
отображение $q\circ F$ является сечением Зейферта и $(q\circ F)''=F$.

Очевидно, сечения Зейферта $f$ и $g$ изотопны над $\pi$
тогда и только тогда, когда соответствующие симметрические
сечения $f''$ и $g''$ симметрично гомотопны (или,
эквивалентно, изотопны). Тогда
$\lessmskips{2mu}X=S(P^*,t)\hm=H_1(P,\partial P;\Z)\cong\Z^{1-\chi(P)}$.
 QED


 \subsection{Применение к интегрируемым гамильтоновым системам}\nopagebreak

В 1994 Болсинов и Фоменко~\cite{BF94} доказали теорему о
топологической траекторной классификации невырожденных
интегрируемых гамильтоновых систем с двумя степенями
свободы на трехмерных многообразиях постоянной энергии.
Мотивировки и краткий обзор см.
в~\cite[\S1]{BFM90},~\cite[\S1]{BF94}. Они доказали, что
две такие системы эквивалентны, если некоторые их
инварианты совпадают. Инвариантом является граф с
некоторыми дополнительными метками на его вершинах и
ребрах. Одним из необходимых ограничений было то, что
рассматриваемая гамильтонова система не имеет неустойчивых
периодических орбит с {\it неориентируемой} сепаратрисой.
Это связано с тем, что в доказательстве использовались {\it сечения локально-тривиальных} расслоений.
Поскольку периодические орбиты с неориентируемой сепаратрисой встречаются в
примерах (например, в волчке Ковалевской), было интересно отбросить это
ограничение. Это оказалось возможным благодаря использованию {\it зейфертовых
сечений зейфертовых расслоений}.

\smallskip
{\bf Теорема.} {\it Пусть $X$ --- множество невырожденных интегрируемых
гамильтоновых систем с двумя степенями свободы на ориентируемых трехмерных
многообразиях постоянной энергии, с точностью до непрерывной траекторной
эквивалентности, сохраняющей ориентацию.
Тогда существует инъекция из $X$ в множество $t$-оснащенных графов  с точностью
до $t$-эквивалентности.}~\cite{RS99'},~\cite{CRS00}, ср.~\cite[Теорема 4.1]{BF94}.

\smallskip
Определения $t$-оснащенных графов и $t$-эквивалентности такое же, как в~\cite{BF94},
см. также~\cite{BFM90},~\cite{CRS00}.
Заметим, что в~\cite[\S12.3 и \S13.5]{BF94} описан образ инъекции из теоремы и зависимость $t$-меток
от ориентации 3-многообразия постоянной энергии.
Кроме того, в~\cite[\S13]{BF94} было построено в некотором смысле более простое оснащение на $W$, названное $t$-молекулой.
Более общая ситуация не требует добавлений и изменений по сравнению с~\cite{BF94}, только $P$-метки могут быть
атомами со звездочками.

Доказательство теоремы основано на следующем общем наблюдении, которое может
быть применено и к другим проблемам.
Бифуркация торов Лиувилля в боттовой интегрируемой гамильтоновой системе может быть описана окрестностью множества $F^{-1}(c)$, где $F$ --- дополнительный интеграл и $c$ --- его критическое значение.
Если критическое подмногообразие дополнительного интеграла $F$, соответствующее $c$, --- окружность, то
эта окрестность является расслоением Зейферта $Q$ над (незамкнутой) 2-поверхностью $P$~\cite{BFM90}.

Чтобы исследовать бифуркацию торов Лиувилля, построим {\it сечение Пуанкаре} потока на $Q$~\cite{BF94}.
Если критическая окружность имеет {\it ориентируемую} сепаратрисную диаграмму
(или, эквивалентно, $P$ не имеет звездочек), то $Q\cong P\times S^1$ и сечения
Пуанкаре могут быть выбраны среди сечений тривиального расслоения.
Тогда сечения Пуанкаре классифицируются методами классической теории препятствий.
Если критическая окружность имеет {\it неориентируемую} сепаратрисную
диаграмму (или, эквивалентно, $P$ не имеет звездочек), то
расслоение Зейферта не является локально-тривиальным расслоением.
Несмотря на это, сечения Пуанкаре являются зейфертовыми аналогами сечений локально-тривиального расслоения.
Доказательство вышеприведенной теоремы с использованием теоремы классификации
зейфертовых сечений аналогично~\cite{BF94}, см. детали в~\cite{CRS00}.







\newpage
\centerline{\bf Литература.}

Звездочками отмечены книги, обзоры и популярные статьи.


\end{document}